\documentclass[12pt]{article}
\usepackage{full page, 
amssymb, amscd, graphicx, 
}
    \title{{\bf Logarithmic
tensor category theory, VI: Expansion condition, associativity of logarithmic
intertwining operators, and the associativity isomorphisms}}
    \author{Yi-Zhi Huang, James Lepowsky and Lin Zhang}
    \date{}

\begin{document}
    \bibliographystyle{alpha}
    \maketitle

    \newtheorem{rema}{Remark}[section]
    \newtheorem{propo}[rema]{Proposition}
    \newtheorem{theo}[rema]{Theorem}
   \newtheorem{defi}[rema]{Definition}
    \newtheorem{lemma}[rema]{Lemma}
    \newtheorem{corol}[rema]{Corollary}
     \newtheorem{exam}[rema]{Example}
\newtheorem{assum}[rema]{Assumption}
     \newtheorem{nota}[rema]{Notation}
        \newcommand{\ba}{\begin{array}}
        \newcommand{\ea}{\end{array}}
        \newcommand{\be}{\begin{equation}}
        \newcommand{\ee}{\end{equation}}
        \newcommand{\bea}{\begin{eqnarray}}
        \newcommand{\eea}{\end{eqnarray}}
        \newcommand{\nno}{\nonumber}
        \newcommand{\nn}{\nonumber\\}
        \newcommand{\lbar}{\bigg\vert}
        \newcommand{\p}{\partial}
        \newcommand{\dps}{\displaystyle}
        \newcommand{\bra}{\langle}
        \newcommand{\ket}{\rangle}
 \newcommand{\res}{\mbox{\rm Res}}
\newcommand{\wt}{\mbox{\rm wt}\;}
\newcommand{\swt}{\mbox{\scriptsize\rm wt}\;}
 \newcommand{\pf}{{\it Proof}\hspace{2ex}}
 \newcommand{\epf}{\hspace{2em}$\square$}
 \newcommand{\epfv}{\hspace{1em}$\square$\vspace{1em}}
        \newcommand{\ob}{{\rm ob}\,}
        \renewcommand{\hom}{{\rm Hom}}
\newcommand{\C}{\mathbb{C}}
\newcommand{\R}{\mathbb{R}}
\newcommand{\Z}{\mathbb{Z}}
\newcommand{\N}{\mathbb{N}}
\newcommand{\A}{\mathcal{A}}
\newcommand{\Y}{\mathcal{Y}}
\newcommand{\Arg}{\mbox{\rm Arg}\;}
\newcommand{\comp}{\mathrm{COMP}}
\newcommand{\lgr}{\mathrm{LGR}}

\newcommand{\dlt}[3]{#1 ^{-1}\delta \bigg( \frac{#2 #3 }{#1 }\bigg) }

\newcommand{\dlti}[3]{#1 \delta \bigg( \frac{#2 #3 }{#1 ^{-1}}\bigg) }

 \makeatletter
\newlength{\@pxlwd} \newlength{\@rulewd} \newlength{\@pxlht}
\catcode`.=\active \catcode`B=\active \catcode`:=\active \catcode`|=\active
\def\sprite#1(#2,#3)[#4,#5]{
   \edef\@sprbox{\expandafter\@cdr\string#1\@nil @box}
   \expandafter\newsavebox\csname\@sprbox\endcsname
   \edef#1{\expandafter\usebox\csname\@sprbox\endcsname}
   \expandafter\setbox\csname\@sprbox\endcsname =\hbox\bgroup
   \vbox\bgroup
  \catcode`.=\active\catcode`B=\active\catcode`:=\active\catcode`|=\active
      \@pxlwd=#4 \divide\@pxlwd by #3 \@rulewd=\@pxlwd
      \@pxlht=#5 \divide\@pxlht by #2
      \def .{\hskip \@pxlwd \ignorespaces}
      \def B{\@ifnextchar B{\advance\@rulewd by \@pxlwd}{\vrule
         height \@pxlht width \@rulewd depth 0 pt \@rulewd=\@pxlwd}}
      \def :{\hbox\bgroup\vrule height \@pxlht width 0pt depth
0pt\ignorespaces}
      \def |{\vrule height \@pxlht width 0pt depth 0pt\egroup
         \prevdepth= -1000 pt}
   }
\def\endsprite{\egroup\egroup}
\catcode`.=12 \catcode`B=11 \catcode`:=12 \catcode`|=12\relax
\makeatother

\def\hboxtr{\FormOfHboxtr} 
\sprite{\FormOfHboxtr}(25,25)[0.5 em, 1.2 ex] 

:BBBBBBBBBBBBBBBBBBBBBBBBB |
:BB......................B |
:B.B.....................B |
:B..B....................B |
:B...B...................B |
:B....B..................B |
:B.....B.................B |
:B......B................B |
:B.......B...............B |
:B........B..............B |
:B.........B.............B |
:B..........B............B |
:B...........B...........B |
:B............B..........B |
:B.............B.........B |
:B..............B........B |
:B...............B.......B |
:B................B......B |
:B.................B.....B |
:B..................B....B |
:B...................B...B |
:B....................B..B |
:B.....................B.B |
:B......................BB |
:BBBBBBBBBBBBBBBBBBBBBBBBB |

\endsprite

\def\shboxtr{\FormOfShboxtr} 
\sprite{\FormOfShboxtr}(25,25)[0.3 em, 0.72 ex] 

:BBBBBBBBBBBBBBBBBBBBBBBBB |
:BB......................B |
:B.B.....................B |
:B..B....................B |
:B...B...................B |
:B....B..................B |
:B.....B.................B |
:B......B................B |
:B.......B...............B |
:B........B..............B |
:B.........B.............B |
:B..........B............B |
:B...........B...........B |
:B............B..........B |
:B.............B.........B |
:B..............B........B |
:B...............B.......B |
:B................B......B |
:B.................B.....B |
:B..................B....B |
:B...................B...B |
:B....................B..B |
:B.....................B.B |
:B......................BB |
:BBBBBBBBBBBBBBBBBBBBBBBBB |

\endsprite


\begin{abstract}
This is the sixth part in a series of papers in which we introduce and
develop a natural, general tensor category theory for suitable module
categories for a vertex (operator) algebra.  In this paper (Part VI),
we construct the appropriate natural associativity isomorphisms
between triple tensor product functors. In fact, 
we establish a ``logarithmic operator product
expansion'' theorem for logarithmic intertwining operators.
In this part, a great deal of
analytic reasoning is needed; the statements of the main theorems
themselves involve convergence assertions.
\end{abstract}


\tableofcontents
\vspace{2em}

In this paper, Part VI of a series of eight papers on logarithmic
tensor category theory, we construct the appropriate natural
associativity isomorphisms between triple tensor product functors.
The sections, equations, theorems and so on are numbered globally in
the series of papers rather than within each paper, so that for
example equation (a.b) is the b-th labeled equation in Section a,
which is contained in the paper indicated as follows: In Part I
\cite{HLZ1}, which contains Sections 1 and 2, we give a detailed
overview of our theory, state our main results and introduce the basic
objects that we shall study in this work.  We include a brief
discussion of some of the recent applications of this theory, and also
a discussion of some recent literature.  In Part II \cite{HLZ2}, which
contains Section 3, we develop logarithmic formal calculus and study
logarithmic intertwining operators.  In Part III \cite{HLZ3}, which
contains Section 4, we introduce and study intertwining maps and
tensor product bifunctors.  In Part IV \cite{HLZ4}, which contains
Sections 5 and 6, we give constructions of the $P(z)$- and
$Q(z)$-tensor product bifunctors using what we call ``compatibility
conditions'' and certain other conditions.  In Part V \cite{HLZ5},
which contains Sections 7 and 8, we study products and iterates of
intertwining maps and of logarithmic intertwining operators and we
begin the development of our analytic approach.  The present paper,
Part VI, contains Sections 9 and 10.  In Part VII \cite{HLZ7}, which
contains Section 11, we give sufficient conditions for the existence
of the associativity isomorphisms.  In Part VIII \cite{HLZ8}, which
contains Section 12, we construct braided tensor category structure.

\paragraph{Acknowledgments}
The authors gratefully
acknowledge partial support {}from NSF grants DMS-0070800 and
DMS-0401302.  Y.-Z.~H. is also grateful for partial support {}from NSF
grant PHY-0901237 and for the hospitality of Institut des Hautes 
\'{E}tudes Scientifiques in the fall of 2007.

\renewcommand{\theequation}{\thesection.\arabic{equation}}
\renewcommand{\therema}{\thesection.\arabic{rema}}
\setcounter{section}{8}
\setcounter{equation}{0}
\setcounter{rema}{0}

\section{The expansion condition for intertwining maps and the
associativity of logarithmic intertwining operators}\label{extsec}

In the present section, we establish results, especially Theorem
\ref{9.7-1}, that form the crucial technical foundation of the
construction of the natural associativity isomorphisms in Section 10.
In Section 7 we have studied the conditions necessary for products and
iterates of certain intertwining maps to exist. In Section 8 we have
proved that products and iterates of such intertwining maps give
elements of the dual space of the vector space tensor product of the
three objects involved that satisfy the $P(z_{1},
z_{2})$-compatibility condition and the $P(z_{1}, z_{2})$-local
grading restriction condition. In Theorem
\ref{characterizationofbackslash}, we have given a characterization of
$W_{1}\hboxtr_{P(z)}W_{2}$ in terms of the $P(z)$-compatibility
condition and the $P(z)$-local grading restriction condition. It is
natural to also try to characterize the subspaces of the dual space
mentioned above by means of the $P(z_{1}, z_{2})$-compatibility
condition and the $P(z_{1}, z_{2})$-local grading restriction
condition, but this time, these two conditions are not enough.  We
have to find additional conditions such that the elements satisfying
all of the appropriate conditions can be obtained from both products
and iterates of suitable intertwining maps.  These additional
conditions constitute what we will call the ``expansion condition.''

Assuming that the convergence condition in Section 7 is satisfied, in
this section we study the conditions for the products of suitable
intertwining maps to be expressible as the iterates of some suitable
intertwining maps, and vice versa.  To do this, following the idea in
\cite{tensor4}, we first study certain properties satisfied by the
product, and separately, the iterate, of two intertwining maps.  Then
we obtain the deeper result that if a product satisfies the properties
naturally satisfied by iterates, it can indeed be expressed as an
iterate, and vice versa for an iterate.  Using these results, we prove
near the end of this section that the condition that products satisfy
the properties for iterates and the condition that iterates satisfy
the properties for products are equivalent to each other, and we
introduce the term ``expansion condition for intertwining maps'' for
either of these equivalent conditions (Definition
\ref{expansion-conditions}).  We show that under the assumption of the
convergence condition and the expansion condition, along with certain
``minor'' conditions, a product or iterate of two intertwining maps
can be expressed as an iterate or product, respectively.  Using the
correspondence between intertwining maps and logarithmic intertwining
operators, these results give the ``associativity of logarithmic
intertwining operators,'' which says that a product of logarithmic
intertwining operators can be expressed as an iterate of logarithmic
intertwining operators.  This associativity of logarithmic
intertwining operators is in fact a strong version of the
``logarithmic operator product expansion,'' which in turn is the
starting point of ``logarithmic conformal field theory,'' studied
extensively by physicists and mathematicians, as we have discussed in
the Introduction.  In this section, this logarithmic operator product
expansion is established as a mathematical theorem.  These results are
also viewed as saying that products or iterates of intertwining maps
or of logarithmic intertwining operators uniquely ``factor through''
suitable tensor product generalized modules.

The results in the present section are generalizations to the
logarithmic setting of the corresponding results in the finitely
reductive case obtained in \cite{tensor4}. See Remark \ref{tensor4}
for a comparison of the main results in the present section with the
corresponding results in \cite{tensor4}, including the corrections of
some minor mistakes. The most crucial results are Theorem \ref{9.7-1},
which essentially constructs the intermediate generalized module, and
Lemma \ref{intertwine-tau}, which essentially constructs the
corresponding intertwining map.  The main difficult aspect of the
proofs of these results is that we have to prove that certain iterated
series converge and that their sums are equal to the sums of iterated
series obtained by changing the order of summation.  These
difficulties, embedding in ``formal calculus,'' are overcome either by
proving the absolute convergence of the associated double or triple
series or by using the convergence of suitable Taylor expansions.  The
reasoning requires considerable use of analytic methods.

We again recall our Assumptions \ref{assum}, \ref{assum-c} and
\ref{assum-exp-set} concerning our category ${\cal C}$.

In this section $z_1$ and $z_2$ will be distinct nonzero complex
numbers, and
\[
z_0=z_1-z_2;
\]
we shall make various assumptions on these numbers below.

For objects $W_1$, $W_2$, $W_3$, $W_4$, $M_1$ and $M_2$ of ${\cal C}$,
let $I_1$, $I_2$, $I^1$ and $I^2$ be $P(z_1)$-, $P(z_2)$-, $P(z_2)$-
and $P(z_0)$-intertwining maps of types ${W_4 \choose W_1\, M_1}$,
${M_1 \choose W_2\, W_3}$, ${W_4 \choose M_2\, W_3}$ and ${M_2\choose
W_1\, W_2}$, respectively.  Then under the assumption of the
convergence condition for intertwining maps in ${\cal C}$ (recall
Proposition \ref{convergence} and Definition \ref{conv-conditions}),
when $|z_1|>|z_2|>|z_0|>0$, both the product $I_1\circ (1_{W_1}\otimes
I_2)$ and the iterate $I^1\circ (I^2\otimes 1_{W_3})$ exist and are
$P(z_1,z_2)$-intertwining maps, by Proposition
\ref{productanditerateareintwmaps}.  In this section we will consider
when such a product can be expressed as such an iterate and vice
versa.  To compare these two types of maps, we shall study some
conditions specific to each type.

Here and below, we let $W_1$, $W_2$ and $W_3$ be arbitrary generalized
$V$-modules.  (Later, these modules will be assumed to be objects of
${\cal C}$.)  We start with the following:

\begin{defi}\label{mudef}{\rm
Let
\[
\lambda \in (W_1\otimes W_2\otimes W_3)^*.
\]
For $w_{(1)}\in W_{1}$, we define the {\it evaluation of $\lambda$ at
$w_{(1)}$} to be the element
\[
\mu^{(1)}_{\lambda, w_{(1)}} \in (W_2\otimes W_3)^{*}
\]
given by
\[
\mu^{(1)}_{\lambda, w_{(1)}}(w_{(2)}\otimes w_{(3)})
=\lambda(w_{(1)}\otimes w_{(2)}\otimes w_{(3)})
\]
for $w_{(2)}\in W_2$ and $w_{(3)}\in W_3$.  For $w_{(3)}\in W_3$, we
define the {\it evaluation of $\lambda$ at $w_{(3)}$} to be the
element
\[
\mu^{(2)}_{\lambda, w_{(3)}} \in (W_1\otimes W_2)^{*}
\]
given by
\[
\mu^{(2)}_{\lambda, w_{(3)}}(w_{(1)}\otimes w_{(2)})
=\lambda(w_{(1)}\otimes w_{(2)}\otimes w_{(3)})
\]
for $w_{(1)}\in W_1$ and $w_{(2)}\in W_2$.
}
\end{defi}

\begin{rema}
{\rm Given $\lambda\in (W_1\otimes W_2\otimes W_3)^{*}$, $w_{(1)}\in
W_{1}$ and $w_{(3)}\in W_3$, it is natural to ask whether the
evaluations $\mu^{(1)}_{\lambda, w_{(1)}}\in (W_{2}\otimes W_{3})^{*}$
and $\mu^{(2)}_{\lambda, w_{(3)}}\in (W_{1}\otimes W_{2})^{*}$ of
$\lambda$ satisfy the $P(z)$-compatibility condition (recall
(\ref{cpb})) for some suitable nonzero complex numbers $z$, under
suitable conditions.  In fact, the formal computations underlying the
next lemma suggest that when $\lambda$ satisfies the $P(z_{1},
z_{2})$-compatibility condition (recall (\ref{zz:cpb})), these
evaluations of $\lambda$ ``almost'' satisfy the
$P(z_{2})$-compatibility condition and the $P(z_{0})$-compatibility
condition, respectively.  However, even when $\lambda$ does satisfy
the $P(z_{1}, z_{2})$-compatibility condition, in general these
evaluations of $\lambda$ do not even satisfy the $P(z_{2})$-lower
truncation condition (Part (a) of the $P(z_{2})$-compatibility
condition) or the $P(z_{0})$-lower truncation condition (Part (a) of
the $P(z_{0})$-compatibility condition).  In particular, when
$\lambda$ in (\ref{cpb}) is replaced by $\mu^{(1)}_{\lambda, w_{(1)}}$
and $z=z_{2}$, the right-hand side of (\ref{cpb}) might not exist in
the usual algebraic sense, and similarly, when $\lambda$ in
(\ref{cpb}) is replaced by $\mu^{(2)}_{\lambda, w_{(3)}}$ and
$z=z_{0}$, the right-hand side of (\ref{cpb}) might not exist
algebraically, and for this reason $\mu^{(1)}_{\lambda, w_{(1)}}$ and
$\mu^{(2)}_{\lambda, w_{(3)}}$ do not in general satisfy the
$P(z_{2})$-compatibility condition or the $P(z_{0})$-compatibility
condition.  But the next result, which generalizes Lemma 14.3 in
\cite{tensor4}, asserts that if $\lambda$ satisfies the $P(z_1,
z_2)$-compatibility condition, then in both cases, under the natural
assumptions on the complex numbers, the right-hand side of (\ref{cpb})
exists {\it analytically} and (\ref{cpb}) holds {\it analytically}, in
the sense of weak absolute convergence, as discussed in Remark
\ref{weakly-abs-conv}.}
\end{rema}

\begin{lemma}\label{mulemma}
Assume that $\lambda\in (W_1\otimes W_2\otimes W_3)^{*}$ satisfies the
$P(z_1, z_2)$-compatibility condition (recall (\ref{zz:cpb})).  If
$|z_2|>|z_0|$ $(>0)$, then for any $v\in V$ and $w_{(1)}\in W_1$,
$w_{(2)}\in W_2$ and $w_{(3)}\in W_3$, 
\begin{eqnarray}\label{rlm4}
\lefteqn{\biggl(Y'_{P(z_0)}(v, x) \mu^{(2)}_{\lambda, w_{(3)}}\biggr)
(w_{(1)} \otimes w_{(2)})}\nn
&&={\rm
Res}_{x_0^{-1}}\Biggl(\tau_{P(z_1,z_2)}\Biggl(
\dlti{x}{x_0^{-1}}{-z_2}\cdot\nn
&&\qquad\qquad\cdot
Y_t((-x_0^2)^{L(0)}e^{-x_0L(1)}e^{xL(1)}(-x^{-2})^{L(0)}v,x_0)
\Biggr)\lambda\Biggr)(w_{(1)}\otimes w_{(2)}\otimes w_{(3)})\nn
&&\quad -{\rm
Res}_{x_0^{-1}}\dlti{x}{-z_2}{+x_0^{-1}}\cdot\nn
&&\qquad\qquad \cdot\lambda (w_{(1)}\otimes w_{(2)}\otimes
Y_3(e^{xL(1)}(-x^{-2})^{L(0)}v,x_0^{-1})w_{(3)}),
\end{eqnarray}
the coefficients of the
monomials in $x$ and $x_{1}$ in
\[
x^{-1}_1 \delta\bigg(\frac{x^{-1}-z_0}{x_1}\bigg)
\biggl(Y'_{P(z_0)}(v, x) \mu^{(2)}_{\lambda, w_{(3)}}\biggr)
(w_{(1)} \otimes w_{(2)})
\]
are absolutely convergent, and we have
\begin{eqnarray}\label{mu12}
\lefteqn{\biggl(\tau_{P(z_0)}\biggl( x^{-1}_1
\delta\bigg(\frac{x^{-1}-z_0}{x_1}\bigg) Y_{t}(v, x)\biggr)
\mu^{(2)}_{\lambda, w_{(3)}}\biggr)(w_{(1)} \otimes w_{(2)})}\nno\\
&&=x^{-1}_1 \delta\bigg(\frac{x^{-1}-z_0}{x_1}\bigg)
\biggl(Y'_{P(z_0)}(v, x) \mu^{(2)}_{\lambda, w_{(3)}}\biggr)
(w_{(1)} \otimes w_{(2)}).
\end{eqnarray}
Analogously, if $|z_1|>|z_2|$ $(>0)$, then for any $v\in V$ and any
$w_{(j)}\in W_j$, 
\begin{eqnarray}\label{rlm7}
\lefteqn{\biggl(Y'_{P(z_2)}(v, x) \mu^{(1)}_{\lambda, w_{(1)}}\biggr)
(w_{(2)} \otimes w_{(3)})}\nn
&&=(Y'_{P(z_1,z_2)}(v,x_0)
\lambda )(w_{(1)}\otimes w_{(2)}\otimes w_{(3)})\nn
&&\quad -\res_{x_1}\dlti{x_0}{z_1}{+x_1}
\lambda (Y_1((-x_0^{-2})^{L(0)}e^{-x_0^{-1}L(1)}v,x_1)
w_{(1)}\otimes w_{(2)}\otimes w_{(3)}),
\end{eqnarray}
the coefficients of the monomials in $x$ and $x_{1}$
in
\[
x^{-1}_1 \delta\bigg(\frac{x^{-1}-z_2}{x_1}\bigg)\biggl(Y'_{P(z_2)}(v, x)
\mu^{(1)}_{\lambda, w_{(1)}}\biggr)(w_{(2)} \otimes w_{(3)})
\]
are absolutely convergent, and we have
\begin{eqnarray}\label{mu23}
\lefteqn{\biggl(\tau_{P(z_2)}\biggl( x^{-1}_1
\delta\bigg(\frac{x^{-1}-z_2}{x_1}\bigg) Y_{t}(v, x)\biggr)
\mu^{(1)}_{\lambda, w_{(1)}}\biggr)(w_{(2)} \otimes w_{(3)})}\nno\\
&&=x^{-1}_1 \delta\bigg(\frac{x^{-1}-z_2}{x_1}\bigg)
\biggl(Y'_{P(z_2)}(v, x) \mu^{(1)}_{\lambda, w_{(1)}}\biggr)
(w_{(2)} \otimes w_{(3)}).
\end{eqnarray}
\end{lemma}
\pf First, for our distinct nonzero complex numbers $z_1$ and $z_2$
with $z_0=z_1-z_2$, by definition of the action $\tau _{P(z_1,z_2)}$
(\ref{tauzzgf}) we have
\begin{eqnarray}\label{lm:1}
\lefteqn{\dlti{x_0}{z_1}{+x_1}\dlt{x_2}{z_0}{+x_1} \lambda
(Y_1((-x_0^{-2})^{L(0)}e^{-x_0^{-1}L(1)}v,x_1)w_{(1)}\otimes
w_{(2)}\otimes w_{(3)})}\nno\\
&&+ \dlti{x_0}{z_2}{+x_2}\dlt{x_1}{-z_0}{+x_2}\lambda (w_{(1)}\otimes
Y_2((-x_0^{-2})^{L(0)}e^{-x_0^{-1}L(1)}v,x_2)w_{(2)}\otimes
w_{(3)})\nno\\
&&= \left(\tau _{P(z_1,z_2)}\left(\dlt{x_1}{x_0^{-1}}{-z_1}\dlt{x_2}{x_0^{-1}}
{-z_2}Y_{t}(v,x_0)\right)\lambda \right)(w_{(1)}\otimes w_{(2)}\otimes
w_{(3)})\nno\\
&&\quad-\dlt{x_1}{-z_1}{+x^{-1}_0}\dlt{x_2}{-z_2}{+x^{-1}_0}\lambda
(w_{(1)}\otimes w_{(2)}\otimes Y_3^o(v,x_0)w_{(3)})
\end{eqnarray}
for any $v\in V$, $w_{(1)}\in W_1$, $w_{(2)}\in W_2$ and $w_{(3)}\in W_3$.
Replacing $v$ by
\[
(-x_0^2)^{L(0)}e^{-x_0L(1)}e^{x_2^{-1}L(1)}(-x_2^2)^{L(0)}v,
\]
using formula (5.3.1) in \cite{FHL}, and then
taking $\res_{x_0^{-1}}$ we get:
\begin{eqnarray}\label{lm4}
\lefteqn{\dlt{x_2}{z_0}{+x_1}\lambda
(Y_1(e^{x_2^{-1}L(1)}(-x_2^2)^{L(0)}v,x_1)w_{(1)}\otimes
w_{(2)}\otimes w_{(3)})}\nno\\
&& +\dlt{x_1}{-z_0}{+x_2}\lambda (w_{(1)}\otimes
Y_2(e^{x_2^{-1}L(1)}(-x_2^2)^{L(0)}v,x_2)w_{(2)}\otimes
w_{(3)})\nno\\
&&= {\rm
Res}_{x_0^{-1}}\Bigg(\tau_{P(z_1,z_2)}\Bigg(\dlt{x_1}{x_0^{-1}}{-z_1}
\dlt{x_2}{x_0^{-1}}{-z_2}\cdot\nno\\
&&\qquad\cdot
Y_{t}((-x_0^2)^{L(0)}e^{-x_0L(1)}e^{x_2^{-1}L(1)}(-x_2^2)^{L(0)}v,x_0)
\Bigg)\lambda
\Bigg)(w_{(1)}\otimes w_{(2)}\otimes w_{(3)}))\nno\\
&&\quad -{\rm
Res}_{x_0^{-1}}\dlt{x_1}{-z_1}{+x^{-1}_0}\dlt{x_2}{-z_2}{+x^{-1}_0}\cdot\nno\\
&&\qquad \cdot\lambda (w_{(1)}\otimes w_{(2)}\otimes
Y_3(e^{x_2^{-1}L(1)}(-x_2^2)^{L(0)}v,x_0^{-1})w_{(3)}).
\end{eqnarray}

By (\ref{taudef}), the left-hand side of (\ref{lm4}) is equal to
\begin{eqnarray}\label{lefthandside}
\left(\tau_{P(z_0)}\left(\dlt{x_1}{x_2}{-z_0}Y_t(v,x_2^{-1})\right)\mu^{(2)}_{\lambda,
w_{(3)}}\right)(w_{(1)}\otimes w_{(2)}).
\end{eqnarray}
Taking $\res_{x_1}$ gives
\begin{eqnarray}\label{resoflefthandside}
(\tau_{P(z_0)}(Y_t(v,x_2^{-1}))\mu^{(2)}_{\lambda,
w_{(3)}})(w_{(1)}\otimes
w_{(2)})=(Y'_{P(z_0)}(v,x_2^{-1})\mu^{(2)}_{\lambda,
w_{(3)}})(w_{(1)}\otimes w_{(2)}).
\end{eqnarray}

By the $P(z_1, z_2)$-compatibility condition and formula
(\ref{consequenceofPz1z2compatformula}) in Remark
\ref{consequenceofPz1z2compat}, the first term on the right-hand side
of (\ref{lm4}) equals
\begin{eqnarray}\label{RHSexpression}
\lefteqn{\dlt{x_1}{x_2}{-z_0}{\rm Res}_{x_0^{-1}}\Bigg(\tau_{P(z_1,z_2)}\Bigg(
\dlt{x_2}{x_0^{-1}}{-z_2}\cdot}\nno\\
&&\qquad\cdot
Y_{t}((-x_0^2)^{L(0)}e^{-x_0L(1)}e^{x_2^{-1}L(1)}(-x_2^2)^{L(0)}v,x_0)
\Bigg)\lambda \Bigg)(w_{(1)}\otimes w_{(2)}\otimes w_{(3)})).
\end{eqnarray}

Now suppose that $|z_2|>|z_0|$ $(>0)$.  Then formula (\ref{l4})
holds. {}From this and (\ref{RHSexpression}), the right-hand side of
(\ref{lm4}) becomes
\begin{eqnarray*}
\lefteqn{\dlt{x_1}{x_2}{-z_0}{\rm
Res}_{x_0^{-1}}\Biggl(\tau_{P(z_1,z_2)}\Biggl(
\dlt{x_2}{x_0^{-1}}{-z_2}\cdot}\nn
&&\qquad\cdot
Y_t((-x_0^2)^{L(0)}e^{-x_0L(1)}e^{x_2^{-1}L(1)}(-x_2^2)^{L(0)}v,x_0)
\Biggr)\lambda\Biggr)(w_{(1)}\otimes w_{(2)}\otimes w_{(3)})\nn
&&-\dlt{x_1}{x_2}{-z_0}{\rm
Res}_{x_0^{-1}}\dlt{x_2}{-z_2}{+x_0^{-1}}\cdot\nn
&&\qquad \cdot\lambda (w_{(1)}\otimes w_{(2)}\otimes
Y_3(e^{x_2^{-1}L(1)}(-x_2^2)^{L(0)}v,x_0^{-1})w_{(3)}).
\end{eqnarray*}
By taking $\res_{x_1}$, we erase the two factors
$\displaystyle\dlt{x_1}{x_2}{-z_0}$, leaving an expression 
which is exactly the right-hand side of 
(\ref{rlm4}) with $x$ replaced by $x_2^{-1}$ (thus proving
(\ref{rlm4})) and,
while not lower truncated in $x_2^{-1}$, can still be multiplied by
$\displaystyle\dlt{x_1}{x_2}{-z_0}$ (when $|z_2|>|z_0|>0$), in the
sense of absolute convergence, yielding this expression again.  That
is, let $X$ be either side of (\ref{lm4}).  Then
\[
X=\dlt{x_1}{x_2}{-z_0}\res_{x_1}X,
\]
in this sense of convergence.  Applying this to (\ref{lefthandside})
and (\ref{resoflefthandside}) gives
\begin{eqnarray}\label{lmu12}
\lefteqn{\left(\tau_{P(z_0)}\left(\dlt{x_1}{x_2}{-z_0}Y_t(v,x_2^{-1})\right)
\mu^{(2)}_{\lambda,w_{(3)}}\right)(w_{(1)}\otimes w_{(2)})}\nno\\
&&\qquad
=\dlt{x_1}{x_2}{-z_0}\biggl(Y'_{P(z_0)}(v,x_2^{-1})\mu^{(2)}_{\lambda,
w_{(3)}}\biggr)(w_{(1)}\otimes w_{(2)}),
\end{eqnarray}
proving (\ref{mu12}) (with $x$ in (\ref{mu12}) replaced by
$x_2^{-1}$).

Analogously, for our distinct nonzero complex numbers $z_1$ and $z_2$
with $z_0=z_1-z_2$, we can also write the definition of
$\tau_{P(z_1,z_2)}$ as
\begin{eqnarray}\label{lm:2}
\lefteqn{\dlti{x_0}{z_2}{+x_2}\dlt{x_1}{-z_0}{+x_2} \lambda (w_{(1)}\otimes
Y_2((-x_0^{-2})^{L(0)}e^{-x_0^{-1}L(1)}v,x_2)w_{(2)}\otimes
w_{(3)})}\nno\\
&&+\dlt{x_1}{-z_1}{+x^{-1}_0}\dlt{x_2}{-z_2}{+x^{-1}_0}\lambda
(w_{(1)}\otimes w_{(2)}\otimes Y_3^o(v,x_0)w_{(3)})\nno\\
&&=\left(\tau_{P(z_1,z_2)}\left(\dlt{x_1}{x_0^{-1}}{-z_1}\dlt{x_2}{x_0^{-1}}
{-z_2}Y_{t}(v,x_0)\right)\lambda \right)(w_{(1)}\otimes w_{(2)}\otimes w_{(3)})\nno\\
&  &\quad- \dlti{x_0}{z_1}{+x_1}\dlt{x_2}{z_0}{+x_1}\lambda
(Y_1((-x_0^{-2})^{L(0)}e^{-x_0^{-1}L(1)}v,x_1)w_{(1)}\otimes
w_{(2)}\otimes w_{(3)})\nno\\
\end{eqnarray}
for $v\in V$, $w_{(1)}\in W_1$, $w_{(2)}\in W_2$ and $w_{(3)}\in W_3$.
Taking $\res_{x_1}$ we get
\begin{eqnarray}\label{lm7}
\lefteqn{\dlti{x_0}{z_2}{+x_2}\lambda (w_{(1)}\otimes
Y_2((-x_0^{-2})^{L(0)}e^{-x_0^{-1}L(1)}v,x_2)w_{(2)}\otimes
w_{(3)})}\nno\\
&  &+ \dlt{x_2}{-z_2}{+x^{-1}_0}\lambda (w_{(1)}\otimes
w_{(2)}\otimes Y^o_3(v,x_0)w_{(3)})\nno\\
& &=\left(\tau
_{P(z_1,z_2)}\left(\dlt{x_2}{x_0^{-1}}{-z_2}Y_{t}(v,x_0)\right)\lambda
\right)(w_{(1)}\otimes w_{(2)}\otimes w_{(3)})\nno\\
& &\quad-\res_{x_1} \dlti{x_0}{z_1}{+x_1}\dlt{x_2}{z_0}{+x_1}\cdot\nno\\
&&\qquad\qquad\cdot\lambda
(Y_1((-x_0^{-2})^{L(0)}e^{-x_0^{-1}L(1)}v,x_1)w_{(1)}\otimes
w_{(2)}\otimes w_{(3)}).
\end{eqnarray}

By the definition of $\tau_{P(z_2)}$ (\ref{taudef}) and formula
(5.3.1) in \cite{FHL}, the left-hand side of (\ref{lm7}) is equal to
\[
\left(\tau_{P(z_2)}\left(\dlt{x_2}{x_0^{-1}}{-z_2}Y_t(v,x_0)\right)\mu^{(1)}_{\lambda,
w_{(1)}}\right)(w_{(2)}\otimes w_{(3)}),
\]
and taking ${\rm Res}_{x_2}$ gives
\[
(\tau_{P(z_2)}(Y_t(v,x_0))\mu^{(1)}_{\lambda,
w_{(1)}})(w_{(2)}\otimes w_{(3)})\\
=(Y'_{P(z_2)}(v,x_0)\mu^{(1)}_{\lambda, w_{(1)}})(w_{(2)}\otimes
w_{(3)}).
\]

Now suppose that $|z_1|>|z_2|>0$.  Then (\ref{l2-1}) holds, and by
(\ref{resofconsequence}), which follows from the
$P(z_1,z_2)$-compatibility condition, the right-hand side of
(\ref{lm7}) becomes
\begin{eqnarray*}
\lefteqn{\dlt{x_2}{x_0^{-1}}{-z_2}(Y'_{P(z_1,z_2)}(v,x_0)
\lambda )(w_{(1)}\otimes w_{(2)}\otimes w_{(3)})}\\
&&-\dlt{x_2}{x_0^{-1}}{-z_2}\res_{x_1}\dlti{x_0}{z_1}{+x_1}\cdot\\
&&\qquad\cdot\lambda (Y_1((-x_0^{-2})^{L(0)}e^{-x_0^{-1}L(1)}v,x_1)
w_{(1)}\otimes w_{(2)}\otimes w_{(3)}).
\end{eqnarray*}
Just as in the proof of (\ref{mu12}), we take $\res_{x_2}$ to obtain 
the right-hand side of (\ref{rlm7}) (thus proving (\ref{rlm7})) 
and then
multiply by $\displaystyle\dlt{x_2}{x_0^{-1}}{-z_2}$, yielding the
same expression, and we obtain
\begin{eqnarray}\label{lmu23}
\lefteqn{\biggl(\tau_{P(z_2)}\biggl(\dlt{x_2}{x_0^{-1}}{-z_2} Y_{t}(v,x_0)\biggr)
\mu^{(1)}_{\lambda, w_{(1)}}\biggr)
(w_{(2)}\otimes w_{(3)})}\nno\\
&&\qquad
=\dlt{x_2}{x_0^{-1}}{-z_2}\biggl(Y'_{P(z_2)}(v, x_0)
\mu^{(1)}_{\lambda,w_{(1)}}\biggr)(w_{(2)}\otimes w_{(3)}),
\end{eqnarray}
proving (\ref{mu23}).
\epf

\begin{rema}{\rm
As we discussed above, Lemma \ref{mulemma} says that under the
appropriate conditions, $\mu^{(2)}_{\lambda, w_{(3)}}$ and
$\mu^{(1)}_{\lambda, w_{(1)}}$ satisfy natural analytic analogues of
the $P(z_0)$-compatibility condition and the $P(z_2)$-compatibility
condition, respectively.  Note that from the proof of (\ref{mu12}),
$(Y'_{P(z_0)}(v,x_2^{-1})\mu^{(2)}_{\lambda, w_{(3)}})(w_{(1)}\otimes
w_{(2)})$ ``behaves qualitatively'' like
$\displaystyle\delta\bigg(\frac{z_2}{-x_2}\bigg)$ and so (\ref{lmu12})
``behaves qualitatively'' like
\[
\dlt{x_1}{x_2}{-z_0}\delta\bigg(\frac{z_2}{-x_2}\bigg),
\]
suggesting the expected convergence when $|z_2|>|z_0|$.  Analogously,
from the proof of (\ref{mu23}), $Y'_{P(z_2)}(v,x_0)\mu^{(1)}_{\lambda,
w_{(1)}}(w_{(2)}\otimes w_{(3)})$ ``behaves qualitatively'' like
$\displaystyle\delta\bigg(\frac{z_1}{x_0^{-1}}\bigg)$ and so
(\ref{lmu23}) ``behaves qualitatively'' like
\[
\dlt{x_2}{x_0^{-1}}{-z_2}\delta\bigg(\frac{z_1}{x_0^{-1}}\bigg),
\]
again suggesting the expected convergence, this time when
$|z_1|>|z_2|$.  (By the $P(z_1 ,z_2)$-lower truncation condition,
${\rm Res}_{x_1}$ of (\ref{lefthandside}) is upper-truncated in $x_2$,
independently of $w_{(1)}$, $w_{(2)}$ and $w_{(3)}$, and ${\rm
Res}_{x_2}$ of the first term on the right-hand side of (\ref{lm7}),
$(Y'_{P(z_1,z_2)}(v,x_0)\lambda)(w_{(1)}\otimes w_{(2)}\otimes
w_{(3)})$, is lower truncated in $x_0$, independently of the
$w_{(j)}$.)}
\end{rema}

\begin{rema}\label{rmk-9.5}
{\rm Given
\[
\lambda\in (W_1\otimes W_2\otimes W_3)^{*},
\]
$w_{(1)}\in
W_{1}$ and $w_{(3)}\in W_3$, it is also natural to ask whether, under
suitable conditions, the evaluations $\mu^{(1)}_{\lambda, w_{(1)}}
\in (W_{2}\otimes W_{3})^{*}$
and $\mu^{(2)}_{\lambda, w_{(3)}}\in (W_{1}\otimes W_{2})^{*}$ 
of $\lambda$ satisfy the
$P(z_{2})$-local grading restriction condition and the $P(z_{0})$-local
grading restriction condition, respectively.  In general, even for
$\lambda$ obtained {}from a product or an iterate of intertwining
maps, these conditions are not satisfied by $\mu^{(1)}_{\lambda,
w_{(1)}}$ or $\mu^{(2)}_{\lambda, w_{(3)}}$, but as we will see below,
for such $\lambda$, the evaluations $\mu^{(1)}_{\lambda, w_{(1)}}$ and
$\mu^{(2)}_{\lambda, w_{(3)}}$ satisfy certain analytic analogues of
these conditions.  These analogues motivate the next four
important conditions on $\lambda\in (W_1\otimes W_2\otimes W_3)^*$.
On the spaces (\ref{W1W2_[C]^Atilde}) and (\ref{W1W2_(C)^Atilde})
for $W_{1}\otimes W_{2}$ and for $W_{2}\otimes W_{3}$, the considerations 
of Remark \ref{set:L(0)s} concerning the semisimple part of $L(0)_{s}$
of $L(0)$ hold, and in particular, on these spaces,
\[
[L'_{P(z)}(0)-L'_{P(z)}(0)_{s}, L'_{P(z)}(0)]=0
\]
and so 
\[
[L'_{P(z)}(0)-L'_{P(z)}(0)_{s}, L'_{P(z)}(0)_{s}]=0.
\]
Hence
\[
e^{yL'_{P(z)}(0)}=e^{yL'_{P(z)}(0)_{s}}e^{y(L'_{P(z)}(0)-L'_{P(z)}(0)_{s})},
\]
where $y$ is a formal variable,
and 
\[
e^{z'L'_{P(z)}(0)}=e^{z'L'_{P(z)}(0)_{s}}e^{z'(L'_{P(z)}(0)-L'_{P(z)}(0)_{s})}
\]
for $z'\in \C$.  (A complex number denoted $z'$ or $-z'$ will 
play this role in the considerations below.) 
Thus $e^{z'L'_{P(z)}(0)}$ maps any $P(z)$-generalized weight vector
$\nu$ in $(W_{2}\otimes W_{3})^{*}$ or $(W_{1}\otimes W_{2})^{*}$
to a $P(z)$-generalized weight vector of the same 
generalized weight. 
An element $\lambda\in (W_1\otimes W_2\otimes W_3)^*$ satisfying one
of the conditions below means essentially that either $\mu^{(1)}_{\lambda,
w_{(1)}}\in (W_2\otimes W_3)^{*}$ or $\mu^{(2)}_{\lambda, w_{(3)}}\in
(W_1\otimes W_2)^{*}$ is the value at $z'=0$ of 
the sum of a series, weakly absolutely convergent in the
sense of Remark \ref{weakly-abs-conv} for $z'$ in a neighborhood of 
$z'=0$, rather than just a finite sum,
of the images under the map
$e^{z'L'_{P(z)}(0)}$ 
of $P(z)$-generalized weight vectors or of ordinary weight vectors
satisfying the $P(z)$-local grading restriction condition or the
$L(0)$-semisimple $P(z)$-local grading restriction condition in
Section 5, and that all of the summands lie in the same subspace whose
grading is restricted.  (For the four stronger conditions, the
restriction on the grading is analogous to (\ref{set:dmltc-1}).)
We shall typically use these conditions for
$z=z_2$ when we consider $\mu^{(1)}_{\lambda, w_{(1)}}$ and for
$z=z_{0}$ when we consider $\mu^{(2)}_{\lambda, w_{(3)}}$.  }
\end{rema}

Recall the
spaces (\ref{W1W2_[C]^Atilde}) and (\ref{W1W2_(C)^Atilde}) and recall that 
\[
e^{yL'_{P(z)}(0)}=e^{yL'_{P(z)}(0)_{s}}e^{y(L'_{P(z)}(0)-L'_{P(z)}(0)_{s})}
\]
and
\[
e^{z'L'_{P(z)}(0)}=e^{z'L'_{P(z)}(0)_{s}}e^{z'(L'_{P(z)}(0)-L'_{P(z)}(0)_{s})}
\]
on these spaces; on the spaces (\ref{W1W2_(C)^Atilde}), 
\[
L'_{P(z)}(0)_{s}=L'_{P(z)}(0). 
\]
Consider formal series 
$\sum_{n\in \C}\lambda_{n}^{(1)}$ and
$\sum_{n\in \C}\lambda_{n}^{(2)}$ with 
\[
\lambda_{n}^{(1)}\in \coprod_{\beta\in \tilde{A}}
((W_{2}\otimes W_{3})^{*})_{[n]}^{(\beta)}
\]
and 
\[
\lambda_{n}^{(2)}
\in\coprod_{\beta\in \tilde{A}}
((W_{1}\otimes W_{2})^{*})_{[n]}^{(\beta)}
\]
for $n\in \C$. Then there exist $K^{(1)}_{n}, K^{(2)}_{n}\in \N$ for $n\in \C$
such that 
\begin{eqnarray}\label{e-y-1}
\sum_{n\in \C}e^{yL'_{P(z)}}\lambda_{n}^{(1)}
&=&\sum_{n\in \C}e^{y(L'_{P(z)}(0)-L'_{P(z)}(0)_{s})}
e^{yL'_{P(z)}(0)_{s}}\lambda^{(1)}_{n}\nn
&=&\sum_{n\in \C}e^{ny}\left(\sum_{i=0}^{K^{(1)}_{n}}\frac{y^{i}}{i!}
(L'_{P(z)}(0)-n)^{i}\lambda^{(1)}_{n}
\right)
\end{eqnarray}
and
\begin{eqnarray}\label{e-y-2}
\sum_{n\in \C}e^{yL'_{P(z)}}\lambda^{(2)}_{n}
&=&\sum_{n\in \C}e^{y(L'_{P(z)}(0)-L'_{P(z)}(0)_{s})}
e^{yL'_{P(z)}(0)_{s}}\lambda^{(2)}_{n}
\nn
&=&\sum_{n\in \C}e^{ny}\left(\sum_{i=0}^{K^{(2)}_{n}}
\frac{y^{i}}{i!}
(L'_{P(z)}(0)-n)^{i}\lambda^{(2)}_{n}
\right).
\end{eqnarray}

\begin{rema}\label{y=>z'}
{\rm The formulas (\ref{e-y-1}) and (\ref{e-y-2}) also hold
with $y$ replaced by $z'\in \C$.}
\end{rema}

We shall be restricting our attention to $n \in \R$; we shall use the
subspace
\begin{equation}\label{W1W2_[R]^Atilde}
((W_1\otimes W_2)^*)_{[\R]}^{( \tilde A )}=
\coprod_{n\in
\R}\coprod_{\beta\in \tilde{A}}((W_1\otimes W_2)^*)_{[n]}^{(\beta)}
\end{equation}
of $((W_1\otimes W_2)^*)_{[\C]}^{( \tilde A )}$ and the
correspondingly defined subspace
\begin{equation}\label{W1W2_(R)^Atilde}
((W_1\otimes W_2)^*)_{(\R)}^{( \tilde A )}=
\coprod_{n\in
\R}\coprod_{\beta\in \tilde{A}}((W_1\otimes W_2)^*)_{(n)}^{(\beta)}
\end{equation}
of $((W_1\otimes W_2)^*)_{(\C)}^{( \tilde A )}$ (recall
(\ref{W1W2_[C]^Atilde}) and (\ref{W1W2_(C)^Atilde})), and similarly
for $W_2\otimes W_3$.

We next introduce the four conditions on
\[
\lambda\in (W_{1}\otimes W_{2}\otimes W_{3})^{*}.
\]
These conditions are motivated by certain properties of such elements
obtained from products or iterates of intertwining maps (see
Proposition \ref{9.7} below).  Essentially there are really only two
conditions, with the designations $P^{(1)}(z)$ and $P^{(2)}(z)$, but
each of them has a non-semisimple version and a semisimple version.  
These four conditions will enter into the formulation of the 
``expansion condition.''

\begin{description}
\item{\bf The $P^{(1)}(z)$-local grading restriction condition}

(a) The {\em $P^{(1)}(z)$-grading condition}: For any $w_{(1)}\in
W_1$, there exists a formal series $\sum_{n\in \R}\lambda^{(1)}_{n}$
with
\[
\lambda^{(1)}_{n}\in \coprod_{\beta\in \tilde{A}}
((W_{2}\otimes W_{3})^{*})_{[n]}^{(\beta)}
\]
for $n\in \R$, an open neighborhood of $z'=0$, and $N \in \N$
such that for $w_{(2)}\in W_{2}$ and $w_{(3)}\in W_{3}$,
the series 
\[
\sum_{n\in \R}(e^{z'L'_{P(z)}(0)}\lambda^{(1)}_{n})(w_{(2)}\otimes w_{(3)})
\]
(recall (\ref{e-y-1}) and Remark \ref{y=>z'}; here we evaluate 
at $w_{(2)}\otimes w_{(3)}$) has the following properties:
\begin{enumerate}

\item[(i)] The series can be written as the iterated series
\[
\sum_{n\in \R}e^{nz'}\left(\left(\sum_{i=0}^{N}\frac{(z')^{i}}{i!}
(L'_{P(z)}(0)-n)^{i}\lambda^{(1)}_{n}
\right)(w_{(2)}\otimes 
w_{(3)})\right)
\]
(recall from Proposition \ref{real-exp-set} that
$\R\times\{0,\dots,N\}$ is a unique expansion set).

\item[(ii)] It is absolutely convergent for $z'\in \C$ in the open
neighborhood of $z'=0$ above.

\item[(iii)] It is absolutely convergent to $\mu^{(1)}_{\lambda,
w_{(1)}}(w_{(2)}\otimes w_{(3)})$ when $z'=0$:
\[
\sum_{n\in \R}\lambda^{(1)}_{n}(w_{(2)}\otimes w_{(3)})
=\mu^{(1)}_{\lambda, w_{(1)}}(w_{(2)}\otimes w_{(3)})
=\lambda(w_{(1)}\otimes w_{(2)}\otimes w_{(3)}).
\]

\end{enumerate}

(b) For any $w_{(1)}\in W_1$, let $W^{(1)}_{\lambda, w_{(1)}}$ be the
smallest doubly graded (or equivalently, $\tilde A$-graded; recall
Remark \ref{singleanddoublegraded}) subspace of $((W_2\otimes
W_3)^{*})_{[\R]}^{(\tilde{A})}$, or equivalently, of $((W_2\otimes
W_3)^{*})_{[\C]}^{(\tilde{A})}$, containing all the terms
$\lambda^{(1)}_{n}$ in the formal series in (a) and stable under the
component operators $\tau_{P(z)}(v\otimes t^{m})$ of the operators
$Y'_{P(z)}(v, x)$ for $v\in V$, $m\in {\mathbb Z}$, and under the
operators $L'_{P(z)}(-1)$, $L'_{P(z)}(0)$ and $L'_{P(z)}(1)$.  (In
view of Remark \ref{stableundercomponentops}, $W^{(1)}_{\lambda,
w_{(1)}}$ indeed exists, just as in the case of the $P(z)$-local
grading restriction condition.)  Then $W^{(1)}_{\lambda, w_{(1)}}$ has
the properties
\begin{eqnarray*}
&\dim(W^{(1)}_{\lambda, w_{(1)}})^{(\beta)}_{[n]}<\infty,&\\
&(W^{(1)}_{\lambda, w_{(1)}})^{(\beta)}_{[n+k]}=0\;\;\mbox{ for
}\;k\in {\mathbb Z} \;\mbox{ sufficiently negative}&
\end{eqnarray*}
for any $n\in \R$ and $\beta\in \tilde A$, where the
subscripts denote the $\R$-grading by
$L'_{P(z)}(0)$-(generalized) eigenvalues and the superscripts denote
the $\tilde A$-grading.

\item{\bf The $L(0)$-semisimple 
$P^{(1)}(z)$-local grading restriction condition}

(a) The {\em $L(0)$-semisimple $P^{(1)}(z)$-grading condition}: For
any $w_{(1)}\in W_1$, there exists  a
formal series $\sum_{n\in \R}\lambda^{(1)}_{n}$ 
with
\[
\lambda^{(1)}_{n}\in \coprod_{\beta\in \tilde{A}}
((W_{2}\otimes W_{3})^{*})_{(n)}^{(\beta)}
\]
for $n\in \R$ and an open neighborhood of $z'=0$
such that for $w_{(2)}\in W_{2}$ and $w_{(3)}\in W_{3}$,
the series 
\[
\sum_{n\in \R}(e^{z'L'_{P(z)}(0)}\lambda^{(1)}_{n})(w_{(2)}\otimes w_{(3)})
\]
has the following properties:
\begin{enumerate}

\item[(i)] It can be written as
\[
\sum_{n\in \R}e^{nz'}\lambda^{(1)}_{n}(w_{(2)}\otimes 
w_{(3)})
\]
(recall from Proposition \ref{real-exp-set} that $\R\times\{0\}$ is a
unique expansion set).

\item[(ii)] It is absolutely convergent for $z'\in \C$ in the
neighborhood of $z'=0$ above.

\item[(iii)] It is absolutely convergent to $\mu^{(1)}_{\lambda,
w_{(1)}}(w_{(2)}\otimes w_{(3)})$ when $z'=0$:
\[
\sum_{n\in \R}\lambda^{(1)}_{n}(w_{(2)}\otimes w_{(3)})
=\mu^{(1)}_{\lambda, w_{(1)}}(w_{(2)}\otimes w_{(3)})
=\lambda(w_{(1)}\otimes w_{(2)}\otimes w_{(3)}).
\]

\end{enumerate}
(Note that such an element $\lambda$ also satisfies the
$P^{(1)}(z)$-grading condition above with the same elements
$\lambda^{(1)}_{n}$.)

(b) For any $w_{(1)}\in W_1$, consider the space $W^{(1)}_{\lambda,
w_{(1)}}$ as above, which in this case is in fact the smallest doubly
graded (or equivalently, $\tilde A$-graded) subspace of $((W_2\otimes
W_3)^{*})_{(\R)}^{(\tilde{A})}$ (or of $((W_2\otimes
W_3)^{*})_{(\C)}^{(\tilde{A})}$) containing all the terms
$\lambda^{(1)}_{n}$ in the formal series in (a) and stable under the
component operators $\tau_{P(z)}(v\otimes t^{m})$ of the operators
$Y'_{P(z)}(v, x)$ for $v\in V$, $m\in {\mathbb Z}$, and under the
operators $L'_{P(z)}(-1)$, $L'_{P(z)}(0)$ and $L'_{P(z)}(1)$. Then
$W^{(1)}_{\lambda, w_{(1)}}$ has the properties
\begin{eqnarray*}
&\dim(W^{(1)}_{\lambda, w_{(1)}})^{(\beta)}_{(n)}<\infty,&\\
&(W^{(1)}_{\lambda, w_{(1)}})^{(\beta)}_{(n+k)}=0\;\;\mbox{ for
}\;k\in {\mathbb Z} \;\mbox{ sufficiently negative}&
\end{eqnarray*}
for any $n\in \R$ and $\beta\in \tilde A$, where the
subscripts denote the $\R$-grading by
$L'_{P(z)}(0)$-eigenvalues and the superscripts denote the $\tilde
A$-grading.

\item{\bf The $P^{(2)}(z)$-local grading restriction condition}

(a) The {\em $P^{(2)}(z)$-grading condition}: For any $w_{(3)}\in
W_3$,  there exists  a formal series $\sum_{n\in \R}\lambda^{(2)}_{n}$ 
with
\[
\lambda^{(2)}_{n}\in \coprod_{\beta\in \tilde{A}}
((W_{1}\otimes W_{2})^{*})_{[n]}^{(\beta)}
\]
for $n\in \R$, an open neighborhood of $z'=0$, and $N \in \N$
such that for $w_{(1)}\in W_{1}$ and $w_{(2)}\in W_{2}$,
the series 
\[
\sum_{n\in \R}(e^{z'L'_{P(z)}(0)}\lambda^{(2)}_{n})(w_{(1)}\otimes w_{(2)})
\]
has the following properties:
\begin{enumerate}

\item[(i)] It can be written as the iterated series
\[
\sum_{n\in \R}e^{nz'}\left(\left(\sum_{i=0}^{N}\frac{(z')^{i}}{i!}
(L'_{P(z)}(0)-n)^{i}\lambda^{(2)}_{n}
\right)(w_{(1)}\otimes 
w_{(2)})\right)
\]
(recall that $\R\times\{0,\dots,N\}$ is a unique expansion set).

\item[(ii)] It is absolutely convergent for $z'\in \C$ in the neighborhood of
$z'=0$ above.

\item[(iii)] It is absolutely convergent to $\mu^{(2)}_{\lambda,
w_{(3)}}(w_{(1)}\otimes w_{(2)})$ when $z'=0$:
\[
\sum_{n\in \R}\lambda^{(2)}_{n}(w_{(1)}\otimes w_{(2)})
=\mu^{(2)}_{\lambda, w_{(3)}}(w_{(1)}\otimes w_{(2)})
=\lambda(w_{(1)}\otimes w_{(2)}\otimes w_{(3)}).
\]

\end{enumerate}

(b) For any $w_{(3)}\in W_3$, let $W^{(2)}_{\lambda, w_{(3)}}$ be the
smallest doubly graded (or equivalently, $\tilde A$-graded) subspace
of $((W_1\otimes W_2)^{*})_{[\R]}^{(\tilde{A})}$, or equivalently, of
$((W_1\otimes W_2)^{*})_{[\C]}^{(\tilde{A})}$, containing all the
terms $\lambda^{(2)}_{n}$ in the formal series in (a) and stable under
the component operators $\tau_{P(z)}(v\otimes t^{m})$ of the operators
$Y'_{P(z)}(v, x)$ for $v\in V$, $m\in {\mathbb Z}$, and under the
operators $L'_{P(z)}(-1)$, $L'_{P(z)}(0)$ and $L'_{P(z)}(1)$.  (As
above, $W^{(2)}_{\lambda, w_{(3)}}$ indeed exists.)  Then
$W^{(2)}_{\lambda, w_{(3)}}$ has the properties
\begin{eqnarray*}
&\dim(W^{(2)}_{\lambda, w_{(3)}})^{(\beta)}_{[n]}<\infty,&\\
&(W^{(2)}_{\lambda, w_{(3)}})^{(\beta)}_{[n+k]}=0\;\;\mbox{ for
}\;k\in {\mathbb Z} \;\mbox{ sufficiently negative}&
\end{eqnarray*}
for any $n\in \R$ and $\beta\in \tilde A$, where the
subscripts denote the $\R$-grading by
$L'_{P(z)}(0)$-(generalized) eigenvalues and the superscripts denote
the $\tilde A$-grading.

\item{\bf The $L(0)$-semisimple 
$P^{(2)}(z)$-local grading restriction condition}

(a) The {\em $L(0)$-semisimple $P^{(2)}(z)$-grading condition}: For
any $w_{(3)}\in W_3$, there exists a formal series $\sum_{n\in
\R}\lambda^{(2)}_{n}$ with
\[
\lambda^{(2)}_{n}\in \coprod_{\beta\in \tilde{A}}
((W_{1}\otimes W_{2})^{*})_{(n)}^{(\beta)}
\]
for $n\in \R$ and an open neighborhood of $z'=0$
such that for $w_{(1)}\in W_{1}$ and $w_{(2)}\in W_{2}$,
the series  
\[
\sum_{n\in \R}(e^{z'L'_{P(z)}(0)}\lambda^{(2)}_{n})(w_{(1)}\otimes w_{(2)})
\]
has the following properties:
\begin{enumerate}

\item[(i)] It can be written as 
\[
\sum_{n\in \R}e^{nz'}\lambda^{(2)}_{n}(w_{(1)}\otimes 
w_{(2)})
\]
(recall that $\R\times\{0\}$ is a unique expansion set).

\item[(ii)] It is absolutely convergent for $z'\in \C$ in the neighborhood of
$z'=0$ above.

\item[(iii)] It is absolutely convergent to $\mu^{(2)}_{\lambda,
w_{(3)}}(w_{(1)}\otimes w_{(2)})$ when $z'=0$:
\[
\sum_{n\in \R}\lambda^{(2)}_{n}(w_{(1)}\otimes w_{(2)})
=\mu^{(2)}_{\lambda, w_{(3)}}(w_{(1)}\otimes w_{(2)})
=\lambda(w_{(1)}\otimes w_{(2)}\otimes w_{(3)}).
\]

\end{enumerate}
(Note that such an element $\lambda$ also satisfies the
$P^{(2)}(z)$-grading condition above with the same elements
$\lambda^{(2)}_{n}$.)

(b) For any $w_{(3)}\in W_3$, consider the space $W^{(2)}_{\lambda,
w_{(3)}}$ as above, which in this case is in fact the smallest doubly
graded (or equivalently, $\tilde A$-graded) subspace of $((W_1\otimes
W_2)^{*})_{(\R)}^{(\tilde{A})}$ (or of $((W_1\otimes
W_2)^{*})_{(\C)}^{(\tilde{A})}$) containing all the terms
$\lambda^{(2)}_{n}$ in the formal series in (a) and stable under the
component operators $\tau_{P(z)}(v\otimes t^{m})$ of the operators
$Y'_{P(z)}(v, x)$ for $v\in V$, $m\in {\mathbb Z}$, and under the
operators $L'_{P(z)}(-1)$, $L'_{P(z)}(0)$ and $L'_{P(z)}(1)$. Then
$W^{(2)}_{\lambda, w_{(3)}}$ has the properties
\begin{eqnarray*}
&\dim(W^{(2)}_{\lambda, w_{(3)}})^{(\beta)}_{(n)}<\infty,&\\
&(W^{(2)}_{\lambda, w_{(3)}})^{(\beta)}_{(n+k)}=0\;\;\mbox{ for
}\;k\in {\mathbb Z} \;\mbox{ sufficiently negative}&
\end{eqnarray*}
for any $n\in \R$ and $\beta\in \tilde A$, where the
subscripts denote the $\R$-grading by
$L'_{P(z)}(0)$-eigenvalues and the superscripts denote the $\tilde
A$-grading.

\end{description}

\begin{rema}\label{part-a}
{\rm Part (a) of each of these conditions says in particular that
there exists $N\in \N$ ($N=0$ in the semisimple case) such that when
applied to an arbitrary element of $W_{2}\otimes W_{3}$ or
$W_{1}\otimes W_{2}$,
$(L'_{P(z)}(0)-L'_{P(z)}(0)_{s})^{N+1}\lambda^{(1)}_{n}$ or
$(L'_{P(z)}(0)-L'_{P(z)}(0)_{s})^{N+1}\lambda^{(2)}_{n}$ becomes $0$
for $n\in \R$.  Thus Part (a) of each of these conditions implies:
\begin{eqnarray*}
\lefteqn{\sum_{n\in \R}(e^{yL'_{P(z)}(0)}\lambda^{(1)}_{n})(w_{(2)}\otimes w_{(3)})}\nn
&&=\sum_{n\in \R}e^{ny}\left(\left(\sum_{i=0}^{N}\frac{y^{i}}{i!}
(L'_{P(z)}(0)-n)^{i}\lambda^{(1)}_{n}
\right)(w_{(2)}\otimes 
w_{(3)})\right)
\end{eqnarray*}
or
\[
\sum_{n\in \R}(e^{yL'_{P(z)}(0)}\lambda^{(1)}_{n})(w_{(2)}\otimes w_{(3)})=
\sum_{n\in \R}e^{ny}\lambda^{(1)}_{n}(w_{(2)}\otimes 
w_{(3)})
\]
or
\begin{eqnarray*}
\lefteqn{\sum_{n\in \R}(e^{yL'_{P(z)}(0)}\lambda^{(2)}_{n})(w_{(1)}\otimes w_{(2)})}\nn
&&=
\sum_{n\in \R}e^{ny}\left(\left(\sum_{i=0}^{N}\frac{y^{i}}{i!}
(L'_{P(z)}(0)-n)^{i}\lambda^{(2)}_{n}
\right)(w_{(1)}\otimes 
w_{(2)})\right)
\end{eqnarray*}
or 
\[
\sum_{n\in \R}(e^{yL'_{P(z)}(0)}\lambda^{(2)}_{n})(w_{(1)}\otimes w_{(2)})=
\sum_{n\in \R}e^{ny}\lambda^{(2)}_{n}(w_{(1)}\otimes 
w_{(2)}).
\]
Part (a) of each of these conditions also asserts in particular, 
in the language of weak absolute
convergence (recall Remark \ref{weakly-abs-conv}),
that, for example, for any $w_{(1)}\in W_1$,
$\mu^{(1)}_{\lambda, w_{(1)}}$ is the sum of a weakly absolutely
convergent series $\sum_{n\in \R}\lambda^{(1)}_{n}$ with
$\lambda^{(1)}_{n}\in \coprod_{\beta\in \tilde{A}} ((W_2\otimes
W_3)^{*})_{[n]}^{(\beta)}$. }
\end{rema}

While the grading restriction conditions above asserting the
existence of the elements $\lambda^{(1)}_{n}$ or $\lambda^{(2)}_{n}$
do not say anything about the uniqueness of these elements, they are
indeed unique:

\begin{propo}\label{unique-lambda-n}
The elements $\lambda^{(1)}_{n}$, $n\in \R$, in the $P^{(1)}(z)$-local
grading restriction condition (or the $L(0)$-semisimple
$P^{(1)}(z)$-local grading restriction condition) and the elements
$\lambda^{(2)}_{n}$, $n\in \R$, in the $P^{(2)}(z)$-local grading
restriction condition (or the $L(0)$-semisimple $P^{(2)}(z)$-local
grading restriction condition) are uniquely determined by the
properties indicated in Part (a) of the conditions.
\end{propo}

Since the proof of this result follows easily from certain facts
established in the proof of Theorem \ref{9.7-1} below, we defer it to
Remark \ref{pf-unique-lambda-n}.  This uniqueness result implies the
following bilinearity result for the elements $\lambda^{(1)}_{n}$ and
$\lambda^{(2)}_{n}$:

\begin{corol}\label{bilincorol}
The set of $\lambda \in (W_1 \otimes W_2 \otimes W_3)^*$ satisfying
any of the four local grading restriction
conditions forms a linear subspace.  The elements $\lambda_{n}^{(1)}$,
$n\in \R$, in the $P^{(1)}(z)$-local grading restriction condition (or
the $L(0)$-semisimple $P^{(1)}(z)$-local grading restriction
condition) are bilinear in $\lambda$ and $w_{(1)}$, and the elements
$\lambda^{(2)}_{n}$, $n\in \R$, in the $P^{(2)}(z)$-local grading
restriction condition (or the $L(0)$-semisimple $P^{(2)}(z)$-local
grading restriction condition) are bilinear in $\lambda$ and
$w_{(3)}$. 
\end{corol}
\pf We prove only the case of the $P^{(1)}(z)$-local grading
restriction condition; the other cases are handled the same way.

We shall use the notation $\lambda_{n}^{(1)}(\lambda, w_{(1)})$, $n\in
\R$, to denote $\lambda_{n}^{(1)}$ in the $P^{(1)}(z)$-local grading
restriction condition to exhibit the dependence of these elements on
$\lambda$ and $w_{(1)}$.  Let $\lambda$ and $\tilde{\lambda}$ be
elements of $(W_{1}\otimes W_{2}\otimes W_{3})^{*}$ satisfying the
$P^{(1)}(z)$-local grading restriction condition, $w_{(1)}$ and
$\tilde{w}_{(1)}$ elements of $W_{1}$, and $a, b, c$ and $d$ complex
numbers.  Then the formal series
\[
\sum_{n\in \R}(ac\lambda_{n}^{(1)}(\lambda, w_{(1)})+
ad\lambda_{n}^{(1)}(\lambda, \tilde{w}_{(1)})+bc\lambda_{n}^{(1)}(\tilde{\lambda}, w_{(1)})
+bd\lambda_{n}^{(1)}(\tilde{\lambda}, \tilde{w}_{(1)}))
\]
satisfies (i) (ii) and (iii) in Part (a) of the $P^{(1)}(z)$-local
grading restriction condition for $a\lambda+b\tilde{\lambda}\in
(W_{1}\otimes W_{2}\otimes W_{3})^{*}$ and
$cw_{(1)}+d\tilde{w}_{(1)}\in W_{1}$, where we use the intersection of
the four open neighborhoods of $z'=0$ and the maximum of the relevant
nonnegative integers $N$.  The summands (for $n\in \R$) also satisfy
Part (b) of the condition.  Thus $a\lambda+b\tilde{\lambda}$ satisfies
the $P^{(1)}(z)$-local grading restriction condition, and by
Proposition \ref{unique-lambda-n} we have
\begin{eqnarray*}
\lefteqn{\lambda_{n}^{(1)}(a\lambda+b\tilde{\lambda}, cw_{(1)}+d\tilde{w}_{(1)})}\nn
&&=ac\lambda_{n}^{(1)}(\lambda, w_{(1)})+
ad\lambda_{n}^{(1)}(\lambda, \tilde{w}_{(1)})+bc\lambda_{n}^{(1)}(\tilde{\lambda}, w_{(1)})
+bd\lambda_{n}^{(1)}(\tilde{\lambda}, \tilde{w}_{(1)})
\end{eqnarray*}
for $n\in \R$.
\epfv

Using the uniqueness and recalling the $\tilde{A}$-homogeneous
subspaces (\ref{W1W2beta}) and (\ref{W1W2W3beta}), we obtain the
following natural $\tilde{A}$-properties of the elements $\lambda_n$
in each of the four conditions:

\begin{propo}\label{lambda-n-Atilde}
Suppose that
\[
\lambda \in ((W_1 \otimes W_2 \otimes W_3)^*)^{(\beta)},
\]
with $\beta \in \tilde{A}$, satisfies the $P^{(1)}(z)$-local grading
restriction condition, and suppose that
\[
w_{(1)}\in W_{1}^{(\beta_{1})},
\]
with $\beta_{1} \in \tilde{A}$.  Then for each $n \in \R$,
\[
\lambda_{n}^{(1)}\in ((W_{2}\otimes W_{3})^{*})^{(\beta+\beta_{1})}_{[n]}.
\]
The analogous statement holds for each of the other three conditions;
for instance, if
\[
\lambda \in ((W_1 \otimes W_2 \otimes W_3)^*)^{(\beta)}
\]
satisfies the $L(0)$-semisimple $P^{(2)}(z)$-local grading restriction
condition and
\[
w_{(3)}\in W_{3}^{(\beta_{3})},
\]
with $\beta_{3} \in \tilde{A}$, then
\[
\lambda_{n}^{(2)}\in ((W_{1}\otimes W_{2})^{*})^{(\beta+\beta_{3})}_{(n)}.
\]
\end{propo}
\pf We prove only the first case mentioned, the proofs in the other
cases being similar.

We have
\[
\mu^{(1)}_{\lambda, w_{(1)}} \in 
((W_{2}\otimes W_{3})^{*})^{(\beta+\beta_{1})},
\]
since for $\beta_{2}, \beta_{3}\in \tilde{A}$ such that
$\beta_{2}+\beta_{3}\ne -\beta-\beta_{1}$ and for $w_{(2)}\in
W_{2}^{(\beta_{2})}$, $w_{(3)}\in W_{3}^{(\beta_{3})}$,
\[
\mu^{(1)}_{\lambda, w_{(1)}}(w_{(2)}\otimes w_{(3)})
=\lambda(w_{(1)}\otimes w_{(2)}\otimes w_{(3)})
=0,
\]
so that we have the absolutely convergent sum
\[
\sum_{n\in \R}\lambda^{(1)}_{n}(w_{(2)}\otimes w_{(3)})
=\mu^{(1)}_{\lambda, w_{(1)}}(w_{(2)}\otimes w_{(3)})
=0.
\]
For each $n \in \R$, let $\tilde{\lambda}_{n}^{(1)}$ be the projection
of ${\lambda}_{n}^{(1)}$ to $((W_{2}\otimes
W_{3})^{*})_{[n]}^{(\beta+\beta_{1})}$:
\[
\tilde{\lambda}_{n}^{(1)}(w_{(2)}\otimes w_{(3)})
=\left\{\begin{array}{ll}\lambda^{(1)}_{n}(w_{(2)}\otimes w_{(3)})&
\;\mbox{ if }\; \beta_{2}+\beta_{3}= -\beta-\beta_{1}\\
0 & \;\mbox{ if }\; \beta_{2}+\beta_{3}\ne -\beta-\beta_{1}.
\end{array}\right.
\]
Then clearly the $\tilde{\lambda}_{n}^{(1)}$ for $n\in \R$ also
satisfy Part (a) of the $P^{(1)}(z)$-local grading restriction
condition, and so by the uniqueness (Proposition
\ref{unique-lambda-n}),
\[
\lambda^{(1)}_{n}=\tilde{\lambda}^{(1)}_{n},
\]
so that
\[
\lambda^{(1)}_{n}\in ((W_{2}\otimes W_{3})^{*})_{[n]}^{(\beta+\beta_{1})}
\]
for $n\in \R$. 
\epfv

In the rest of this section, we shall focus on the case that the
convergence condition for intertwining maps in $\mathcal{C}$ holds and
that the generalized $V$-modules that we start with are objects of
$\mathcal{C}$. Recall the categories $\mathcal{M}_{sg}$ and
$\mathcal{GM}_{sg}$ from Notation $\ref{MGM}$, and recall Assumptions
\ref{assum}, \ref{assum-c} and \ref{assum-exp-set} on the category
$\mathcal{C}$.

\begin{rema}\label{I1I2'}
{\rm Let $W_{1}$, $W_{2}$, $W_{3}$ and $W_{4}$ be generalized
$V$-modules.  Given an $\tilde{A}$-compatible map
\[
F: W_{1}\otimes W_{2}\otimes W_{3}\to \overline{W}_{4}
\]
as in Remark \ref{Atildecompatcorrespondence}, there is a canonical
$\tilde{A}$-compatible map $G$ from $W'_{4}$ to $(W_{1}\otimes
W_{2}\otimes W_{3})^{*}$ corresponding to $F$ under the indicated
canonical isomorphism between the spaces of such maps. We shall denote
the map $G$ corresponding to $F$ by $F'$:
\[
F'=G: W'_{4} \to (W_{1}\otimes W_{2}\otimes W_{3})^{*}.
\]
Assume the convergence condition for intertwining maps in
$\mathcal{C}$ and that all generalized $V$-modules considered are
objects of $\mathcal{C}$.  Let $I_{1}$, $I_{2}$, $I^1$ and $I^2$ be
$P(z_1)$-, $P(z_2)$-, $P(z_2)$- and $P(z_0)$-intertwining maps of
types ${W_4}\choose {W_1M_1}$, ${M_1}\choose {W_2W_3}$, ${W_4}\choose
{M_2W_3}$ and ${M_2}\choose {W_1W_2}$, respectively.  Then the maps
\begin{equation}\label{adj-prod}
(I_1\circ (1_{W_1}\otimes I_2))': W_{4}'\to (W_{1}\otimes W_{2}\otimes
W_{3})^{*}
\end{equation}
for $|z_{1}|>|z_{2}|>0$ and
\begin{equation}\label{adj-iter}
(I^1\circ (I^2\otimes 1_{W_3}))': W_{4}'\to (W_{1}\otimes
W_{2}\otimes W_{3})^{*}
\end{equation}
for $|z_{2}|>|z_{0}|>0$ are well-defined $\tilde{A}$-compatible maps.
In particular, for $w'_{(4)}\in W'_{4}$, we have the elements
\begin{equation}\label{lambda1}
(I_1\circ (1_{W_1}\otimes I_2))'(w'_{(4)})\in (W_{1}\otimes
W_{2}\otimes W_{3})^{*}
\end{equation}
and
\begin{equation}\label{lambda2}
(I^1\circ (I^2\otimes 1_{W_3}))'(w'_{(4)})\in (W_{1}\otimes
W_{2}\otimes W_{3})^{*}.
\end{equation}
}
\end{rema}

Proposition \ref{lambda-n-Atilde} applies to such elements $\lambda$
when $w'_{(4)}$ is $\tilde{A}$-homogeneous, since for
\[
w'_{(4)} \in (W_{4}')^{(\beta_{4})},
\]
with $\beta_{4}\in \tilde{A}$, we have
\[
\lambda \in ((W_1 \otimes W_2 \otimes W_3)^*)^{(\beta_{4})}
\]
for $\lambda$ either of the elements (\ref{lambda1}), (\ref{lambda2}),
by the $\tilde{A}$-compatibility.  This yields the natural
$\tilde{A}$-properties of the corresponding elements $\lambda_n$ in
each of the four conditions (the complex numbers $z$ being chosen in
the ways that they arise naturally in the theory):

\begin{propo}\label{lambda-n-a-tilde}
Assume that the convergence condition for intertwining maps in
$\mathcal{C}$ holds.  Let $W_{1}$, $W_{2}$, $W_{3}$, $W_{4}$, $M_{1}$
and $M_{2}$ be objects of $\mathcal{C}$ and let $I_{1}$, $I_{2}$,
$I^1$ and $I^2$ be $P(z_1)$-, $P(z_2)$-, $P(z_2)$- and
$P(z_0)$-intertwining maps of types ${W_4}\choose {W_1M_1}$,
${M_1}\choose {W_2W_3}$, ${W_4}\choose {M_2W_3}$ and ${M_2}\choose
{W_1W_2}$, respectively.  Let
\[
w'_{(4)}\in (W_{4}')^{(\beta_{4})},
\]
with
\[
\beta_{4}\in \tilde{A}.
\]
Assume that $|z_{1}|>|z_{2}|>0$ and let
\[
\lambda = (I_1\circ (1_{W_1}\otimes I_2))'(w'_{(4)}).
\]
Suppose that $\lambda$ satisifies the $P^{(1)}(z_2)$-local grading
restriction condition or, respectively, the $P^{(2)}(z_0)$-local
grading restriction condition.  Then for $w_{(1)}\in
W_{1}^{(\beta_{1})}$ or, respectively, $w_{(3)}\in
W_{3}^{(\beta_{3})}$, with $\beta_{j} \in \tilde{A}$, we have
\[
\lambda_{n}^{(1)}\in ((W_{2}\otimes W_{3})^{*})^{(\beta_{4}+\beta_{1})}_{[n]}
\]
or, respectively,
\[
\lambda_{n}^{(2)}\in ((W_{1}\otimes
W_{2})^{*})^{(\beta_{4}+\beta_{3})}_{[n]}
\]
for each $n\in \R$.  When $\mathcal{C}$ is in $\mathcal{M}_{sg}$,
suppose instead that $\lambda$ satisfies the $L(0)$-semisimple
$P^{(1)}(z_2)$-local grading restriction condition or, respectively,
the $L(0)$-semisimple $P^{(2)}(z_0)$-local grading restriction
condition.  Then for $w_{(1)}$ and $w_{(3)}$ as above, we have
\[
\lambda_{n}^{(1)}\in ((W_{2}\otimes W_{3})^{*})^{(\beta_{4}+\beta_{1})}_{(n)}
\]
or, respectively,
\[
\lambda_{n}^{(2)}\in ((W_{1}\otimes W_{2})^{*})^{(\beta_{4}+\beta_{3})}_{(n)}
\]
for each $n\in \R$.  The analogous conclusions hold if, instead,
$|z_{2}|>|z_{0}|>0$ and $(I^1\circ (I^2\otimes 1_{W_3}))'(w'_{(4)})$
satisfies the $P^{(1)}(z_2)$-local grading restriction condition
(or the $L(0)$-semisimple $P^{(1)}(z_2)$-local grading restriction
condition when $\mathcal{C}$ is in $\mathcal{M}_{sg}$) or the
$P^{(2)}(z_0)$-local grading restriction condition (or the
$L(0)$-semisimple $P^{(2)}(z_0)$-local grading restriction condition
when $\mathcal{C}$ is in $\mathcal{M}_{sg}$).  \epf
\end{propo}

In the next result, we prove that for the product of a
$P(z_{1})$-intertwining map $I_1$ and a $P(z_{2})$-intertwining map
$I_2$, each element (\ref{lambda1}) of the image of the map
(\ref{adj-prod}) satisfies the $P^{(1)}(z_{2})$-local grading
restriction condition and that for the iterate of a
$P(z_{2})$-intertwining map $I^1$ and a $P(z_{0})$-intertwining map
$I^2$, each element (\ref{lambda2}) of the image of the map
(\ref{adj-iter}) satisfies the $P^{(2)}(z_{0})$-local grading
restriction condition.

\begin{propo}\label{9.7}
Assume that the convergence condition for intertwining maps in
$\mathcal{C}$ holds.  Let $W_{1}$, $W_{2}$, $W_{3}$, $W_{4}$, $M_{1}$
and $M_{2}$ be objects of $\mathcal{C}$ and let $I_{1}$, $I_{2}$,
$I^1$ and $I^2$ be $P(z_1)$-, $P(z_2)$-, $P(z_2)$- and
$P(z_0)$-intertwining maps of types ${W_4}\choose {W_1M_1}$,
${M_1}\choose {W_2W_3}$, ${W_4}\choose {M_2W_3}$ and ${M_2}\choose
{W_1W_2}$, respectively.  Let $w'_{(4)}\in W'_4$.  If $|z_1|>|z_2|>0$,
then
\[
(I_1\circ (1_{W_1}\otimes I_2))'(w'_{(4)}) \in
(W_{1}\otimes W_{2}\otimes W_{3})^{*}
\]
satisfies the $P^{(1)}(z_2)$-local grading restriction condition (or
the $L(0)$-semisimple $P^{(1)}(z_2)$-local grading restriction
condition when $\mathcal{C}$ is in $\mathcal{M}_{sg}$), and if
$|z_2|>|z_0|>0$, then
\[
(I^1\circ (I^2\otimes 1_{W_3}))'(w'_{(4)}) \in
(W_{1}\otimes W_{2}\otimes W_{3})^{*}
\]
satisfies the $P^{(2)}(z_0)$-local grading restriction condition (or
the $L(0)$-semisimple $P^{(2)}(z_0)$-local grading restriction
condition when $\mathcal{C}$ is in $\mathcal{M}_{sg}$).  
Moreover, suppose that
$\mathcal{C}$ is closed under images (recall Definition
\ref{closedunderimages}).  Let $w_{(1)}\in W_1$ and $w_{(3)}\in W_3$.
Take
\[
\lambda_{n}^{(1)} \in (W_{2}\otimes W_{3})^{*}
\]
and
\[
\lambda_{n}^{(2)} \in (W_{1}\otimes W_{2})^{*},
\]
$n\in \R$, to be the elements constructed in the proof below.  Then
the corresponding spaces
\[
W^{(1)}_{(I_1\circ (1_{W_1}\otimes I_2))'(w'_{(4)}), w_{(1)}} \subset
(W_{2}\otimes W_{3})^{*}
\]
and 
\[
W^{(2)}_{(I^1\circ (I^2\otimes 1_{W_3}))'(w'_{(4)}), w_{(3)}} \subset
(W_{1}\otimes W_{2})^{*}
\]
(constructed using these elements $\lambda_{n}^{(1)}$ and
$\lambda_{n}^{(2)}$), equipped with the vertex operator maps given by
$Y'_{P(z_{2})}$ and $Y'_{P(z_{0})}$, respectively, and the operators
$L'_{P(z_{2})}(j)$ and $L'_{P(z_{0})}(j)$, respectively, for $j=-1, 0,
1$, are generalized $V$-submodules of objects of $\mathcal{C}$
included in $(W_{2}\otimes W_{3})^{*}$ and $(W_{1}\otimes W_{2})^{*}$,
respectively.  In particular, for any $n\in \R$, the doubly-graded
subspaces
\[
W_{\lambda_{n}^{(1)}} \subset
W^{(1)}_{(I_1\circ (1_{W_1}\otimes I_2))'(w'_{(4)}), w_{(1)}}
\]
and
\[
W_{\lambda_{n}^{(2)}} \subset
W^{(2)}_{(I^1\circ (I^2\otimes 1_{W_3}))'(w'_{(4)}), w_{(3)}}
\]
(recall the notation $W_{\lambda}$ in the $P(z)$-local grading
restriction condition in Section 5) are also generalized
$V$-submodules of objects of $\mathcal{C}$ included in $(W_{2}\otimes
W_{3})^{*}$ and $(W_{1}\otimes W_{2})^{*}$, respectively;
$W_{\lambda_{n}^{(1)}}$ is the generalized $V$-submodule generated by
$\lambda_{n}^{(1)}$ (the smallest generalized $V$-submodule containing
$\lambda_{n}^{(1)}$), and analogously for $\lambda_{n}^{(2)}$.
\end{propo}
\pf 
Let $w_{(1)}\in W_{1}$. 
For $n\in \R$, let $m_{(1), n}'\in M_{1}^{*}$ be defined by 
\[
m_{(1), n}'(m_{(1)})=\langle w'_{(4)}, I_1(w_{(1)}\otimes  
\pi_{n}(m_{(1)}))\rangle
\]
for $m_{(1)}\in M_{1}$. Since $w_{(1)}$ and $w'_{(4)}$ are finite sums
of $\tilde{A}$-homogeneous elements and $I_{1}$ is
$\tilde{A}$-compatible, $m_{(1), n}'$ is also a finite sum of
$\tilde{A}$-homogeneous elements.  Since by definition, for
$m_{(1)}\in (M_{1})_{[m]}$, $m_{(1), n}'(m_{(1)})=0$ when $m\ne n$, we
see that $m_{(1), n}'\in (M'_{1})_{[n]}$.

Let 
\[
\lambda^{(1)}_{n}=m_{(1), n}'\circ I_{2} \in (W_{2}\otimes
W_{3})^{*}.
\]
{}From Notation \ref{scriptN}, we have
$\lambda^{(1)}_{n}= I_{2}'(m_{(1), n}')$, where $I_{2}': M_{1}' \to
(W_{2}\otimes W_{3})^{*}$ is as indicated in Notation \ref{scriptN}.
Since $m_{(1), n}'\in (M'_{1})_{[n]}$, by Proposition \ref{im:abc}(b),
\[
\lambda^{(1)}_{n}= I_{2}'(m_{(1), n}') \in ((W_{2}\otimes
W_{3})^{*})_{[n]}
\]
for $n\in \R$.  In addition, since $I_{2}'$ is $\tilde{A}$-compatible
and $m_{(1), n}'$ is a finite sum of $\tilde{A}$-homogeneous elements,
$\lambda^{(1)}_{n}$ is also a finite sum of $\tilde{A}$-homogeneous
elements.  By Proposition \ref{exp-set}, the set $\{(n, i)\in \C\times
\N\;|\; (L(0)-n)^{i}(M'_{1})_{[n]}\ne 0\}$ is included in a (unique
expansion) set of the form $\R\times \{0, \dots, N\}$, and its subset
$\{(n, i)\in \C\times \N\;|\; (L'_{P(z_{2})}(0)-n)^{i}
\lambda^{(1)}_{n} \ne 0\}$ is included in the same set (recall that
$I_{2}'$ intertwines the various actions, including those of $L(0)$
and $L'_{P(z_{2})}(0)$).

When $|z_1|>|z_2|>0$, the product of $I_1$ and $I_2$ exists.
For $w_{(2)}\in W_2$ and $w_{(3)}\in W_3$ we have
\begin{eqnarray*}
\lefteqn{\mu^{(1)}_{(I_1\circ (1_{W_1}\otimes I_2))'(w'_{(4)}),
w_{(1)}}(w_{(2)}\otimes w_{(3)})}\nno\\
&&=\langle w'_{(4)}, I_1(w_{(1)}\otimes I_2(w_{(2)}\otimes w_{(3)}))
\rangle\nno\\
&&=\sum_{n\in \R}\langle w'_{(4)}, 
I_1(w_{(1)}\otimes \pi_{n}(I_2(w_{(2)}\otimes w_{(3)})))
\rangle\nno\\
&&=\sum_{n\in \R}m_{(1), n}'(I_2(w_{(2)}\otimes w_{(3)}))\nno\\
&&=\sum_{n\in \R}\lambda^{(1)}_{n}(w_{(2)}\otimes w_{(3)}),
\end{eqnarray*}
an absolutely convergent series. 

Let $\mathcal{Y}_{1}=\mathcal{Y}_{I_{1},0}$ and
$\mathcal{Y}_{2}=\mathcal{Y}_{I_{2},0}$ (recall Proposition
\ref{im:correspond}) so that 
\begin{eqnarray*}
I_{1}(w_{(1)}\otimes w)&=&\Y_{1}(w_{(1)}, z_{1})w,\\
I_{2}(w_{(2)}\otimes w_{(3)})&=&\Y_{2}(w_{(2)}, z_{2})w_{(3)}
\end{eqnarray*}
for $w_{(1)}\in W_{1}$, $w_{(2)}\in W_{2}$,
$w_{(3)}\in W_{3}$ and $w\in M_{1}$ (recall the ``substitution'' notation {}from
(\ref{im:f(z)}), where we choose $p=0$).  
By Proposition \ref{im:abc}(b), the map $I_{2}'$ preserves generalized
weights. For $z'\in \C$, we also have, using (\ref{log:p2}) (recall 
(\ref{log:not3}) and Remark \ref{3.33}),
\begin{eqnarray*}
\lefteqn{e^{z'L(0)}I_2(w_{(2)}\otimes w_{(3)})}\nn
&&=y^{L(0)}\Y_{2}(w_{(2)}, x_{2}) w_{(3)}
\lbar_{y^{m}=e^{mz'},\; \log y=z',\; x_{2}^{m}=e^{m\log z_{2}},\; \log x_{2}=\log z_{2}}\nn
&&=\Y_{2}(y^{L(0)}w_{(2)}, x_{2}y) y^{L(0)}w_{(3)}
\lbar_{y^{m}=e^{mz'},\; \log y=z',\; x_{2}^{m}=e^{m\log z_{2}},\; \log x_{2}=\log z_{2}}\nn
&&=\Y_{2}(e^{z'L(0)}w_{(2)}, x) e^{z'L(0)}w_{(3)}\lbar_{x^{m}=e^{m((\log z_{2})+z')},\;
\log x=(\log z_{2})+z'}.
\end{eqnarray*}
Thus
\begin{eqnarray*}
\lefteqn{\sum_{n\in \R}(e^{z'L'_{P(z_{2})}(0)}
\lambda^{(1)}_{n})(w_{(2)}\otimes w_{(3)})}\nn
&&=\sum_{n\in \R}(e^{z'L'_{P(z_{2})}(0)}(I_{2}'(m_{(1), n}')))
(w_{(2)}\otimes w_{(3)})\nno\\
&&=\sum_{n\in \R}(I_{2}'(e^{z'L(0)}m_{(1), n}'))
(w_{(2)}\otimes w_{(3)})\nno\\
&&=\sum_{n\in \R}(e^{z'L(0)}m_{(1), n}')
(I_2(w_{(2)}\otimes w_{(3)}))\nno\\
&&=\sum_{n\in \R}m_{(1), n}'
(e^{z'L(0)}I_2(w_{(2)}\otimes w_{(3)}))\nno\\
&&=\sum_{n\in \R}\langle w'_{(4)}, 
I_1(w_{(1)}\otimes \pi_{n}(e^{z'L(0)}I_2(w_{(2)}\otimes w_{(3)})))
\rangle\nno\\
&&=\sum_{n\in \R}\langle w'_{(4)}, 
\Y_1(w_{(1)}, x_{1})\cdot\nn
&&\quad\quad\quad
\pi_{n}(\Y_{2}(e^{z'L(0)}w_{(2)}, x) e^{z'L(0)}
w_{(3)})
\rangle\lbar_{x_{1}^{m}=e^{m\log z_{1}},\; \log x_{1}=\log z_{1},\;
x^{m}=e^{m((\log z_{2})+z')},\;
\log x=(\log z_{2})+z'}
\end{eqnarray*}
and when $z'$ is in a small open neighborhood of $0$ 
such that in particular $|z_{1}|>|e^{z'} z_{2}|>0$, 
this series is absolutely convergent for all $w_{(1)}\in W_{1}$,
$w_{(2)}\in W_{2}$, $w_{(3)}\in W_{3}$, and $w_{(4)}'\in W'_{4}$
(note that there exists $p\in \Z$ 
such that $(\log z_{2})+z'=l_{p}(z_{2}e^{z'})$; in fact, 
if $z_{2}$ is not a positive real number, then $(\log z_{2})+z'=\log(z_{2}e^{z'})$
for $z'$ sufficiently near $0$ and if $z_{2}$ is a positive real number,
then either $(\log z_{2})+z'=\log(z_{2}e^{z'})$ or 
$(\log z_{2})+z'=l_{-1}(z_{2}e^{z'})$ (recall (\ref{branch2}))).
Hence $(I_1\circ
(1_{W_1}\otimes I_2))'(w'_{(4)})$ satisfies the $P^{(1)}(z_2)$-grading
condition.

We know that the map $I_{2}'$ preserves generalized weights, and
$I_{2}'$ is also $\tilde{A}$-compatible.  Thus the image under
$I_{2}'$ of the generalized $V$-submodule of the generalized
$V$-module $M'_{1}$ generated by the elements $m_{(1), n}'$ for $n\in
\R$ (that is, the smallest (strongly $\tilde{A}$-graded) generalized
$V$-submodule containing these elements) satisfies the two grading
restriction conditions (\ref{lgrc1}) and (\ref{lgrc2}).  Since
$W^{(1)}_{(I_1\circ (1_{W_1}\otimes I_2))'(w'_{(4)}), w_{(1)}}$ is
this image, Part (b) of the $P^{(1)}(z_2)$-local grading restriction
condition holds.  

When $\mathcal{C}$ is in $\mathcal{M}_{sg}$, the same arguments show
that $(I_1\circ (1_{W_1}\otimes I_2))'(w'_{(4)})$ satisfies the
$L(0)$-semisimple $P^{(1)}(z_2)$-local grading restriction condition.

Moreover, $I_{2}'(M'_{1})$ and $W^{(1)}_{(I_1\circ (1_{W_1}\otimes
I_2))'(w'_{(4)}), w_{(1)}}$ are generalized $V$-modules, and so is
$W_{\lambda_{n}^{(1)}}$ for each $n \in \R$.  If $\mathcal{C}$ is
closed under images, then $I_{2}'(M'_{1})$ is an object of
$\mathcal{C}$ included in $(W_{2}\otimes W_{3})^{*}$, since $M'_{1}
\in \ob \mathcal{C}$, and so $W^{(1)}_{(I_1\circ (1_{W_1}\otimes
I_2))'(w'_{(4)}), w_{(1)}}$ and $W_{\lambda_{n}^{(1)}}$ for each $n
\in \R$ are generalized submodules of objects of $\mathcal{C}$
included in $(W_{2}\otimes W_{3})^{*}$.

Now we handle the other case analogously.  Let $w_{(3)}\in
W_{3}$.  For $n\in \R$, let $m_{(2), n}'\in M_{2}^{*}$ be
defined by
\[
m_{(2), n}'(m_{(2)})=\langle w'_{(4)}, I^1(\pi_{n}(m_{(2)})
\otimes w_{(3)})\rangle
\]
for $m_{(2)}\in M_{2}$. Then $m_{(2), n}'$ is a finite sum of
$\tilde{A}$-homogeneous elements and is an element of
$(M'_{2})_{[n]}$.  Let 
\[
\lambda^{(2)}_{n}=m_{(2), n}'\circ I^{2} \in
(W_{1}\otimes W_{2})^{*}.
\]
Then 
\[
\lambda^{(2)}_{n}= (I^{2})'(m_{(2),
n}') \in ((W_{1}\otimes W_{2})^{*})_{[n]}
\]
for $n\in \R$ and is a finite sum of $\tilde{A}$-homogeneous elements,
where $(I^{2})': M'_{2} \to (W_{1}\otimes W_{2})^{*}$ is as indicated
in Notation \ref{scriptN}.  By Proposition \ref{exp-set}, the set
$\{(n, i)\in \C\times \N\;|\; (L(0)-n)^{i}(M'_{2})_{[n]}\ne 0\}$ is
included in a (unique expansion) set of the form $\R\times \{0, \dots,
N\}$, and its subset $\{(n, i)\in \C\times \N\;|\;
(L'_{P(z_{0})}(0)-n)^{i} \lambda^{(2)}_{n} \ne 0\}$ is included in the
same set.

When $|z_2|>|z_0|>0$, the iterate of $I^1$ and $I^2$ exists.
For $w_{(1)}\in W_1$ and $w_{(2)}\in W_2$, 
\begin{eqnarray*}
\lefteqn{\mu^{(2)}_{(I^1\circ (I^{2}\otimes 1_{W_3}))'(w'_{(4)}),
w_{(3)}}(w_{(1)}\otimes w_{(2)})}\nno\\
&&=\langle w'_{(4)}, I^1(I^2(w_{(1)}\otimes w_{(2)})\otimes w_{(3)})
\rangle\nno\\
&&=\sum_{n\in \R}\langle w'_{(4)}, 
I^1(\pi_{n}(I^2(w_{(1)}\otimes w_{(2)}))\otimes w_{(3)})
\rangle\nno\\
&&=\sum_{n\in \R}m_{(2), n}'(I^2(w_{(1)}\otimes w_{(2)}))\nno\\
&&=\sum_{n\in \R}\lambda^{(2)}_{n}(w_{(1)}\otimes w_{(2)}).
\end{eqnarray*}

Let $\mathcal{Y}^{1}=\mathcal{Y}_{I^{1},0}$ and
$\mathcal{Y}^{2}=\mathcal{Y}_{I^{2},0}$ so that 
\begin{eqnarray*}
I^{1}(w\otimes w_{(3)})&=&\Y^{1}(w, z_{2})w_{(3)},\\
I^{2}(w_{(1)}\otimes w_{(2)})&=&\Y^{2}(w_{(1)}, z_{0})w_{(2)}
\end{eqnarray*}
for $w_{(1)}\in W_{1}$, $w_{(2)}\in W_{2}$,
$w_{(3)}\in W_{3}$ and $w\in M_{2}$.
By Proposition \ref{im:abc}(b), the map $(I^{2})'$ preserves
generalized weights. For $z'\in \C$, we also have 
\begin{eqnarray*}
\lefteqn{e^{z'L(0)}I^2(w_{(1)}\otimes w_{(2)})}\nn
&&=y^{L(0)}\Y^{2}(w_{(1)}, x_{0}) w_{(2)}
\lbar_{y^{m}=e^{mz'},\; \log y=z',\; x_{0}^{m}=e^{m\log z_{0}},\; \log x_{0}=\log z_{0}}\nn
&&=\Y^{2}(y^{L(0)}w_{(1)}, x_{0}y) y^{L(0)}w_{(2)}
\lbar_{y^{m}=e^{mz'},\; \log y=z',\; x_{0}^{m}=e^{m\log z_{0}},\; \log x_{0}=\log z_{0}}\nn
&&=\Y^{2}(e^{z'L(0)}w_{(2)}, x) e^{z'L(0)}w_{(2)}\lbar_{x^{m}=e^{m((\log z_{0})+z')},\;
\log x=(\log z_{0})+z'}.
\end{eqnarray*}
Hence
\begin{eqnarray*}
\lefteqn{\sum_{n\in \R}(e^{z'L'_{P(z_{0})}(0)}
\lambda^{(2)}_{n})(w_{(1)}\otimes w_{(2)})}\nn
&&=\sum_{n\in \R}(e^{z'L'_{P(z_{0})}(0)}((I^{2})'(m_{(2), n}')))
(w_{(1)}\otimes w_{(2)})\nno\\
&&=\sum_{n\in \R}((I^{2})'(e^{z'L(0)}m_{(2), n}'))
(w_{(1)}\otimes w_{(2)})\nno\\
&&=\sum_{n\in \R}(e^{z'L(0)}m_{(2), n}')
(I^2(w_{(1)}\otimes w_{(2)}))\nno\\
&&=\sum_{n\in \R}m_{(2), n}'
(e^{z'L(0)}I^2(w_{(1)}\otimes w_{(2)}))\nno\\
&&=\sum_{n\in \R}\langle w'_{(4)}, 
I^1(\pi_{n}(e^{z'L(0)}I^2(w_{(1)}\otimes w_{(2)}))\otimes w_{(3)})
\rangle\nno\\
&&=\sum_{n\in \R}\langle w'_{(4)}, 
\Y^1(\pi_{n}(\Y^{2}(e^{z'L(0)}w_{(1)}, x) e^{z'L(0)}
w_{(2)}), x_{2})\cdot\nn
&&\quad\quad\quad\quad\quad\quad\cdot
w_{(3)}
\rangle\lbar_{x_{2}^{m}=e^{m\log z_{2}},\; \log x_{2}=\log z_{2},\;
x^{m}=e^{m((\log z_{0})+z')},\;
\log x=(\log z_{0})+z'},
\end{eqnarray*}
which is absolutely convergent for all $w_{(1)}\in W_{1}$,
$w_{(2)}\in W_{2}$, $w_{(3)}\in W_{3}$, and $w_{(4)}'\in W'_{4}$,
when $z'$ is in a small open neighborhood of $0$ such that in particular,
$|z_{2}|>|e^{z'} z_{0}|>0$.
Thus $(I^1\circ
(I^2\otimes 1_{W_3}))'(w'_{(4)})$ satisfies the $P^{(1)}(z_0)$-grading
condition.

Since the map $(I^{2})'$ preserves generalized weights and is also
$\tilde{A}$-compatible, the image under $(I^{2})'$ of the (strongly
$\tilde{A}$-graded) generalized $V$-submodule of the generalized
$V$-module $M_{2}'$ generated by the elements $m_{(2), n}'$ for $n\in
\R$ satisfies the two grading restriction conditions (\ref{lgrc1}) and
(\ref{lgrc2}).  Since $W^{(2)}_{(I^1\circ (I^2\otimes
1_{W_3}))'(w'_{(4)}), w_{(3)}}$ is this image, Part (b) of the
$P^{(2)}(z_0)$-local grading restriction condition holds.  

The rest of the proof proceeds as above.
\epfv

We will often need to prove that certain generalized $V$-modules (or
ordinary $V$-modules), in the {\it original} sense of Definitions
\ref{cvamodule}, \ref{moduleMobius} and
\ref{definitionofgeneralizedmodule} ({\it without} the assumption of
being strongly graded) are indeed strongly graded, and the main
nontrivial properties to verify will often be the grading-restriction
conditions (\ref{set:dmltc}) and (\ref{set:dmfin}).  Thus we shall
find the following definition useful, in particular in the next
result:

\begin{defi}\label{doublygraded}
{\rm In the setting of Definition \ref{def:dgw} (the definition of
``strongly graded''), a generalized $V$-module or a $V$-module (not
necessarily strongly graded, of course) is {\it doubly graded} if it
satisfies all the conditions in Definition \ref{def:dgw} except
perhaps for (\ref{set:dmltc}) and (\ref{set:dmfin}).  The
doubly-graded generalized $V$-submodule {\it generated by} given
elements of a doubly-graded generalized $V$-module is (of course) the
smallest doubly-graded (or equivalently, $\tilde{A}$-graded)
generalized $V$-submodule containing the elements; similarly for
doubly-graded $V$-modules.}
\end{defi}

\begin{rema}\label{submodstrgraded}
{\rm A doubly-graded generalized $V$-submodule of a generalized
$V$-module is of course strongly graded; similarly for $V$-modules.
(Recall that a generalized $V$-module and a $V$-module are
$\tilde{A}$-graded, and in addition strongly graded, by our
conventions.)}
\end{rema}

\begin{rema}
{\rm Such structures have arisen in Propositions \ref{im:abc} and
\ref{im-q:abc}.}
\end{rema}

In general, for the product of a $P(z_{1})$-intertwining map $I_{1}$
and a $P(z_{2})$-intertwining map $I_{2}$, the elements of the image
of the map (\ref{adj-prod}) might not satisfy the
$P^{(2)}(z_{0})$-local grading restriction condition and for the
iterate of a $P(z_{2})$-intertwining map $I^{1}$ and a
$P(z_{0})$-intertwining map $I^{2}$, the elements of the image of the
map (\ref{adj-iter}) might not satisfy the $P^{(1)}(z_{2})$-local
grading restriction condition.  But if they do, we have important
consequences.  In the theorem below, we shall prove a fundamental
consequence, which plays an essential role in the rest of this section
and in our construction of the associativity isomorphisms in the next
section.  The content of this theorem is essentially this: Given a
product, and assuming the relevant condition, we construct a certain
generalized $V$-module which will become (by virtue of Lemma
\ref{intertwine-tau} below) an intermediate module for a suitable
iterate.  This will allow us to express the product as an iterate, and
vice versa when we start with an iterate.  The hard part of the proof
is to show that each term in the series given by the $P^{(2)}(z_{0})$-
or $P^{(1)}(z_{2})$-local grading restriction condition satisfies the
$P(z_{0})$- or $P(z_{2})$-compatibility condition, respectively.

\begin{theo}\label{9.7-1}
Assume that the convergence condition for intertwining maps in
$\mathcal{C}$ holds and that 
\[
|z_1|>|z_2|>|z_{0}|>0.
\]
(Recall that $z_0 = z_1 - z_2$.)  Let $W_{1}$, $W_{2}$, $W_{3}$,
$W_{4}$, $M_{1}$ and $M_{2}$ be objects of $\mathcal{C}$ and let
$I_{1}$, $I_{2}$, $I^1$ and $I^2$ be $P(z_1)$-, $P(z_2)$-, $P(z_2)$-
and $P(z_0)$-intertwining maps of types ${W_4}\choose {W_1M_1}$,
${M_1}\choose {W_2W_3}$, ${W_4}\choose {M_2W_3}$ and ${M_2}\choose
{W_1W_2}$, respectively.  Let $w'_{(4)}\in W'_4$.
\begin{enumerate}

\item
Suppose that $(I_1\circ (1_{W_1}\otimes I_2))'(w'_{(4)})$ satisfies
Part (a) of the $P^{(2)}(z_0)$-local grading restriction condition,
that is, the $P^{(2)}(z_0)$-grading condition (or the
$L(0)$-semisimple $P^{(2)}(z_0)$-grading condition when $\mathcal{C}$
is in $\mathcal{M}_{sg}$).  For any $w_{(3)}\in W_{3}$, let
$\sum_{n\in \R}\lambda_{n}^{(2)}$ be a series weakly absolutely
convergent to
\[
\mu^{(2)}_{(I_1\circ (1_{W_1}\otimes I_2))'(w'_{(4)}), w_{(3)}} \in
(W_1 \otimes W_2)^*
\]
as indicated in the $P^{(2)}(z_0)$-grading condition (or the
$L(0)$-semisimple $P^{(2)}(z_0)$-grading condition), and suppose in
addition that the elements $\lambda_{n}^{(2)} \in (W_{1}\otimes
W_{2})^{*}$, $n\in \R$, satisfy the $P(z_{0})$-lower truncation
condition (Part (a) of the $P(z_{0})$-compatibility condition in
Section 5).  Then each $\lambda_{n}^{(2)}$ satisfies the (full)
$P(z_{0})$-compatibility condition.  Moreover, the corresponding space
\[
W^{(2)}_{(I_1\circ (1_{W_1}\otimes I_2))'(w'_{(4)}), w_{(3)}} \subset
(W_{1}\otimes W_{2})^{*},
\]
equipped with the vertex operator map given by $Y'_{P(z_{0})}$ and the
operators $L'_{P(z_{0})}(j)$ for $j=-1, 0, 1$, is a doubly-graded
generalized $V$-module, and when $\mathcal{C}$ is in
$\mathcal{M}_{sg}$, a doubly-graded $V$-module.  In particular, if
$(I_1\circ (1_{W_1}\otimes I_2))'(w'_{(4)})$ satisfies the full
$P^{(2)}(z_0)$-local grading restriction condition (or the
$L(0)$-semisimple $P^{(2)}(z_0)$-local grading restriction condition
when $\mathcal{C}$ is in $\mathcal{M}_{sg}$), then $W^{(2)}_{(I_1\circ
(1_{W_1}\otimes I_2))'(w'_{(4)}), w_{(3)}}$ is a generalized
$V$-module, that is, an object of
$\mathcal{GM}_{sg}$ (or a $V$-module, that is, an object of
$\mathcal{M}_{sg}$, when $\mathcal{C}$ is in
$\mathcal{M}_{sg}$); in this case, the assumption that each
$\lambda_{n}^{(2)}$ satisfies the $P(z_0)$-lower truncation condition
is redundant.

\item
Analogously, suppose that $(I^1\circ (I^2\otimes 1_{W_3}))'(w'_{(4)})$
satisfies Part (a) of the $P^{(1)}(z_2)$-local grading restriction
condition, that is, the $P^{(1)}(z_2)$-grading condition (or the
$L(0)$-semisimple $P^{(1)}(z_2)$-grading condition when $\mathcal{C}$
is in $\mathcal{M}_{sg}$).  For any $w_{(1)}\in W_{1}$, let
$\sum_{n\in \R}\lambda_{n}^{(1)}$ be a series weakly absolutely
convergent to
\[
\mu^{(1)}_{(I^1\circ (I^2\otimes 1_{W_3}))'(w'_{(4)}), w_{(1)}} \in
(W_2 \otimes W_3)^*
\]
as indicated in the $P^{(1)}(z_2)$-grading condition (or the
$L(0)$-semisimple $P^{(1)}(z_2)$-grading condition), and suppose in
addition that the elements $\lambda_{n}^{(1)} \in (W_{2}\otimes
W_{3})^{*}$, $n\in \R$, satisfy the $P(z_{2})$-lower truncation
condition (Part (a) of the $P(z_{2})$-compatibility condition).  Then
each $\lambda_{n}^{(1)}$ satisfies the (full) $P(z_{2})$-compatibility
condition.  Moreover, the corresponding space
\[
W^{(1)}_{(I^1\circ (I^2\otimes 1_{W_3}))'(w'_{(4)}), w_{(1)}} \subset
(W_{2}\otimes W_{3})^{*},
\]
equipped with the vertex operator map given by $Y'_{P(z_{2})}$ and the
operators $L'_{P(z_{2})}(j)$ for $j=-1, 0, 1$, is a doubly-graded
generalized $V$-module, and when $\mathcal{C}$ is in
$\mathcal{M}_{sg}$, a doubly-graded $V$-module.  In particular, if
$(I^1\circ (I^2\otimes 1_{W_3}))'(w'_{(4)})$ satisfies the full
$P^{(1)}(z_2)$-local grading restriction condition (or the
$L(0)$-semisimple $P^{(1)}(z_2)$-local grading restriction condition
when $\mathcal{C}$ is in $\mathcal{M}_{sg}$), then $W^{(1)}_{(I^1\circ
(I^2\otimes 1_{W_3}))'(w'_{(4)}), w_{(1)}}$ is a generalized
$V$-module, that is, an object of
$\mathcal{GM}_{sg}$ (or a $V$-module, that is, an object of
$\mathcal{M}_{sg}$, when $\mathcal{C}$ is in
$\mathcal{M}_{sg}$); in this case, the assumption that each
$\lambda_{n}^{(1)}$ satisfies the $P(z_2)$-lower truncation condition
is redundant.
\end{enumerate}
\end{theo}
\pf We will prove only Part 1 of the theorem, involving $I_1\circ
(1_{W_1}\otimes I_2)$; Part 2 is proved entirely analogously.

To prove that $W^{(2)}_{(I_1\circ (1_{W_1}\otimes I_2))'(w'_{(4)}),
w_{(3)}}$ is a doubly-graded generalized $V$-module (and when
$\mathcal{C}$ is in $\mathcal{M}_{sg}$, a doubly-graded $V$-module),
we claim that it is sufficient to prove that each $\lambda_{n}^{(2)}$,
$n\in \R$, satisfies Part (b) of the $P(z_{0})$-compatibility
condition (and hence the $P(z_{0})$-compatibility condition itself,
since these elements are assumed to satisfy Part (a)).  Indeed by
Lemma \ref{a-tilde-comp}, the space
\[
(\comp_{P(z_0)}((W_1\otimes W_2)^*)) \cap ((W_1\otimes
W_2)^*)^{(\tilde A)}
\]
is $\tilde{A}$-graded, and hence so is its intersection
\[
M=(\comp_{P(z_0)}((W_1\otimes W_2)^*)) \cap ((W_1\otimes
W_2)^*)_{[{\mathbb C}]}^{(\tilde A)}
\] 
with the $\tilde A$-graded space $((W_1\otimes W_2)^*)_{[{\mathbb
C}]}^{(\tilde A)}$.  By Theorem \ref{stable}, this space $M$ is also
$L'_{P(z_{0})}(0)$-stable and hence $\C$-graded and therefore doubly
graded.  By Theorem \ref{wk-mod} and Remark
\ref{stableundercomponentops}, $M$ is a weak $V$-module and hence in
fact a doubly-graded generalized $V$-module; when $\mathcal{C}$ is in
$\mathcal{M}_{sg}$, we replace the subscript $[\C]$ by $(\C)$, and $M$
is a doubly-graded $V$-module.  By our hypothesis that each
$\lambda_{n}^{(2)}$ satisfies the $P(z_{0})$-compatibility condition,
we have that $\lambda_{n}^{(2)} \in M$, and so
\[
W^{(2)}_{(I_1\circ (1_{W_1}\otimes I_2))'(w'_{(4)}), w_{(3)}} \subset
M,
\]
proving our claim.

The proof below of Part (b) of the $P(z_{0})$-compatibility condition
is a generalization of the proof of (14.51) in \cite{tensor4}.  The
proof here is (necessarily) much more elaborate.  When ${\cal C}$ is
in ${\cal M}_{sg}$, the proof below of course simplifies to a certain
extent, but even in this case, our setting is more general than that
in \cite{tensor4}, and the proof here is correspondingly more
delicate.

Let $\mathcal{Y}_{1}=\mathcal{Y}_{I_{1}, 0}$ and
$\mathcal{Y}_{2}=\mathcal{Y}_{I_{2},0}$ (recall Proposition
\ref{im:correspond}) so that 
\begin{eqnarray}
I_{1}(w_{(1)}\otimes w)&=&\Y_{1}(w_{(1)}, z_{1})w,\label{9.7-1--4}\\
I_{2}(w_{(2)}\otimes w_{(3)})&=&\Y_{2}(w_{(2)}, z_{2})w_{(3)}\label{9.7-1--3}
\end{eqnarray}
for $w_{(1)}\in W_{1}$, $w_{(2)}\in W_{2}$, $w_{(3)}\in W_{3}$ and
$w\in M_{1}$ (recall the ``substitution'' notation {}from
(\ref{im:f(z)}), where we choose $p=0$).  For $z\in \C^{\times}$, let
$I_{1}^{z}$ and $I_{2}^{z}$ be the $P(z_{0}+zz_{2})$- and
$P(zz_{2})$-intertwining maps $I_{\mathcal{Y}_{1}, 0}$ and
$I_{\mathcal{Y}_{2}, 0}$, respectively (assuming that $z_{0}+zz_{2}
\ne 0$), so that
\begin{eqnarray}
I^{z}_{1}(w_{(1)}\otimes w)&=&\Y_{1}(w_{(1)}, z_{0}+zz_{2})w,\label{9.7-1--2}\\
I^{z}_{2}(w_{(2)}\otimes w_{(3)})&=&\Y_{2}(w_{(2)}, zz_{2})w_{(3)}
\label{9.7-1--1}
\end{eqnarray}
for $w_{(1)}\in W_{1}$, $w_{(2)}\in W_{2}$, $w_{(3)}\in W_{3}$ and
$w\in M_{1}$; these maps are ``deformations'' of (\ref{9.7-1--4}) and
(\ref{9.7-1--3}), which correspond to $z=1$.  Since $|z_{1}|>|z_{2}|>|z_{0}|>0$,
there exists a sufficiently small neighborhood of 
$z=1$ such that 
\[
|z_{0}+zz_{2}|>|zz_{2}|>|z_{0}|>0
\]
(recall that 
\[
z_{1}=z_{0}+z_{2}).
\]
Since
\begin{equation}\label{z-prod}
\sum_{n\in \R}\langle w_{(4)}', I_{1}^{z}
(w_{(1)}\otimes \pi_{n}(I_{2}^{z}(w_{(2)}\otimes w_{(3)})))\rangle
=\sum_{n\in \R}\langle w_{(4)}', \Y_{1}(w_{(1)}, z_{0}+zz_{2})
\pi_{n}(\Y_{2}(w_{(2)}, zz_{2})w_{(3)})\rangle
\end{equation}
is absolutely convergent when $|z_{0}+zz_{2}|>|zz_{2}|>0$
for $w_{(1)}\in W_{1}$, 
$w_{(2)}\in W_{2}$, $w_{(3)}\in W_{3}$ and $w_{(4)}'\in W_{4}'$,
the product $I_{1}^{z}\circ (1_{W_{2}}\otimes I_{2}^{z})$ exists 
for $z$ in a sufficiently small neighborhood of $z=1$. 

We shall establish a relationship between 
(\ref{z-prod}) and a certain Taylor series expansion in $\log z$. 

The case $j=0$, $z=z_{0}$ and 
$\lambda=\mu^{(2)}_{(I_1\circ (1_{W_1}\otimes
I_2))'(w'_{(4)}), w_{(3)}}$
of (\ref{LP'(j)})
gives 
\begin{eqnarray}\label{9.7-1-0}
\lefteqn{(L'_{P(z_{0})}(0)\mu^{(2)}_{(I_1\circ (1_{W_1}\otimes
I_2))'(w'_{(4)}), w_{(3)}})(w_{(1)}\otimes w_{(2)})}\nno\\
&&=\mu^{(2)}_{(I_1\circ (1_{W_1}\otimes
I_2))'(w'_{(4)}), w_{(3)}}(w_{(1)}\otimes L(0)w_{(2)}
+(L(0)+z_{0}L(-1))w_{(1)}\otimes w_{(2)}).
\end{eqnarray}
Let $x$ be a formal variable. Then recalling the notation from
(\ref{im:f(z)}) with $p=0$ and Remark \ref{I1I2'}, and using
Definition \ref{productanditerateexisting}, (\ref{9.7-1--4}) and
(\ref{9.7-1--3}), we have
\begin{eqnarray}\label{14.43}
\lefteqn{((1-x)^{-L'_{P(z_{0})}(0)}
\mu^{(2)}_{(I_1\circ (1_{W_1}\otimes
I_2))'(w'_{(4)}), w_{(3)}})(w_{(1)}\otimes w_{(2)})}\nno\\
&&=(e^{-\log(1-x)L'_{P(z_{0})}(0)}
\mu^{(2)}_{(I_1\circ (1_{W_1}\otimes
I_2))'(w'_{(4)}), w_{(3)}})(w_{(1)}\otimes w_{(2)})\nno\\
&&=\mu^{(2)}_{(I_1\circ (1_{W_1}\otimes
I_2))'(w'_{(4)}), w_{(3)}}(e^{\log(1-x)(-z_{0}L(-1)-L(0))}w_{(1)}\otimes 
e^{-\log(1-x)L(0)}w_{(2)})\nno\\
&&=\mu^{(2)}_{(I_1\circ (1_{W_1}\otimes
I_2))'(w'_{(4)}), w_{(3)}}((1-x)^{-z_{0}L(-1)-L(0)}w_{(1)}\otimes 
(1-x)^{-L(0)}w_{(2)})\nno\\
&&=((I_1\circ (1_{W_1}\otimes
I_2))'(w'_{(4)}))((1-x)^{-z_{0}L(-1)-L(0)}
w_{(1)}\otimes 
(1-x)^{-L(0)}w_{(2)}\otimes w_{(3)})\nno\\
&&=\langle w_{(4)}', (I_1\circ (1_{W_1}\otimes
I_2))((1-x)^{-z_{0}L(-1)-L(0)}
w_{(1)}\otimes 
(1-x)^{-L(0)}w_{(2)}\otimes w_{(3)})\rangle\nno\\
&&=\langle w_{(4)}', I_1((1-x)^{-z_{0}L(-1)-L(0)}
w_{(1)}\otimes I_2((1-x)^{-L(0)}w_{(2)}\otimes w_{(3)})\rangle\nno\\
&&=\langle w'_{(4)}, {\cal Y}_{1}((1-x)^{-z_{0}L(-1)-L(0)}
w_{(1)}, x_{1}){\cal Y}_{2}(
(1-x)^{-L(0)}w_{(2)}, x_{2})
w_{(3)}\rangle_{W_{4}}\lbar_{x_{1}=z_{1}, \; x_{2}=z_{2}},\nno\\
&&
\end{eqnarray}
where, because of Proposition \ref{formal=proj}, the coefficient of
each power of $x$ on the right-hand side of (\ref{14.43}) has any of
the meanings discussed in Remark \ref{4notations}, and in particular,
means an absolutely convergent multisum or an analytic function of
$z_{1}$ and $z_{2}$.  The equality (\ref{14.43}) says that the
coefficient of each power of $x$ on the left-hand side of
(\ref{14.43}) is equal to such an absolutely convergent multisum or
such an analytic function obtained from the coefficient of the same
power of $x$ on the right-hand side.

Using Remark \ref{log:Lj2rema}, we have 
\begin{eqnarray*}
\lefteqn{{\cal Y}_{1}((1-x)^{-(x_{1}-x_{2})L(-1)-L(0)}
w_{(1)}, x_{1}){\cal Y}_{2}(
(1-x)^{-L(0)}w_{(2)}, x_{2})}\nn
&&=(1-x)^{-(L(0)-x_{2}L(-1))}{\cal Y}_{1}(
w_{(1)}, x_{1})
{\cal Y}_{2}(
w_{(2)}, x_{2})(1-x)^{L(0)-x_{2}L(-1)}.
\end{eqnarray*}
Then by Proposition \ref{formal=proj}, we see that
the right-hand side of (\ref{14.43}) is equal to
\begin{eqnarray}\label{14.44}
\lefteqn{\langle w'_{(4)}, (1-x)^{-(L(0)-x_{2}L(-1))}{\cal Y}_{1}(
w_{(1)}, x_{1})\cdot}\nno\\
&&\hspace{3em}\cdot{\cal Y}_{2}(
w_{(2)}, x_{2})(1-x)^{L(0)-x_{2}L(-1)}
w_{(3)}\rangle_{W_{4}}\lbar_{x_{1}=z_{1}, \; x_{2}=z_{2}}.
\end{eqnarray}
Lemma 9.3 in \cite{tensor2}, which used only the bracket formula for
$L(0)$ and $L(-1)$, gives the formula
\[
(1-x)^{L(0)-x_{2}L(-1)}=e^{x_{2}xL(-1)}(1-x)^{L(0)},
\]
and so by (\ref{log:p1}), (\ref{14.44}) is equal to
\begin{eqnarray}\label{14.45}
\lefteqn{\langle w'_{(4)}, (1-x)^{-L(0)}e^{-x_{2}xL(-1)}{\cal Y}_{1}(
w_{(1)}, x_{1})\cdot}\nno\\
&&\hspace{3em}\cdot {\cal Y}_{2}(
w_{(2)}, x_{2})e^{x_{2}xL(-1)}(1-x)^{L(0)}
w_{(3)}\rangle_{W_{4}}\lbar_{x_{1}=z_{1}, \; x_{2}=z_{2}}\nno\\
&&=\langle (1-x)^{-L'(0)}w'_{(4)}, {\cal Y}_{1}(
w_{(1)}, x_{1}-x_{2}x)\cdot\nno\\
&&\hspace{3em}\cdot{\cal Y}_{2}(
w_{(2)}, x_{2}-x_{2}x)(1-x)^{L(0)}
w_{(3)}\rangle_{W_{4}}\lbar_{x_{1}=z_{1}, \; x_{2}=z_{2}}
.
\end{eqnarray}
Thus the left-hand side of (\ref{14.43})
is equal to the right-hand side of (\ref{14.45}) (as formal power
series in $x$). 

Let 
\begin{equation}\label{l0z}
l^{0}(z)=\left\{\begin{array}{ll}\log z&0\le \arg z <\pi\\
\log z-2\pi i&\pi \le \arg z<2\pi\end{array}\right.,
\end{equation}
which is a single-valued branch of the logarithm of $z$ in the 
complex plane with a cut along the negative real line. We use this region
to choose a branch because we will need a single-valued branch such that 
$z=1$ is in the interior of the region.
By Proposition
\ref{analytic} and (\ref{eaL0-general}), 
\begin{equation}\label{g-zeta}
g(\zeta_{1}, \zeta_{2}, z)=\langle e^{-l^{0}(z)L'(0)}w'_{(4)}, {\cal Y}_{1}(
w_{(1)}, x_{1}){\cal Y}_{2}(
w_{(2)}, x_{2})e^{l^{0} (z)L(0)}
w_{(3)}\rangle_{W_{4}}\lbar_{x_{1}=\zeta_{1}, \; x_{2}=\zeta_{2}}
\end{equation}
is a single-valued function 
of $\zeta_{1}, \zeta_{2}$ and $z$ 
defined on  $|\zeta_{1}|>|\zeta_{2}|>0$, $\arg z\ne \pi$
and analytic when  $\arg \zeta_{1}, \zeta_{2}\ne 0$. 
The restriction of $g(\zeta_{1}, \zeta_{2}, z)$
to the subregion 
\[
|\zeta_{1}|>|\zeta_{2}|>0,\; 0\le \arg \zeta_{1},
\arg \zeta_{2}<\pi, \; \arg z\ne \pi
\]
can be analytically extended
(i) to a single-valued function $g_{1}(\zeta_{1}, \zeta_{2}, z)$
defined on the regions given by $|\zeta_{1}|>|\zeta_{2}|>0$, $\arg \zeta_{1}\ne \pi$,
$\arg z\ne \pi$ and analytic when 
$\arg \zeta_{2}\ne 0$;
(ii) to a single-valued function $g_{2}(\zeta_{1}, \zeta_{2}, z)$ 
defined on the region given by 
$|\zeta_{1}|>|\zeta_{2}|>0$, 
$\arg \zeta_{2}\ne \pi$, $\arg z\ne \pi$ and analytic 
when $\arg \zeta_{1}\ne 0$; and (iii) to a single-valued analytic function
$g_{3}(\zeta_{1}, \zeta_{2}, z)$
defined on the region given by 
$|\zeta_{1}|>|\zeta_{2}|>0$, $\arg \zeta_{1}, \arg \zeta_{2}\ne \pi$, 
$\arg z\ne \pi$ (recall Proposition
\ref{analytic}). 
For convenience, we shall use $h(\zeta_{1}, \zeta_{2}, z)$ to
denote $g(\zeta_{1}, \zeta_{2}, z)$ when $\arg z_{1}, \arg z_{2}\ne 0$,
to denote $g_{1}(\zeta_{1}, \zeta_{2}, z)$ when $\arg z_{1}=0$,
$\arg z_{2}\ne 0$, 
to denote $g_{2}(\zeta_{1}, \zeta_{2}, z)$ when $\arg z_{1}\ne 0$,
$\arg z_{2}=0$ and to denote $g_{3}(\zeta_{1}, \zeta_{2}, z)$
when $\arg z_{1}=\arg z_{2}= 0$.
Then in particular, $h(\zeta_{1}, \zeta_{2}, z)$ is analytic near
$\zeta_{1}=z_{1}$, $\zeta_{2}=z_{2}$ and $z=1$ and
we see that there exists a sufficiently small open
neighborhood of $z=1$ such that in this neighborhood, as the 
composition of the single-valued 
analytic function $h(\zeta_{1}, \zeta_{2}, z)$
of $\zeta_{1}, \zeta_{2}$ and $z$  with the analytic functions 
\[
\zeta_{1}=z_{0}+zz_{2}
\]
and 
\[
\zeta_{2}=zz_{2}
\]
of $z$, $h(z_{0}+zz_{2}, zz_{2}, z)$
is a single-valued analytic function 
of $z$.  For $z$ satisfying
\[
|z_{0}+zz_{2}|>|zz_{2}|>0,
\]
\begin{equation}\label{lessthanpi}
0 \le \arg (z_{0}+zz_{2}), \arg (zz_{2}) < \pi,
\end{equation}
\[
\arg z\ne \pi,
\]
by definition, we have
\begin{eqnarray}\label{9.7-1-1}
\lefteqn{h(z_{0}+zz_{2}, zz_{2}, z)}\nn
&&=g(z_{0}+zz_{2}, zz_{2}, z)\nno\\
&&=\langle e^{-l^{0}(z)L'(0)}w'_{(4)}, {\cal Y}_{1}(
w_{(1)}, x_{1}){\cal Y}_{2}(
w_{(2)}, x_{2})e^{l^{0}(z)L(0)}
w_{(3)}\rangle_{W_{4}}\lbar_{x_{1}=z_{0}+zz_{2}, \; x_{2}=zz_{2}};
\end{eqnarray}
moreover, $\pi$ can be increased to $2\pi$ in (\ref{lessthanpi}) if
$\arg z_{1}$ and/or $\arg z_{2}$ is positive, according to the cases
discussed above.

In the region $\arg z\ne \pi$, the analytic function 
\[
z'=l^{0}(z)
\]
is single-valued and univalent and $z'=0$ is in the image
of the region. So the composition of the function $h(z_{0}+zz_{2}, zz_{2}, z)$  
with the inverse function 
\[
z=e^{z'}
\]
of the function $z'=l^{0}(z)$ gives us a single-valued analytic
function
\begin{equation}\label{def-f-z'}
f(z')=h(z_{0}+e^{z'}z_{2}, e^{z'}z_{2}, e^{z'})
\end{equation}
of $z'$ in a sufficiently small open neighborhood of $z'=0$.  In
particular, we can expand $f(z')$ as a power series in $z'$.  Since
the power series expansion of any function analytic at $z'=0$ is
uniquely determined by its derivatives at $z'=0$, we can find the
power series expansion of $f(z')$ as follows: Since $h(z_{0}+zz_{2},
zz_{2}, z)$ is analytic at $z=1$, we first expand it as a power series
in $z-1$ in a sufficiently small open neighborhood of $z=1$.  Then the
power series expansion of $f(z')$ in a sufficiently small open
neighborhood $U$ of $z'=0$ is obtained by replacing each nonnegative
integral power of $z-1$ by the corresponding power of $\sum_{k\in
\Z_{+}}\frac{(z')^{k}}{k!}$. Since the convergence of the power series
expansion of $f(z')$ is independent of $w_{(1)}\in W_{1}$, $w_{(2)}\in
W_{2}$, $w_{(3)}\in W_{3}$ and $w'_{(4)}\in W_{4}'$, we can choose $U$
to be independent of these elements.

We now want to give this power series explicitly (first, in powers of
$z-1$) using the right-hand side of (\ref{9.7-1-1}). To do this, we
have to restrict $z$ to be in a subset such that $h(z_{0}+zz_{2},
zz_{2}, z)$ is equal to the right-hand side of (\ref{9.7-1-1}).

Let 
\[
O=e^{U}.
\]
Then $O$ is an open subset, containing $1$, of the domain of the
analytic function $h(z_{0}+zz_{2}, zz_{2}, z)$ of $z$, and is
independent of $w_{(1)}$, $w_{(2)}$, $w_{(3)}$ and $w'_{(4)}$.  We
choose $U$ to be small enough so that the power series expansion of
$h(z_{0}+zz_{2}, zz_{2}, z)$ near $z=1$ is absolutely convergent for
$z\in O$.

Let $P$ be (i) the set of $z \in \C$ with
\begin{equation}\label{0arg2pi}
0\le \arg (z_{0}+zz_{2}), 
\arg (zz_{2})< 2\pi
\end{equation}
when $\arg z_{1}, \arg z_{2}\ne 0$ (so that in this case, $P$ simply
equals $\C$); or (ii) the set of $z$ with
\[
0\le \arg (z_{0}+zz_{2})<\pi, \;
0\le \arg (zz_{2})<2\pi
\]
when $\arg z_{1}=0$ but $\arg z_{2}\ne 0$; 
or (iii) the set of $z$ with
\[
0\le \arg (z_{0}+zz_{2})<2\pi, \;
0\le \arg (zz_{2})<\pi
\]
when $\arg z_{1}\ne 0$ but $\arg z_{2}= 0$;
or (iv) the set of $z$ with
\[
0\le \arg (z_{0}+zz_{2}), \arg (zz_{2})<\pi
\]
when $\arg z_{1}=\arg z_{2}= 0$.
Then by definition, for
\[
z \in O\cap P,
\]
we have
\begin{equation}\label{h=g}
h(z_{0}+zz_{2}, zz_{2}, z)=g(z_{0}+zz_{2}, zz_{2}, z),
\end{equation}
which is given by the right-hand side of (\ref{9.7-1-1}).  Note that
\[
1\in O\cap P.
\]

For
\[
z'=l^{0}(z)\in U
\]
($z\in O$), the power series expansion of $f(z')$ can be obtained
explicitly from the right-hand side of (\ref{9.7-1-1}) as follows:

By Propositions \ref{exp-set} and \ref{formal=proj} (and in
particular, (\ref{triple-sum})), (\ref{eaL0-general}) and Remark
\ref{4notations}, one of the equivalent meanings of (\ref{g-zeta}) is
an absolutely convergent series in the region
$|\zeta_{1}|>|\zeta_{2}|>0$ of the form
\begin{equation}\label{6-tuple-series-1}
\sum_{s=0}^{S}\sum_{t=0}^{T}\sum_{p, q\in \R}\sum_{k=0}^{K}\sum_{l=0}^{L}
b_{s,t,p,q,k,l}e^{p\log \zeta_{1}}(\log \zeta_{1})^{k}
e^{q\log \zeta_{2}}(\log \zeta_{2})^{l}e^{c_{s}l^{0}(z)}(l^{0}(z))^{t},
\end{equation}
where $b_{s,t,p,q,k,l}\in \C$ and $c_{s}\in \R$.  By
Lemma \ref{po-ser-an}, the derivatives of (\ref{6-tuple-series-1}) are given by 
the absolutely convergent series 
\begin{equation}\label{6-tuple-series-1.3}
\sum_{s=0}^{S}\sum_{t=0}^{T}\sum_{p, q\in \R}\sum_{k=0}^{K}\sum_{l=0}^{L}
b_{s,t,p,q,k,l}
\frac{\partial^{i}}{\partial \zeta_{1}^{i}}\frac{\partial^{j}}{\partial \zeta_{2}^{j}}
\frac{\partial^{n}}{\partial z^{n}}\left(e^{p\log \zeta_{1}}(\log \zeta_{1})^{k}
e^{q\log \zeta_{2}}(\log \zeta_{2})^{l}e^{c_{s}l^{0}(z)}(l^{0}(z))^{t}\right)
\end{equation}
for $i, j, n\in \N$. 

On the other hand, the coefficients of the expansion of
$h(z_{0}+zz_{2}, zz_{2}, z)$ as a power series in $z-1$ are given by
its derivatives at $z=1$. By the chain rule, there exist $\alpha_{m,
i, j, n}\in \C$ (depending on $z_{2}$) such that for any analytic
function $F(\zeta_{1}, \zeta_{2}, z)$ of $\zeta_{1}, \zeta_{2}, z$
near $\zeta_{1}=z_{1}$, $\zeta_{2}=z_{2}$ and $z=1$ (in particular,
for $F(\zeta_{1}, \zeta_{2}, z) =h(\zeta_{1}, \zeta_{2}, z)$),
\begin{equation}\label{6-tuple-series-1.7}
\frac{\partial^{m}}{\partial z^{m}}
F(z_{0}+zz_{2}, zz_{2}, z)\lbar_{z=1}
=\sum_{i+j+n=m, \;i, j, n\in \N}\alpha_{m, i, j, n}
\frac{\partial^{i}}{\partial \zeta_{1}^{i}}\frac{\partial^{j}}{\partial \zeta_{2}^{j}}
\frac{\partial^{n}}{\partial z^{n}}
F(\zeta_{1}, \zeta_{2}, z)\lbar_{\zeta_{1}=z_{1},\;\zeta_{2}=z_{2},\; z=1}.
\end{equation}
For $z\in O\cap P$, $h(z_{0}+zz_{2}, zz_{2}, z)$ is equal to the
right-hand side of (\ref{9.7-1-1}). But one of the meanings of the
right-hand side of (\ref{9.7-1-1}) is the absolutely convergent series 
\begin{equation}\label{6-tuple-series}
\sum_{s=0}^{S}\sum_{t=0}^{T}\sum_{p, q\in \R}\sum_{k=0}^{K}\sum_{l=0}^{L}
b_{s, t, p, q, k, l}e^{p\log(z_{0}+zz_{2})}(\log(z_{0}+zz_{2}))^{k}
e^{q\log(zz_{2})}(\log(zz_{2}))^{l}e^{c_{s}l^{0}(z)}(l^{0}(z))^{t}
\end{equation}
(see (\ref{g-zeta}) and (\ref{6-tuple-series-1})).
Using (\ref{6-tuple-series-1.3}), (\ref{6-tuple-series-1.7}) and 
(\ref{6-tuple-series}), we see that 
the $m$-th derivative of $h(z_{0}+zz_{2}, zz_{2}, z)$
at $z=1$ is equal to the absolutely convergent series
\begin{eqnarray}\label{6-tuple-series-2}
\lefteqn{\sum_{s=0}^{S}\sum_{t=0}^{T}\sum_{p, q\in \R}\sum_{k=0}^{K}\sum_{l=0}^{L}
b_{s, t, p, q, k, l}\sum_{i+j+n=m, \;i, j, n\in \N}\alpha_{m, i, j, n}
\cdot}\nn
&&\cdot\frac{\partial^{i}}{\partial \zeta_{1}^{i}}\frac{\partial^{j}}{\partial \zeta_{2}^{j}}
\frac{\partial^{n}}{\partial z^{n}}\left(e^{p\log \zeta_{1}}(\log \zeta_{1})^{k}
e^{q\log \zeta_{2}}(\log \zeta_{2})^{l}e^{c_{s}l^{0}(z)}(l^{0}(z))^{t}\right)
\lbar_{\zeta_{1}=z_{1},\;\zeta_{2}=z_{2},\;z=1}
\end{eqnarray}
for $m\in \N$. (Note that when $\arg z_{1}$ or $\arg z_{2}$ is $0$,
$1$ is not in the interior of $O\cap P$ and we have to calculate the
derivatives above using only $z\in O\cap P$. But the result is the
same.)  Thus we see that the coefficients of the expansion of
$h(z_{0}+zz_{2}, zz_{2}, z)$ as a power series in $z-1$ are given by
(\ref{6-tuple-series-2}) divided by $m!$.

By (\ref{6-tuple-series-1.7}) and (\ref{6-tuple-series}), we see that
the coefficients of this power series in $z-1$ are also equal to the
coefficients of the power series in $z-1$ obtained from
(\ref{6-tuple-series}) by replacing
\begin{eqnarray}
e^{p\log(z_{0}+zz_{2})}(\log(z_{0}+zz_{2}))^{k}
&=&e^{p\log(z_{1}+(z-1)z_{2})}(\log(z_{1}+(z-1)z_{2}))^{k},\label{monomia-1}\nn
e^{q\log(zz_{2})}(\log(zz_{2}))^{l}
&=&e^{q\log(z_{2}+(z-1)z_{2})}(\log(z_{2}+(z-1)z_{2}))^{l},\nn
e^{c_{s}l^{0}(z)}(l^{0}(z))^{t}&=&e^{c_{s}l^{0}(1+(z-1))}(l^{0}(1+(z-1)))^{t}\label{monomia-2}
\end{eqnarray}
by their power series expansions near $z=1$. We have shown that the
power series obtained in this way has the indicated sums of absolutely
convergent series as coefficients and is absolutely convergent to
$h(z_{0}+zz_{2}, zz_{2}, z)$ when $z\in O$ (since we chose $U$ and $O$
small enough, above).  We have succeeded in giving this power series
in $z-1$ explicitly.  Finally, as above, we replace each nonnegative
integral power of $z-1$ by the corresponding power of $\sum_{k\in
\Z_{+}}\frac{(z')^{k}}{k!}$ to obtain the power series expansion of
$f(z')$ for $z'\in U$.

Note that the constant terms of both the power series expansion of 
\[
l^{0}(z)=l^{0}(1+(z-1))
\]
near $z=1$ and the formal power series $\log (1+x)$ are $0$, and in
fact, this expansion of $l^{0}(z)$ is obtained by substituting $z-1$
for $x$ in the formal series $\log (1+x)$.  (This is of course a
reflection of the fact that the formal power series notation ``$\log
(1+x)$'' is in effect choosing a branch of a multivalued function.)
Thus from the explicit expansion procedure obtaining the power series
in $z-1$ above and the precise meaning of the right-hand side of
(\ref{14.45}), we see that, as sums of absolutely convergent series,
the coefficient of the $n$-th power of $x$ in the formal power series
in $x$ given by the right-hand side of (\ref{14.45}) is exactly the
same as $(-1)^{n}$ times the coefficient of the $n$-th power of $z-1$
in the power series in $z-1$ obtained above. Thus, if we substitute
$-(z-1)$ for $x$ in the right-hand side of (\ref{14.45}), we obtain
the explicit expansion above of the right-hand side of (\ref{9.7-1-1})
as a power series in $z-1$.

Because of the explicit calculations and discussions above, when
$z'\in U$, the power series expansion of $f(z')$ can be obtained using
the right-hand side of (\ref{14.45}) and the following two steps:
(i) Substitute $1-e^{y}\in -y+y^{2}\C[[y]]$ for $x$ in the right-hand
side of (\ref{14.45}) and
(ii) substitute $z'$ for $y$ in the resulting series.  This power
series in $z'$ as the expansion of $f(z')$ must be absolutely
convergent in the neighborhood $U$ of $z'=0$ and its sum is equal to
the single-valued analytic function $h(z_{0}+e^{z'}z_{2}, e^{z'}z_{2},
e^{z'})$.  In particular, for $z'\in U$ and $z\in P$, this power
series in $z'=l^{0}(z)$ is absolutely convergent to the right-hand
side of (\ref{9.7-1-1}).
 
Applying the same steps (i) and (ii) above to the left-hand side of
(\ref{14.43}), we also obtain a power series $S(z')$ in $z'$.  Since
the left-hand side of (\ref{14.43}) is equal to the right-hand side of
(\ref{14.45}) as formal power series in $x$, we see that the power
series expansion of $f(z')$ and the power series $S(z')$ are the
same. In particular, in the neighborhood $U$ of $z'=0$, $S(z')$ is
absolutely convergent to $f(z')$. Since $f(z')$ is equal to the
right-hand side of (\ref{9.7-1-1}) when $z'=l^{0}(z)\in U$ and $z\in
P$, $S(z')$ is absolutely convergent to the right-hand side of
(\ref{9.7-1-1}) when $z'=l^{0}(z)\in U$ and $z\in P$, that is, when
$z\in O\cap P$. (Recall from (\ref{0arg2pi}) that in case $\arg
z_{1}\ne 0$ and $\arg z_{2}\ne 0$, $P=\C$, so that in this case,
$S(z')$ is absolutely convergent to the right-hand side of
(\ref{9.7-1-1}) whenever $z'\in U$.)
 
By assumption, for $w_{(1)}\in W_{1}$ and $w_{(2)}\in W_{2}$, the
series
\begin{equation}\label{sumlambda}
\sum_{n\in \R}\lambda^{(2)}_{n}(w_{(1)}\otimes w_{(2)})
\end{equation}
converges absolutely to
\[
\mu^{(2)}_{(I_1\circ (1_{W_1}\otimes
I_2))'(w'_{(4)}), w_{(3)}}(w_{(1)}\otimes w_{(2)}),
\]
and 
\[
\sum_{n\in \R}(e^{z'L'_{P(z_{0})}(0)}
\lambda^{(2)}_{n})(w_{(1)}\otimes w_{(2)})
\]
is absolutely convergent for $z'$ in an open neighborhood of $z'=0$
independent of $w_{(1)}$ and  $w_{(2)}$.  (Note that the elements
$\lambda^{(2)}_{n}$ for $n \in \R$ depend on $w'_{(4)} \in W'_4$ and
$w_{(3)} \in W_3$.)

Let
\[
Q=\{z\in \C\;|\;0\le \arg z<\pi\}.
\]
Note that for  $z\in Q$, 
\[
l^{0}(z)=\log z. 
\]
We will show that
the series 
\begin{equation}\label{iter-sum-0}
\sum_{n\in \R}(e^{-(l^{0}(z))L'_{P(z_{0})}(0)}\lambda^{(2)}_{n})(w_{(1)}\otimes 
w_{(2)}),
\end{equation}
which is absolutely convergent in an open neighborhood of $z=1$ independent of
$w_{(1)}$ and  $w_{(2)}$ and absolutely convergent to
\[
\mu^{(2)}_{(I_1\circ (1_{W_1}\otimes I_2))'(w'_{(4)}),
w_{(3)}}(w_{(1)}\otimes w_{(2)})
\]
for $z=1$, gives a double series of the form $\sum_{n\in
\R}\sum_{i=0}^{N}a_{n, i}e^{-n\log z}(-\log z)^{i}$ absolutely
convergent to $f(\log z)$ for $z$ in a nonempty open subset of $O\cap
P\cap Q$, and since by Proposition \ref{real-exp-set} $\R\times \{0,
\dots, N\}$ is a unique expansion set, the coefficients $a_{n, i}$ and
related numbers are uniquely determined.

Since 
\[
\lambda_{n}^{(2)}\in \coprod_{\beta\in \tilde{A}}((W_{1}\otimes
W_{2})^{*})_{[n]}^{(\beta)},
\]
we have 
\[
L'_{P(z_{0})}(0)_{s}\lambda_{n}^{(2)}=n\lambda_{n}^{(2)}
\] 
(recall Remarks \ref{set:L(0)s} and \ref{rmk-9.5}).

In the case that $\mathcal{C}$ is in $\mathcal{M}_{sg}$,
$(I_{1}\circ (1_{W_{1}}\otimes I_{2}))'(w'_{(4)})$ satisfies 
the $L(0)$-semisimple $P^{(2)}(z_{0})$-grading condition, so that
\[
\lambda_{n}^{(2)}\in \coprod_{\beta\in \tilde{A}}((W_{1}\otimes
W_{2})^{*})_{(n)}^{(\beta)},
\]
and we have
\[
L'_{P(z_{0})}(0)_{s}=L'_{P(z_{0})}(0)
\]
and 
\[
L'_{P(z_{0})}(0)\lambda_{n}^{(2)}
=n\lambda_{n}^{(2)}.
\] 
In particular, 
the proof below will give the desired result,
and is in fact simpler, in this case.

{}From the $P^{(2)}(z_{0})$-grading condition, we have, for $z$ in an
open neighborhood of $1$ independent of $w_{(1)}\in W_{1}$ and
$w_{(2)}\in W_{2}$, and with $N$ independent of $w_{(1)}$ and
$w_{(2)}$,
\begin{eqnarray}\label{iter-sum}
\lefteqn{\sum_{n\in \R}(e^{-(l^{0}(z))L'_{P(z_{0})}(0)}\lambda^{(2)}_{n})
(w_{(1)}\otimes 
w_{(2)})}\nn
&&=\sum_{n\in \R}e^{-n(l^{0}(z))}\left(\left(\sum_{i=0}^{N}\frac{(-l^{0}(z))^{i}}{i!}
(L'_{P(z_{0})}(0)-L'_{P(z_{0})}(0)_{s})^{i}\lambda^{(2)}_{n}
\right)(w_{(1)}\otimes 
w_{(2)})\right).\nn
\end{eqnarray}

The derivative with respect to $z'=l^{0}(z)$ of the iterated series 
(\ref{iter-sum}) is
\begin{eqnarray}\label{iter-sum-der}
&&\sum_{n\in \R}\frac{\partial}{\partial l^{0}(z)}
(e^{-(l^{0}(z))L'_{P(z_{0})}(0)}\lambda^{(2)}_{n})
(w_{(1)}\otimes w_{(2)})\nn
&&=\sum_{n\in \R}\frac{\partial}{\partial l^{0}(z)}
\left(e^{-n(l^{0}(z))}\left(\left(\sum_{i=0}^{N}\frac{(-l^{0}(z))^{i}}{i!}
(L'_{P(z_{0})}(0)-L'_{P(z_{0})}(0)_{s})^{i}\lambda^{(2)}_{n}
\right)(w_{(1)}\otimes 
w_{(2)})\right)\right)\nn
&&=\sum_{n\in \R}
(-n)e^{-n(l^{0}(z))}\left(\left(\sum_{i=0}^{N}\frac{(-l^{0}(z))^{i}}{i!}
(L'_{P(z_{0})}(0)-L'_{P(z_{0})}(0)_{s})^{i}\lambda^{(2)}_{n}\right)
(w_{(1)}\otimes 
w_{(2)})\right)\nn
&&\quad+\sum_{n\in \R}
e^{-n(l^{0}(z))}\left(\left(\sum_{i=1}^{N}\frac{-(-l^{0}(z))^{i-1}}{(i-1)!}
(L'_{P(z_{0})}(0)-L'_{P(z_{0})}(0)_{s})^{i}\lambda^{(2)}_{n}\right)
(w_{(1)}\otimes 
w_{(2)})\right)\nn
&&=-\sum_{n\in \R}
e^{-n(l^{0}(z))}\cdot\nn
&&\quad\quad\quad\quad\cdot 
\left(\left(\sum_{i=0}^{N}\frac{(-l^{0}(z))^{i}}{i!}
(L'_{P(z_{0})}(0)-L'_{P(z_{0})}(0)_{s})^{i}L'_{P(z_{0})}(0)_{s}\lambda^{(2)}_{n}
\right)(w_{(1)}\otimes 
w_{(2)})\right)\nn
&&\quad-\sum_{n\in \R}
e^{-n(l^{0}(z))}\left(\left(\sum_{i=0}^{N}\frac{(-l^{0}(z))^{i}}{i!}
(L'_{P(z_{0})}(0)-L'_{P(z_{0})}(0)_{s})^{i+1}
\lambda^{(2)}_{n}
\right)(w_{(1)}\otimes 
w_{(2)})\right)\nn
&&=-\sum_{n\in \R}
e^{-n(l^{0}(z))}\cdot\nn
&&\quad\quad\quad\quad\cdot \left(\left(\sum_{i=0}^{N}\frac{(-l^{0}(z))^{i}}{i!}
L'_{P(z_{0})}(0)(L'_{P(z_{0})}(0)-L'_{P(z_{0})}(0)_{s})^{i}\lambda^{(2)}_{n}
\right)(w_{(1)}\otimes 
w_{(2)})\right)\nn
&&=-\sum_{n\in \R}(L'_{P(z_{0})}(0)e^{-(l^{0}(z))L'_{P(z_{0})}(0)}
\lambda^{(2)}_{n})
(w_{(1)}\otimes w_{(2)})\nn
&&=-\sum_{n\in \R}(e^{-(l^{0}(z))L'_{P(z_{0})}(0)}
\lambda^{(2)}_{n})
(w_{(1)}\otimes L(0)w_{(2)}+(L(0)+z_{0}L(-1))w_{(1)}\otimes w_{(2)})
\end{eqnarray}
(recall (\ref{LP'(j)})). Since (\ref{iter-sum-0}) is absolutely
convergent for $z'=l^{0}(z)$ in an open neighborhood of $z'=0$
independent of $w_{(1)}\in W_{1}$ and $w_{(2)}\in W_{2}$, so is the
left-hand side of (\ref{iter-sum}).  Thus the right-hand side of
(\ref{iter-sum-der}) and consequently the left-hand side of
(\ref{iter-sum-der}) is absolutely convergent for $z'=l^{0}(z)$ in the
same neighborhood of $z'=0$.  Since the map $l^{0}$ is univalent in a
neighborhood of $z=1$, we see that there exists an open neighborhood
$\Pi$ of $z=1$ independent of $w_{(1)}\in W_{1}$ and $w_{(2)}\in
W_{2}$ such that both sides of (\ref{iter-sum}) and of
(\ref{iter-sum-der}) are absolutely convergent.  For later use, we may
and do choose $\Pi$ to be a small open disk centered at $1$.  The same
calculation and argument show that all the higher derivatives with
respect to $z'=l^{0}(z)$ of the iterated series (\ref{iter-sum}) are
also absolutely convergent for $z\in \Pi$.

Since $l^{0}(z)=\log z$
for $z\in Q$, 
we see by Proposition 
\ref{log-coeff-conv<=>iterate-conv} that 
\[
\sum_{n\in \R}e^{-n(\log z)}
((L'_{P(z_{0})}(0)-L'_{P(z_{0})}(0)_{s})^{i}\lambda^{(2)}_{n})(w_{(1)}\otimes w_{(2)})
\]
is absolutely convergent for $z\in \Pi\cap Q$, for each
$i=0, \dots, N$. Since $n \in \R$, we have
\[
|e^{-n(\log z)}|=|e^{-n(l^{0}(z))}|=|e^{-n(l^{0}(\overline{z}))}|
\]
for $z\in \C^{\times}$, and since $\overline{\Pi}=\Pi$, for $z\in \Pi$, 
\begin{equation}\label{coeff-log-sum}
\sum_{n\in \R}e^{-n(l^{0}(z))}
((L'_{P(z_{0})}(0)-L'_{P(z_{0})}(0)_{s})^{i}\lambda^{(2)}_{n})(w_{(1)}\otimes w_{(2)})
\end{equation}
is absolutely convergent for $i=0, \dots, N$. Thus
the double series
\begin{equation}\label{double-sum}
\sum_{n\in \R}\sum_{i=0}^{N}e^{-n(l^{0}(z))}\frac{(-l^{0}(z))^{i}}{i!}
((L'_{P(z_{0})}(0)-L'_{P(z_{0})}(0)_{s})^{i}\lambda^{(2)}_{n})(w_{(1)}\otimes w_{(2)})
\end{equation}
is absolutely convergent for $z\in \Pi$. Since for any $z\in \C^{\times}$,
\[
\sum_{n\in \R}e^{-n(l^{0}(z))}
((L'_{P(z_{0})}(0)-L'_{P(z_{0})}(0)_{s})^{i}\lambda^{(2)}_{n})
(w_{(1)}\otimes w_{(2)})
\]
can be written as a series of the form of $\sum_{n\in
\R}a_{n}e^{-n(\log z)}$, by Lemma \ref{po-ser-an} the sums of
(\ref{coeff-log-sum}) give analytic functions of $z\in \Pi$ for $i=0,
\dots, N$. Thus the sum of (\ref{iter-sum}), or equivalently, the sum
of (\ref{double-sum}), gives an analytic function of $z\in \Pi$. Using
(\ref{iter-sum-der}) repeatedly, we see that for $k\in \N$, the $k$-th
derivative with respective to $z'=l^{0}(z)$ of this analytic function
is given by the absolutely convergent series
\begin{equation}\label{k-th-der}
(-1)^{k}\sum_{n\in
\R}((L'_{P(z_{0})}(0))^{k}e^{-(l^{0}(z))L'_{P(z_{0})}(0)}\lambda^{(2)}_{n})
(w_{(1)}\otimes w_{(2)})
\end{equation}
for $z \in \Pi$, and its $k$-th derivative with respect to $z'$ at
$z'=0$ or equivalently at $z=1$ is given by the absolutely convergent
series
\begin{equation}\label{k-th-der-at-0}
(-1)^{k}\sum_{n\in \R}((L'_{P(z_{0})}(0))^{k}\lambda^{(2)}_{n})
(w_{(1)}\otimes w_{(2)}).
\end{equation}
 
Since 
\[
\sum_{n\in \R}\lambda^{(2)}_{n}
\]
is weakly absolutely convergent to 
\[
\mu^{(2)}_{(I_1\circ (1_{W_1}\otimes
I_2))'(w'_{(4)}), w_{(3)}},
\]
using (\ref{9.7-1-0}) we have, starting as in (\ref{14.43}),
\begin{eqnarray}\label{9.7-1-2}
\lefteqn{((1-x)^{-L'_{P(z_{0})}(0)}
\mu^{(2)}_{(I_1\circ (1_{W_1}\otimes
I_2))'(w'_{(4)}), w_{(3)}})(w_{(1)}\otimes w_{(2)})}\nno\\
&&=(e^{-\log (1-x) L'_{P(z_{0})}(0)}
\mu^{(2)}_{(I_1\circ (1_{W_1}\otimes
I_2))'(w'_{(4)}), w_{(3)}})(w_{(1)}\otimes w_{(2)})\nno\\
&&=\mu^{(2)}_{(I_1\circ (1_{W_1}\otimes
I_2))'(w'_{(4)}), w_{(3)}}(e^{\log (1-x)(-z_{0}L(-1)-L(0))}w_{(1)}\otimes 
e^{-\log (1-x)L(0)}w_{(2)})\nno\\
&&=\sum_{n\in \R}\lambda^{(2)}_{n}
(e^{\log (1-x)(-z_{0}L(-1)-L(0))}w_{(1)}\otimes 
e^{-\log (1-x)L(0)}w_{(2)})\nno\\
&&=\sum_{n\in \R}(e^{-\log (1-x)L'_{P(z_{0})}(0)}\lambda^{(2)}_{n})
(w_{(1)}\otimes w_{(2)})\nno\\
&&=\sum_{n\in \R}((1-x)^{-L'_{P(z_{0})}(0)}\lambda^{(2)}_{n})
(w_{(1)}\otimes w_{(2)}),
\end{eqnarray}
where the absolute convergence holds for the coefficient of each power
of $x$ in the formal power series in $x$ in (\ref{9.7-1-2}).
We also have
\begin{eqnarray}\label{9.7-1-3}
\lefteqn{(1-x)^{-L'_{P(z_{0})}(0)}\lambda^{(2)}_{n}}\nno\\
&&=
e^{-\log(1-x)L'_{P(z_{0})}(0)}\lambda^{(2)}_{n}\nno\\
&&=
e^{-\log (1-x)(L'_{P(z_{0})}(0)-L'_{P(z_{0})}(0)_{s})}
e^{-\log (1-x)L'_{P(z_{0})}(0)_{s}}\lambda^{(2)}_{n}\nno\\
&&=e^{-n\log (1-x)}\sum_{i=0}^{K_{n}}\frac{(-\log (1-x))^{i}}{i!}
(L'_{P(z_{0})}(0)-L'_{P(z_{0})}(0)_{s})^{i}\lambda^{(2)}_{n}\nno\\
&&=(1-x)^{-n}\sum_{i=0}^{K_{n}}\frac{(-\log (1-x))^{i}}{i!}
(L'_{P(z_{0})}(0)-n)^{i}\lambda^{(2)}_{n}
\end{eqnarray}
where $K_n \in \N$; cf. (\ref{e-y-2}).
{}From  (\ref{9.7-1-2}) and (\ref{9.7-1-3}), we obtain
\begin{eqnarray}\label{9.7-1-4}
\lefteqn{((1-x)^{-L'_{P(z_{0})}(0)}
\mu^{(2)}_{(I_1\circ (1_{W_1}\otimes
I_2))'(w'_{(4)}), w_{(3)}})(w_{(1)}\otimes w_{(2)})}\nn
&&=\sum_{n\in \R}((1-x)^{-L'_{P(z_{0})}(0)}\lambda^{(2)}_{n})
(w_{(1)}\otimes w_{(2)})\nno\\
&&=\sum_{n\in \R}(1-x)^{-n}\sum_{i=0}^{K_{n}}\frac{(-\log(1-x))^{i}}{i!}
((L'_{P(z_{0})}(0)-n)^{i}\lambda^{(2)}_{n})
(w_{(1)}\otimes w_{(2)}),\nno\\
&&
\end{eqnarray}
with absolute convergence for each power of $x$, as in (\ref{9.7-1-2}).

Recall that we have proved that if we substitute $1-e^{y}$ for $x$ in
the left-hand side of (\ref{14.43}), which is the same as the
left-hand side of (\ref{9.7-1-4}), and then substitute $l^{0}(z)$ for
$y$, we obtain an absolutely convergent power series $S(l^{0}(z))$ in
$l^{0}(z)$ for $z\in O$, and for $z\in O\cap Q$, so that
$l^{0}(z)=\log z$, the sum of this series $S(\log z)$ is equal to
$f(\log z)$, that is, if we also use $S(\log z)$ to denote its sum,
then
\[
S(\log z)=f(\log z)
\]
(recall (\ref{def-f-z'})).
Moreover, for  
\[
z\in O\cap P\cap Q,
\]
this is also equal to the right-hand side of (\ref{9.7-1-1}) (recall
(\ref{h=g})).

The same substitution steps in the right-hand side of (\ref{9.7-1-4})
give the same absolutely convergent series $S(l^{0}(z))$.
Substituting $1-e^{y}$ for $x$ in (\ref{9.7-1-4}) and using Remark
\ref{part-a}, we obtain
\begin{eqnarray}\label{9.7-1-4-y}
\lefteqn{(e^{-yL'_{P(z_{0})}(0)}
\mu^{(2)}_{(I_1\circ (1_{W_1}\otimes
I_2))'(w'_{(4)}), w_{(3)}})(w_{(1)}\otimes w_{(2)})}\nn
&&=\sum_{n\in \R}(e^{-yL'_{P(z_{0})}(0)}\lambda^{(2)}_{n})
(w_{(1)}\otimes w_{(2)})\nno\\
&&=\sum_{n\in \R}e^{-ny}\sum_{i=0}^{N}\frac{(-y)^{i}}{i!}
((L'_{P(z_{0})}(0)-n)^{i}\lambda^{(2)}_{n})
(w_{(1)}\otimes w_{(2)}),
\end{eqnarray}
where the absolute convergence holds for the coefficient of each power
of $y$ in the formal power series in $y$ in (\ref{9.7-1-4-y}).  Thus
$S(l^{0}(z))$ is equal to the series obtained by substituting
$l^{0}(z)$ for $y$ in (\ref{9.7-1-4-y}).  In particular, for $k\in
\N$, the $k$-th derivative with respect to $l^{0}(z)$ of $S(l^{0}(z))$
at $l^{0}(z)=0$ (or equivalently at $z=1$) is equal to the constant
term of the $k$-th derivative with respect to $y$ of
(\ref{9.7-1-4-y}), and this is equal to (\ref{k-th-der-at-0}).  We
know that $S(l^{0}(z))$ is an absolutely convergent power series in
$l^{0}(z)$ for $z\in O$.  In particular, for $k\in \N$, the sum of the
$k$-th derivative at $l^{0}(z)=0$ of the series $S(l^{0}(z))$ is equal
to the $k$-th derivative at $l^{0}(z)=0$ of the analytic function
given by the sum of $S(l^{0}(z))$.  Since for $k\in \N$, the $k$-th
derivative at $l^{0}(z)=0$ of the analytic function given by the sum
of (\ref{iter-sum}) and the $k$-th derivative at $l^{0}(z)=0$ of the
analytic function given by the sum of $S(l^{0}(z))$ are equal, these
two analytic functions must be equal on an open neighborhood $\Gamma$
of $z=1$ in the intersection of their domains.  Clearly we can choose
$\Gamma$ to be independent of $w_{(1)}\in W_{1}$ and $w_{(2)}\in
W_{2}$.  Then for
\[
z\in O\cap P\cap Q\cap \Pi\cap \Gamma,
\]
\begin{equation}\label{9.7-1-5}
\sum_{n\in \R}e^{-n\log z}\left(\sum_{i=0}^{N}\frac{(-\log z)^{i}}{i!}
((L'_{P(z_{0})}(0)-n)^{i}\lambda^{(2)}_{n})
(w_{(1)}\otimes w_{(2)})\right)
\end{equation}
is absolutely convergent to the right-hand side of (\ref{9.7-1-1}) and
in fact, we have proved that the corresponding double series, namely,
\begin{equation}\label{9.7-1-5.0}
\sum_{n\in \R}\sum_{i=0}^{N}e^{-n\log z}\frac{(-\log z)^{i}}{i!}
((L'_{P(z_{0})}(0)-n)^{i}\lambda^{(2)}_{n})
(w_{(1)}\otimes w_{(2)}),
\end{equation}
is also absolutely convergent to the right-hand side of (\ref{9.7-1-1}).

But using (\ref{9.7-1--2})--(\ref{z-prod}) and Definition
\ref{productanditerateexisting}, and recalling 
Remark \ref{I1I2'}, we obtain, when $|z_{0}+zz_{2}|>|zz_{2}|>0$,
\begin{eqnarray}\label{9.7-1-5.1}
\lefteqn{\langle e^{-(\log z)L'(0)}w'_{(4)}, 
(I^{z}_1\circ (1_{W_1}\otimes
I^{z}_2))(
w_{(1)}\otimes 
w_{(2)}\otimes e^{(\log z)L(0)}
w_{(3)})\rangle}\nn
&&=
((I^{z}_1\circ (1_{W_1}\otimes
I^{z}_2))'(e^{-(\log z)L'(0)}w'_{(4)}))(
w_{(1)}\otimes 
w_{(2)}\otimes e^{(\log z)L(0)}
w_{(3)})\nn
&&=(\mu^{(2)}_{(I^{z}_1\circ (1_{W_1}\otimes
I^{z}_2))'(e^{-(\log z)L'(0)}w'_{(4)}), e^{(\log z)L(0)}w_{(3)}})
(w_{(1)}\otimes w_{(2)}),
\end{eqnarray}
and for
\[
z\in O\cap P\cap Q\cap \Pi\cap\Gamma,
\]
the left-hand side equals (\ref{9.7-1-1}).  Thus for
\[
z\in O\cap P\cap Q\cap \Pi\cap\Gamma
\]
such that $|z_{0}+zz_{2}|>|zz_{2}|$,
\begin{eqnarray}\label{9.7-1-6}
\lefteqn{(\mu^{(2)}_{(I^{z}_1\circ (1_{W_1}\otimes
I^{z}_2))'(e^{-(\log z)L'(0)}w'_{(4)}), e^{(\log z)L(0)}w_{(3)}})
(w_{(1)}\otimes w_{(2)})}\nno\\
&&=\sum_{n\in \R}\sum_{i=0}^{N}e^{-n\log z}\frac{(-\log z)^{i}}{i!}
((L'_{P(z_{0})}(0)-n)^{i}\lambda^{(2)}_{n})
(w_{(1)}\otimes w_{(2)}).
\end{eqnarray}
Since $O$, $P$, $Q$, $\Pi$ and $\Gamma$ are all independent of
$w_{(1)}\in W_{1}$ and $w_{(2)}\in W_{2}$, (\ref{9.7-1-6}) holds for
$z\in O\cap P\cap Q\cap \Pi\cap\Gamma$ such that
$|z_{0}+zz_{2}|>|zz_{2}|$ and for all $w_{(1)}\in W_{1}$ and
$w_{(2)}\in W_{2}$. Thus for $z\in O\cap P\cap Q\cap \Pi\cap\Gamma$
such that $|z_{0}+zz_{2}|>|zz_{2}|$, we have
\begin{eqnarray}\label{9.7-1-6-1}
\lefteqn{\mu^{(2)}_{(I^{z}_1\circ (1_{W_1}\otimes
I^{z}_2))'(e^{-(\log z)L'(0)}w'_{(4)}), e^{(\log z)L(0)}w_{(3)}}}\nno\\
&&=\sum_{n\in \R}\sum_{i=0}^{N}e^{-n\log z}\frac{(-\log z)^{i}}{i!}
((L'_{P(z_{0})}(0)-n)^{i}\lambda^{(2)}_{n}),
\end{eqnarray}
where the right-hand side is understood as the sum of the weakly
absolutely convergent (double) series denoted by the same notation
(recall Remark \ref{weakly-abs-conv}).

Now $I^{z}_1$ and $I^{z}_2$ are $P(z_{0}+zz_{2})$- and
$P(zz_{2})$-intertwining maps, and when
\[
|z_{0}+zz_{2}|>|zz_{2}|>0,
\]
$I^{z}_1\circ (1_{W_1}\otimes I^{z}_2)$ is a
$P(z_0+zz_2,zz_2)$-intertwining map, by Proposition
\ref{productanditerateareintwmaps}.  Thus by Proposition \ref{8.12},
\begin{equation}\label{Izcompat}
(I^{z}_1\circ (1_{W_1}\otimes I^{z}_2))'(e^{-(\log z)L'(0)}w'_{(4)})
\in (W_1 \otimes W_2 \otimes W_3)^*
\end{equation}
satisfies the $P(z_{0}+zz_{2}, zz_{2})$-compatibility condition, and
\[
(z_{0}+zz_{2})-zz_{2}=z_{0}.
\]
Then by Lemma \ref{mulemma}, when
\[
|z_{0}+zz_{2}|>|zz_{2}|>|z_{0}|>0,
\]
for $v \in V$ the coefficients of the monomials in $x$ and $x_{1}$ in
\[
x^{-1}_1 \delta\left(\frac{x^{-1}-z_{0}}{x_1}\right)
\left(Y'_{P(z_{0})}(v, x)\mu^{(2)}_{(I^{z}_1\circ (1_{W_1}\otimes
I^{z}_2))'(e^{-(\log z)L'(0)}w'_{(4)}), e^{(\log z)L(0)}w_{(3)}}\right)
(w_{(1)}\otimes w_{(2)})
\]
are absolutely convergent and 
we have
\begin{eqnarray}\label{mu12-1}
\lefteqn{\biggl(\tau_{P(z_0)}\biggl( x^{-1}_1
\delta\bigg(\frac{x^{-1}-z_0}{x_1}\bigg) Y_{t}(v, x)\biggr)
\mu^{(2)}_{(I^{z}_1\circ (1_{W_1}\otimes
I^{z}_2))'(e^{-(\log z)L'(0)}w'_{(4)}), e^{(\log z)L(0)}w_{(3)}}\biggr)
(w_{(1)} \otimes w_{(2)})}\nno\\
&&=x^{-1}_1 \delta\bigg(\frac{x^{-1}-z_0}{x_1}\bigg)
\biggl(Y'_{P(z_0)}(v, x) \mu^{(2)}_{(I^{z}_1\circ (1_{W_1}\otimes
I^{z}_2))'(e^{-(\log z)L'(0)}w'_{(4)}), e^{(\log z)L(0)}w_{(3)}}\biggr)
(w_{(1)} \otimes w_{(2)}).\nn
\end{eqnarray}

As we previewed in (\ref{z-prod}), let $R$ be a sufficiently small
open neighborhood of $z=1$ such that
\begin{equation}\label{domainR}
|z_{0}+zz_{2}|>|zz_{2}|>|z_{0}|>0
\end{equation}
for $z\in R$. Then for 
\[
z \in O\cap P\cap Q\cap \Pi \cap \Gamma \cap R,
\] 
(\ref{mu12-1}) holds.  
{}From (\ref{9.7-1-6-1}) and (\ref{mu12-1}), we
obtain
\begin{eqnarray}\label{mu12-2}
\lefteqn{\Biggl(\tau_{P(z_0)}\biggl( x^{-1}_1
\delta\biggl(\frac{x^{-1}-z_0}{x_1}\bigg) Y_{t}(v, x)\biggr)}\nno\\
&&\quad\quad\quad
\Biggl(\sum_{n\in \R}\sum_{i=0}^{N}e^{-n\log z}\frac{(-\log z)^{i}}{i!}
((L'_{P(z_{0})}(0)-n)^{i}\lambda^{(2)}_{n})\Biggr)\Biggr)
(w_{(1)}\otimes w_{(2)})\nno\\
&&=x^{-1}_1 \delta\bigg(\frac{x^{-1}-z_0}{x_1}\bigg)\cdot\nno\\
&&\quad\quad\quad\cdot 
\Biggl(Y'_{P(z_0)}(v, x) \Biggl(\sum_{n\in \R}\sum_{i=0}^{N}e^{-n\log z}
\frac{(-\log z)^{i}}{i!}
((L'_{P(z_{0})}(0)-n)^{i}\lambda^{(2)}_{n})\Biggr)\Biggr)
(w_{(1)}\otimes w_{(2)})\nn
\end{eqnarray}
for
\[
z \in O\cap P\cap Q\cap \Pi \cap \Gamma \cap R,
\] 
a set that is independent of $w_{(1)}$ and $w_{(2)}$, and as in
(\ref{mu12-1}), the meaning of the right-hand side is that the
coefficient of each monomial in $x$ and $x_1$ is the sum of an
absolutely convergent series, each term of which now involves the
weakly absolutely convergent double sum over $n \in \R$ and $i=0,
\dots, N$.

We shall need to bring the double sums over $n$ and $i$ to the outside,
on both sides of (\ref{mu12-2}).

First we do this for the left-hand side of (\ref{mu12-2}).  Using the
definition (\ref{taudef}) of $\tau_{P(z_0)}$ and the definition
(\ref{yo}) of the opposite vertex operator $Y^{o}$, we can write the
definition of $\tau_{P(z_0)}$ more explicitly as
\begin{eqnarray}\label{mu12-1.0}
\lefteqn{\Biggl(\tau_{P(z_0)}\biggl( x^{-1}_1
\delta\biggl(\frac{x^{-1}-z_0}{x_1}\bigg) Y_{t}(v, x)\biggr)\lambda\Biggr)
(w_{(1)}\otimes w_{(2)})}\nn
&&=z_{0}^{-1}\delta\left(\frac{x^{-1}-x_1}{z_0}\right)\lambda
(Y_1(e^{xL(1)}(-x^{-2})^{L(0)}v,x_1)w_{(1)}\otimes
w_{(2)})\nno\\
&& \quad +x_{1}^{-1}\delta\left(\frac{z_0-x^{-1}}{-x_1}\right)\lambda(w_{(1)}\otimes
Y_2(e^{xL(1)}(-x^{-2})^{L(0)}v,x^{-1})w_{(2)})
\end{eqnarray}
for $v\in V$, $w_{(1)}\in W_{1}$, $w_{(2)}\in W_{2}$
and $\lambda\in (W_{1}\otimes W_{2})^{*}$.
In particular, the left-hand side of (\ref{mu12-1}) is equal to
\begin{eqnarray}\label{mu12-1.1}
\lefteqn{z_{0}^{-1}\delta\left(\frac{x^{-1}-x_1}{z_0}\right)\cdot}\nn
&&\cdot\mu^{(2)}_{(I^{z}_1\circ (1_{W_1}\otimes
I^{z}_2))'(e^{-(\log z)L'(0)}w'_{(4)}), e^{(\log z)L(0)}w_{(3)}}
(Y_1(e^{xL(1)}(-x^{-2})^{L(0)}v,x_1)w_{(1)}\otimes
w_{(2)})\nno\\
&& +x_{1}^{-1}\delta\left(\frac{z_0-x^{-1}}{-x_1}\right)\cdot\nn
&&\quad\cdot\mu^{(2)}_{(I^{z}_1\circ (1_{W_1}\otimes
I^{z}_2))'(e^{-(\log z)L'(0)}w'_{(4)}), e^{(\log z)L(0)}w_{(3)}}(w_{(1)}\otimes
Y_2(e^{xL(1)}(-x^{-2})^{L(0)}v,x^{-1})w_{(2)}).\nn
\end{eqnarray}

For 
\[
z \in O\cap P\cap Q\cap \Pi \cap \Gamma\cap R,
\]
by (\ref{9.7-1-6-1}), the coefficients of the monomials in $x$ and $x_{1}$ in
\begin{eqnarray}\label{mu12-1.1-1}
\lefteqn{\sum_{n\in \R}\sum_{i=0}^{N}e^{-n\log z}
\frac{(-\log z)^{i}}{i!}
\cdot}\nno\\
&&\quad\quad\quad\quad\quad\quad\cdot 
((L'_{P(z_{0})}(0)-n)^{i}\lambda^{(2)}_{n})
(Y_1(e^{xL(1)}(-x^{-2})^{L(0)}v,x_1)w_{(1)}\otimes
w_{(2)})
\end{eqnarray}
and 
\begin{eqnarray}\label{mu12-1.1-2}
\lefteqn{\sum_{n\in \R}\sum_{i=0}^{N}e^{-n\log z}
\frac{(-\log z)^{i}}{i!}
\cdot}\nno\\
&&\quad\quad\quad\quad\quad\quad\cdot 
((L'_{P(z_{0})}(0)-n)^{i}\lambda^{(2)}_{n})(w_{(1)}\otimes
Y_2(e^{xL(1)}(-x^{-2})^{L(0)}v,x^{-1})w_{(2)})\nn
\end{eqnarray}
are absolutely convergent to the 
corresponding coefficients of the monomials in $x$ and $x_{1}$ in
\[
\mu^{(2)}_{(I^{z}_1\circ (1_{W_1}\otimes
I^{z}_2))'(e^{-(\log z)L'(0)}w'_{(4)}), e^{(\log z)L(0)}w_{(3)}}
(Y_1(e^{xL(1)}(-x^{-2})^{L(0)}v,x_1)w_{(1)}\otimes
w_{(2)})
\]
and 
\[
\mu^{(2)}_{(I^{z}_1\circ (1_{W_1}\otimes
I^{z}_2))'(e^{-(\log z)L'(0)}w'_{(4)}), e^{(\log z)L(0)}w_{(3)}}(w_{(1)}\otimes
Y_2(e^{xL(1)}(-x^{-2})^{L(0)}v,x^{-1})w_{(2)}),
\]
respectively. Then, as finite linear combinations of the coefficients
of the monomials in $x$ and $x_{1}$ in (\ref{mu12-1.1-1}) and
(\ref{mu12-1.1-2}), the coefficients of the monomials in $x$ and
$x_{1}$ in
\begin{eqnarray}\label{mu12-1.1-3}
\lefteqn{z_{0}^{-1}\delta\left(\frac{x^{-1}-x_1}{z_0}\right)\Biggl(\sum_{n\in \R}
\sum_{i=0}^{N}e^{-n\log z}
\frac{(-\log z)^{i}}{i!}
\cdot}\nno\\
&&\quad\quad\quad\quad\quad\quad\cdot 
((L'_{P(z_{0})}(0)-n)^{i}\lambda^{(2)}_{n})
(Y_1(e^{xL(1)}(-x^{-2})^{L(0)}v,x_1)w_{(1)}\otimes
w_{(2)})\Biggr)\nno\\
&& \quad+x_{1}^{-1}\delta\left(\frac{z_0-x^{-1}}{-x_1}\right)\Biggl(\sum_{n\in \R}
\sum_{i=0}^{N}e^{-n\log z}\frac{(-\log z)^{i}}{i!}
\cdot\nno\\
&&\quad\quad\quad\quad\quad\quad\cdot 
((L'_{P(z_{0})}(0)-n)^{i}\lambda^{(2)}_{n})(w_{(1)}\otimes
Y_2(e^{xL(1)}(-x^{-2})^{L(0)}v,x^{-1})w_{(2)})\Biggr)\nn
\end{eqnarray}
are absolutely convergent to the corresponding 
coefficients of the monomials in $x$ and $x_{1}$ in
(\ref{mu12-1.1}). 

Now for 
\[
z \in O\cap P\cap Q\cap \Pi \cap \Gamma\cap R,
\]
we consider the coefficients of the monomials in $x$ and $x_{1}$ in
\begin{eqnarray}\label{mu12-1.1-4}
\lefteqn{\sum_{n\in \R}\sum_{i=0}^{N}e^{-n\log z}
\frac{(-\log z)^{i}}{i!}\cdot}\nno\\
&&\quad\quad\cdot 
\Biggl(\tau_{P(z_0)}\biggl( x^{-1}_1
\delta\biggl(\frac{x^{-1}-z_0}{x_1}\bigg) Y_{t}(v, x)\biggr)
((L'_{P(z_{0})}(0)-n)^{i}\lambda^{(2)}_{n})
(w_{(1)}\otimes w_{(2)})\Biggr)\nn
&&=\sum_{n\in \R}\sum_{i=0}^{N}e^{-n\log z}
\frac{(-\log z)^{i}}{i!}
z_{0}^{-1}\delta\left(\frac{x^{-1}-x_1}{z_0}\right)
\cdot\nno\\
&&\quad\quad\quad\quad\quad\quad\cdot 
((L'_{P(z_{0})}(0)-n)^{i}\lambda^{(2)}_{n})
(Y_1(e^{xL(1)}(-x^{-2})^{L(0)}v,x_1)w_{(1)}\otimes
w_{(2)})\nno\\
&& \quad+\sum_{n\in \R}\sum_{i=0}^{N}e^{-n\log z}
\frac{(-\log z)^{i}}{i!}
x_{1}^{-1}\delta\left(\frac{z_0-x^{-1}}{-x_1}\right)
\cdot\nno\\
&&\quad\quad\quad\quad\quad\quad\cdot 
((L'_{P(z_{0})}(0)-n)^{i}\lambda^{(2)}_{n})(w_{(1)}\otimes
Y_2(e^{xL(1)}(-x^{-2})^{L(0)}v,x^{-1})w_{(2)}),\nn
\end{eqnarray}
where we have used (\ref{mu12-1.0}). 
The coefficient of each monomial in $x$ and $x_{1}$ in the
right-hand side of (\ref{mu12-1.1-4}) is the sum over $n$ and $i$
of the sum of the corresponding monomial in $x$ and $x_{1}$ in
\begin{eqnarray}
\lefteqn{e^{-n\log z}
\frac{(-\log z)^{i}}{i!}
z_{0}^{-1}\delta\left(\frac{x^{-1}-x_1}{z_0}\right)
\cdot}\nno\\
&&\quad\quad\cdot 
((L'_{P(z_{0})}(0)-n)^{i}\lambda^{(2)}_{n})
(Y_1(e^{xL(1)}(-x^{-2})^{L(0)}v,x_1)w_{(1)}\otimes
w_{(2)})
\end{eqnarray}
and in
\begin{eqnarray}
\lefteqn{e^{-n\log z}
\frac{(-\log z)^{i}}{i!}
x_{1}^{-1}\delta\left(\frac{z_0-x^{-1}}{-x_1}\right)
\cdot}\nno\\
&&\quad\quad\cdot 
((L'_{P(z_{0})}(0)-n)^{i}\lambda^{(2)}_{n})(w_{(1)}\otimes
Y_2(e^{xL(1)}(-x^{-2})^{L(0)}v,x^{-1})w_{(2)}).
\end{eqnarray}
Thus since finite linear combinations of the coefficients of the
monomials in $x$ and $x_{1}$ in (\ref{mu12-1.1-1}) and
(\ref{mu12-1.1-2}) are absolutely convergent, the coefficients of the
monomials in $x$ and $x_{1}$ in the right-hand side of
(\ref{mu12-1.1-4}) are also absolutely convergent.  Moreover, these
(absolutely convergent) coefficients are equal to the (absolutely
convergent) coefficients of the monomials in $x$ and $x_{1}$ in
(\ref{mu12-1.1-3}).  Thus by (\ref{mu12-1.1-4}), for
\[
z \in O\cap P\cap Q\cap \Pi \cap \Gamma\cap R,
\]
the coefficient of each monomial in $x$ and $x_{1}$ in 
\begin{eqnarray}\label{mu12-1.2}
\lefteqn{\sum_{n\in \R}\sum_{i=0}^{N}e^{-n\log z}
\frac{(-\log z)^{i}}{i!}\cdot}\nno\\
&&\quad\quad\cdot 
\Biggl(\tau_{P(z_0)}\biggl( x^{-1}_1
\delta\biggl(\frac{x^{-1}-z_0}{x_1}\bigg) Y_{t}(v, x)\biggr)
((L'_{P(z_{0})}(0)-n)^{i}\lambda^{(2)}_{n})
(w_{(1)}\otimes w_{(2)})\Biggr)\nn
\end{eqnarray}
is an absolutely convergent (double) series and converges to the
coefficient of the corresponding monomial in $x$ and $x_{1}$ in the
left-hand side of (\ref{mu12-1}) (or of (\ref{mu12-2})).

Now we need to bring the double sum over $n$ and $i$ on the right-hand
side of (\ref{mu12-2}) to the outside, and in the process, we shall
need to increase $N$ and restrict the range of $z$.

Taking ${\rm Res}_{x_1}$ in (\ref{mu12-1.2}) and using (\ref{3.18-1}), 
we see that the coefficient of each monomial
in $x$ in
\begin{equation}\label{mu12-1.2.5}
\sum_{n\in \R}\sum_{i=0}^{N}e^{-n\log z}
\frac{(-\log z)^{i}}{i!}Y'_{P(z_0)}(v, x) 
((L'_{P(z_{0})}(0)-n)^{i}\lambda^{(2)}_{n})
(w_{(1)}\otimes w_{(2)})
\end{equation}
is an absolutely convergent series and that it converges to the
coefficient of the corresponding monomial in $x$ in the result of
applying ${\rm Res}_{x_1}$ to the left-hand side of (\ref{mu12-1})
(or of (\ref{mu12-2})), namely,
\begin{equation}\label{mu12-1.2.6}
\Biggl(Y'_{P(z_0)}(v, x) \Biggl(\sum_{n\in \R}\sum_{i=0}^{N}e^{-n\log z}
\frac{(-\log z)^{i}}{i!}
((L'_{P(z_{0})}(0)-n)^{i}\lambda^{(2)}_{n})\Biggr)\Biggr)
(w_{(1)}\otimes w_{(2)})
\end{equation}
(recall (\ref{9.7-1-6-1})).

By (\ref{rlm4}), with the continuing assumption
on $z$, we thus have 
\begin{eqnarray}\label{mu12-1.4}
\lefteqn{\sum_{n\in \R}\sum_{i=0}^{N}e^{-n\log z}
\frac{(-\log z)^{i}}{i!}Y'_{P(z_0)}(v, x) 
((L'_{P(z_{0})}(0)-n)^{i}\lambda^{(2)}_{n})
(w_{(1)}\otimes w_{(2)})}\nn
&&+{\rm
Res}_{x_0^{-1}}\dlti{x}{-zz_2}{+x_0^{-1}}((I^{z}_1\circ (1_{W_1}\otimes
I^{z}_2))'(e^{-(\log z)L'(0)}w'_{(4)})) (w_{(1)}\otimes w_{(2)}\nn
&&\qquad\qquad \qquad\otimes
Y_3(e^{xL(1)}(-x^{-2})^{L(0)}v,x_0^{-1})e^{(\log z)L(0)}w_{(3)})\nn
&&=\Biggl(Y'_{P(z_0)}(v, x) \Biggl(\sum_{n\in \R}\sum_{i=0}^{N}e^{-n\log z}
\frac{(-\log z)^{i}}{i!}
((L'_{P(z_{0})}(0)-n)^{i}\lambda^{(2)}_{n})\Biggr)\Biggr)
(w_{(1)}\otimes w_{(2)})\nn
&&\quad +{\rm
Res}_{x_0^{-1}}\dlti{x}{-zz_2}{+x_0^{-1}}((I^{z}_1\circ (1_{W_1}\otimes
I^{z}_2))'(e^{-(\log z)L'(0)}w'_{(4)})) (w_{(1)}\otimes w_{(2)}\nn
&&\qquad\qquad \qquad\otimes
Y_3(e^{xL(1)}(-x^{-2})^{L(0)}v,x_0^{-1})e^{(\log z)L(0)}w_{(3)})\nn
&&=(Y'_{P(z_0)}(v, x) \mu^{(2)}_{(I^{z}_1\circ (1_{W_1}\otimes
I^{z}_2))'(e^{-(\log z)L'(0)}w'_{(4)}), e^{(\log z)L(0)}w_{(3)}})
(w_{(1)}\otimes w_{(2)})\nn
&&\quad +{\rm
Res}_{x_0^{-1}}\dlti{x}{-zz_2}{+x_0^{-1}}((I^{z}_1\circ (1_{W_1}\otimes
I^{z}_2))'(e^{-(\log z)L'(0)}w'_{(4)})) (w_{(1)}\otimes w_{(2)}\nn
&&\qquad\qquad \qquad\otimes
Y_3(e^{xL(1)}(-x^{-2})^{L(0)}v,x_0^{-1})e^{(\log z)L(0)}w_{(3)})\nn
&&={\rm
Res}_{x_0^{-1}}\Biggl(\tau_{P(z_0+zz_2,zz_2)}\Biggl(
\dlti{x}{x_0^{-1}}{-zz_2}
Y_t((-x_0^2)^{L(0)}e^{-x_0L(1)}e^{xL(1)}(-x^{-2})^{L(0)}v,x_0)
\Biggr)\cdot\nn
&&\qquad\qquad\qquad\cdot ((I^{z}_1\circ (1_{W_1}\otimes
I^{z}_2))'(e^{-(\log z)L'(0)}w'_{(4)}))\Biggr)(w_{(1)}\otimes w_{(2)}\otimes 
e^{(\log z)L(0)}w_{(3)}),\nn
\end{eqnarray}
where the equalities include  the information that 
the coefficient of each monomial in each expression involving a double sum
is absolutely convergent. 

Recall from the $P^{(2)}(z_0)$-grading condition that the elements
$\lambda^{(2)}_{n}$ for $n \in \R$ depend on $w_{(3)}$ (and on
$w'_{(4)}$) and that the sets $\Pi$ and $\Gamma$ cannot be assumed
independent of $w_{(3)}$, nor can our integer $N \in \N$.  Since we
now need to use $\lambda^{(2)}_{n}$, $n\in \R$, for different
$w_{(3)}$, we denote $\lambda^{(2)}_{n}$ by
$\lambda^{(2)}_{n}[w_{(3)}]$ for $n\in \R$ and we denote $\Pi$,
$\Gamma$ and $N$ by $\Pi[w_{(3)}]$, $\Gamma[w_{(3)}]$ and
$N[w_{(3)}]$, respectively.  Then with these notations, the
convergence properties and formulas that we have proved hold for all
$w_{(3)}\in W_{3}$.

In order to handle the second term in the left-hand side of
(\ref{mu12-1.4}), we shall need to consider the following application
of a certain conjugated operator to $w_{(3)}$, and to treat this
element as an analogue of $w_{(3)}$:
\begin{eqnarray}\label{def-X}
X&=&e^{-(\log z)L(0)}Y_3(e^{xL(1)}(-x^{-2})^{L(0)}v,x_0^{-1}) e^{(\log
z)L(0)}w_{(3)}\nn
&=&Y_3(e^{-(\log z)L(0)}e^{xL(1)}(-x^{-2})^{L(0)}v, e^{-\log z}x_0^{-1})w_{(3)},
\end{eqnarray}
where we obtain the second expression by using the conjugation formula
(\ref{710}). Since $e^{xL(1)}(-x^{-2})^{L(0)}v\in V[x, x^{-1}]$
and $Y_{3}(u, y)w_{(3)}\in W_{3}((y))$ for $u\in V$ and $y$ a formal variable, 
the right-hand side of 
(\ref{def-X}) is of the form
\[
Y_3(z^{-L(0)}e^{xL(1)}(-x^{-2})^{L(0)}v, (z x_0)^{-1})w_{(3)}
= \sum_{l \le L}\sum_{m=-M}^{M}\sum_{s=-S}^{S}
z^{s+l}w^{l,m}_{s}x_{0}^{l}x^{m}
\]
for certain integers $L$, $M \ge 0$ and $S \ge 0$ and certain
(determined) elements $w^{l, m}_{s} \in W_3$. That is,
\[
X = \sum_{l \le L}\sum_{m=-M}^{M}\sum_{s=-S}^{S}
z^{s+l}w^{l,m}_{s}x_{0}^{l}x^{m}.
\]

With this notation, for $z \in R$ (recall (\ref{domainR}) and
(\ref{9.7-1-5.1})), the second term in the left-hand side of
(\ref{mu12-1.4}) equals
\begin{equation}\label{ResZ}
{\rm Res}_{x_0^{-1}}\dlti{x}{-zz_2}{+x_0^{-1}}Z,
\end{equation}
where
\begin{equation}\label{Z}
Z=(I^{z}_1\circ (1_{W_1}\otimes I^{z}_2))'(e^{-(\log
z)L'(0)}w'_{(4)})) (w_{(1)}\otimes w_{(2)} \otimes e^{(\log z)L(0)}X).
\end{equation}

The coefficient of $x_{0}^{l}x^{m}$ in (\ref{Z}) is the precisely
determined expression
\begin{eqnarray*}
\lefteqn{(I^{z}_1\circ (1_{W_1}\otimes I^{z}_2))'(e^{-(\log
z)L'(0)}w'_{(4)})) \left(w_{(1)}\otimes w_{(2)} \otimes e^{(\log
z)L(0)} \left(\sum_{s=-S}^{S}z^{s+l}w^{l,m}_{s}\right)\right)}\nno\\
&&=\sum_{s=-S}^{S}z^{s+l} (I^{z}_1\circ
(1_{W_1}\otimes I^{z}_2))'(e^{-(\log z)L'(0)}w'_{(4)}))
(w_{(1)}\otimes w_{(2)} \otimes e^{(\log z)L(0)}w^{l,m}_{s}),\nn
\end{eqnarray*}
which by (\ref{9.7-1-5.1}) equals
\begin{equation}\label{coeffofx0lxm}
\sum_{s=-S}^{S}z^{s+l}
(\mu^{(2)}_{(I^{z}_1\circ (1_{W_1}\otimes
I^{z}_2))'(e^{-(\log z)L'(0)}w'_{(4)}), e^{(\log z)L(0)}w^{l,m}_{s}})
(w_{(1)}\otimes w_{(2)})
\end{equation}
(with $z \in R$, as we have assumed above).  In particular,
\begin{equation}\label{Zexpansion}
Z= \sum_{l \le L}\sum_{m=-M}^{M}
\Biggl(\sum_{s=-S}^{S}z^{s+l}
(\mu^{(2)}_{(I^{z}_1\circ (1_{W_1}\otimes
I^{z}_2))'(e^{-(\log z)L'(0)}w'_{(4)}), e^{(\log z)L(0)}w^{l,m}_{s}})
(w_{(1)}\otimes w_{(2)})\Biggr)x_{0}^{l}x^{m}.
\end{equation}

Now the residue (\ref{ResZ}) involves only finitely many powers of
$x_0$ in (\ref{Zexpansion}), and in particular, the coefficient of
each monomial in $x$ in (\ref{ResZ}) is a finite linear combination of
expressions of the form (\ref{coeffofx0lxm}), involving only finitely
many $l$, independently of the power of $x$ in (\ref{ResZ}).  We apply
our results above to the corresponding finite family of elements
$w^{l,m}_{s}$ that arise in this way, and we use
$\lambda^{(2)}_{n}[w^{l,m}_{s}]$, $\Pi[w^{l,m}_{s}]$,
$\Gamma[w^{l,m}_{s}]$ and $N[w^{l,m}_{s}]$, as defined above.  Let
\[
\widetilde{N}=\max\{N[w^{l,m}_{s}]\},
\]
\[
\widetilde{\Pi}= \bigcap(\Pi[w^{l,m}_{s}]),
\]
\[
\widetilde{\Gamma}= \bigcap(\Gamma[w^{l,m}_{s}]),
\]
and let
\begin{equation}\label{rangeofz}
z \in O\cap P\cap Q\cap \widetilde{\Pi} \cap \widetilde{\Gamma}\cap R.
\end{equation}
Then by (\ref{9.7-1-6-1}), 
\begin{eqnarray}
\lefteqn{\mu^{(2)}_{(I^{z}_1\circ (1_{W_1}\otimes
I^{z}_2))'(e^{-(\log z)L'(0)}w'_{(4)}), e^{(\log z)L(0)}w^{l,m}_{s}}}\nno\\
&&=\sum_{n\in \R}\sum_{i=0}^{\widetilde{N}} e^{-n\log z}\frac{(-\log z)^{i}}{i!}
((L'_{P(z_{0})}(0)-n)^{i}\lambda^{(2)}_{n}[w^{l,m}_{s}]),
\end{eqnarray}
where the right-hand side is understood as the sum of a weakly
absolutely convergent double series.  For $n \in \R$, form the formal
Laurent series
\begin{equation}\label{Lambda}
\Lambda^{(2)}_{n}(x_{0},x;z) = \sum_{l \le L}\sum_{m=-M}^{M}
\left(\sum_{s=-S}^{S}z^{s+l}
\lambda^{(2)}_{n}[w^{l,m}_{s}]\right)x_{0}^{l}x^{m}.
\end{equation}
Then the coefficient of each power of $x$ in the second term in the
left-hand side of (\ref{mu12-1.4}) equals its coefficient in
\begin{eqnarray}\label{mu12-1.5}
\lefteqn{{\rm Res}_{x_0^{-1}}\dlti{x}{-zz_2}{+x_0^{-1}}Z}\nno\\
&&={\rm Res}_{x_0^{-1}}\dlti{x}{-zz_2}{+x_0^{-1}}
\sum_{l\le L}\sum_{m=-M}^{M} \sum_{s=-S}^{S}z^{s+l}
\cdot\nno\\
&&\quad\quad\quad\cdot
\left(\sum_{n\in \R}\sum_{i=0}^{\widetilde{N}} e^{-n\log
z}\frac{(-\log z)^{i}}{i!}
((L'_{P(z_{0})}(0)-n)^{i}\lambda^{(2)}_{n}[w^{l,m}_{s}])
(w_{(1)}\otimes w_{(2)})\right)x_{0}^{l}x^{m}\nno\\
&&={\rm Res}_{x_0^{-1}}\dlti{x}{-zz_2}{+x_0^{-1}}\cdot\nno\\
&&\quad\quad\quad\cdot
\sum_{n\in\R}\sum_{i=0}^{\widetilde{N}}e^{-n\log z}\frac{(-\log
z)^{i}}{i!}((L'_{P(z_{0})}(0)-n)^{i} \Lambda^{(2)}_{n}(x_{0},x;z))
(w_{(1)}\otimes w_{(2)})\nno\\
&&=\sum_{n\in\R}\sum_{i=0}^{\widetilde{N}}
e^{-n\log z}\frac{(-\log z)^{i}}{i!}
{\rm Res}_{x_0^{-1}}\dlti{x}{-zz_2}{+x_0^{-1}}\cdot\nno\\
&&\quad\quad\quad\cdot
((L'_{P(z_{0})}(0)-n)^{i} \Lambda^{(2)}_{n}(x_{0},x;z))
(w_{(1)}\otimes w_{(2)}),
\end{eqnarray}
where we have double absolute convergence; recall that $\widetilde{N}$
and the range (\ref{rangeofz}) of $z$ are independent of the power of
$x$.

In order to reach our goal of bringing the double sum over $n$ and $i$
to the outside on the right-hand side of (\ref{mu12-2}), we need to
multiply (\ref{mu12-1.2.5}), (\ref{mu12-1.2.6}) and the other
expressions in (\ref{mu12-1.4}) by
\begin{equation}\label{dltx1x-1-z0}
\dlt{x_1}{x^{-1}}{-z_0};
\end{equation}
in particular, we need to show that we can do this for each of the
expressions.

First we do this for the expressions in (\ref{mu12-1.5}), which are
equal to the second term in the left-hand side of (\ref{mu12-1.4}).
Since $|zz_{2}|>|z_{0}|$, by Lemma \ref{deltalemma}, formula
(\ref{l4}),
\begin{equation}\label{prod-delta}
\dlt{x_1}{x^{-1}}{-z_0}\dlti{x}{-zz_2}{+x_0^{-1}}
\end{equation}
is a formal 
Laurent series in $x$, $x_{1}$ and $x_{0}$ each of whose coefficients
is an absolutely convergent series of the form $\sum_{j\in \N}a_{j}$ 
($a_{j}\in \C$). 
Thus for 
\[
z \in O\cap P\cap Q\cap \widetilde{\Pi} \cap \widetilde{\Gamma}\cap R,
\]
the coefficient of each monomial in $x$ and $x_1$ in
\begin{eqnarray}\label{mu12-1.5.0}
\lefteqn{\dlt{x_1}{x^{-1}}{-z_0}\sum_{n\in\R}\sum_{i=0}^{\widetilde{N}}
e^{-n\log z}\frac{(-\log z)^{i}}{i!}
{\rm Res}_{x_0^{-1}}\dlti{x}{-zz_2}{+x_0^{-1}}\cdot}\nno\\
&&\quad\quad\quad\cdot
((L'_{P(z_{0})}(0)-n)^{i} \Lambda^{(2)}_{n}(x_{0},x;z))
(w_{(1)}\otimes w_{(2)})\nn
&&=\dlt{x_1}{x^{-1}}{-z_0}{\rm Res}_{x_0^{-1}}
\dlti{x}{-zz_2}{+x_0^{-1}}
\cdot\nn
&&\quad\quad\cdot
\left(\sum_{n\in \R}\sum_{i=0}^{\widetilde{N}}e^{-n\log z}\frac{(-\log z)^{i}}{i!}
((L'_{P(z_{0})}(0)-n)^{i}\Lambda^{(2)}_{n}(x_{0},x;z))
(w_{(1)}\otimes w_{(2)})\right)\nn
&&={\rm Res}_{x_0^{-1}}\dlt{x_1}{x^{-1}}{-z_0}
\dlti{x}{-zz_2}{+x_0^{-1}}\cdot\nn
&&\quad\quad\cdot
\left(\sum_{n\in \R}\sum_{i=0}^{\widetilde{N}}e^{-n\log z}\frac{(-\log z)^{i}}{i!}
((L'_{P(z_{0})}(0)-n)^{i}\Lambda^{(2)}_{n}(x_{0},x;z))
(w_{(1)}\otimes w_{(2)})\right)\nn
\end{eqnarray}
is a finite linear combination of products of pairs of absolutely
convergent series and hence is a finite linear combination of
absolutely convergent triple series (over $j\in \N$, $n \in \R$ and
$i=0,\dots, {\widetilde{N}}$).  In particular, the second term in the
left-hand side of (\ref{mu12-1.4}) can be multiplied by
(\ref{dltx1x-1-z0}), in the sense of absolute convergence, and
moreover, for
\[
z \in O\cap P\cap Q\cap \widetilde{\Pi} \cap \widetilde{\Gamma}\cap R,
\]
the coefficient of each monomial in $x$ and $x_1$ in
\begin{eqnarray}\label{mu12-1.5.00}
\lefteqn{\sum_{n\in \R}\sum_{i=0}^{\widetilde{N}}e^{-n\log z}\frac{(-\log z)^{i}}{i!}
\cdot}\nn
&& \cdot
{\rm Res}_{x_0^{-1}}\dlt{x_1}{x^{-1}}{-z_0}\dlti{x}{-zz_2}{+x_0^{-1}}
((L'_{P(z_{0})}(0)-n)^{i}\Lambda^{(2)}_{n}(x_{0},x;z))(w_{(1)}\otimes w_{(2)})\nn
\end{eqnarray}
is absolutely convergent to the coefficient of the corresponding monomial 
in (\ref{mu12-1.5.0}).

We now re-express (\ref{mu12-1.5.0}) and (\ref{mu12-1.5.00}) by using
the explicit dependence of (\ref{Lambda}) on $z$.  Since the
coefficient of each power of $x_{0}$ in (\ref{Lambda}) is a finite
sum, we have
\begin{eqnarray*}
\lefteqn{\dlt{x_1}{x^{-1}}{-z_0}\dlti{x}{-zz_2}{+x_0^{-1}}
((L'_{P(z_{0})}(0)-n)^{i}\Lambda^{(2)}_{n}(x_{0},x;z))(w_{(1)}\otimes w_{(2)})}\nn
&&=\dlt{x_1}{x^{-1}}{-z_0}\dlti{x}{-zz_2}{+x_0^{-1}}
\cdot\nn
&&\quad\quad\quad \quad \cdot \sum_{l \le L}\sum_{m=-M}^{M}
\sum_{s=-S}^{S}z^{s+l}((L'_{P(z_{0})}(0)-n)^{i}
\lambda^{(2)}_{n}[w^{l,m}_{s}])(w_{(1)}\otimes w_{(2)})
x_{0}^{l}x^{m}\nn
&&=\sum_{l \le L}\sum_{m=-M}^{M}
\sum_{s=-S}^{S}z^{s+l}\dlt{x_1}{x^{-1}}{-z_0}\dlti{x}{-zz_2}{+x_0^{-1}}\cdot\nn
&&\quad\quad\quad \quad\quad\quad\quad \quad \cdot 
((L'_{P(z_{0})}(0)-n)^{i}
\lambda^{(2)}_{n}[w^{l,m}_{s}])(w_{(1)}\otimes w_{(2)})x_{0}^{l}x^{m}.
\end{eqnarray*}
Thus for
\[
z \in O\cap P\cap Q\cap \widetilde{\Pi} \cap \widetilde{\Gamma}\cap R,
\]
the coefficient of each monomial in $x$ and $x_1$ in
(\ref{mu12-1.5.00}), written as
\begin{eqnarray}\label{mu12-1.5.01}
\lefteqn{\sum_{n\in \R}\sum_{i=0}^{\widetilde{N}}e^{-n\log z}\frac{(-\log z)^{i}}{i!}
\cdot}\nn
&&\quad\quad \cdot {\rm Res}_{x_0^{-1}}\sum_{l \le L}\sum_{m=-M}^{M}
\sum_{s=-S}^{S}z^{s+l}\dlt{x_1}{x^{-1}}{-z_0}\dlti{x}{-zz_2}{+x_0^{-1}}\cdot\nn
&&\quad\quad\quad \quad\cdot 
((L'_{P(z_{0})}(0)-n)^{i}
\lambda^{(2)}_{n}[w^{l,m}_{s}])(w_{(1)}\otimes w_{(2)})x_{0}^{l}x^{m}\nn
&&=\sum_{n\in \R}\sum_{i=0}^{\widetilde{N}}
\sum_{l \le L}\sum_{m=-M}^{M}
\sum_{s=-S}^{S}e^{-(n-s-l)\log z}\frac{(-\log z)^{i}}{i!}
\cdot\nn
&&\quad\quad \cdot 
{\rm Res}_{x_0^{-1}}\dlt{x_1}{x^{-1}}{-z_0}\dlti{x}{-zz_2}{+x_0^{-1}}\cdot\nn
&&\quad\quad\quad\quad \cdot ((L'_{P(z_{0})}(0)-n)^{i}
\lambda^{(2)}_{n}[w^{l,m}_{s}])(w_{(1)}\otimes w_{(2)})x_{0}^{l}x^{m}\nn
&&=\sum_{p\in \R}\sum_{i=0}^{\widetilde{N}}
\sum_{l \le L}\sum_{m=-M}^{M}
\sum_{s=-S}^{S}e^{-p\log z}\frac{(-\log z)^{i}}{i!}
\cdot\nn
&&\quad\quad \cdot 
{\rm Res}_{x_0^{-1}}\dlt{x_1}{x^{-1}}{-z_0}\dlti{x}{-zz_2}{+x_0^{-1}}\cdot\nn
&&\quad\quad\quad\quad \cdot 
((L'_{P(z_{0})}(0)-(p+s+l))^{i}
\lambda^{(2)}_{p+s+l}[w^{l,m}_{s}])(w_{(1)}\otimes w_{(2)})x_{0}^{l}x^{m}\nn
&&=\sum_{n\in \R}\sum_{i=0}^{\widetilde{N}}
e^{-n\log z}\frac{(-\log z)^{i}}{i!}
\cdot\nn
&&\quad\quad \cdot 
{\rm Res}_{x_0^{-1}}\dlt{x_1}{x^{-1}}{-z_0}\dlti{x}{-zz_2}{+x_0^{-1}}\cdot\nn
&&\quad\quad\quad\quad \cdot 
\sum_{l \le L}\sum_{m=-M}^{M}
\sum_{s=-S}^{S}((L'_{P(z_{0})}(0)-(n+s+l))^{i}
\lambda^{(2)}_{n+s+l}[w^{l,m}_{s}])(w_{(1)}\otimes w_{(2)})x_{0}^{l}x^{m}\nn
\end{eqnarray}
is absolutely convergent to the coefficient of the corresponding
monomial in (\ref{mu12-1.5.0}), written as
\begin{eqnarray}\label{mu12-1.5.02}
\lefteqn{\dlt{x_1}{x^{-1}}{-z_0}{\rm Res}_{x_0^{-1}} \dlti{x}{-zz_2}{+x_0^{-1}}
\cdot}\nn
&&\quad\quad \cdot 
\sum_{n\in \R}\sum_{i=0}^{\widetilde{N}}e^{-n\log z}\frac{(-\log z)^{i}}{i!}
\cdot\nn
&&\quad\quad\quad \quad\cdot \sum_{l \le L}\sum_{m=-M}^{M}
\sum_{s=-S}^{S}z^{s+l}
((L'_{P(z_{0})}(0)-n)^{i}
\lambda^{(2)}_{n}[w^{l,m}_{s}])(w_{(1)}\otimes w_{(2)})x_{0}^{l}x^{m}\nn
&&=\dlt{x_1}{x^{-1}}{-z_0}{\rm Res}_{x_0^{-1}}\dlti{x}{-zz_2}{+x_0^{-1}}
\cdot\nn
&&\quad\quad \cdot \sum_{n\in \R}\sum_{i=0}^{\widetilde{N}}
e^{-n\log z}\frac{(-\log z)^{i}}{i!}
\cdot\nn
&&\quad\quad\quad\quad \cdot 
\sum_{l \le L}\sum_{m=-M}^{M}
\sum_{s=-S}^{S}((L'_{P(z_{0})}(0)-(n+s+l))^{i}
\lambda^{(2)}_{n+s+l}[w^{l,m}_{s}])(w_{(1)}\otimes w_{(2)})x_{0}^{l}x^{m}.\nn
\end{eqnarray}
Moreover, recall that (\ref{mu12-1.5.0}), and hence
(\ref{mu12-1.5.02}), is equal to the product, in the sense of absolute
convergence, of (\ref{dltx1x-1-z0}) and the second term in the
left-hand side of (\ref{mu12-1.4}).

Next we show that we can multiply the right-hand side of
(\ref{mu12-1.4}) by (\ref{dltx1x-1-z0}), by using the $P(z_0+zz_2,
zz_2)$-compatibility condition.  Since the expression
(\ref{Izcompat}) satisfies this condition, we have
\begin{eqnarray*}
\lefteqn{\Biggl(\tau_{P(z_0+zz_2,zz_2)}\Biggl(
\dlt{x_{1}}{x_0^{-1}}{-(z_0+zz_2)}
\dlti{x}{x_0^{-1}}{-zz_2}
Y_t((-x_0^2)^{L(0)}e^{-x_0L(1)}e^{xL(1)}(-x^{-2})^{L(0)}v,x_0)
\Biggr)\cdot}\nn
&&\qquad\qquad\qquad\cdot ((I^{z}_1\circ (1_{W_1}\otimes
I^{z}_2))'(e^{-(\log z)L'(0)}w'_{(4)}))\Biggr)(w_{(1)}\otimes w_{(2)}\otimes 
e^{(\log z)L(0)}w_{(3)})\nn
&&=\dlt{x_{1}}{x_0^{-1}}{-(z_0+zz_2)}
\dlti{x}{x_0^{-1}}{-zz_2}\biggl(Y'_{P(z_0+zz_2,zz_2)}
((-x_0^2)^{L(0)}e^{-x_0L(1)}e^{xL(1)}(-x^{-2})^{L(0)}v,x_0)\cdot\nn
&&\qquad\qquad\qquad\cdot ((I^{z}_1\circ (1_{W_1}\otimes
I^{z}_2))'(e^{-(\log z)L'(0)}w'_{(4)}))\biggr)(w_{(1)}\otimes w_{(2)}\otimes 
e^{(\log z)L(0)}w_{(3)}),\nn
\end{eqnarray*}
and taking ${\rm Res}_{x_1}$ and using Remark
\ref{consequenceofPz1z2compat} (and in particular,
(\ref{resofconsequence})), and then applying $\res_{x_{0}^{-1}}$, we
obtain
\begin{eqnarray}\label{mu12-1.5-1}
\lefteqn{\res_{x_{0}^{-1}}\Biggl(\tau_{P(z_0+zz_2,zz_2)}\Biggl(
\dlti{x}{x_0^{-1}}{-zz_2}
Y_t((-x_0^2)^{L(0)}e^{-x_0L(1)}e^{xL(1)}(-x^{-2})^{L(0)}v,x_0)
\Biggr)\cdot}\nn
&&\qquad\qquad\qquad\cdot ((I^{z}_1\circ (1_{W_1}\otimes
I^{z}_2))'(e^{-(\log z)L'(0)}w'_{(4)}))\Biggr)(w_{(1)}\otimes w_{(2)}\otimes 
e^{(\log z)L(0)}w_{(3)})\nn
&&=\res_{x_{0}^{-1}}\dlti{x}{x_0^{-1}}{-zz_2}\biggl(Y'_{P(z_0+zz_2,zz_2)}
((-x_0^2)^{L(0)}e^{-x_0L(1)}e^{xL(1)}(-x^{-2})^{L(0)}v,x_0)\cdot\nn
&&\qquad\qquad\qquad\cdot ((I^{z}_1\circ (1_{W_1}\otimes
I^{z}_2))'(e^{-(\log z)L'(0)}w'_{(4)}))\biggr)(w_{(1)}\otimes w_{(2)}\otimes 
e^{(\log z)L(0)}w_{(3)})\nn
&&=\res_{x_{0}^{-1}}\dlti{x_0}{x^{-1}}{+zz_2}\biggl(Y'_{P(z_0+zz_2,zz_2)}
((-x_0^2)^{L(0)}e^{-x_0L(1)}e^{xL(1)}(-x^{-2})^{L(0)}v,x_0)\cdot\nn
&&\qquad\qquad\qquad\cdot ((I^{z}_1\circ (1_{W_1}\otimes
I^{z}_2))'(e^{-(\log z)L'(0)}w'_{(4)}))\biggr)(w_{(1)}\otimes w_{(2)}\otimes 
e^{(\log z)L(0)}w_{(3)}).\nn
\end{eqnarray}
The right-hand side of (\ref{mu12-1.5-1}) and thus the left-hand side,
which is the right-hand side of (\ref{mu12-1.4}), involves only
finitely many negative powers of $x$.  In particular, the sum of the
two terms on the left-hand side of (\ref{mu12-1.4}) involves only
finitely many negative powers of $x$ and hence lies in
$x^{m_{0}}\C[[x]]$ for some $m_{0}\in \Z$, so that the coefficients of
$x^{m}$ for $m < m_{0}$ in the left-hand side of (\ref{mu12-1.4})
cancel.  Moreover, we can multiply the right-hand side of
(\ref{mu12-1.4}) by (\ref{dltx1x-1-z0}), to obtain
\begin{eqnarray}\label{mu12-1.5-4}
\lefteqn{\dlt{x_1}{x^{-1}}{-z_0}{\rm
Res}_{x_0^{-1}}\Biggl(\tau_{P(z_0+zz_2,zz_2)}\Biggl(
\dlti{x}{x_0^{-1}}{-zz_2}
Y_t((-x_0^2)^{L(0)}e^{-x_0L(1)}e^{xL(1)}(-x^{-2})^{L(0)}v,x_0)
\Biggr)\cdot}\nn
&&\qquad\qquad\qquad\cdot ((I^{z}_1\circ (1_{W_1}\otimes
I^{z}_2))'(e^{-(\log z)L'(0)}w'_{(4)}))\Biggr)(w_{(1)}\otimes w_{(2)}\otimes 
e^{(\log z)L(0)}w_{(3)}).\qquad\;\;\;\nn
\end{eqnarray}

The first term in the left-hand side of (\ref{mu12-1.4}) is equal to 
\[
\sum_{n\in \R}\sum_{i=0}^{N}e^{-n\log z}
\frac{(-\log z)^{i}}{i!}g_{n, i}(x),
\]
where
\begin{equation}\label{g}
g_{n, i}(x)=Y'_{P(z_0)}(v, x) 
((L'_{P(z_{0})}(0)-n)^{i}\lambda^{(2)}_{n}[w_{(3)}])
(w_{(1)}\otimes w_{(2)})
\end{equation}
for $n\in \R$ and $i=0, \dots, N$, where we have double absolute
convergence for
\[
z \in O\cap P\cap Q\cap \Pi[w_{(3)}] \cap \Gamma[w_{(3)}] \cap R.
\]
By (\ref{Lambda}) and (\ref{mu12-1.5}), which we rewrite as in
(\ref{mu12-1.5.02}), the second term in the left-hand side of
(\ref{mu12-1.4}) is equal to
\begin{equation}\label{sumh}
\sum_{n\in \R}\sum_{i=0}^{\widetilde{N}}e^{-n\log z}
\frac{(-\log z)^{i}}{i!}h_{n, i}(x),
\end{equation}
where 
\begin{eqnarray}\label{h}
h_{n, i}(x)&=&{\rm Res}_{x_0^{-1}}\dlti{x}{-zz_2}{+x_0^{-1}}\cdot\nn
&&\quad\cdot \sum_{l \le L}\sum_{m=-M}^{M}
\sum_{s=-S}^{S}((L'_{P(z_{0})}(0)-(n+s+l))^{i}
\lambda^{(2)}_{n+s+l}[w^{l,m}_{s}])(w_{(1)}\otimes
w_{(2)})x_{0}^{l}x^{m}\nn
\end{eqnarray}
for $n\in \R$ and $i=0, \dots, \widetilde{N}$; here we have double
absolute convergence for 
\[
z \in O\cap P\cap Q\cap \widetilde{\Pi} \cap \widetilde{\Gamma}\cap R.
\]

For $N<i\le \max(N, \widetilde{N})$, set $g_{n, i}(x)=0$ and 
for $\widetilde{N}<i\le \max(N, \widetilde{N})$, set $h_{n, i}(x)=0$.
Let 
\[
f_{n, i}(x)=g_{n, i}(x)+h_{n, i}(x)
\]
for $n\in \R$ and $i=0, \dots, \max(N, \widetilde{N})$.
Then the left-hand side of (\ref{mu12-1.4}) is equal to 
\begin{equation}\label{mu12-1.5-2}
\sum_{n\in \R}\sum_{i=0}^{\max(N, \widetilde{N})}e^{-n\log z}
\frac{(-\log z)^{i}}{i!}f_{n, i}(x),
\end{equation}
with double absolute convergence for each power of $x$ for
\[
z \in O\cap P\cap Q\cap \Pi[w_{(3)}] \cap\widetilde{\Pi} \cap
\Gamma[w_{(3)}] \cap \widetilde{\Gamma}\cap R.
\]
That is, for each $m\in \Z$, 
\begin{equation}\label{mu12-1.5-2.1}
\sum_{n\in \R}\sum_{i=0}^{\widetilde{N}}e^{-n\log z}
\frac{(-\log z)^{i}}{i!}\res_{x}x^{-m-1}f_{n, i}(x)
\end{equation}
is absolutely convergent for $z$ in this set. 

Recall that the left-hand side of (\ref{mu12-1.4}) in fact lies in
$x^{m_{0}}\C[[x]]$.  Thus for $m<m_{0}$,
\[
\sum_{n\in \R}\sum_{i=0}^{\widetilde{N}}e^{-n\log z}
\frac{(-\log z)^{i}}{i!}\res_{x}x^{-m-1}f_{n, i}(x)=0
\]
for 
\[
z \in O\cap P\cap Q\cap \Pi[w_{(3)}] \cap\widetilde{\Pi} \cap \Gamma[w_{(3)}]
\cap \widetilde{\Gamma}\cap R.
\]
Since by Proposition \ref{real-exp-set} $\R\times \{0, \dots, \max(N,
\widetilde{N})\}$ is a unique expansion set,
\[
\res_{x}x^{-m-1}f_{n, i}(x)=0
\]
for $m<m_{0}$, $n\in \R$ and $i=0, \dots, \max(N, \widetilde{N})$,
and so
\[
f_{n, i}(x)\in x^{m_{0}}\C[[x]]
\]
for $n\in \R$ and $i=0, \dots, \max(N, \widetilde{N})$.  Thus for
\[
z \in O\cap P\cap Q\cap \Pi[w_{(3)}] \cap\widetilde{\Pi} \cap \Gamma[w_{(3)}]
\cap \widetilde{\Gamma}\cap R,
\]
\begin{equation}\label{sumf}
\sum_{n\in \R}\sum_{i=0}^{\max(N, \widetilde{N})}e^{-n\log z}
\frac{(-\log z)^{i}}{i!}
\dlt{x_1}{x^{-1}}{-z_0}f_{n, i}(x)
\end{equation}
is a formal Laurent series in $x$ and $x_{1}$ whose coefficients, as
finite linear combinations of the coefficients of powers of $x$ in
(\ref{mu12-1.5-2}), are equal to the corresponding (absolutely
convergent) coefficients in
\[
\dlt{x_1}{x^{-1}}{-z_0}\sum_{n\in \R}\sum_{i=0}^{\max(N, \widetilde{N})}e^{-n\log z}
\frac{(-\log z)^{i}}{i!}f_{n, i}(x).
\]
This expression equals the product of (\ref{dltx1x-1-z0}) and the
left-hand side of (\ref{mu12-1.4}), and hence (\ref{sumf}) is
absolutely convergent to (\ref{mu12-1.5-4}).

Moreover, from (\ref{mu12-1.5.0})--(\ref{mu12-1.5.02}),
(\ref{mu12-1.5-2}) and (\ref{mu12-1.5-2.1}), for
\[
z \in O\cap P\cap Q\cap \widetilde{\Pi} \cap \widetilde{\Gamma}\cap R,
\]
the coefficient of each monomial in $x$ and $x_1$ in
\begin{eqnarray*}
\lefteqn{\sum_{n\in \R}\sum_{i=0}^{\widetilde{N}}e^{-n\log z}
\frac{(-\log z)^{i}}{i!}\dlt{x_1}{x^{-1}}{-z_0}h_{n, i}(x)}\nn
&&=\sum_{n\in \R}\sum_{i=0}^{\max(N, \widetilde{N})}e^{-n\log z}
\frac{(-\log z)^{i}}{i!}\dlt{x_1}{x^{-1}}{-z_0}h_{n, i}(x)
\end{eqnarray*}
is absolutely convergent to the corresponding coefficient in
\begin{eqnarray*}
\lefteqn{\dlt{x_1}{x^{-1}}{-z_0}\sum_{n\in \R}\sum_{i=0}^{\widetilde{N}}e^{-n\log z}
\frac{(-\log z)^{i}}{i!}h_{n, i}(x)}\nn
&&=\dlt{x_1}{x^{-1}}{-z_0}\sum_{n\in \R}\sum_{i=0}^{\max(N, \widetilde{N})}e^{-n\log z}
\frac{(-\log z)^{i}}{i!}h_{n, i}(x).
\end{eqnarray*}

Assume that
\[
z \in O\cap P\cap Q\cap \Pi[w_{(3)}] \cap\widetilde{\Pi} \cap
\Gamma[w_{(3)}] \cap \widetilde{\Gamma}\cap R.
\]
We are now ready to bring the double sum over $n$ and $i$ in the
right-hand side of (\ref{mu12-2}) to the outside.  This right-hand
side equals
\begin{eqnarray}\label{righthandside}
\lefteqn{x^{-1}_1 \delta \bigg(\frac{x^{-1}-z_0}{x_1}\bigg)
\cdot}\nno\\
&&\cdot
\Biggl(Y'_{P(z_0)}(v, x) \Biggl(\sum_{n\in \R}\sum_{i=0}^{N}e^{-n\log
z} \frac{(-\log z)^{i}}{i!}
((L'_{P(z_{0})}(0)-n)^{i}\lambda^{(2)}_{n})[w_{(3)}]\Biggr)\Biggr)
(w_{(1)}\otimes w_{(2)});\nn
\end{eqnarray}
we recall that the coefficient of each monomial in $x$ and $x_1$ in
(\ref{righthandside}) is the sum of an absolutely convergent series,
each term of which involves the weakly absolutely convergent double
sum over $n$ and $i$.  Using the absolute convergence of the
coefficients in (\ref{mu12-1.2.5}) to those in (\ref{mu12-1.2.6}), we
rewrite (\ref{righthandside}) as
\begin{equation}\label{righthandsiderewrite}
x^{-1}_1 \delta \bigg(\frac{x^{-1}-z_0}{x_1}\bigg)
\sum_{n\in \R}\sum_{i=0}^{N}
e^{-n\log z} \frac{(-\log z)^{i}}{i!} g_{n, i}(x).
\end{equation}
What we need to show is that the coefficient of each monomial in $x$
and $x_1$ in
\begin{eqnarray}\label{mu12-1.3}
\lefteqn{\sum_{n\in \R}\sum_{i=0}^{N}e^{-n\log z}
\frac{(-\log z)^{i}}{i!}\cdot}\nno\\
&&\quad\quad\cdot 
x^{-1}_1\delta\biggl(\frac{x^{-1}-z_0}{x_1}\bigg)\Biggl(Y'_{P(z_0)}(v, x)
((L'_{P(z_{0})}(0)-n)^{i}\lambda^{(2)}_{n}[w_{(3)}])\Biggr)
(w_{(1)}\otimes w_{(2)})\nno\\
&&=\sum_{n\in \R}\sum_{i=0}^{N} e^{-n\log z} \frac{(-\log z)^{i}}{i!}
x^{-1}_1 \delta \bigg(\frac{x^{-1}-z_0}{x_1}\bigg)
g_{n, i}(x) 
\end{eqnarray}
is absolutely convergent and that it converges to the corresponding
(absolutely convergent) coefficient in (\ref{righthandsiderewrite}).
We have:
\begin{eqnarray*}
\lefteqn{\sum_{n\in \R}\sum_{i=0}^{N}e^{-n\log z}
\frac{(-\log z)^{i}}{i!}\dlt{x_1}{x^{-1}}{-z_0}g_{n, i}(x)}\nn
&&=\sum_{n\in \R}\sum_{i=0}^{\max(N, \widetilde{N})}e^{-n\log z}
\frac{(-\log z)^{i}}{i!}\dlt{x_1}{x^{-1}}{-z_0}g_{n, i}(x)\nn
&&=\sum_{n\in \R}\sum_{i=0}^{\max(N, \widetilde{N})}e^{-n\log z}
\frac{(-\log z)^{i}}{i!}\dlt{x_1}{x^{-1}}{-z_0}f_{n, i}(x)\nn
&&\quad -\sum_{n\in \R}\sum_{i=0}^{\max(N, \widetilde{N})}e^{-n\log z}
\frac{(-\log z)^{i}}{i!}\dlt{x_1}{x^{-1}}{-z_0}h_{n, i}(x),
\end{eqnarray*}
and the coefficient of each monomial in $x$ and $x_{1}$ is absolutely
convergent to the corresponding coefficient in
\begin{eqnarray*}
\lefteqn{\dlt{x_1}{x^{-1}}{-z_0}\sum_{n\in
\R}\sum_{i=0}^{\max(N, \widetilde{N})}e^{-n\log z}
\frac{(-\log z)^{i}}{i!}f_{n, i}(x)}\nn
&&\quad -\dlt{x_1}{x^{-1}}{-z_0}\sum_{n\in \R}\sum_{i=0}^{\max(N,
\widetilde{N})}e^{-n\log z}
\frac{(-\log z)^{i}}{i!}h_{n, i}(x)\nn
&&=\dlt{x_1}{x^{-1}}{-z_0}\sum_{n\in \R}\sum_{i=0}^{\max(N,
\widetilde{N})}e^{-n\log z}
\frac{(-\log z)^{i}}{i!}g_{n, i}(x)\nn
&&=\dlt{x_1}{x^{-1}}{-z_0}\sum_{n\in \R}\sum_{i=0}^{N}e^{-n\log z}
\frac{(-\log z)^{i}}{i!}g_{n, i}(x),
\end{eqnarray*}
as desired.  (Note that in particular, the coefficient of each
monomial in the inner expression $\dlt{x_1}{x^{-1}}{-z_0}g_{n, i}(x)$
in (\ref{mu12-1.3}) is absolutely convergent, since it is essentially
a sub-sum of the relevant absolutely convergent series.)

We have succeeded in bringing the double sums on both sides of
(\ref{mu12-2}) to the outside: For
\[
z \in O\cap P\cap Q\cap \Pi[w_{(3)}] \cap\widetilde{\Pi} \cap
\Gamma[w_{(3)}] \cap \widetilde{\Gamma}\cap R,
\]
\begin{eqnarray}\label{mu12-2-1}
\lefteqn{\sum_{n\in \R}\sum_{i=0}^{N}e^{-n\log z}
\frac{(-\log z)^{i}}{i!}\cdot}\nno\\
&&\quad\quad\cdot 
\Biggl(\tau_{P(z_0)}\biggl( x^{-1}_1
\delta\biggl(\frac{x^{-1}-z_0}{x_1}\bigg) Y_{t}(v, x)\biggr)
((L'_{P(z_{0})}(0)-n)^{i}\lambda^{(2)}_{n}[w_{(3)}])\Biggr)
(w_{(1)}\otimes w_{(2)})\nn
&&=\sum_{n\in \R}\sum_{i=0}^{N}e^{-n\log z}
\frac{(-\log z)^{i}}{i!}\cdot\nno\\
&&\quad\quad\cdot 
x^{-1}_1
\delta\biggl(\frac{x^{-1}-z_0}{x_1}\bigg)
\Biggl(Y'_{P(z_0)}(v, x)
((L'_{P(z_{0})}(0)-n)^{i}\lambda^{(2)}_{n}[w_{(3)}])\Biggr)
(w_{(1)}\otimes w_{(2)}),
\end{eqnarray}
where the coefficients of the monomials in $x$ and $x_{1}$ in the
double sums on both sides of (\ref{mu12-2-1}) are absolutely
convergent and are equal.

Now we are able to apply Proposition \ref{real-exp-set}.  Since
$\R\times \{0, \dots, N\}$ is a unique expansion set,
we conclude that the expansion coefficients of the double sums on 
the left- and right-hand
sides of (\ref{mu12-2-1}) are equal. In particular, taking $i=0$, we obtain
\begin{eqnarray*}
\lefteqn{\Biggl(\tau_{P(z_0)}\biggl( x^{-1}_1
\delta\biggl(\frac{x^{-1}-z_0}{x_1}\bigg) Y_{t}(v, x)\biggr)
\lambda^{(2)}_{n}[w_{(3)}]\Biggr)
(w_{(1)}\otimes w_{(2)})}\nno\\
&&=x^{-1}_1 \delta\Bigg(\frac{x^{-1}-z_0}{x_1}\bigg)
(Y'_{P(z_0)}(v, x)
\lambda^{(2)}_{n}[w_{(3)}])
(w_{(1)}\otimes w_{(2)})
\end{eqnarray*}
for each $n\in \R$, and so we have proved that each
$\lambda^{(2)}_{n}=\lambda^{(2)}_{n}[w_{(3)}]$ satisfies Part (b) of
the $P(z_0)$-compatibility condition.

This completes the proof.  \epfv

\begin{rema}\label{tensor4}
{\rm Here we relate the proof of Theorem \ref{9.7-1} to the
corresponding analysis for the special case treated in \cite{tensor4}.
The proof above of Part (b) of the $P(z_{0})$-compatibility condition
is a generalization of the proof of (14.51) in \cite{tensor4}.  In the
proof of (14.51) in \cite{tensor4}, for a series of the form
\begin{equation}\label{sumoverninD}
\sum_{n\in D}a_{n}z^{-n}
\end{equation}
where $D$ is a strictly increasing sequence in $\R$, in order to
determine the coefficients $a_{n} \in \C$ uniquely from the sum of the
series, the series is required to be absolutely convergent in an open
set of the form $0<|z^{-1}|<r$ because Lemma 14.5 in \cite{tensor4}
was proved in \cite{tensor4} only for such a series. (Here the first
author would like to correct some minor mistakes in \cite{tensor4}:
First, in the $P(z_{2})$-local grading-restriction condition
(respectively, in the $P(z_{1}-z_{2})$-local grading-restriction
condition) in \cite{tensor4}, we should require that the series
depending on $z'$ obtained by applying $e^{z'L'_{P(z_{2})}(0)}$
(respectively, $e^{z'L'_{P(z_{1}-z_{2})}(0)}$) to each term of the
weakly absolutely convergent series of $P(z_{2})$-weight vectors in
$(W_{2}\otimes W_{3})^{*}$ (respectively, in $(W_{1}\otimes
W_{2})^{*}$) be weakly absolutely convergent for $z'$ in a
neighborhood of $z'=0$.  This is implicitly used in the proof of
(14.51) in \cite{tensor4} and follows easily from the convergence and
extension property in \cite{tensor4}.  But it is not clear to the
first author whether this can be proved by assuming the
$P(z_{2})$-local grading-restriction condition or the
$P(z_{1}-z_{2})$-local grading-restriction condition in \cite{tensor4}
for all $z_{1}$ and $z_{2}$ satisfying
$|z_{1}|>|z_{2}|>|z_{1}-z_{2}|>0$.  Second, in the proof of (14.51) in
\cite{tensor4}, the domain
\[
0<|z|<\frac{|z_{0}|}{2|z_{2}|}
\]
should be replaced by the intersection of 
\[
0<|z^{-1}|<\frac{|z_{2}|}{|z_{0}|} \; \; \; (>1)
\]
and 
\[
|z_{0}+zz_{2}|>|zz_{2}|>0;
\]
since the expansion (14.48) and consequently the series in (14.49) is
of the form (\ref{sumoverninD}) rather than $\sum_{n\in D}a_{n}z^{n}$,
the correct domain is
\[
0<|z^{-1}|<\frac{|z_{2}|}{|z_{0}|}.
\]
The reason why the right-hand side of (14.48) is absolutely convergent
in this domain and not just in its intersection with
\[
|z_{0}+zz_{2}|>|zz_{2}|>0
\]
is that when $D$ is a strictly increasing sequence in $\R$, the
absolute convergence of (\ref{sumoverninD}) at one particular $z$ such
that
\[
|z^{-1}|=r\ne 0
\]
implies that it is also absolutely convergent at any $z$ satisfying 
\[
0<|z^{-1}|\le r.)
\]
However, in our proof of Theorem \ref{9.7-1} above, because $\R\times
\{0, \dots, N\}$ is a unique expansion set, the double absolute
convergence of
\[
\sum_{n\in \R}\sum_{i=0}^{N}a_{n, i}z^{-n}(-\log z)^{i}
\]
for $z$ in {\it any} nonempty open subset of $\C^{\times}$, not
necessarily containing an open subset of the form $0<|z^{-1}|<r$,
implies that the coefficients $a_{n, i}$ are uniquely determined by
the sum of the series; here we do not need the absolute convergence of
the series for $z^{-1}$ near $0$.  But our convergence-condition
assumption gives only the absolute convergence of {\it iterated}
series of the form
\[
\sum_{n\in \R}\left(\sum_{i=0}^{N}a_{n, i}z^{-n}(-\log z)^{i}\right)
\]
in a nonempty open set, and we had to prove the absolute convergence
of the corresponding double series
\[
\sum_{n\in \R}\sum_{i=0}^{N}a_{n, i}z^{-n}(-\log z)^{i}
\]
in the same open set using Proposition
\ref{log-coeff-conv<=>iterate-conv} (or Corollary
\ref{double-conv<=>iterate-conv}). This proof of the absolute
convergence of these double series is one of the hard parts of the
proof of Theorem \ref{9.7-1}, and this was not needed in the proof of
(14.51) in \cite{tensor4} because there we have no finite sum over
powers of $\log z$.  Another hard part of the proof above is to show
that (\ref{mu12-2}) implies that each $\lambda_{n}^{(2)}$ satisfies
the $P(z_{0})$-compatibility condition.  This part of the proof
amounts to a proof that certain triple series are absolutely
convergent, so that suitable iterated sums exist and are equal.  (In
fact, this part of the proof was also needed in the proof of (14.51)
in \cite{tensor4} (with double rather than triple sums, since there
are no finite sums over powers of $\log z$) but was not given
there. The last part of the proof of Theorem \ref{9.7-1} above gives
this missing detail, in our present much more general case.)  Also,
even in the case considered in \cite{tensor4}, the proof of Theorem
\ref{9.7-1} above establishes a stronger statement than (14.51) in
\cite{tensor4}: For each $n\in \R$, $\lambda_{n}^{(2)}$ satisfies the
$P(z_{0})$-compatibility condition even if $(I_{1} \circ
(1_{W_{1}}\otimes I_{2}))'(w'_{(4)})$ is not assumed to satisfy Part
(b) of the $P(z_{0})$-local grading restriction condition.  Finally,
we comment that the proof above certainly also proves (14.51) in
\cite{tensor4} as a special case.}
\end{rema}

\begin{rema}\label{pf-unique-lambda-n}
{\rm We now use the part of the proof of Theorem \ref{9.7-1} from
(\ref{iter-sum}) to (\ref{k-th-der-at-0}) to prove the part of
Proposition \ref{unique-lambda-n} on the uniqueness of the elements
$\lambda^{(2)}_{n}$, $n\in \R$, with the properties indicated in Part
(a) of the $P^{(2)}(z)$-local grading restriction condition; the other
three cases are handled the same way.  Using the proof from
(\ref{iter-sum}) to (\ref{k-th-der-at-0}) with $z_{0}$ and $-l^{0}(z)$
replaced by $z$ and $z'$, respectively, we have that the sum of
\begin{eqnarray}\label{unique-lambda-n-1}
\lefteqn{\sum_{n\in \R}(e^{z'L'_{P(z)}(0)}\lambda^{(2)}_{n})(w_{(1)}\otimes w_{(2)})}\nn
&&=\sum_{n\in \R}e^{nz'}\left(\left(\sum_{i=0}^{N}\frac{(z')^{i}}{i!}
(L'_{P(z)}(0)-L'_{P(z)}(0)_{s})^{i}\lambda^{(2)}_{n}
\right)(w_{(1)}\otimes 
w_{(2)})\right)
\end{eqnarray}
is an analytic function of $z'$ for $z'$ in an open neighborhood of
$0$, that its $k$-th derivative with respect to $z'$ at $z'=0$ is the
sum of the absolutely convergent series
\begin{equation}\label{unique-lambda-n-2}
\sum_{n\in \R}(L'_{P(z)}(0)^{k}
\lambda^{(2)}_{n})(w_{(1)}\otimes w_{(2)}),
\end{equation}
and that the iterated sum on the right-hand side of
(\ref{unique-lambda-n-1}) equals the corresponding double sum,
absolutely convergent in a suitably small neighborhood of $z'=0$
independent of $w_{(1)}$ and $w_{(2)}$.  Using (\ref{LP'(j)})
repeatedly and then using (iii) of Part (a) of the $P^{(2)}(z)$-local
grading restriction condition, we see that (\ref{unique-lambda-n-2})
is equal to
\begin{eqnarray}\label{unique-lambda-n-3}
\lefteqn{\sum_{i=0}^{k}{k\choose i}\sum_{n\in \R}\lambda^{(2)}_{n}
((L(0)+zL(-1))^{k-i}w_{(1)}\otimes L(0)^{i}w_{(2)})}\nn
&&=\sum_{i=0}^{k}{k\choose i}\mu^{(2)}_{\lambda,
w_{(3)}}((L(0)+zL(-1))^{k-i}w_{(1)}\otimes L(0)^{i}w_{(2)}).
\end{eqnarray}
Since the right-hand side of (\ref{unique-lambda-n-3}) is independent
of $\lambda_{n}^{(2)}$, $n\in \R$, the analytic function obtained from
the double sum corresponding to (\ref{unique-lambda-n-1}) is also
independent of $\lambda_{n}^{(2)}$, $n\in \R$, that is, if the formal
series $\sum_{n\in \R}\lambda^{(2)}_{n}$ and $\sum_{n\in
\R}\tilde{\lambda}^{(2)}_{n}$ both satisfy Part (a) of the
$P^{(2)}(z)$-local grading restriction condition, then
\[
\sum_{n\in \R}(e^{z'L'_{P(z)}(0)}\lambda^{(2)}_{n})(w_{(1)}\otimes w_{(2)})
\]
and 
\[
\sum_{n\in \R}(e^{z'L'_{P(z)}(0)}\tilde{\lambda}^{(2)}_{n})(w_{(1)}\otimes w_{(2)})
\]
are analytic functions equal to each other in a suitably small open
neighborhood of $z'=0$.  Thus in this neighborhood,
\begin{eqnarray}\label{unique-lambda-n-4}
\lefteqn{\sum_{n\in \R}\sum_{i=0}^{N}e^{nz'}\frac{(z')^{i}}{i!}
((L'_{P(z)}(0)-L'_{P(z)}(0)_{s})^{i}(\lambda^{(2)}_{n}-\tilde{\lambda}^{(2)}_{n}))
(w_{(1)}\otimes w_{(2)})}
\nn
&&\quad\quad\quad\quad\quad
 =\sum_{n\in \R}(e^{z'L'_{P(z)}(0)}(\lambda^{(2)}_{n}-\tilde{\lambda}^{(2)}_{n}))
(w_{(2)}\otimes w_{(3)})\quad\quad\quad\quad\quad\quad\nn
&&\quad\quad\quad\quad\quad=0,
\end{eqnarray}
and so by Proposition \ref{real-exp-set}, $\R\times \{0, \dots, N\}$
being a unique expansion set, we have
\[((L'_{P(z)}(0)-L'_{P(z)}(0)_{s})^{i}(\lambda^{(2)}_{n}-\tilde{\lambda}^{(2)}_{n}))
(w_{(1)}\otimes w_{(2)})=0
\]
for $n\in \R$ and $i=0, \dots, N$. In particular (for $i=0$),
\[
\lambda^{(2)}_{n}-\tilde{\lambda}^{(2)}_{n}=0
\]
for $n\in \R$, proving the uniqueness.}
\end{rema}

We will be invoking the uniqueness (Proposition \ref{unique-lambda-n})
and the bilinearity (Corollary \ref{bilincorol}) of the elements
$\lambda^{(1)}_{n}$ and $\lambda^{(2)}_{n}$ below.

We now relate Proposition \ref{9.7} and Theorem \ref{9.7-1} to
$\hboxtr_{P(z)}$ and $\boxtimes_{P(z)}$ for suitable $z\in
\C^{\times}$; recall Definitions \ref{def-hboxtr} and \ref{pz-tp}.  We
will sometimes use Definition \ref{doublygraded} and Remark
\ref{submodstrgraded}.  First we relate Proposition \ref{9.7} to
$\hboxtr_{P(z)}$, and this will serve as motivation for Corollary
\ref{bar-boxbackslash}, in which we relate Theorem \ref{9.7-1} to
$\hboxtr_{P(z)}$.

\begin{rema}
{\rm Assume that $\mathcal{C}$ is closed under images.  In the setting
and under all the assumptions of Proposition \ref{9.7}, we have
(according to this result): If $\lambda=(I_1\circ (1_{W_1}\otimes
I_2))'(w'_{(4)})$ (respectively, $\lambda=(I^1\circ (I^2\otimes
1_{W_3}))'(w'_{(4)})$), then $W^{(1)}_{\lambda, w_{(1)}}$
(respectively, $W^{(2)}_{\lambda, w_{(3)}}$) is a generalized
$V$-submodule of an object of $\mathcal{C}$ included in $(W_{2}\otimes
W_{3})^{*}$ (respectively, included in $(W_{1}\otimes W_{2})^{*}$),
and in particular, for each $n\in \R$ the $P(z_{2})$-generalized
weight vector $\lambda_{n}^{(1)}$ (respectively, the
$P(z_{0})$-generalized weight vector $\lambda^{(2)}_{n}$), which
generates a generalized $V$-submodule of $W^{(1)}_{\lambda, w_{(1)}}$
(respectively, of $W^{(2)}_{\lambda, w_{(3)}}$), also generates a
generalized $V$-submodule of an object of $\mathcal{C}$ included in
$(W_2\otimes W_3)^{*}$ (respectively, $(W_1\otimes W_2)^{*}$).  Hence
by Proposition \ref{backslash=union},
\[
\lambda_{n}^{(1)} \in W_2\hboxtr_{P(z_2)} W_3
\]
and
\[
\lambda_{n}^{(2)} \in W_1\hboxtr_{P(z_0)} W_2
\]
for each $n\in \R$.}
\end{rema}

Invoking the last assertion of Proposition \ref{backslash=union}, we
have the following corollary of Theorem \ref{9.7-1}:

\begin{corol}\label{bar-boxbackslash}
Assume that the convergence condition for
intertwining maps in $\mathcal{C}$ holds and that
\[
|z_1|>|z_2|>|z_{0}|>0.
\]
Let $W_{1}$, $W_{2}$, $W_{3}$, $W_{4}$, $M_{1}$ and $M_{2}$ be objects
of $\mathcal{C}$ and let $I_{1}$, $I_{2}$, $I^1$ and $I^2$ be
$P(z_1)$-, $P(z_2)$-, $P(z_2)$- and $P(z_0)$-intertwining maps of
types ${W_4}\choose {W_1M_1}$, ${M_1}\choose {W_2W_3}$, ${W_4}\choose
{M_2W_3}$ and ${M_2}\choose {W_1W_2}$, respectively.  Let $w'_{(4)}\in
W'_4$.
\begin{enumerate}

\item
Suppose that
\[
\lambda=(I_1\circ (1_{W_1}\otimes I_2))'(w'_{(4)})
\]
satisfies the (full) $P^{(2)}(z_0)$-local grading restriction
condition (or the $L(0)$-semisimple $P^{(2)}(z_0)$-local grading
restriction condition when $\mathcal{C}$ is in $\mathcal{M}_{sg}$).
For any $w_{(3)}\in W_{3}$, let $\sum_{n\in \R}\lambda_{n}^{(2)}$ be
the (unique) series weakly absolutely convergent to
\[
\mu^{(2)}_{\lambda, w_{(3)}} \in (W_{1}\otimes W_{2})^{*}
\]
as indicated in the $P^{(2)}(z_0)$-grading condition (or the
$L(0)$-semisimple $P^{(2)}(z_0)$-grading condition).  If for each $n
\in \R$ the generalized $V$-submodule of the generalized $V$-module
$W^{(2)}_{\lambda, w_{(3)}}$ (given by Theorem \ref{9.7-1}) generated
by $\lambda^{(2)}_{n}$ is a generalized $V$-submodule of some object
of $\mathcal{C}$ (depending on $n$) included in $(W_{1}\otimes
W_{2})^{*}$, then
\[
\lambda^{(2)}_{n} \in W_1\hboxtr_{P(z_0)} W_2.
\]

\item
Analogously, suppose that
\[
\lambda=(I^1\circ (I^2\otimes 1_{W_3}))'(w'_{(4)})
\]
satisfies the (full) $P^{(1)}(z_2)$-local grading restriction
condition (or the $L(0)$-semisimple $P^{(1)}(z_2)$-local grading
restriction condition when $\mathcal{C}$ is in $\mathcal{M}_{sg}$).
For any $w_{(1)}\in W_{1}$, let $\sum_{n\in \R}\lambda_{n}^{(1)}$ be
the (unique) series weakly absolutely convergent to
\[
\mu^{(1)}_{\lambda, w_{(1)}} \in (W_{2}\otimes W_{3})^{*}
\]
as indicated in the $P^{(1)}(z_2)$-grading condition (or the
$L(0)$-semisimple $P^{(1)}(z_2)$-grading condition).  If for each $n
\in \R$ the generalized $V$-submodule of the generalized $V$-module
$W^{(1)}_{\lambda, w_{(1)}}$ (given by Theorem \ref{9.7-1}) generated
by $\lambda^{(1)}_{n}$ is a generalized $V$-submodule of some object
of $\mathcal{C}$ (depending on $n$) included in $(W_{2}\otimes
W_{3})^{*}$, then
\[
\lambda^{(1)}_{n} \in W_2\hboxtr_{P(z_2)} W_3.
\]
\epfv

\end{enumerate}

\end{corol}

Next we shall express a product of suitable intertwining maps as an
iterate and vice versa.  This is accomplished in Theorem
\ref{lgr=>asso} below.  We shall actually carry this out only for the
case of expressing a product as an iterate, which is Part 1 of Theorem
\ref{lgr=>asso}; this case is based on Lemma \ref{intertwine-tau}
below.  Expresing an iterate as a product (Part 2 of the theorem) is
proved analogously.  We start with hypotheses for the lemma and the
theorem.

Assume that $\mathcal{C}$ is closed under images, that the convergence
condition for intertwining maps in $\mathcal{C}$ holds, and that
\[
|z_1|>|z_2|>|z_{0}|>0.
\]
Let $W_{1}$, $W_{2}$, $W_{3}$,
$W_{4}$ and $M_{1}$ be objects of $\mathcal{C}$ and let
$I_{1}$ and $I_{2}$ be $P(z_1)$- and $P(z_2)$-intertwining maps of 
types ${W_4}\choose {W_1M_1}$ and
${M_1}\choose {W_2W_3}$, respectively.  Set
\[
G=(I_1\circ (1_{W_1}\otimes I_2))'\in
\hom (W'_4, (W_1\otimes W_2\otimes W_3)^*)
\]
(cf. Remark \ref{I1I2'}).  Suppose that $W_1\boxtimes_{P(z_0)} W_2$
exists in $\mathcal{C}$ and that for each $w'_{(4)}\in W_{4}'$,
$G(w'_{(4)})$, as in Corollary \ref{bar-boxbackslash}, satisfies the
$P^{(2)}(z_0)$-local grading restriction condition (or the
$L(0)$-semisimple $P^{(2)}(z_0)$-local grading restriction condition
when $\mathcal{C}$ is in $\mathcal{M}_{sg}$). For $w_{(3)}\in W_{3}$,
let
\begin{equation}\label{lambdan2w'w}
\sum_{n\in \R}\lambda_{n}^{(2)}(w'_{(4)}, w_{(3)})
\end{equation}
be the (unique) series weakly absolutely convergent to
\begin{equation}\label{mu2G}
\mu^{(2)}_{G(w'_{(4)}), w_{(3)}} \in (W_1\otimes W_2)^*
\end{equation}
as indicated in the $P^{(2)}(z_0)$-grading condition (or the
$L(0)$-semisimple $P^{(2)}(z_0)$-grading condition).  Suppose further
that for each $n\in \R$, $w'_{(4)}\in W_{4}'$ and $w_{(3)}\in W_{3}$,
the generalized $V$-submodule of $W^{(2)}_{G(w'_{(4)}), w_{(3)}}$
generated by $\lambda_{n}^{(2)}(w'_{(4)}, w_{(3)})$ is a generalized
$V$-submodule of some object of $\mathcal{C}$ included in
$(W_{1}\otimes W_{2})^{*}$.  Using Part 1 of Corollary
\ref{bar-boxbackslash}, which is based on and follows from Part 1 of
Theorem \ref{9.7-1}, we shall now prove that the product $I_1\circ
(1_{W_1}\otimes I_2)$ of the intertwining maps $I_{1}$ and $I_{2}$ can
be written as an iterate of suitable intertwining maps, which is Part
1 of Theorem \ref{lgr=>asso} below. First we formulate and prove a
lemma under these assumptions.  This lemma is the core of the proof of
the theorem.

Recall from Proposition \ref{tensor1-13.7} that since $\mathcal{C}$ is
closed under images, the existence of $W_1\boxtimes_{P(z_0)} W_2$ in
$\mathcal{C}$ implies that $W_1\hboxtr_{P(z_0)} W_2$ is an object of
$\mathcal{C}$ and that
\begin{equation}\label{boxtensor=backslash'}
W_1\boxtimes_{P(z_0)} W_2 = (W_1\hboxtr_{P(z_0)} W_2)'.
\end{equation}
By Corollary \ref{bar-boxbackslash}, for any $w'_{(4)}\in W_{4}'$ and
$w_{(3)}\in W_{3}$, $\lambda_{n}^{(2)}(w'_{(4)}, w_{(3)})$ is an element (of generalized
weight $n$) of
\[
W_1\hboxtr_{P(z_0)} W_2 \subset (W_1\otimes W_2)^*
\]
for $n \in \R$.  Thus we have the element
\begin{equation}\label{mutilde2G}
\widetilde{\mu}^{(2)}_{G(w'_{(4)}), w_{(3)}}\in 
\overline{W_1\hboxtr_{P(z_0)} W_2},
\end{equation}
the formal completion of the generalized $V$-module (and object of
$\mathcal{C}$) $W_1\hboxtr_{P(z_0)} W_2$, whose homogeneous components
are the $P(z_{0})$-generalized weight vectors
\[
\lambda_{n}^{(2)}(w'_{(4)}, w_{(3)})\in W_1\hboxtr_{P(z_0)} W_2,
\]
$n\in \R$.  Recalling the notation
\[
\pi_n: \overline{W_1\hboxtr_{P(z_0)} W_2} \to W_1\hboxtr_{P(z_0)} W_2
\]
from Definition \ref{Wbardef}, we have
\begin{equation}\label{lambda-tilde-mu}
\lambda_{n}^{(2)}(w'_{(4)}, w_{(3)})=
\pi_n(\widetilde{\mu}^{(2)}_{G(w'_{(4)}), w_{(3)}}),
\end{equation}
and so from (\ref{lambdan2w'w}) and (\ref{mu2G}) we have the weakly
absolutely convergent series
\begin{equation}\label{musumpi}
\mu^{(2)}_{G(w'_{(4)}), w_{(3)}}=
\sum_{n\in \R}\pi_n(\widetilde{\mu}^{(2)}_{G(w'_{(4)}), w_{(3)}})
=\sum_{n\in \R}\lambda_{n}^{(2)}(w'_{(4)}, w_{(3)}).
\end{equation}
Note the distinction between the different sums (\ref{mutilde2G}) and
(\ref{musumpi}) of the same elements $\lambda_{n}^{(2)}(w'_{(4)},
w_{(3)})$; they take place in different spaces.  Recall from
Definition \ref{defofWprime} and (\ref{boxtensor=backslash'}) that we
have a canonical pairing
\[
\langle \cdot,\cdot \rangle_{W_1\hboxtr_{P(z_0)} W_2}
\]
between $W_1 \boxtimes_{P(z_0)} W_2$ and $\overline{W_1\hboxtr_{P(z_0)}
W_2}$.  By Corollary \ref{bilincorol}, the element (\ref{mutilde2G})
depends bilinearly on $w'_{(4)}$ and $w_{(3)}$, so that we have a
linear map
\[
\widetilde{G}: W_1\boxtimes_{P(z_0)} W_2\to (W'_4\otimes
W_3)^*
\]
determined by the condition
\begin{equation}\label{tildeG}
\widetilde{G}(w)(w'_{(4)}\otimes w_{(3)})=\langle
w,\widetilde{\mu}^{(2)}_{G(w'_{(4)}),w_{(3)}}\rangle_{W_1\hboxtr_{P(z_0)} W_2}
\end{equation}
for $w\in W_1\boxtimes_{P(z_0)} W_2$.  Moreover, generalizing the
corresponding lemma in \cite{tensor4}, we have the following lemma, in
which $\tau_{Q(z_2)}$ (recall Section 5.1, in particular, Definition
\ref{deftauQ}) appears naturally:

\begin{lemma}\label{intertwine-tau}
Under these assumptions,  the linear map
\[
\widetilde{G}\in\hom(W_1\boxtimes_{P(z_0)}
W_2, (W'_4\otimes W_3)^*)
\]
intertwines the actions $\tau_{W_1\boxtimes_{P(z_0)} W_2}$ and
$\tau_{Q(z_2)}$ of $V\otimes \iota_+{\mathbb C}[t,t^{-1},
(z_2+t)^{-1}]$, and also the corresponding
$\mathfrak{s}\mathfrak{l}(2)$ actions, on $W_1\boxtimes_{P(z_0)} W_2$
and on $(W'_4\otimes W_3)^*$ (recall Section 5.1 for these
actions).
\end{lemma}
\pf We shall prove the assertion about $\tau_{W_1\boxtimes_{P(z_0)}
W_2}$ and $\tau_{Q(z_2)}$, and at the end of this proof we shall
briefly comment that one can similarly prove the assertion about
$\mathfrak{s}\mathfrak{l}(2)$ by considering appropriate aspects of
the case $v=\omega$.

As in Proposition \ref{tensor1-13.7}, we shall denote the vertex
operator map for $W_1\boxtimes_{P(z_0)} W_2$ by $Y_{P(z_{0})}$.  By
(\ref{tildeG}) and the definition (\ref{(5.1)}) of the
$\tau_{Q(z_2)}$-action (see also (\ref{5.2})), we need to show that
\begin{eqnarray}\label{needtoshow}
\lefteqn{z^{-1}_2\delta\left(\frac{x_1-x_0}{z_2}\right)
\widetilde{G}(Y_{P(z_0)}(v, x_0)w)}\nno\\
&&=\tau_{Q(z_2)}\Biggl(
z^{-1}_2\delta\left(\frac{x_1-x_0}{z_2}\right) Y_{t}(v,
x_0)\Biggr)\widetilde{G}(w)
\end{eqnarray}
for $v\in V$ and $w\in W_1\boxtimes_{P(z_0)} W_2$, or equivalently,
that
\begin{eqnarray}\label{as:need0}
\lefteqn{\left\langle z^{-1}_2\delta\left(\frac{x_1-x_0}{z_2}\right)
Y_{P(z_0)}(v, x_0)w, \widetilde{\mu}^{(2)}_{G(w_{(4)}'),
w_{(3)}} \right\rangle_{W_1\shboxtr_{P(z_0)}W_2}}\nno\\
&&=\left(\tau_{Q(z_2)}\Biggl(
z^{-1}_2\delta\left(\frac{x_1-x_0}{z_2}\right) Y_{t}(v,
x_0)\Biggr)\widetilde{G}(w)\right)(w_{(4)}'\otimes w_{(3)})
\end{eqnarray}
for $v\in V$, $w\in W_1\boxtimes_{P(z_0)} W_2$, $w_{(4)}'\in W_{4}'$,
$w_{(3)}\in W_{3}$.  Note that the left-hand sides of
(\ref{needtoshow}) and of (\ref{as:need0}) involve only finitely many
negative powers of $x_0$.

By Proposition \ref{span}, we need only prove our assertion for
\[
w=\pi_{n}(w_{(1)}\boxtimes_{P(z_0)} w_{(2)})
\]
for any $n\in \R$ and
$w_{(1)}\in W_{1}$, $w_{(2)}\in W_{2}$.  (Again recall the notation
$\pi_n$ from Definition \ref{Wbardef}.)

By (\ref{boxpair}), for $n\in \R$ we have 
\begin{eqnarray}\label{*-to-box-0}
\lefteqn{(\lambda_{n}^{(2)}(w'_{(4)}, w_{(3)}))(w_{(1)}\otimes w_{(2)})}\nn
&&=\langle 
\lambda_{n}^{(2)}(w'_{(4)}, w_{(3)}), 
w_{(1)}\boxtimes_{P(z_{0})} w_{(2)}\rangle_{W_1\boxtimes_{P(z_0)} W_2}\nn
&&=\langle \lambda_{n}^{(2)}(w'_{(4)}, w_{(3)}), 
\pi_{n}(w_{(1)}\boxtimes_{P(z_{0})} w_{(2)})\rangle\nn
&&=\langle \widetilde{\mu}^{(2)}_{G(w'_{(4)}),
w_{(3)}}, 
\pi_{n}(w_{(1)}\boxtimes_{P(z_{0})} w_{(2)})\rangle,
\end{eqnarray}
where the last pairing is between $\overline{W_1\hboxtr_{P(z_0)} W_2}$
and $W_1\boxtimes_{P(z_0)} W_2$.

Recalling (\ref{3.21}), (\ref{yo}), (\ref{y'}) and the definition
(\ref{y'-p-z}) of $Y'_{P(z_0)}$ (see also (\ref{Y'def})), we define
\[
Y'^{o}_{P(z_{0})}(v, x_{0}): (W_{1}\otimes W_{2})^{*}\to 
(W_{1}\otimes W_{2})^{*}[[x_{0}, x_{0}^{-1}]]
\]
by
\begin{eqnarray}\label{y'o-ext}
Y'^{o}_{P(z_{0})}(v, x_{0})\mu&=&\tau_{P(z_{0})}(Y^{o}_{t}(v, x_{0}))\mu\nn
&=&Y'_{P(z_0)}(e^{x_0L(1)}(-x_0^{-2})^{L(0)}v,
x_0^{-1})\mu
\end{eqnarray}
for $\mu\in (W_{1}\otimes W_{2})^{*}$.
Then 
\begin{equation}\label{yio-y-boxbs}
Y'^{o}_{P(z_{0})}(v, x_{0})\lbar_{W_{1}\hboxtr_{P(z_{0})}W_{2}}
=Y^{o}_{W_{1}\hboxtr_{P(z_{0})}W_{2}}(v, x_{0}).
\end{equation}
We have natural maps
\[
Y'_{P(z_{0})}(v, x_{0}): \overline{W_1\hboxtr_{P(z_0)} W_2} \to 
\overline{W_1\hboxtr_{P(z_0)} W_2}[[x_{0}, x_{0}^{-1}]]
\]
and
\[
Y'^{o}_{P(z_{0})}(v, x_{0}): \overline{W_1\hboxtr_{P(z_0)} W_2} \to 
\overline{W_1\hboxtr_{P(z_0)} W_2}[[x_{0}, x_{0}^{-1}]];
\]
and 
\begin{eqnarray}\label{*-to-box}
\lefteqn{\pi_{n}(Y'^{o}_{P(z_{0})}(v, x_{0})\widetilde{\mu}^{(2)}_{G(w'_{(4)}),
w_{(3)}})
(w_{(1)}\otimes w_{(2)})}\nn
&&=\langle \pi_{n}(Y'^{o}_{P(z_0)}(v,x_0)\widetilde{\mu}^{(2)}_{G(w'_{(4)}),
w_{(3)}})), w_{(1)}\boxtimes_{P(z_{0})} w_{(2)}\rangle\nn
&&=\langle Y'^{o}_{P(z_{0})}(v, x_{0})
\widetilde{\mu}^{(2)}_{G(w'_{(4)}),
w_{(3)}}, 
\pi_{n}(w_{(1)}\boxtimes_{P(z_{0})} w_{(2)})\rangle
\end{eqnarray}
for $n\in \R$ (cf. (\ref{*-to-box-0})).

Now taking  
\[
w=\pi_{n}(w_{(1)}\boxtimes_{P(z_0)} w_{(2)})
\]
in the left-hand side of (\ref{as:need0}) and using  the definition 
of $Y_{P(z_{0})}$
(see Proposition \ref{tensor1-13.7}),
(\ref{*-to-box}) and
(\ref{yio-y-boxbs}) 
(and, as we have done above, 
dropping the subscripts for the
pairings), we see that the left-hand side of (\ref{as:need0}) becomes
\begin{eqnarray}\label{as:l}
\lefteqn{\left\langle z^{-1}_2\delta\left(\frac{x_1-x_0}{z_2}\right)
Y_{P(z_0)}(v, x_0)(\pi_{n}(w_{(1)}\boxtimes_{P(z_0)} w_{(2)})), 
\widetilde{\mu}^{(2)}_{G(w'_{(4)}),w_{(3)}} \right\rangle}\nno\\
&&=z^{-1}_2\delta\left(\frac{x_1-x_0}{z_2}\right)\langle
Y'^o_{P(z_0)}(v, x_0)\widetilde{\mu}^{(2)}_{G(w'_{(4)}),
w_{(3)}}, \pi_{n}(w_{(1)}\boxtimes_{P(z_0)}
w_{(2)})\rangle\nno\\
&&=z^{-1}_2\delta\left(\frac{x_1-x_0}{z_2}\right)
\pi_{n}(Y'^{o}_{P(z_{0})}(v, x_{0})\widetilde{\mu}^{(2)}_{G(w'_{(4)}),
w_{(3)}})(w_{(1)}\otimes
w_{(2)}).\nn
\end{eqnarray}
Recall that this expression involves only finitely many negative
powers of $x_0$.

Taking 
\[
w=\pi_{n}(w_{(1)}\boxtimes_{P(z_0)} w_{(2)})
\]
in the right-hand side of 
(\ref{as:need0}) and using the definition 
(\ref{5.2}) of
$\tau_{Q(z_2)}$, (\ref{tildeG}) and (\ref{*-to-box-0}), we obtain
\begin{eqnarray}\label{as:o}
\lefteqn{\bigg(\tau_{Q(z_2)}\biggl(z^{-1}_2\delta\left(\frac{x_1-x_0}
{z_2}\right) Y_{t}(v,x_0)\biggr)\widetilde{G}(\pi_{n}(w_{(1)}\boxtimes_{P(z_{0})}
w_{(2)}))\bigg)(w'_{(4)}\otimes
w_{(3)})}\nn
&&=x_0^{-1}\delta\bigg(\frac{x_1-z_2}{x_0}\bigg)
(\widetilde{G}(\pi_{n}(w_{(1)}\boxtimes_{P(z_{0})}
w_{(2)})))(Y'^o_4(v,x_1)w'_{(4)}\otimes
w_{(3)})\nn
&&\qquad-x_0^{-1}\delta\bigg(\frac{z_2-x_1}{-x_0}\bigg)
(\widetilde{G}(\pi_{n}(w_{(1)}\boxtimes_{P(z_{0})}
w_{(2)})))(w'_{(4)}\otimes
Y_3(v,x_1)w_{(3)})\nn
&&=x_0^{-1}\delta\bigg(\frac{x_1-z_2}{x_0}\bigg)\langle \pi_{n}(w_{(1)}\boxtimes_{P(z_{0})}
w_{(2)}),
\widetilde{\mu}^{(2)}_{G(Y'^o_4(v,x_1)w'_{(4)}),w_{(3)}}
\rangle\nn
&&\qquad-x_0^{-1}\delta\bigg(\frac{z_2-x_1}{-x_0}\bigg)\langle 
\pi_{n}(w_{(1)}\boxtimes_{P(z_{0})}
w_{(2)}),
\widetilde{\mu}^{(2)}_{G(w'_{(4)}),Y_3(v,x_1)w_{(3)}}
\rangle\nn
&&=x_0^{-1}\delta\bigg(\frac{x_1-z_2}{x_0}\bigg)
(\lambda_{n}^{(2)}(Y'^o_4(v,x_1)w'_{(4)}, w_{(3)}))(w_{(1)}
\otimes w_{(2)})\nn
&&\qquad-x_0^{-1}\delta\bigg(\frac{z_2-x_1}{-x_0}\bigg)
(\lambda_{n}^{(2)}(w'_{(4)}, Y_3(v,x_1)w_{(3)}))(w_{(1)}
\otimes w_{(2)}).
\end{eqnarray}

In order to prove the equality of (\ref{as:l}) and (\ref{as:o}), we
will consider the sums of both expressions over $n \in \R$.

By the definitions (\ref{y'o-ext}) 
and (\ref{Y'def}) of $Y^{\prime o}_{P(z_0)}$ and  $Y'_{P(z_0)}$, we
have, using (\ref{Zoo}),
\begin{eqnarray}\label{as:p}
\lefteqn{(Y'^{o}_{P(z_{0})}(v, x_{0})\mu)(w_{(1)}\otimes
w_{(2)})}\nn
&&=\big(Y'_{P(z_0)}(e^{x_0L(1)}(-x_0^{-2})^{L(0)}v,
x_0^{-1})\mu\big)(w_{(1)}\otimes w_{(2)})\nn
&&=
\mu(w_{(1)}\otimes Y_2(v,x_0)w_{(2)})\nno\\
&&\quad+\res_{x_2}z_0^{-1}
\delta\bigg(\frac{x_0-x_2}{z_0}\bigg)\mu
(Y_1(v, x_2)w_{(1)}\otimes w_{(2)})
\end{eqnarray}
for $\mu\in (W_{1}\otimes W_{2})^{*}$, $w_{(1)}\in W_{1}$ and $w_{(2)}\in W_{2}$.
Taking
\[
\mu=\mu^{(2)}_{G(w'_{(4)}),w_{(3)}},
\]
we thus have
\begin{eqnarray*}
\lefteqn{(Y'^{o}_{P(z_{0})}(v, x_{0})\mu^{(2)}_{G(w'_{(4)}),w_{(3)}})
(w_{(1)}\otimes w_{(2)})}\nn
&&=
\mu^{(2)}_{G(w'_{(4)}),w_{(3)}}(w_{(1)}\otimes Y_2(v,x_0)w_{(2)})\nno\\
&&\quad+\res_{x_2}z_0^{-1}
\delta\bigg(\frac{x_0-x_2}{z_0}\bigg)\mu^{(2)}_{G(w'_{(4)}),w_{(3)}}
(Y_1(v, x_2)w_{(1)}\otimes w_{(2)}).
\end{eqnarray*}
{}From (\ref{musumpi}), whose right-hand side is weakly absolutely
convergent,
\begin{eqnarray*}
\lefteqn{(Y'^{o}_{P(z_{0})}(v, x_{0})
\sum_{n\in \R}\pi_n(\widetilde{\mu}^{(2)}_{G(w'_{(4)}), w_{(3)}}))
(w_{(1)}\otimes w_{(2)})}\nn
&&=
\sum_{n\in \R}(\pi_n(\widetilde{\mu}^{(2)}_{G(w'_{(4)}), w_{(3)}})
(w_{(1)}\otimes Y_2(v,x_0)w_{(2)}))\nno\\
&&\quad+\res_{x_2}z_0^{-1}
\delta\bigg(\frac{x_0-x_2}{z_0}\bigg)
\sum_{n\in \R}(\pi_n(\widetilde{\mu}^{(2)}_{G(w'_{(4)}), w_{(3)}})
(Y_1(v, x_2)w_{(1)}\otimes w_{(2)}))\nno\\
&&=
\sum_{n\in \R}(\pi_n(\widetilde{\mu}^{(2)}_{G(w'_{(4)}), w_{(3)}})
(w_{(1)}\otimes Y_2(v,x_0)w_{(2)}))\nno\\
&&\quad+\sum_{n\in \R}\res_{x_2}z_0^{-1}
\delta\bigg(\frac{x_0-x_2}{z_0}\bigg)
(\pi_n(\widetilde{\mu}^{(2)}_{G(w'_{(4)}), w_{(3)}})
(Y_1(v, x_2)w_{(1)}\otimes w_{(2)})),
\end{eqnarray*}
since in the second term, the coefficient of each monomial in $x_2$
involves only the single infinite sum $\sum_{n\in \R}$, and so we can
switch $z_0^{-1}\delta\left(\frac{x_0-x_2}{z_0}\right)$ and
$\sum_{n\in \R}$.  Thus
\begin{eqnarray}\label{Y'oswitch}
\lefteqn{(Y'^{o}_{P(z_{0})}(v, x_{0})
\sum_{n\in \R}\pi_n(\widetilde{\mu}^{(2)}_{G(w'_{(4)}), w_{(3)}}))
(w_{(1)}\otimes w_{(2)})}\nn
&&=
\sum_{n\in \R}(Y'^{o}_{P(z_{0})}(v, x_{0})
\pi_n(\widetilde{\mu}^{(2)}_{G(w'_{(4)}), w_{(3)}}))
(w_{(1)}\otimes w_{(2)}),
\end{eqnarray}
with absolute convergence for the coefficient of each monomial in
$x_0$.

Now the product of the first term in the right-hand side of (\ref{as:p})
with $z_2^{-1}\delta\left(\frac{x_1-x_0}{z_2}\right)$ exists algebraically,
and since the product
\begin{equation}\label{as:p-3}
z_2^{-1}\delta\left(\frac{x_1-x_0}{z_2}\right)
z_0^{-1}\delta\bigg(\frac{x_0-x_2}{z_0}\bigg)
\end{equation}
exists in the sense of absolute convergence, by (\ref{l2-2}) in
Proposition \ref{deltalemma} (since $|z_2|>|z_{0}|>0$), the product of
the second term in the right-hand side of (\ref{as:p}) with
$z_2^{-1}\delta\left(\frac{x_1-x_0}{z_2}\right)$ exists in the sense
of absolute convergence.  (That is, the sum over the integral powers
of $x_0$ obtained from extracting the coefficient of any monomial in
$x_0$, $x_1$ and $x_2$ is absolutely convergent.)  Thus the product of
the left-hand side of (\ref{as:p}) with
$z_2^{-1}\delta\left(\frac{x_1-x_0}{z_2}\right)$ also exists in the
sense of absolute convergence.  Again taking
\[
\mu=\mu^{(2)}_{G(w'_{(4)}),
w_{(3)}},
\]
we thus have
\begin{eqnarray}\label{as:p-1}
\lefteqn{z_2^{-1}\delta\left(\frac{x_1-x_0}{z_2}\right)
(Y'^{o}_{P(z_{0})}(v, x_{0})\mu^{(2)}_{G(w'_{(4)}),
w_{(3)}})(w_{(1)}\otimes
w_{(2)})}\nn
&&=z_2^{-1}\delta\left(\frac{x_1-x_0}{z_2}\right)
\mu^{(2)}_{G(w'_{(4)}),
w_{(3)}}(w_{(1)}\otimes Y_2(v,x_0)w_{(2)})\nno\\
&&\quad+z_2^{-1}\delta\left(\frac{x_1-x_0}{z_2}\right)
\res_{x_2}z_0^{-1}\delta\bigg(\frac{x_0-x_2}{z_0}\bigg)\mu^{(2)}_{G(w'_{(4)}),
w_{(3)}}
(Y_1(v, x_2)w_{(1)}\otimes w_{(2)}),
\end{eqnarray}
in the sense of absolute convergence.  We will need a variant of this
formula, with the left-hand side replaced by $\sum_{n\in \R}$ applied
to (\ref{as:l}) (see (\ref{as:p-6}) below).

{}From (\ref{musumpi}), the right-hand side of (\ref{as:p-1}) is equal
to
\begin{eqnarray}\label{as:p-2}
\lefteqn{z_2^{-1}\delta\left(\frac{x_1-x_0}{z_2}\right)
\sum_{n\in \R}\pi_{n}(\widetilde{\mu}^{(2)}_{G(w'_{(4)}),
w_{(3)}})(w_{(1)}\otimes Y_2(v,x_0)w_{(2)})}\nno\\
&&\quad+z_2^{-1}\delta\left(\frac{x_1-x_0}{z_2}\right)
\res_{x_2}z_0^{-1}\delta\bigg(\frac{x_0-x_2}{z_0}\bigg)
\sum_{n\in \R}\pi_{n}(\widetilde{\mu}^{(2)}_{G(w'_{(4)}),
w_{(3)}})
(Y_1(v, x_2)w_{(1)}\otimes w_{(2)})\nn
&&=z_2^{-1}\delta\left(\frac{x_1-x_0}{z_2}\right)
\sum_{n\in \R}\pi_{n}(\widetilde{\mu}^{(2)}_{G(w'_{(4)}),
w_{(3)}})(w_{(1)}\otimes Y_2(v,x_0)w_{(2)})\nno\\
&&\quad+
\res_{x_2}z_2^{-1}\delta\left(\frac{x_1-x_0}{z_2}\right)
z_0^{-1}\delta\bigg(\frac{x_0-x_2}{z_0}\bigg)\sum_{n\in \R}
\pi_{n}(\widetilde{\mu}^{(2)}_{G(w'_{(4)}),
w_{(3)}})
(Y_1(v, x_2)w_{(1)}\otimes w_{(2)}),\nn
\end{eqnarray}
where the sums over $n \in \R$ are absolutely convergent.

In the first term in the right-hand side of (\ref{as:p-2}), the
coefficient of each monomial in $x_{0}$ and $x_{1}$ involves only the
single infinite sum $\sum_{n\in \R}$, and so we can switch
$z_2^{-1}\delta\left(\frac{x_1-x_0}{z_2}\right)$ and $\sum_{n\in \R}$
in this term.  In the second term, both (\ref{as:p-3}) and
\[
\sum_{n\in \R}
\pi_{n}(\widetilde{\mu}^{(2)}_{G(w'_{(4)}),
w_{(3)}})
(Y_1(v, x_2)w_{(1)}\otimes w_{(2)})
\]
are formal series in $x_{0}$, $x_{1}$ and $x_{2}$ each of whose
coefficients is an absolutely convergent series, and both series are
truncated from below in powers of $x_{2}$.  Thus the double sums
obtained from the coefficients of the product of these two formal
series in $x_{0}$, $x_{1}$ and $x_{2}$ are also absolutely convergent
and in particular, we can switch (\ref{as:p-3}) and $\sum_{n\in \R}$
in the second term.  So we see, using (\ref{as:p}), that the
right-hand side of (\ref{as:p-2}) is equal to
\begin{eqnarray}\label{as:p-4}
\lefteqn{\sum_{n\in \R}z_2^{-1}\delta\left(\frac{x_1-x_0}{z_2}\right)
\pi_{n}(\widetilde{\mu}^{(2)}_{G(w'_{(4)}),
w_{(3)}})(w_{(1)}\otimes Y_2(v,x_0)w_{(2)})}\nno\\
&&\quad+\sum_{n\in \R}
z_2^{-1}\delta\left(\frac{x_1-x_0}{z_2}\right)\res_{x_2}
z_0^{-1}\delta\bigg(\frac{x_0-x_2}{z_0}\bigg)
\pi_{n}(\widetilde{\mu}^{(2)}_{G(w'_{(4)}),
w_{(3)}})
(Y_1(v, x_2)w_{(1)}\otimes w_{(2)})\nn
&&=\sum_{n\in \R}z_2^{-1}\delta\left(\frac{x_1-x_0}{z_2}\right)
(Y'^{o}_{P(z_{0})}(v, x_{0})\pi_{n}(\widetilde{\mu}^{(2)}_{G(w'_{(4)}),
w_{(3)}}))(w_{(1)}\otimes
w_{(2)}),
\end{eqnarray}
and the corresponding double sums in the two terms in the left-hand
side and thus the corresponding double sums in the right-hand side are
all absolutely convergent.  Hence, using the fact that the coefficient
of each monomial in the right-hand side of (\ref{Y'oswitch}) is
absolutely convergent, we obtain that the right-hand side of
(\ref{as:p-4}) is equal to
\begin{eqnarray}\label{as:p-5}
\lefteqn{\sum_{n\in \R}z_2^{-1}\delta\left(\frac{x_1-x_0}{z_2}\right)
(Y'^{o}_{P(z_{0})}(v, x_{0})\pi_{n}(\widetilde{\mu}^{(2)}_{G(w'_{(4)}),
w_{(3)}}))(w_{(1)}\otimes
w_{(2)})}\nn
&&=z_2^{-1}\delta\left(\frac{x_1-x_0}{z_2}\right)\sum_{n\in \R}
(Y'^{o}_{P(z_{0})}(v, x_{0})\pi_{n}(\widetilde{\mu}^{(2)}_{G(w'_{(4)}),
w_{(3)}}))(w_{(1)}\otimes
w_{(2)})\nn
&&=z_2^{-1}\delta\left(\frac{x_1-x_0}{z_2}\right)\sum_{n\in \R}
\pi_{n}(Y'^{o}_{P(z_{0})}(v, x_{0})\widetilde{\mu}^{(2)}_{G(w'_{(4)}),
w_{(3)}})(w_{(1)}\otimes
w_{(2)})\nn
&&=\sum_{n\in \R}z_2^{-1}\delta\left(\frac{x_1-x_0}{z_2}\right)
\pi_{n}(Y'^{o}_{P(z_{0})}(v, x_{0})\widetilde{\mu}^{(2)}_{G(w'_{(4)}),
w_{(3)}})(w_{(1)}\otimes
w_{(2)}),
\end{eqnarray}
and we continue to have multiple absolute convergence.  Note that on
the right-hand side, for each $n \in \R$ the sum over the integral
powers of $x_0$ is a finite sum, in view of our comment after
(\ref{as:l}).

By the results from (\ref{as:p-2}) to (\ref{as:p-5}) we obtain
\begin{eqnarray}\label{as:p-6}
\lefteqn{\sum_{n\in \R}
z_2^{-1}\delta\left(\frac{x_1-x_0}{z_2}\right)
\pi_{n}(Y'^{o}_{P(z_{0})}(v, x_{0})
\widetilde{\mu}^{(2)}_{G(w'_{(4)}),
w_{(3)}})(w_{(1)}\otimes
w_{(2)})}\nn
&&=z_2^{-1}\delta\left(\frac{x_1-x_0}{z_2}\right)
\mu^{(2)}_{G(w'_{(4)}),
w_{(3)}}(w_{(1)}\otimes Y_2(v,x_0)w_{(2)})\nno\\
&&\quad+z_2^{-1}\delta\left(\frac{x_1-x_0}{z_2}\right)
\res_{x_2}z_0^{-1}\delta\bigg(\frac{x_0-x_2}{z_0}\bigg)\mu^{(2)}_{G(w'_{(4)}),
w_{(3)}}
(Y_1(v, x_2)w_{(1)}\otimes w_{(2)})\nno\\
\end{eqnarray}
(with absolute convergence).  Note that in this variant of
(\ref{as:p-1}), the left-hand side is the sum over $n \in \R$ of the
right-hand side of (\ref{as:l}).)

We now relate the right-hand side of (\ref{as:p-6}) to (\ref{as:o}).
By the definitions of $\mu^{(2)}_{G(Y'^o_4(v,x_1)w'_{(4)}),w_{(3)}}$
and $\mu^{(2)}_{G(w'_{(4)}),Y_3(v,x_1)w_{(3)}}$ (cf. (\ref{mu2G});
recall (\ref{yo}) and (\ref{y'})), we have 
\begin{eqnarray}\label{as:q}
\lefteqn{x_0^{-1}\delta\bigg(\frac{x_1-z_2}{x_0}\bigg)
\mu^{(2)}_{G(Y'^o_4(v,x_1)w'_{(4)}),w_{(3)}}(w_{(1)}
\otimes w_{(2)})}\nn
&&\qquad-x_0^{-1}\delta\bigg(\frac{z_2-x_1}{-x_0}\bigg)
\mu^{(2)}_{G(w'_{(4)}),Y_3(v,x_1)w_{(3)}}(w_{(1)}
\otimes w_{(2)})\nn
&&=x_0^{-1}\delta\bigg(\frac{x_1-z_2}{x_0}\bigg)\langle
Y'^o_4(v,x_1)w'_{(4)},
(I_1\circ (1_{W_1}\otimes I_2))(w_{(1)} \otimes w_{(2)} \otimes w_{(3)})
\rangle\nn
&&\qquad-x_0^{-1}\delta\bigg(\frac{z_2-x_1}{-x_0}\bigg)\langle
w'_{(4)},
(I_1\circ (1_{W_1}\otimes I_2))(w_{(1)} \otimes w_{(2)} \otimes 
Y_3(v,x_1)w_{(3)})
\rangle\nn
&&=x_0^{-1}\delta\bigg(\frac{x_1-z_2}{x_0}\bigg)\langle
w'_{(4)},Y_4(v,x_1)
(I_1\circ (1_{W_1}\otimes I_2))(w_{(1)} \otimes w_{(2)} \otimes w_{(3)})
\rangle\nn
&&\qquad-x_0^{-1}\delta\bigg(\frac{z_2-x_1}{-x_0}\bigg)\langle
w'_{(4)},
(I_1\circ (1_{W_1}\otimes I_2))(w_{(1)} \otimes w_{(2)} \otimes 
Y_3(v,x_1)w_{(3)})
\rangle.
\end{eqnarray}
Now using the formula obtained by taking $\res_{x_1}$ of (\ref{F12})
then replacing $x_0$ by $x_1$ and $x_2$ by $x_0$, we see that the right-hand side of 
(\ref{as:q}) is equal to
\begin{eqnarray}\label{as:r}
\lefteqn{z_2^{-1}\delta\bigg(\frac{x_1-x_0}{z_2}\bigg)\langle
w'_{(4)},
(I_1\circ (1_{W_1}\otimes I_2))(w_{(1)} \otimes Y_2(v,x_0)w_{(2)} \otimes 
w_{(3)})
\rangle}\nno\\
&&\quad+x_0^{-1}\delta\bigg(\frac{x_1-z_2}{x_0}\bigg)\res_{x_2}
z_1^{-1}\delta\bigg(\frac{x_1-x_2}{z_1}\bigg)\langle w'_{(4)},
(I_1\circ (1_{W_1}\otimes I_2))(Y_1(v,x_2)w_{(1)} \otimes w_{(2)} \otimes 
w_{(3)})
\rangle\nn
&&=z_2^{-1}\delta\bigg(\frac{x_1-x_0}{z_2}\bigg)
\mu^{(2)}_{G(w'_{(4)}),w_{(3)}}(w_{(1)}
\otimes Y_2(v,x_0)w_{(2)})\nno\\
&&\quad+x_0^{-1}\delta\bigg(\frac{x_1-z_2}{x_0}\bigg)\res_{x_2}
z_1^{-1}\delta\bigg(\frac{x_1-x_2}{z_1}\bigg)
\mu^{(2)}_{G(w'_{(4)}),w_{(3)}}(Y_1(v,x_2)w_{(1)}
\otimes w_{(2)}).
\end{eqnarray}
Using (\ref{l2-2}) and (\ref{l2-1}), we obtain
\begin{eqnarray*}
z^{-1}_2\delta\left(\frac{x_1-x_0}{z_2}\right)
z_0^{-1}\delta\bigg(\frac{x_0-x_2}{z_0}\bigg)
&=&x^{-1}_1\delta\left(\frac{z_1+x_2}{x_1}\right)
x_0^{-1}\delta\bigg(\frac{z_0+x_2}{x_0}\bigg)
\nn
&=&z_1^{-1}\delta\bigg(\frac{x_1-x_2}{z_1}\bigg)
x_0^{-1}\delta\bigg(\frac{x_1-z_2}{x_0}\bigg),
\end{eqnarray*}
so that the right-hand side of (\ref{as:r}) is equal to 
\begin{eqnarray}\label{as:s}
\lefteqn{z_2^{-1}\delta\bigg(\frac{x_1-x_0}{z_2}\bigg)
\mu^{(2)}_{G(w'_{(4)}),w_{(3)}}(w_{(1)}
\otimes Y_2(v,x_0)w_{(2)})}\nno\\
&&\quad+z^{-1}_2\delta\left(\frac{x_1-x_0}{z_2}\right)\res_{x_2}z_0^{-1}
\delta\bigg(\frac{x_0-x_2}{z_0}\bigg)
\mu^{(2)}_{G(w'_{(4)}),w_{(3)}}(Y_1(v,x_2)w_{(1)}
\otimes w_{(2)})\nno\\
\end{eqnarray}

{}From (\ref{as:q}), (\ref{as:r}) and (\ref{as:s}), we obtain that the
right-hand side of (\ref{as:p-6}) equals
\begin{eqnarray*}
\lefteqn{x_0^{-1}\delta\bigg(\frac{x_1-z_2}{x_0}\bigg)
\mu^{(2)}_{G(Y'^o_4(v,x_1)w'_{(4)}),w_{(3)}}(w_{(1)}
\otimes w_{(2)})}\nn
&&\quad-x_0^{-1}\delta\bigg(\frac{z_2-x_1}{-x_0}\bigg)
\mu^{(2)}_{G(w'_{(4)}),Y_3(v,x_1)w_{(3)}}(w_{(1)}
\otimes w_{(2)}),
\end{eqnarray*}
so that by (\ref{as:p-6}),
\begin{eqnarray}\label{as:t}
\lefteqn{\sum_{n\in \R}
z_2^{-1}\delta\left(\frac{x_1-x_0}{z_2}\right)
\pi_{n}(Y'^{o}_{P(z_{0})}(v, x_{0})
\widetilde{\mu}^{(2)}_{G(w'_{(4)}),
w_{(3)}})(w_{(1)}\otimes
w_{(2)})}\nn
&&=x_0^{-1}\delta\bigg(\frac{x_1-z_2}{x_0}\bigg)
\mu^{(2)}_{G(Y'^o_4(v,x_1)w'_{(4)}),w_{(3)}}(w_{(1)}
\otimes w_{(2)})\nn
&& \quad-x_0^{-1}\delta\bigg(\frac{z_2-x_1}{-x_0}\bigg)
\mu^{(2)}_{G(w'_{(4)}),Y_3(v,x_1)w_{(3)}}(w_{(1)}
\otimes w_{(2)})
\end{eqnarray}
for all $w_{(1)}\in W_{1}$ and $w_{(2)}\in W_{2}$. The right-hand side of (\ref{as:t})
is equal to 
\begin{eqnarray}\label{as:t-1}
\lefteqn{x_0^{-1}\delta\bigg(\frac{x_1-z_2}{x_0}\bigg)
\sum_{n\in \R}\lambda_{n}^{(2)}(Y'^o_4(v,x_1)w'_{(4)},w_{(3)})(w_{(1)}
\otimes w_{(2)})}\nn
&&\quad-x_0^{-1}\delta\bigg(\frac{z_2-x_1}{-x_0}\bigg)
\sum_{n\in \R}\lambda_{n}^{(2)}(w'_{(4)},Y_3(v,x_1)w_{(3)})(w_{(1)}
\otimes w_{(2)})
\end{eqnarray}
(recall (\ref{musumpi})).  Since the only infinite sums in
(\ref{as:t-1}) are those over $n \in \R$, we can move $\sum_{n \in
\R}$ to the left to obtain from (\ref{as:t})
\begin{eqnarray}\label{as:u}
\lefteqn{\sum_{n\in \R}z^{-1}_2\delta\left(\frac{x_1-x_0}{z_2}\right)
\pi_{n}(Y'^{o}_{P(z_{0})}(v, x_{0})
\widetilde{\mu}^{(2)}_{G(w'_{(4)}),
w_{(3)}})(w_{(1)}\otimes
w_{(2)})}\nn
&&=\sum_{n\in \R}x_0^{-1}\delta\bigg(\frac{x_1-z_2}{x_0}\bigg)
\lambda_{n}^{(2)}(Y'^o_4(v,x_1)w'_{(4)},w_{(3)})(w_{(1)}
\otimes w_{(2)})\nn
&& \quad-
\sum_{n\in \R}x_0^{-1}\delta\bigg(\frac{z_2-x_1}{-x_0}\bigg)
\lambda_{n}^{(2)}(w'_{(4)},Y_3(v,x_1)w_{(3)})(w_{(1)}
\otimes w_{(2)})
\end{eqnarray}
for $w_{(1)}\in W_{1}$ and $w_{(2)}\in W_{2}$.  That is, the sums over
$n \in \R$ of (\ref{as:l}) and (\ref{as:o}) are equal.

We now set up an application of Proposition \ref{real-exp-set}, by
first establishing from (\ref{as:u}) an equality of formal power
series in $y$ (see (\ref{as:v-1}) below) and then specializing $y$ to
$z'$ and proving and using certain convergence assertions.  For $k\in
\N$,
\begin{eqnarray*}
\lefteqn{
\sum_{n\in \R}x_0^{-1}\delta\bigg(\frac{x_1-z_2}{x_0}\bigg)\cdot}\nn
&&\quad\quad\quad\cdot
\sum_{i=0}^{k}{k\choose i}\lambda_{n}^{(2)}(Y'^o_4(v,x_1)w'_{(4)},w_{(3)})(
(L(0)+z_{0}L(-1))^{k-i}w_{(1)}
\otimes L(0)^{i}w_{(2)})\nn
&&\quad-
\sum_{n\in \R}x_0^{-1}\delta\bigg(\frac{z_2-x_1}{-x_0}\bigg)\cdot\nn
&&\quad\quad\quad\cdot
\sum_{i=0}^{k}{k\choose i}\lambda_{n}^{(2)}(w'_{(4)},Y_3(v,x_1)w_{(3)})
((L(0)+z_{0}L(-1))^{k-i}w_{(1)}
\otimes L(0)^{i}w_{(2)})\nn
&&=
\sum_{n\in \R}z^{-1}_2\delta\left(\frac{x_1-x_0}{z_2}\right)\cdot\nn
&&\quad\quad\quad\cdot
\sum_{i=0}^{k}{k\choose i}\pi_{n}(Y'^{o}_{P(z_{0})}(v, x_{0})
\widetilde{\mu}^{(2)}_{G(w'_{(4)}),
w_{(3)}})((L(0)+z_{0}L(-1))^{k-i}w_{(1)}
\otimes L(0)^{i}w_{(2)}),
\end{eqnarray*}
that is (by (\ref{LP'(j)})),
\begin{eqnarray*}
\lefteqn{
\sum_{n\in \R}x_0^{-1}\delta\bigg(\frac{x_1-z_2}{x_0}\bigg)
L_{P(z_{0})}'(0)^{k}
(\lambda_{n}^{(2)}(Y'^o_4(v,x_1)w'_{(4)},w_{(3)}))(w_{(1)}
\otimes w_{(2)})}\nn
&&\quad-
\sum_{n\in \R}x_0^{-1}\delta\bigg(\frac{z_2-x_1}{-x_0}\bigg)
L_{P(z_{0})}'(0)^{k}
(\lambda_{n}^{(2)}(w'_{(4)},Y_3(v,x_1)w_{(3)}))
(w_{(1)}
\otimes w_{(2)})\nn
&&=
\sum_{n\in \R}z^{-1}_2\delta\left(\frac{x_1-x_0}{z_2}\right)
L_{P(z_{0})}'(0)^{k}
(\pi_{n}(Y'^{o}_{P(z_{0})}(v, x_{0})
\widetilde{\mu}^{(2)}_{G(w'_{(4)}),
w_{(3)}}))(w_{(1)}
\otimes w_{(2)}),
\end{eqnarray*}
which gives
\begin{eqnarray}\label{as:v}
\lefteqn{
\sum_{n\in \R}x_0^{-1}\delta\bigg(\frac{x_1-z_2}{x_0}\bigg)
e^{yL_{P(z_{0})}'(0)}
(\lambda_{n}^{(2)}(Y'^o_4(v,x_1)w'_{(4)},w_{(3)}))(
(w_{(1)}
\otimes w_{(2)})}\nn
&&\quad-
\sum_{n\in \R}x_0^{-1}\delta\bigg(\frac{z_2-x_1}{-x_0}\bigg)
e^{yL_{P(z_{0})}'(0)}
(\lambda_{n}^{(2)}(w'_{(4)},Y_3(v,x_1)w_{(3)}))
(w_{(1)}
\otimes w_{(2)})\nn
&&=
\sum_{n\in \R}z^{-1}_2\delta\left(\frac{x_1-x_0}{z_2}\right)
e^{yL_{P(z_{0})}'(0)}
(\pi_{n}(Y'^{o}_{P(z_{0})}(v, x_{0})
\widetilde{\mu}^{(2)}_{G(w'_{(4)}),
w_{(3)}}))(w_{(1)}
\otimes w_{(2)}).\nn
\end{eqnarray}

Now both sides of (\ref{as:v}) are formal Laurent series in $x_{0}$
and $x_{1}$ with coefficients in $\C[[y]]$, and on the left-hand side,
the coefficient of each monomial in $x_0$ and $x_1$ involves only
finitely many pairs of vectors in $W'_4$ and $W_3$. Also, since
$W_{1}\hboxtr_{P(z_{0})}W_{2}$ is an object of $\mathcal{C}$, by
(\ref{mutilde2G}) and Assumption \ref{assum-exp-set} there exists
$K\in \N$ such that
\begin{equation}\label{as:v-0}
(L'_{P(z_{0})}(0)-n)^{K+1}(\pi_{n}(Y'^{o}_{P(z_{0})}(v, x_{0})
\widetilde{\mu}^{(2)}_{G(w'_{(4)}),
w_{(3)}}))=0
\end{equation}
for $n\in \R$.  Thus by Remark \ref{part-a} and (\ref{as:v-0}), for
each pair $p, q\in \Z$ there exists $N_{p, q}\in \N$ with
\begin{equation}\label{NgeK}
N_{p, q} \ge K
\end{equation}
such that 
\begin{eqnarray}\label{as:v-0-1}
\lefteqn{\res_{x_{0}}\res_{x_{1}}x_{0}^{p}x_{1}^{q}
\sum_{n\in \R}x_0^{-1}\delta\bigg(\frac{x_1-z_2}{x_0}\bigg)
e^{yL_{P(z_{0})}'(0)}
(\lambda_{n}^{(2)}(Y'^o_4(v,x_1)w'_{(4)},w_{(3)}))(
(w_{(1)}
\otimes w_{(2)})}\nn
&&\quad-\res_{x_{0}}\res_{x_{1}}x_{0}^{p}x_{1}^{q}
\sum_{n\in \R}x_0^{-1}\delta\bigg(\frac{z_2-x_1}{-x_0}\bigg)
e^{yL_{P(z_{0})}'(0)}
(\lambda_{n}^{(2)}(w'_{(4)},Y_3(v,x_1)w_{(3)}))
(w_{(1)}
\otimes w_{(2)})\nn
&&=\res_{x_{0}}\res_{x_{1}}x_{0}^{p}x_{1}^{q}
\sum_{n\in \R}x_0^{-1}\delta\bigg(\frac{x_1-z_2}{x_0}\bigg)
e^{ny}\cdot\nn
&&\quad\quad\quad\cdot \sum_{i=0}^{N_{p, q}}\frac{y^{i}}{i!}
((L'_{P(z_{0})}(0)-n)^{i}(\lambda_{n}^{(2)}(Y'^o_4(v,x_1)w'_{(4)},w_{(3)})))
(w_{(1)}
\otimes w_{(2)})\nn
&&\quad-\res_{x_{0}}\res_{x_{1}}x_{0}^{p}x_{1}^{q}
\sum_{n\in \R}x_0^{-1}\delta\bigg(\frac{z_2-x_1}{-x_0}\bigg)
e^{ny}\cdot\nn
&&\quad\quad\quad\cdot \sum_{i=0}^{N_{p, q}}\frac{y^{i}}{i!}
((L'_{P(z_{0})}(0)-n)^{i}(\lambda_{n}^{(2)}(w'_{(4)},Y_3(v,x_1)w_{(3)})))
(w_{(1)}
\otimes w_{(2)})\quad\quad\quad\quad\quad\quad\quad
\end{eqnarray}
and 
\begin{eqnarray}\label{as:v-0-2}
\lefteqn{\res_{x_{0}}\res_{x_{1}}x_{0}^{p}x_{1}^{q}
\sum_{n\in \R}z^{-1}_2\delta\left(\frac{x_1-x_0}{z_2}\right)
e^{yL_{P(z_{0})}'(0)}
(\pi_{n}(Y'^{o}_{P(z_{0})}(v, x_{0})
\widetilde{\mu}^{(2)}_{G(w'_{(4)}),
w_{(3)}}))(w_{(1)}
\otimes w_{(2)})}\nn
&&=\res_{x_{0}}\res_{x_{1}}x_{0}^{p}x_{1}^{q}
\sum_{n\in \R}z^{-1}_2\delta\left(\frac{x_1-x_0}{z_2}\right)
e^{ny}\cdot\nn
&&\quad\quad\quad\cdot \sum_{i=0}^{N_{p, q}}\frac{y^{i}}{i!}
((L'_{P(z_{0})}(0)-n)^{i}(\pi_{n}(Y'^{o}_{P(z_{0})}(v, x_{0})
\widetilde{\mu}^{(2)}_{G(w'_{(4)}),
w_{(3)}})))(w_{(1)}
\otimes w_{(2)}).\quad\quad\quad\quad\quad
\end{eqnarray}
In particular, we obtain from (\ref{as:v}), (\ref{as:v-0-1})
and (\ref{as:v-0-2})
\begin{eqnarray}\label{as:v-1}
\lefteqn{\res_{x_{0}}\res_{x_{1}}x_{0}^{p}x_{1}^{q}
\sum_{n\in \R}x_0^{-1}\delta\bigg(\frac{x_1-z_2}{x_0}\bigg)
e^{ny}\cdot}\nn
&&\quad\quad\quad\cdot \sum_{i=0}^{N_{p, q}}\frac{y^{i}}{i!}
((L'_{P(z_{0})}(0)-n)^{i}(\lambda_{n}^{(2)}(Y'^o_4(v,x_1)w'_{(4)},w_{(3)})))
(w_{(1)}
\otimes w_{(2)})\nn
&&\quad-\res_{x_{0}}\res_{x_{1}}x_{0}^{p}x_{1}^{q}
\sum_{n\in \R}x_0^{-1}\delta\bigg(\frac{z_2-x_1}{-x_0}\bigg)
e^{ny}\cdot\nn
&&\quad\quad\quad\cdot \sum_{i=0}^{N_{p, q}}\frac{y^{i}}{i!}
((L'_{P(z_{0})}(0)-n)^{i}(\lambda_{n}^{(2)}(w'_{(4)},Y_3(v,x_1)w_{(3)})))
(w_{(1)}
\otimes w_{(2)})\nn
&&=\res_{x_{0}}\res_{x_{1}}x_{0}^{p}x_{1}^{q}
\sum_{n\in \R}z^{-1}_2\delta\left(\frac{x_1-x_0}{z_2}\right)
e^{ny}\cdot\nn
&&\quad\quad\quad\cdot \sum_{i=0}^{N_{p, q}}\frac{y^{i}}{i!}
((L'_{P(z_{0})}(0)-n)^{i}(\pi_{n}(Y'^{o}_{P(z_{0})}(v, x_{0})
\widetilde{\mu}^{(2)}_{G(w'_{(4)}),
w_{(3)}})))(w_{(1)}
\otimes w_{(2)}).
\end{eqnarray}

We shall substitute $z'$ for $y$ in the two sides of (\ref{as:v-1}),
thus obtaining a common power series in $z'$, and we shall show that
this common power series is the power series expansion of two analytic
functions, which must then be equal.  We start with the left-hand
side.

Using the part of the proof of Theorem \ref{9.7-1} from
(\ref{iter-sum}) to (\ref{k-th-der-at-0}) with $-l^{0}(z)$ replaced by
$z'$ (as we did earlier in Remark \ref{pf-unique-lambda-n}), we see
using Proposition \ref{log-coeff-conv<=>iterate-conv} that for each
$i=0, \dots, N_{p, q}$, the series
\begin{eqnarray}\label{as:v-2.3}
\lefteqn{
\sum_{n\in \R}
\res_{x_{0}}\res_{x_{1}}x_{0}^{p}x_{1}^{q}
x_0^{-1}\delta\bigg(\frac{x_1-z_2}{x_0}\bigg)
e^{nz'}\cdot}\nn
&&\quad\quad\quad\cdot 
((L'_{P(z_{0})}(0)-n)^{i}(\lambda_{n}^{(2)}(Y'^o_4(v,x_1)w'_{(4)},w_{(3)})))
(w_{(1)}
\otimes w_{(2)})\nn
&&\quad-
\sum_{n\in \R}
\res_{x_{0}}\res_{x_{1}}x_{0}^{p}x_{1}^{q}
x_0^{-1}\delta\bigg(\frac{z_2-x_1}{-x_0}\bigg)
e^{nz'}\cdot\nn
&&\quad\quad\quad\cdot 
((L'_{P(z_{0})}(0)-n)^{i}(\lambda_{n}^{(2)}(w'_{(4)},Y_3(v,x_1)w_{(3)})))
(w_{(1)}
\otimes w_{(2)})
\end{eqnarray}
is absolutely convergent in an open neighborhood of $z'=0$, and that
by Lemma \ref{po-ser-an}, (\ref{as:v-2.3}) is in fact absolutely
convergent to an analytic function of $z'$ in this neighborhood.

The sum of (\ref{as:v-2.3}) as an analytic function of $z'$ has an
expansion as a power series in $z'$ in a small disk centered at $z'=0$
and the coefficients of the expansion are determined by its
derivatives at $z'=0$. By Lemma \ref{po-ser-an}, for each $k\in \N$
the $k$-th derivative at $z'=0$ of the sum of (\ref{as:v-2.3}) is the
sum of the absolutely convergent series
\begin{eqnarray*}
\lefteqn{
\sum_{n\in \R}
\res_{x_{0}}\res_{x_{1}}x_{0}^{p}x_{1}^{q}
x_0^{-1}\delta\bigg(\frac{x_1-z_2}{x_0}\bigg)
n^{k}\cdot}\nn
&&\quad\quad\quad\cdot 
((L'_{P(z_{0})}(0)-n)^{i}(\lambda_{n}^{(2)}(Y'^o_4(v,x_1)w'_{(4)},w_{(3)})))
(w_{(1)}
\otimes w_{(2)})\nn
&&\quad-
\sum_{n\in \R}
\res_{x_{0}}\res_{x_{1}}x_{0}^{p}x_{1}^{q}
x_0^{-1}\delta\bigg(\frac{z_2-x_1}{-x_0}\bigg)
n^{k}\cdot\nn
&&\quad\quad\quad\cdot 
((L'_{P(z_{0})}(0)-n)^{i}(\lambda_{n}^{(2)}(w'_{(4)},Y_3(v,x_1)w_{(3)})))
(w_{(1)}
\otimes w_{(2)})
\end{eqnarray*}
for each $i=0, \dots, N_{p, q}$. Thus we see that the expansion of the
sum of (\ref{as:v-2.3}) as a power series in $z'$ is 
\begin{eqnarray}\label{as:v-2.7}
\lefteqn{
\sum_{k\in \N}\Biggl(\sum_{n\in \R}
\res_{x_{0}}\res_{x_{1}}x_{0}^{p}x_{1}^{q}
x_0^{-1}\delta\bigg(\frac{x_1-z_2}{x_0}\bigg)
\frac{n^{k}}{k!}\cdot}\nn
&&\quad\quad\quad\cdot 
((L'_{P(z_{0})}(0)-n)^{i}(\lambda_{n}^{(2)}(Y'^o_4(v,x_1)w'_{(4)},w_{(3)})))
(w_{(1)}
\otimes w_{(2)})\Biggr)(z')^{k}\nn
&&\quad-
\sum_{k\in \N}\Biggl(\sum_{n\in \R}
\res_{x_{0}}\res_{x_{1}}x_{0}^{p}x_{1}^{q}
x_0^{-1}\delta\bigg(\frac{z_2-x_1}{-x_0}\bigg)
\frac{n^{k}}{k!}\cdot\nn
&&\quad\quad\quad\cdot 
((L'_{P(z_{0})}(0)-n)^{i}(\lambda_{n}^{(2)}(w'_{(4)},Y_3(v,x_1)w_{(3)})))
(w_{(1)}
\otimes w_{(2)})\Biggr)(z')^{k}
\end{eqnarray}
for each $i=0, \dots, N_{p, q}$. 
In particular, the power series obtained by substituting $z'$ for $y$ 
in the left-hand side of (\ref{as:v-1}) is absolutely convergent to
the sum of the doubly absolutely convergent series
\begin{eqnarray}\label{as:v-2.8}
\lefteqn{
\sum_{n\in \R}
\res_{x_{0}}\res_{x_{1}}x_{0}^{p}x_{1}^{q}
x_0^{-1}\delta\bigg(\frac{x_1-z_2}{x_0}\bigg)
e^{nz'}\cdot}\nn
&&\quad\quad\quad\cdot \sum_{i=1}^{N_{p, q}}\frac{(z')^{i}}{i!}
((L'_{P(z_{0})}(0)-n)^{i}(\lambda_{n}^{(2)}(Y'^o_4(v,x_1)w'_{(4)},w_{(3)})))
(w_{(1)}
\otimes w_{(2)})\nn
&&\quad-
\sum_{n\in \R}
\res_{x_{0}}\res_{x_{1}}x_{0}^{p}x_{1}^{q}
x_0^{-1}\delta\bigg(\frac{z_2-x_1}{-x_0}\bigg)
e^{nz'}\cdot\nn
&&\quad\quad\quad\cdot \sum_{i=1}^{N_{p, q}}\frac{(z')^{i}}{i!}
((L'_{P(z_{0})}(0)-n)^{i}(\lambda_{n}^{(2)}(w'_{(4)},Y_3(v,x_1)w_{(3)})))
(w_{(1)}
\otimes w_{(2)})
\end{eqnarray}
for $z'$ in 
the small disk above. 

We now consider the right-hand side of (\ref{as:v-1}) analogously.
Since $G(w'_{(4)})$ satisfies the $P^{(2)}(z_{0})$-local grading
restriction condition, for $z'$ in a neighborhood of $z'=0$, the
series
\[
\sum_{n\in \R}(e^{z'L'_{P(z_{0})}(0)}\pi_{n}(\widetilde{\mu}^{(2)}_{G(w'_{(4)}),
w_{(3)}}))(w_{(1)}
\otimes w_{(2)})=
\sum_{n\in \R}(e^{z'L'_{P(z_{0})}(0)}\lambda_{n}^{(2)}(w'_{(4)}, w_{(3)}))
(w_{(1)}\otimes w_{(2)})
\]
(recall (\ref{lambda-tilde-mu}) and (\ref{musumpi})) is absolutely
convergent for $w_{(1)}\in W_{1}$ and $w_{(2)}\in W_{2}$, and for $z'$
in this neighborhood, the sums of these series for $w_{(1)}\in W_{1}$
and $w_{(2)}\in W_{2}$ give an element $\mu^{(2)}_{G(w'_{(4)}),
w_{(3)}}(z')$ of $(W_{1}\otimes W_{2})^{*}$. Taking
\[
\mu=\mu^{(2)}_{G(w'_{(4)}),
w_{(3)}}(z')
\]
in (\ref{as:p}), we obtain 
\begin{eqnarray}\label{as:p-z'}
\lefteqn{(Y'^{o}_{P(z_{0})}(v, x_{0})\mu^{(2)}_{G(w'_{(4)}),
w_{(3)}}(z'))(w_{(1)}\otimes
w_{(2)})}\nn
&&=
(\mu^{(2)}_{G(w'_{(4)}),
w_{(3)}}(z'))(w_{(1)}\otimes Y_2(v,x_0)w_{(2)})\nno\\
&&\quad+\res_{x_2}z_0^{-1}
\delta\bigg(\frac{x_0-x_2}{z_0}\bigg)(\mu^{(2)}_{G(w'_{(4)}),
w_{(3)}}(z'))
(Y_1(v, x_2)w_{(1)}\otimes w_{(2)})\nn
&&=\sum_{n\in \R}(e^{z'L'_{P(z_{0})}(0)}\lambda_{n}^{(2)}(w'_{(4)}, w_{(3)}))
(w_{(1)}\otimes Y_2(v,x_0)w_{(2)})\nno\\
&&\quad+\res_{x_2}z_0^{-1}
\delta\bigg(\frac{x_0-x_2}{z_0}\bigg)
\sum_{n\in \R}(e^{z'L'_{P(z_{0})}(0)}\lambda_{n}^{(2)}(w'_{(4)}, w_{(3)}))
(Y_1(v, x_2)w_{(1)}\otimes w_{(2)})\nn
&&=\sum_{n\in \R}(e^{z'L'_{P(z_{0})}(0)}\lambda_{n}^{(2)}(w'_{(4)}, w_{(3)}))
(w_{(1)}\otimes Y_2(v,x_0)w_{(2)})\nno\\
&&\quad+\res_{x_2}
\sum_{n\in \R}z_0^{-1}
\delta\bigg(\frac{x_0-x_2}{z_0}\bigg)
(e^{z'L'_{P(z_{0})}(0)}\lambda_{n}^{(2)}(w'_{(4)}, w_{(3)}))
(Y_1(v, x_2)w_{(1)}\otimes w_{(2)}),\nn
\end{eqnarray}
where in the second term, the coefficient of each monomial in $x_{0}$
involves only the single infinite
sum $\sum_{n\in \R}$, so that $\sum_{n\in \R}$ can be switched with 
$z_0^{-1}
\delta\bigg(\frac{x_0-x_2}{z_0}\bigg)$.
Thus we obtain
\begin{eqnarray}\label{as:p-z'-1}
\lefteqn{\left(Y'^{o}_{P(z_{0})}(v, x_{0})\sum_{n\in \R}(e^{z'L'_{P(z_{0})}(0)}
\lambda_{n}^{(2)}(w'_{(4)}, w_{(3)}))
\right)(w_{(1)}\otimes
w_{(2)})}\nn
&&=\sum_{n\in \R}(Y'^{o}_{P(z_{0})}(v, x_{0})e^{z'L'_{P(z_{0})}(0)}
\lambda_{n}^{(2)}(w'_{(4)}, w_{(3)}))
(w_{(1)}\otimes
w_{(2)}),
\end{eqnarray}
with absolute convergence for the coefficient of each monomial in
$x_{0}$.

Let 
\[
\widetilde{\mu}^{(2)}_{G(w'_{(4)}),
w_{(3)}}(z')\in \overline{W_{1}\hboxtr_{P(z_{0})}W_{2}}
\]
be the element whose homogeneous components of generalized weight $n\in \R$
are
\[
e^{z'L'_{P(z_{0})}(0)}
\lambda_{n}^{(2)}(w'_{(4)}, w_{(3)}).
\]
In particular, we have
\begin{equation}\label{as:p-z'-1.1}
\widetilde{\mu}^{(2)}_{G(w'_{(4)}),
w_{(3)}}(z')=e^{z'L'_{P(z_{0})}(0)}\widetilde{\mu}^{(2)}_{G(w'_{(4)}),
w_{(3)}}.
\end{equation}

Since $|z_{2}|>|z_{0}|>0$, $|e^{-z'}z_{2}|>|z_{0}|>0$ for $|z'|$
sufficiently small.  Then the exact same arguments from (\ref{as:p-3})
to (\ref{as:p-5}) with $z_{2}$ replaced by $e^{-z'}z_{2}$, $v$ by
$e^{-z'L(0)}v$ and $\mu^{(2)}_{G(w'_{(4)}), w_{(3)}}$ by
\[
\mu=\mu^{(2)}_{G(w'_{(4)}),
w_{(3)}}(z')
\]
show that the  coefficient of each monomial in $x_{0}$ and $x_{1}$ in
\begin{eqnarray}\label{as:p-z'-2}
\lefteqn{\sum_{n\in \R}(e^{-z'}z_2)^{-1}\delta\left(\frac{x_1-x_0}{e^{-z'}z_2}\right)
(Y'^{o}_{P(z_{0})}(e^{-z'L(0)}v, x_{0})
\pi_{n}(\widetilde{\mu}^{(2)}_{G(w'_{(4)}),
w_{(3)}}(z')))
(w_{(1)}\otimes
w_{(2)})}\nn
&&=\sum_{n\in \R}(e^{-z'}z_2)^{-1}\delta\left(\frac{x_1-x_0}{e^{-z'}z_2}\right)
(Y'^{o}_{P(z_{0})}(e^{-z'L(0)}v, x_{0})e^{z'L'_{P(z_{0})}(0)}
\lambda_{n}^{(2)}(w'_{(4)}, w_{(3)}))
(w_{(1)}\otimes
w_{(2)})\nn
\end{eqnarray}
is doubly absolutely convergent for $z'$ in a neighborhood of $z'=0$
independent of $w_{(1)}$ and $w_{(2)}$ and that
\begin{eqnarray}\label{as:p-z'-3}
\lefteqn{\sum_{n\in \R}(e^{-z'}z_2)^{-1}\delta\left(\frac{x_1-x_0}{e^{-z'}z_2}\right)
\pi_{n}(Y'^{o}_{P(z_{0})}(e^{-z'L(0)}v, x_{0})\widetilde{\mu}^{(2)}_{G(w'_{(4)}),
w_{(3)}}(z'))
(w_{(1)}\otimes
w_{(2)})}\nn
&&=\sum_{n\in \R}(e^{-z'}z_2)^{-1}\delta\left(\frac{x_1-x_0}{e^{-z'}z_2}\right)
(Y'^{o}_{P(z_{0})}(e^{-z'L(0)}v, x_{0})\pi_{n}(\widetilde{\mu}^{(2)}_{G(w'_{(4)}),
w_{(3)}}(z')))
(w_{(1)}\otimes
w_{(2)})\nn
\end{eqnarray}
for $z'$ in this same neighborhood, again with double absolute
convergence.  Replacing $x_{1}$ and $x_{0}$ in (\ref{as:p-z'-2}) and
(\ref{as:p-z'-3}) by $e^{-z'}x_{1}$ and $e^{-z'}x_{0}$, respectively,
and dividing both sides by $e^{z'}$, we see that the coefficient of
each monomial in $x_{0}$ and $x_{1}$ in
\begin{equation}\label{as:p-z'-3.2}
\sum_{n\in \R}z_2^{-1}\delta\left(\frac{x_1-x_0}{z_2}\right)
(Y'^{o}_{P(z_{0})}(e^{-z'L(0)}v, e^{-z'}x_{0})e^{z'L'_{P(z_{0})}(0)}
\lambda_{n}^{(2)}(w'_{(4)}, w_{(3)}))
(w_{(1)}\otimes
w_{(2)})
\end{equation}
is doubly absolutely convergent in the same neighborhood of $z'=0$ and that
\begin{eqnarray}\label{as:p-z'-3.3}
\lefteqn{\sum_{n\in \R}z_2^{-1}\delta\left(\frac{x_1-x_0}{z_2}\right)
\pi_{n}(Y'^{o}_{P(z_{0})}(e^{-z'L(0)}v, e^{-z'}x_{0})\widetilde{\mu}^{(2)}_{G(w'_{(4)}),
w_{(3)}}(z'))
(w_{(1)}\otimes
w_{(2)})}\nn
&&=\sum_{n\in \R}z_2^{-1}\delta\left(\frac{x_1-x_0}{z_2}\right)
(Y'^{o}_{P(z_{0})}(e^{-z'L(0)}v, e^{-z'}x_{0})\pi_{n}(\widetilde{\mu}^{(2)}_{G(w'_{(4)}),
w_{(3)}}(z')))
(w_{(1)}\otimes
w_{(2)})\nn
\end{eqnarray}
for $z'$ in this same neighborhood, with double absolute convergence.

Using (\ref{yio-y-boxbs}), (\ref{yo-l-1}) and (\ref{sl2opp-2}), we
have, as in (\ref{710}) but for $Y^{o}$ (or by invoking (\ref{710}),
(\ref{y'}) and (\ref{L'(n)}))
\begin{eqnarray}\label{as:p-z'-3.5}
\lefteqn{
e^{z'L'_{P(z_{0})}(0)}Y'^{o}_{P(z_{0})}(v, x_{0})
\lambda_{n}^{(2)}(w'_{(4)}, w_{(3)})}\nn
&&=Y'^{o}_{P(z_{0})}(e^{-z'L(0)}v, e^{-z'}x_{0})e^{z'L'_{P(z_{0})}(0)}
\lambda_{n}^{(2)}(w'_{(4)}, w_{(3)})
\end{eqnarray}
for each $n\in \R$, and so
\begin{eqnarray}\label{as:p-z'-3.6}
\lefteqn{
e^{z'L'_{P(z_{0})}(0)}Y'^{o}_{P(z_{0})}(v, x_{0})\widetilde{\mu}^{(2)}_{G(w'_{(4)}),
w_{(3)}}}\nn
&&=Y'^{o}_{P(z_{0})}(e^{-z'L(0)}v, e^{-z'}x_{0})e^{z'L'_{P(z_{0})}(0)}
\widetilde{\mu}^{(2)}_{G(w'_{(4)}),
w_{(3)}}.
\end{eqnarray}
Also, by definition, $e^{z'L'_{P(z_{0})}(0)}$ commutes with $\pi_{n}$ for $n\in \R$. 
{}From these formulas, we obtain
\begin{eqnarray}\label{as:p-z'-4}
\lefteqn{\sum_{n\in \R}z_2^{-1}\delta\left(\frac{x_1-x_0}{z_2}\right)
(e^{z'L'_{P(z_{0})}(0)}\pi_{n}(Y'^{o}_{P(z_{0})}(v, x_{0})
\widetilde{\mu}^{(2)}_{G(w'_{(4)}),
w_{(3)}}))
(w_{(1)}\otimes
w_{(2)})}\nn
&&=\sum_{n\in \R}z_2^{-1}\delta\left(\frac{x_1-x_0}{z_2}\right)
\pi_{n}(e^{z'L'_{P(z_{0})}(0)}Y'^{o}_{P(z_{0})}(v, x_{0})
\widetilde{\mu}^{(2)}_{G(w'_{(4)}),
w_{(3)}})
(w_{(1)}\otimes
w_{(2)})\nn
&&=\sum_{n\in \R}z_2^{-1}\delta\left(\frac{x_1-x_0}{z_2}\right)
\pi_{n}(Y'^{o}_{P(z_{0})}(e^{-z'L(0)}v, e^{-z'}x_{0})e^{z'L'_{P(z_{0})}(0)}
\widetilde{\mu}^{(2)}_{G(w'_{(4)}),
w_{(3)}})
(w_{(1)}\otimes
w_{(2)})\nn
&&=\sum_{n\in \R}z_2^{-1}\delta\left(\frac{x_1-x_0}{z_2}\right)
\pi_{n}(Y'^{o}_{P(z_{0})}(e^{-z'L(0)}v, e^{-z'}x_{0})
\widetilde{\mu}^{(2)}_{G(w'_{(4)}),
w_{(3)}}(z'))
(w_{(1)}\otimes
w_{(2)}),
\end{eqnarray}
with double absolute convergence in the same neighborhood of $z'=0$
for each coefficient in $x_0$ and $x_1$ for each sum, since we know
that the right-hand side has double absolute convergence. (We recall
that the double sums are over $n\in \R$ and over the integral powers
of $x_{0}$. The operator $e^{z'L'_{P(z_{0})}(0)}$ is applied to the
indicated vectors, and in each case, it acts as a convergent sum of
operators on a suitable {\it finite-}dimensional vector space because
$W_{1}\hboxtr_{P(z_{0})} W_{2}$ is a generalized module.)

{}From (\ref{as:v-0}), (\ref{NgeK}) and (\ref{as:v-0-2}) and this
double absolute convergence for the left-hand side of
(\ref{as:p-z'-4}), we have, for each pair $p, q\in \Z$,  
\begin{eqnarray}\label{as:p-z'-5}
\lefteqn{\res_{x_{0}}\res_{x_{1}}x_{0}^{p}x_{1}^{q}
\sum_{n\in \R}z_2^{-1}\delta\bigg(\frac{x_1-x_0}{z_2}\bigg)
(e^{z'L'_{P(z_{0})}(0)}\pi_{n}(Y'^{o}_{P(z_{0})}(v, x_{0})\widetilde{\mu}^{(2)}_{G(w'_{(4)}),
w_{(3)}}))(w_{(1)}
\otimes w_{(2)})}\nn
&&=\res_{x_{0}}\res_{x_{1}}x_{0}^{p}x_{1}^{q}
\sum_{n\in \R}z^{-1}_2\delta\left(\frac{x_1-x_0}{z_2}\right)
e^{nz'}\cdot\nn
&&\quad\quad\quad\cdot \left(\sum_{i=0}^{N_{p, q}}\frac{(z')^{i}}{i!}
((L'_{P(z_{0})}(0)-n)^{i}(\pi_{n}(Y'^{o}_{P(z_{0})}(v,
x_{0})\widetilde{\mu}^{(2)}_{G(w'_{(4)}),
w_{(3)}})))(w_{(1)}
\otimes w_{(2)})\right);\nn
\end{eqnarray}
on the right-hand side, each summand of the doubly absolutely
convergent sum (in the neighborhood of $z'=0$ above) has been replaced
by a finite sum over $i$.

Since both sides of (\ref{as:p-z'-5}) are doubly absolutely
convergent, we can write (\ref{as:p-z'-5}) with the sums over $n \in
\R$ on the outside and the sums over the integral powers of $x_0$ on
the inside:
\begin{eqnarray}\label{as:p-z'-6}
\lefteqn{\sum_{n\in \R}\Biggl(\res_{x_{0}}\res_{x_{1}}x_{0}^{p}x_{1}^{q}
z_2^{-1}\delta\bigg(\frac{x_1-x_0}{z_2}\bigg)
\cdot}\nn
&&\quad\quad\quad\quad\quad\quad\cdot
(e^{z'L'_{P(z_{0})}(0)}\pi_{n}(Y'^{o}_{P(z_{0})}(v, x_{0})
\widetilde{\mu}^{(2)}_{G(w'_{(4)}),
w_{(3)}}))(w_{(1)}
\otimes w_{(2)})\Biggr)\nn
&&=\sum_{n\in \R}\Biggl(\res_{x_{0}}\res_{x_{1}}x_{0}^{p}x_{1}^{q}
z_2^{-1}\delta\bigg(\frac{x_1-x_0}{z_2}\bigg)e^{nz'}\cdot\nn
&&\quad\quad\quad\cdot
\Biggl(\sum_{i=0}^{N_{p, q}}\frac{(z')^{i}}{i!}
((L'_{P(z_{0})}(0)-n)^{i}(\pi_{n}(Y'^{o}_{P(z_{0})}(v,
x_{0})\widetilde{\mu}^{(2)}_{G(w'_{(4)}),
w_{(3)}})))(w_{(1)}
\otimes w_{(2)})\Biggr)\Biggr)\nn
&&=
\sum_{n\in \R}
e^{nz'}\Biggl(\sum_{i=0}^{N_{p, q}}\frac{(z')^{i}}{i!}
\Biggl((L'_{P(z_{0})}(0)-n)^{i}\cdot\nn
&&\quad\quad\quad\cdot \Biggl(\res_{x_{0}}\res_{x_{1}}x_{0}^{p}x_{1}^{q}
z^{-1}_2\delta\left(\frac{x_1-x_0}{z_2}\right)\pi_{n}(Y'^{o}_{P(z_{0})}(v,
x_{0})\widetilde{\mu}^{(2)}_{G(w'_{(4)}),
w_{(3)}})\Biggr)\Biggr)(w_{(1)}
\otimes w_{(2)})\Biggr),\nn
\end{eqnarray}
where the last equality follows from the finiteness of the sum over
integral powers of $x_0$ on the right-hand side of (\ref{as:p-5}) for
each $n \in \R$; in particular, on the right-hand side of
(\ref{as:p-z'-6}), for each $n \in \R$ the sum is finite, and thus the
same is true of the left-hand side of (\ref{as:p-z'-6}), the powers of
$x_0$ entering into the inner sum being the same on the two sides.
The outer sums (over $n \in \R$) are of course absolutely convergent
in our neighborhood of $z'=0$, which, we recall, is independent of
$w_{(1)}\in W_{1}$ and $w_{(2)}\in W_{2}$.

We again use the part of the proof of Theorem \ref{9.7-1} from
(\ref{iter-sum}) to (\ref{k-th-der-at-0}) with $-l^{0}(z)$ replaced by
$z'$ (as in Remark \ref{pf-unique-lambda-n} and (\ref{as:v-2.3}); what
follows is a variant of the argument in (\ref{as:v-2.3})--(\ref{as:v-2.8})):
Since the left-hand side of (\ref{as:p-z'-6}) is absolutely convergent
in a neighborhood of $z'=0$ independent of $w_{(1)}\in W_{1}$ and
$w_{(2)}\in W_{2}$, from (\ref{LP'(j)}), for $k\in \N$ the series
\begin{eqnarray}\label{as:p-z'-7}
\lefteqn{
\sum_{n\in \R}\Biggl(\res_{x_{0}}\res_{x_{1}}x_{0}^{p}x_{1}^{q}
z_2^{-1}\delta\bigg(\frac{x_1-x_0}{z_2}\bigg)
\Biggl((L'_{P(z_{0})}(0))^{k}e^{z'L'_{P(z_{0})}(0)}\cdot}\nn
&&\quad\quad\quad\quad\quad\quad\cdot
\pi_{n}(Y'^{o}_{P(z_{0})}(v, x_{0})\widetilde{\mu}^{(2)}_{G(w'_{(4)}),
w_{(3)}})\Biggr)(w_{(1)}
\otimes w_{(2)})\Biggr)\nn
&&=\sum_{n\in \R}\Biggl(\res_{x_{0}}\res_{x_{1}}x_{0}^{p}x_{1}^{q}
z_2^{-1}\delta\bigg(\frac{x_1-x_0}{z_2}\bigg)
\Biggl(e^{z'L'_{P(z_{0})}(0)}
\pi_{n}(Y'^{o}_{P(z_{0})}(v, x_{0})\widetilde{\mu}^{(2)}_{G(w'_{(4)}),
w_{(3)}})\Biggr)\cdot\nn
&&\quad\quad\quad \cdot\Biggl(\sum_{i=0}^{k}
{k\choose i}(L(0)+z_{0}L(-1))^{i}w_{(1)}
\otimes (L(0))^{k-i}w_{(2)}\Biggr)\Biggr)
\end{eqnarray}
obtained by summing the term-by-term $k$-th derivatives with respect
to $z'$ on the left-hand side of (\ref{as:p-z'-6}) is absolutely
convergent in the same neighborhood. By (\ref{as:p-z'-6}),
equivalently, the series of term-by-term $k$-th derivatives with
respect to $z'$ for the right-hand side of (\ref{as:p-z'-6}) is
absolutely convergent (as an iterated series) in the same
neighborhood, for $k\in \N$.  Thus by Proposition
\ref{log-coeff-conv<=>iterate-conv}, for each $i=0, \dots, N_{p. q}$,
\begin{eqnarray}\label{as:p-z'-8}
\lefteqn{\sum_{n\in \R} e^{nz'}\Biggl((L'_{P(z_{0})}(0)-n)^{i}\cdot}\nn
&&\quad\quad\quad\cdot \Biggl(\res_{x_{0}}\res_{x_{1}}x_{0}^{p}x_{1}^{q}
z^{-1}_2\delta\left(\frac{x_1-x_0}{z_2}\right)\pi_{n}(Y'^{o}_{P(z_{0})}(v,
x_{0})\widetilde{\mu}^{(2)}_{G(w'_{(4)}),
w_{(3)}})\Biggr)\Biggr)(w_{(1)}
\otimes w_{(2)})\nn
\end{eqnarray}
is absolutely convergent in the same neighborhood.  From Lemma
\ref{po-ser-an}, (\ref{as:p-z'-8}) is an analytic function in this
neighborhood, and thus so is the right-hand side of (\ref{as:p-z'-6}),
which equals the sum of the absolutely convergent double series
\begin{eqnarray}\label{as:p-z'-9}
\lefteqn{
\sum_{n\in \R}\sum_{i=0}^{N_{p, q}}
e^{nz'}\frac{(z')^{i}}{i!}
\Biggl((L'_{P(z_{0})}(0)-n)^{i}\cdot}\nn
&&\quad\quad\quad\cdot \Biggl(\res_{x_{0}}\res_{x_{1}}x_{0}^{p}x_{1}^{q}
z^{-1}_2\delta\left(\frac{x_1-x_0}{z_2}\right)\pi_{n}(Y'^{o}_{P(z_{0})}(v,
x_{0})\widetilde{\mu}^{(2)}_{G(w'_{(4)}),
w_{(3)}})\Biggr)\Biggr)(w_{(1)}
\otimes w_{(2)}).\nn
\end{eqnarray}
Moreover, the $k$-th derivative of (\ref{as:p-z'-9}) with respect to
with respect to $z'$ at $z'=0$ is the sum of the absolutely convergent
series obtained by setting $z'=0$ in the left-hand side of
(\ref{as:p-z'-7}), namely,
\begin{equation}\label{as:p-z'-10}
\sum_{n\in \R}\res_{x_{0}}\res_{x_{1}}x_{0}^{p}x_{1}^{q}
z_2^{-1}\delta\bigg(\frac{x_1-x_0}{z_2}\bigg)
((L'_{P(z_{0})}(0))^{k}
\pi_{n}(Y'^{o}_{P(z_{0})}(v, x_{0})\widetilde{\mu}^{(2)}_{G(w'_{(4)}),
w_{(3)}}))(w_{(1)}
\otimes w_{(2)}).\nn
\end{equation}

This information determines a power series expansion of the analytic
function (\ref{as:p-z'-9}) in a small disk centered at $z'=0$, and in
fact we know that it is obtained as follows: The right-hand side of
(\ref{as:v-1}), which is a formal power series in $y$ whose
coefficients are absolutely convergent sums over $n \in \R$, is
obtained by applying $\res_{x_{0}}\res_{x_{1}}x_{0}^{p}x_{1}^{q}$ to
the right-hand side of (\ref{as:v}), and the coefficient of $y^k/k!$
in this formal power series is exactly (\ref{as:p-z'-10}).  Hence in a
small disk centered at $z'=0$, the substitution of $z'$ for $y$ in the
right-hand side of (\ref{as:v-1}) gives a convergent power series
expansion of the analytic function (\ref{as:p-z'-9}), or equivalently,
(\ref{as:p-z'-6}).

We have shown that the left- and right-hand sides of (\ref{as:v-1})
with $y$ replaced by $z'$ are absolutely convergent to the sums of
(\ref{as:v-2.8}) and of (\ref{as:p-z'-9}), respectively, for $z'$ in
small disks centered at $z'=0$.  Thus the analytic functions
(\ref{as:v-2.8}) and (\ref{as:p-z'-9}) must be equal in the
intersection of these disks.  That is, for $p, q\in \Z$,
\begin{eqnarray}\label{as:v-3}
\lefteqn{\sum_{n\in \R}\sum_{i=0}^{N_{p, q}}e^{nz'}\frac{(z')^{i}}{i!}
\res_{x_{0}}\res_{x_{1}}x_{0}^{p}x_{1}^{q}
x_0^{-1}\delta\bigg(\frac{x_1-z_2}{x_0}\bigg)\cdot}\nn
&&\quad\quad\quad \cdot
((L'_{P(z_{0})}(0)-n)^{i}(\lambda_{n}^{(2)}(Y'^o_4(v,x_1)w'_{(4)},w_{(3)})))
(w_{(1)}
\otimes w_{(2)})\nn
&&\quad-\sum_{n\in \R}\sum_{i=0}^{N_{p, q}}e^{nz'}\frac{(z')^{i}}{i!}
\res_{x_{0}}\res_{x_{1}}x_{0}^{p}x_{1}^{q}
x_0^{-1}\delta\bigg(\frac{z_2-x_1}{-x_0}\bigg)
\cdot\nn
&&\quad\quad\quad\cdot 
((L'_{P(z_{0})}(0)-n)^{i}(\lambda_{n}^{(2)}(w'_{(4)},Y_3(v,x_1)w_{(3)})))
(w_{(1)}
\otimes w_{(2)})\nn
&&=\sum_{n\in \R}\sum_{i=0}^{N_{p, q}}e^{nz'}\frac{(z')^{i}}{i!}
\Biggl((L'_{P(z_{0})}(0)-n)^{i}\cdot\nn
&&\quad\quad\quad\cdot \Biggl(\res_{x_{0}}\res_{x_{1}}x_{0}^{p}x_{1}^{q}
z^{-1}_2\delta\left(\frac{x_1-x_0}{z_2}\right)\pi_{n}(Y'^{o}_{P(z_{0})}(v,
x_{0})\widetilde{\mu}^{(2)}_{G(w'_{(4)}),
w_{(3)}})\Biggr)\Biggr)(w_{(1)}
\otimes w_{(2)}),\nn
\end{eqnarray}
with double absolute convergence for $z'$ in a small disk centered at $z'=0$. 

We can now apply Proposition \ref{real-exp-set}.  Since $\R\times \{0,
\dots, N_{p, q}\}$ is a unique expansion set for any $p, q\in \Z$, we
obtain from (\ref{as:v-3}), taking $i=0$, that
\begin{eqnarray*}
\lefteqn{\res_{x_{0}}\res_{x_{1}}x_{0}^{p}x_{1}^{q}
x_0^{-1}\delta\bigg(\frac{x_1-z_2}{x_0}\bigg)
(\lambda_{n}^{(2)}(Y'^o_4(v,x_1)w'_{(4)}, w_{(3)}))(w_{(1)}\otimes w_{(2)})}\nn
&&\quad-\res_{x_{0}}\res_{x_{1}}x_{0}^{p}x_{1}^{q}
x_0^{-1}\delta\bigg(\frac{z_2-x_1}{-x_0}\bigg)
(\lambda_{n}^{(2)}(w'_{(4)}, Y_3(v,x_1)w_{(3)}))(w_{(1)}\otimes w_{(2)})\nn
&&=\res_{x_{0}}\res_{x_{1}}x_{0}^{p}x_{1}^{q}
z^{-1}_2\delta\left(\frac{x_1-x_0}{z_2}\right)
\pi_{n}(Y'^{o}_{P(z_{0})}(v, x_{0})
\widetilde{\mu}^{(2)}_{G(w'_{(4)}),
w_{(3)}})(w_{(1)}\otimes
w_{(2)})
\end{eqnarray*}
for $n\in \R$ and $p, q\in \Z$, or equivalently, that
\begin{eqnarray*}
\lefteqn{x_0^{-1}\delta\bigg(\frac{x_1-z_2}{x_0}\bigg)
(\lambda_{n}^{(2)}(Y'^o_4(v,x_1)w'_{(4)}, w_{(3)}))(w_{(1)}\otimes w_{(2)})}\nn
&&\quad-x_0^{-1}\delta\bigg(\frac{z_2-x_1}{-x_0}\bigg)
(\lambda_{n}^{(2)}(w'_{(4)}, Y_3(v,x_1)w_{(3)}))(w_{(1)}\otimes w_{(2)})\nn
&&=z^{-1}_2\delta\left(\frac{x_1-x_0}{z_2}\right)
\pi_{n}(Y'^{o}_{P(z_{0})}(v, x_{0})
\widetilde{\mu}^{(2)}_{G(w'_{(4)}),
w_{(3)}})(w_{(1)}\otimes
w_{(2)})
\end{eqnarray*}
for $n\in \R$.

We have thus proved the equality of (\ref{as:l}) and (\ref{as:o}), and
this proves (\ref{as:need0}) and hence (\ref{needtoshow}).

Similarly, one can prove that $\widetilde{G}$ intertwines the corresponding
$\mathfrak{s}\mathfrak{l}(2)$ actions; in fact, the appropriate
argument arises from setting $v=\omega$ above, and taking the relevant
three components at each step.
\epfv

Recalling the assumptions given before Lemma \ref{intertwine-tau}, we
see that the map $\widetilde{G}$ is $\widetilde{A}$-compatible (recall
Definition \ref{defJAtildecompat} and (\ref{W1W2beta})), from the
definitions, Remark \ref{I1I2'} and Proposition
\ref{lambda-n-a-tilde}.  Thus by Lemma \ref{intertwine-tau} and
Proposition \ref{qz} there is a unique $Q(z_2)$-intertwining map
$\widetilde{I}$ of type ${W_1\hboxtr_{P(z_0)} W_2\choose W'_4\,\,W_3}$
such that
\[
\widetilde{G}(w)(w'_{(4)}\otimes w_3)=\langle w, \widetilde{I}(w'_{(4)}\otimes 
w_{(3)})\rangle
\]
for $w\in W_1\boxtimes_{P(z_0)} W_2$, $w'_{(4)}\in W'_{4}$ and
$w_{(3)}\in W_{3}$.  By Corollary \ref{Q(z)P(z)iso}, there exists a
unique $P(z_{2})$-intertwining map $I$ of type ${W_4\choose
W_1\boxtimes_{P(z_0)} W_2\,\,W_3}$ such that
\[
\langle w, \widetilde{I}(w'_{(4)}\otimes 
w_{(3)})\rangle=\langle w'_{(4)}, I(w\otimes 
w_{(3)})\rangle,
\]
or equivalently,
\[
\widetilde{G}(w)(w'_{(4)}\otimes w_3)=\langle w'_{(4)}, I(w\otimes 
w_{(3)})\rangle,
\]
or equivalently (by (\ref{tildeG})),
\[
\langle w, \widetilde{\mu}^{(2)}_{G(w'_{(4)}),w_{(3)}}\rangle_{W_1\hboxtr_{P(z_0)}
W_2}= \langle w'_{(4)},I(w \otimes w_{(3)})\rangle
\]
for $w\in W_1\boxtimes_{P(z_0)} W_2$, $w'_{(4)}\in W'_{4}$ and
$w_{(3)}\in W_{3}$.  Taking
\[
w=\pi_{n}(w_{(1)}\boxtimes_{P(z_0)} w_{(2)})
\]
for $w_{(1)}\in W_{1}$, $w_{(2)}\in W_{2}$ and $n \in \R$ and invoking
Proposition \ref{span}, we see that $I$ is unique such that
\[
\langle \pi_{n}(w_{(1)}\boxtimes_{P(z_0)} w_{(2)}),
\widetilde{\mu}^{(2)}_{G(w'_{(4)}),w_{(3)}}\rangle_{W_1\hboxtr_{P(z_0)}
W_2}= \langle w'_{(4)},I(\pi_{n}(w_{(1)}\boxtimes_{P(z_0)} w_{(2)})
\otimes w_{(3)})\rangle,
\]
or equivalently (by (\ref{*-to-box-0})),
\begin{equation}\label{lambda=I}
(\lambda_{n}^{(2)}(w'_{(4)}, w_{(3)}))(w_{(1)}\otimes w_{(2)})
= \langle w'_{(4)},I(\pi_{n}(w_{(1)}\boxtimes_{P(z_0)} w_{(2)})
\otimes w_{(3)})\rangle
\end{equation}
for all $w_{(j)}\in W_{j}$ and $w'_{(4)}\in W'_{4}$.

Now we sum (\ref{lambda=I}) over $n \in \R$ to obtain the equality
\[
\widetilde{\mu}^{(2)}_{G(w'_{(4)}),w_{(3)}}(w_{(1)}\otimes w_{(2)})
= \langle w'_{(4)},I((w_{(1)}\boxtimes_{P(z_0)} w_{(2)})
\otimes w_{(3)})\rangle
\]
of absolutely convergent sums; for the left-hand side we recall
(\ref{lambdan2w'w}) and (\ref{mu2G}) and for the right-hand side we
invoke the convergence condition for intertwining maps in
$\mathcal{C}$ for the $P(z_0)$-intertwining map $\boxtimes_{P(z_0)}$.
That is, from Definition \ref{productanditerateexisting}, Remarks
\ref{Atildecompatcorrespondence} and \ref{I1I2'} and Definition
\ref{mudef},
\begin{eqnarray*}
\lefteqn{\langle w'_{(4)},I_1(w_{(1)}\otimes I_2(w_{(2)}\otimes
w_{(3)}))\rangle}\nn
&&=\langle w'_{(4)}, (I_1\circ (1_{W_1}\otimes I_2))
(w_{(1)}\otimes w_{(2)}\otimes w_{(3)})\rangle\nn
&&=(G(w'_{(4)}))(w_{(1)}\otimes w_{(2)}\otimes w_{(3)})
\nn
&&=\langle w'_{(4)},I((w_{(1)}\boxtimes_{P(z_0)} w_{(2)})\otimes 
w_{(3)})\rangle.
\end{eqnarray*}
Moreover, by Proposition \ref{intermediate}, this equality,
\begin{equation}\label{I1I2=I}
\langle w'_{(4)},I_1(w_{(1)}\otimes I_2(w_{(2)}\otimes w_{(3)}))\rangle
=\langle w'_{(4)},I((w_{(1)}\boxtimes_{P(z_0)} w_{(2)})\otimes 
w_{(3)})\rangle
\end{equation}
for all $w_{(j)}\in W_{j}$ and $w'_{(4)}\in W'_{4}$ determines the
$P(z_2)$-intertwining map $I$ uniquely. 

We have now proved Part 1 of the theorem below, which states in
particular that under the assumptions above, this product of
intertwining maps can be written as an iterate of certain intertwining
maps. Moreover, as is guaranteed by Proposition \ref{intermediate} and
proved directly above, the intermediate module can always be taken as
$W_{1}\boxtimes_{P(z_0)} W_{2}$.  Part 2 of this theorem is proved
analogously.

\begin{theo}\label{lgr=>asso}
Assume that $\mathcal{C}$ is closed under images, that the convergence
condition for intertwining maps in $\mathcal{C}$ holds and that
\[
|z_1|>|z_2|>|z_{0}|>0.
\]
Let $W_{1}$, $W_{2}$, $W_{3}$, $W_{4}$, $M_{1}$ and $M_{2}$ be objects
of $\mathcal{C}$.  Assume also that $W_1\boxtimes_{P(z_0)} W_2$ and
$W_2\boxtimes_{P(z_2)} W_3$ exist in $\mathcal{C}$.

\begin{enumerate}

\item  Let
$I_{1}$ and $I_{2}$ be $P(z_1)$- and $P(z_2)$-intertwining maps of 
types ${W_4}\choose {W_1M_1}$ and
${M_1}\choose {W_2W_3}$, respectively. 
Suppose that for each $w'_{(4)} \in W'_{4}$,
\[
\lambda=(I_1\circ (1_{W_1}\otimes I_2))'(w'_{(4)}) \in
(W_{1}\otimes W_{2} \otimes W_{3})^{*}
\]
satisfies the $P^{(2)}(z_0)$-local grading restriction condition (or
the $L(0)$-semisimple $P^{(2)}(z_0)$-local grading restriction
condition when $\mathcal{C}$ is in $\mathcal{M}_{sg}$). For
$w'_{(4)}\in W'_{4}$ and $w_{(3)}\in W_{3}$, let $\sum_{n\in
\R}\lambda_{n}^{(2)}$ be the (unique) series weakly absolutely
convergent to $\mu^{(2)}_{\lambda, w_{(3)}}$ as indicated in the
$P^{(2)}(z_0)$-grading condition (or the $L(0)$-semisimple
$P^{(2)}(z_0)$-grading condition).  Suppose also that for each $n \in
\R$, $w'_{(4)} \in W'_4$ and $w_{(3)} \in W_3$, the generalized
$V$-submodule of the generalized $V$-module $W^{(2)}_{\lambda,
w_{(3)}}$ (given by Theorem \ref{9.7-1}) generated by
$\lambda_{n}^{(2)}$ is a generalized $V$-submodule of some object of
$\mathcal{C}$ included in $(W_1 \otimes W_2)^*$.  Then the product
\[
I_1\circ (1_{W_1}\otimes I_2)
\]
can be expressed as an iterate, and in fact, there exists a unique
$P(z_2)$-intertwining map $I^{1}$ of type ${W_4\choose
W_1\boxtimes_{P(z_0)} W_2\,\,W_3}$ such that
\[
\langle w'_{(4)},I_1(w_{(1)} \otimes I_2(w_{(2)} \otimes w_{(3)}))\rangle
=\langle w'_{(4)}, I^{1}((w_{(1)}\boxtimes_{P(z_0)} w_{(2)})\otimes 
w_{(3)})\rangle
\]
for all $w_{(1)}\in W_1$, $w_{(2)}\in W_2$, $w_{(3)}\in W_3$ and
$w'_{(4)}\in W'_4$.

\item Analogously, let $I^1$ and $I^2$ be $P(z_2)$- and
$P(z_0)$-intertwining maps of types ${W_4}\choose {M_2W_3}$ and
${M_2}\choose {W_1W_2}$, respectively. Suppose that for each
$w'_{(4)} \in W'_{4}$,
\[
\lambda=(I^1\circ (I^2 \otimes 1_{W_3}))'(w'_{(4)}) \in
(W_{1}\otimes W_{2} \otimes W_{3})^{*}
\]
satisfies the $P^{(1)}(z_2)$-local grading restriction condition (or
the $L(0)$-semisimple $P^{(1)}(z_2)$-local grading restriction
condition when $\mathcal{C}$ is in $\mathcal{M}_{sg}$). For
$w'_{(4)}\in W'_{4}$ and $w_{(1)}\in W_{1}$, let $\sum_{n\in
\R}\lambda_{n}^{(1)}$ be the (unique) series weakly absolutely
convergent to $\mu^{(1)}_{\lambda, w_{(1)}}$ as indicated in the
$P^{(1)}(z_2)$-grading condition (or the $L(0)$-semisimple
$P^{(1)}(z_2)$-grading condition).  Suppose also that for each $n \in
\R$, $w'_{(4)} \in W'_4$ and $w_{(1)} \in W_1$, the generalized
$V$-submodule of the generalized $V$-module $W^{(1)}_{\lambda,
w_{(1)}}$ (given by Theorem \ref{9.7-1}) generated by
$\lambda_{n}^{(1)}$ is a generalized $V$-submodule of some object of
$\mathcal{C}$ included in $(W_2 \otimes W_3)^*$.  Then the iterate
\[
I^1\circ (I^2 \otimes 1_{W_3})
\]
can be expressed as a product, and in fact, there exists a unique
$P(z_1)$-intertwining map $I_{1}$ of type ${W_4\choose
W_1\,\,W_2\boxtimes_{P(z_2)} W_3}$ such that
\[
\langle w'_{(4)}, I^1(I^2(w_{(1)}\otimes w_{(2)})\otimes w_{(3)})\rangle
=\langle w'_{(4)}, I_{1}(w_{(1)}\otimes (w_{(2)}\boxtimes_{P(z_2)}w_{(3)}))
\rangle
\]
for all $w_{(1)}\in W_1$, $w_{(2)}\in W_2$, $w_{(3)}\in W_3$ and
$w'_{(4)}\in W'_4$. \epf

\end{enumerate}
\end{theo}

We know from Section 4, in particular, Proposition
\ref{im:correspond}, in which we shall take $p=0$, that
$P(z)$-intertwining maps are equivalent to suitable evaluations of
logarithmic intertwining operators or ordinary intertwining operators
at $z$.  Thus Theorem \ref{lgr=>asso} in fact says that under all of
the assumptions in the theorem, the following associativity of
logarithmic and of ordinary intertwining operators holds:

\begin{corol}\label{lgr=>asso-op}
Assume that $\mathcal{C}$ is closed under images, that the convergence
condition for intertwining maps in $\mathcal{C}$ holds and that
\[
|z_1|>|z_2|>|z_{0}|>0.
\]
Let $W_{1}$, $W_{2}$, $W_{3}$, $W_{4}$, $M_{1}$ and $M_{2}$ be objects
of $\mathcal{C}$.  Assume also that $W_1\boxtimes_{P(z_0)} W_2$ and
$W_2\boxtimes_{P(z_2)} W_3$ exist in $\mathcal{C}$.

\begin{enumerate}

\item Let $\Y_{1}$ and $\Y_{2}$ be logarithmic intertwining operators
(ordinary intertwining operators in the 
case that $\mathcal{C}$ is in $\mathcal{M}_{sg}$) 
of types ${W_4}\choose {W_1M_1}$ and ${M_1}\choose {W_2W_3}$,
respectively.  Suppose that for each $w'_{(4)} \in W'_{4}$, the
element $\lambda\in (W_{1}\otimes W_{2}\otimes W_{3})^{*}$ given by
\[
\lambda(w_{(1)}\otimes w_{(2)}\otimes w_{(3)})
=\langle w'_{(4)}, \Y_1(w_{(1)}, x_{1})\Y_2(w_{(2)}, x_{2})w_{(3)}\rangle
\lbar_{x_{1}=z_{1},\;x_{2}=z_{2}}
\]
(recalling (\ref{prodabbr})) for $w_{(1)}\in W_{1}$, $w_{(2)}\in
W_{2}$ and $w_{(3)}\in W_{3}$ satisfies the $P^{(2)}(z_0)$-local
grading restriction condition (or the $L(0)$-semisimple
$P^{(2)}(z_0)$-local grading restriction condition when $\mathcal{C}$
is in $\mathcal{M}_{sg}$). For $w'_{(4)}\in W'_{4}$ and $w_{(3)}\in
W_{3}$, let $\sum_{n\in \R}\lambda_{n}^{(2)}$ be the (unique) series
weakly absolutely convergent to $\mu^{(2)}_{\lambda, w_{(3)}}$ as
indicated in the $P^{(2)}(z_0)$-grading condition (or the
$L(0)$-semisimple $P^{(2)}(z_0)$-grading condition).  Suppose also
that for each $n \in \R$, $w'_{(4)} \in W'_4$ and $w_{(3)} \in W_3$,
the generalized $V$-submodule of the generalized $V$-module
$W^{(2)}_{\lambda, w_{(3)}}$ (given by Theorem \ref{9.7-1}) generated
by $\lambda_{n}^{(2)}$ is a generalized $V$-submodule of some object
of $\mathcal{C}$ included in $(W_1 \otimes W_2)^*$.  Then there exists
a unique logarithmic intertwining operator (a unique ordinary
intertwining operator in the case that $\mathcal{C}$ is in
$\mathcal{M}_{sg}$) $\Y^{1}$ of type ${W_4\choose
W_1\boxtimes_{P(z_0)} W_2\,\,W_3}$ such that
\begin{eqnarray}\label{prod=>iter}
\lefteqn{\langle w'_{(4)},\Y_1(w_{(1)}, x_{1}) \Y_2(w_{(2)}, x_{2}) w_{(3)}\rangle
\lbar_{x_{1}=z_{1},\;x_{2}=z_{2}}}\nn
&&=\langle w'_{(4)}, \Y^{1}(\Y_{\boxtimes_{P(z_{0})}, 0}
(w_{(1)}, x_{0})w_{(2)}, x_{2})w_{(3)})\rangle
\lbar_{x_{0}=z_{0},\;x_{2}=z_{2}}
\end{eqnarray}
(recalling (\ref{recover}) and (\ref{iterabbr})) for all $w_{(1)}\in
W_1$, $w_{(2)}\in W_2$, $w_{(3)}\in W_3$ and $w'_{(4)}\in W'_4$.  In
particular, the product of the logarithmic intertwining operators
(ordinary intertwining operators in the case that $\mathcal{C}$ is in
$\mathcal{M}_{sg}$) $\Y_{1}$ and $\Y_{2}$ evaluated at $z_{1}$ and
$z_{2}$, respectively, can be expressed as an iterate (with the
intermediate generalized $V$-module $W_1\boxtimes_{P(z_0)} W_2$) of
logarithmic intertwining operators (ordinary intertwining operators in
the case that $\mathcal{C}$ is in $\mathcal{M}_{sg}$) evaluated at
$z_{2}$ and $z_{0}$.

\item Analogously, let $\Y^1$ and $\Y^2$ be logarithmic intertwining
operators (ordinary intertwining operators in the case that
$\mathcal{C}$ is in $\mathcal{M}_{sg}$) of types ${W_4}\choose
{M_2W_3}$ and ${M_2}\choose {W_1W_2}$, respectively. Suppose that for
each $w'_{(4)}\in W'_{4}$, the element $\lambda\in (W_{1}\otimes
W_{2}\otimes W_{3})^{*}$ given by
\[
\lambda(w_{(1)}\otimes w_{(2)}\otimes w_{(3)})
=\langle w'_{(4)}, \Y^{1}(\Y^2(w_{(1)}, x_{0})w_{(2)}, x_{2})w_{(3)}\rangle
\lbar_{x_{0}=z_{0},\;x_{2}=z_{2}}
\]
(recalling (\ref{iterabbr})) satisfies the $P^{(1)}(z_2)$-local
grading restriction condition (or the $L(0)$-semisimple
$P^{(1)}(z_2)$-local grading restriction condition when $\mathcal{C}$
is in $\mathcal{M}_{sg}$). For $w'_{(4)}\in W'_{4}$ and $w_{(1)}\in
W_{1}$, let $\sum_{n\in \R}\lambda_{n}^{(1)}$ be the (unique) series
weakly absolutely convergent to $\mu^{(1)}_{\lambda, w_{(1)}}$ as
indicated in the $P^{(1)}(z_2)$-grading condition (or the
$L(0)$-semisimple $P^{(1)}(z_2)$-grading condition).  Suppose also
that for each $n \in \R$, $w'_{(4)} \in W'_4$ and $w_{(1)} \in W_1$,
the generalized $V$-submodule of the generalized $V$-module
$W^{(1)}_{\lambda, w_{(1)}}$ (given by Theorem \ref{9.7-1}) generated
by $\lambda_{n}^{(1)}$ is a generalized $V$-submodule of some object
of $\mathcal{C}$ included in $(W_2 \otimes W_3)^*$.  Then there exists
a unique logarithmic intertwining operator (a unique ordinary
intertwining operator in the case that $\mathcal{C}$ is in
$\mathcal{M}_{sg}$) $\Y_{1}$ of type ${W_4\choose
W_1\,\,W_2\boxtimes_{P(z_2)} W_3}$ such that
\begin{eqnarray}\label{iter=>prod}
\lefteqn{\langle w'_{(4)}, \Y^1(\Y^2(w_{(1)}, x_{0})w_{(2)}, x_{2})w_{(3)}\rangle
\lbar_{x_{0}=z_{0},\;x_{2}=z_{2}}}\nn
&&=\langle w'_{(4)}, \Y_{1}(w_{(1)}, x_{1})\Y_{\boxtimes_{P(z_2)}, 0}(w_{(2)}, x_{2})w_{(3)}
\rangle\lbar_{x_{1}=z_{1},\;x_{2}=z_{2}}
\end{eqnarray}
(again recalling (\ref{recover}) and (\ref{prodabbr})) for all
$w_{(1)}\in W_1$, $w_{(2)}\in W_2$, $w_{(3)}\in W_3$ and $w'_{(4)}\in
W'_4$.  In particular, the iterate of the logarithmic intertwining
operators (ordinary intertwining operators in the case that
$\mathcal{C}$ is in $\mathcal{M}_{sg}$) $\Y^{1}$ and $\Y^{2}$
evaluated at $z_{2}$ and $z_{0}$, respectively, can be expressed as a
product (with the intermediate generalized $V$-module
$W_2\boxtimes_{P(z_2)} W_3$) of logarithmic intertwining operators
(ordinary intertwining operators in the case that $\mathcal{C}$ is in
$\mathcal{M}_{sg}$) evaluated at $z_{1}$ and $z_{2}$.

\end{enumerate}
\end{corol}
\pf We prove only Part 1, the proof of Part 2 being analogous.

By (\ref{log:IYp}), we have
\[
\langle w'_{(4)}, \Y_1(w_{(1)}, x_{1})\Y_2(w_{(2)}, x_{2})w_{(3)}\rangle
\lbar_{x_{1}=z_{1},\;x_{2}=z_{2}}
=\langle w'_{(4)}, I_{\Y_1, 0}(w_{(1)} \otimes 
I_{\Y_2, 0}(w_{(2)}\otimes w_{(3)}))\rangle
\]
for $w_{(1)}\in W_{1}$, $w_{(2)}\in W_{2}$, $w_{(3)}\in W_{3}$ and
$w'_{(4)} \in W'_{4}$. By Part 1 of Theorem \ref{lgr=>asso}, There
exists a unique $P(z_{2})$-intertwining map $I^{1}$ of type
${W_{4}\choose W_{1}\boxtimes_{P(z_{0})}W_{2}\; W_{3}}$ such that
\[
\langle w'_{(4)}, I_{\Y_1, 0}(w_{(1)} \otimes 
I_{\Y_2, 0}(w_{(2)}\otimes w_{(3)}))\rangle=
\langle w'_{(4)}, I^{1}((w_{(1)} \boxtimes_{P(z_{0})} 
w_{(2)})\otimes w_{(3)})\rangle.
\]
By Proposition \ref{im:correspond}, we have
\[
\langle w'_{(4)}, I^{1}((w_{(1)} \boxtimes_{P(z_{0})} 
w_{(2)})\otimes w_{(3)})\rangle
=\langle w'_{(4)}, \Y_{I^{1}, 0}(\Y_{\boxtimes_{P(z_{0})}, 0}(w_{(1)}, x_{0})
w_{(2)}, x_{2}) w_{(3)})\rangle\lbar_{x_{0}=z_{0},\;x_{2}=z_{2}}.
\]
Taking $\Y^{1}=\Y_{I^{1}, 0}$, we obtain (\ref{prod=>iter}).  Since
$I^{1}$ is unique, $\Y^{1}=\Y_{I^{1}, 0}$ is also unique, by
Proposition \ref{im:correspond} (as in Corollary \ref{intermediate2}).
In the case that $\mathcal{C}$ is in $\mathcal{M}_{sg}$, $\Y^{1}$ is
an ordinary intertwining operator, by Remark \ref{log:ordi}.  \epfv

Theorem \ref{lgr=>asso} and Corollary \ref{lgr=>asso-op} both have two
parts, each with a major assumption, involving the
$P^{(2)}(z_{0})$-local grading restriction condition (or the
$L(0)$-semisimple $P^{(2)}(z_{0})$-local grading restriction
condition) in Part 1 and the $P^{(1)}(z_{2})$-local grading
restriction condition (or the $L(0)$-semisimple $P^{(1)}(z_{2})$-local
grading restriction condition) in Part 2, and the resulting pair of
assumptions essentially form most of what we will call the ``expansion
condition'' (see Definition \ref{expansion-conditions} below).  We
would now like to show that these two major assumptions are actually
equivalent to each other. For this, we need the equivalence (Theorem
\ref{asso-io} below) of two versions of the associativity of
logarithmic or ordinary intertwining operators, namely, that every
product can be expressed as an iterate, and on the other hand, that
every iterate can be expressed as a product (recall the conclusions of
Part 1 and Part 2 of Corollary \ref{lgr=>asso-op}).  Theorem
\ref{asso-io} and the lemma below used in its proof do not use any
results in Section 8 or any of the results in the present section that
we have obtained so far.

Recall that in Section \ref{convsec} we proved two formulas,
(\ref{nosub}) and (\ref{nosub2}), using the maps $\Omega_r$ (recall
(\ref{Omega_r})), on writing products of intertwining operators
satisfying certain conditions in terms of iterates, and vice versa.
In the next lemma, we shall use (\ref{nosub}) and (\ref{nosub2}) to
prove analogues of these two formulas.  In the statement and proof of
this lemma, we shall use the analyticity, Proposition \ref{analytic},
and Proposition \ref{formal=proj} and Remark \ref{4notations} to
rewrite the consequences (\ref{i2p}) and (\ref{p2i}) of (\ref{nosub})
and (\ref{nosub2}), respectively, and to write analogous expressions.

\begin{lemma}\label{123=321}
Assume that the convergence condition for intertwining maps in
$\mathcal{C}$ holds.  Let $W_{1}$, $W_{2}$, $W_{3}$, $W_{4}$, $M_{1}$
and $M_{2}$ be objects of $\mathcal{C}$.  Then:
\begin{enumerate}

\item 
For any nonzero complex numbers $z_1$, $z_2$ such that
\[
|z_1|>|z_0|>0, \;\; |z_2|>|z_0|>0
\]
(with $z_0=z_1 - z_2$ as usual),  there exist $p, q\in \Z$ such that 
for any logarithmic (in particular, ordinary) intertwining operators 
${\cal Y}^1$ and
${\cal Y}^2$ of types ${W_4}\choose {M_2 W_3}$ and ${M_2}\choose {W_1
W_2}$, respectively, we have
\begin{eqnarray}\label{(12)3-3(21)}
\lefteqn{\langle w'_{(4)}, {\cal Y}^1({\cal
Y}^2(w_{(1)},x_0)w_{(2)},x_2)w_{(3)}\rangle_{W_4}\lbar_{x_0=z_0,\;
x_2=z_2}}\nno\\
&&=\langle e^{z_1L'(1)}w'_{(4)}, \Omega_{-1}({\cal Y}^1) (w_{(3)},
y_1)\Omega_{-1}({\cal Y}^2)(w_{(2)}, y_2)
\cdot\nn
&&\quad\quad\quad\quad\quad\quad\quad\quad\cdot
w_{(1)}\rangle_{W_4}\lbar_{y_1^{n}=e^{nl_{p}(-z_1)},\; 
\log y_{1}=l_{p}(-z_1),\;y^{n}_2=e^{nl_{q}(-z_{0})},\;
\log y_{2}=l_{q}(-z_{0})}\nno\\
\end{eqnarray}
for all $w_{(1)}\in W_1$, $w_{(2)}\in W_2$, $w_{(3)}\in W_3$ and
$w'_{(4)}\in W'_4$.

\item  For any nonzero complex numbers $z_1$, $z_2$ such that
\[
|z_1|>|z_2|>0,\;\; |z_0|>|z_2|>0,
\]
there exist $\tilde{p}, \tilde{q}\in \Z$ such that for any logarithmic
(in particular, ordinary) intertwining operators ${\cal Y}_1$ and
${\cal Y}_2$ of types ${W_4}\choose {W_1 M_1}$ and ${M_1}\choose {W_2
W_3}$, respectively, we have
\begin{eqnarray}\label{1(23)-(32)1}
\lefteqn{\langle w'_{(4)}, {\cal Y}_1(w_{(1)},x_1){\cal
Y}_2(w_{(2)},x_2)w_{(3)}\rangle_{W_4}\lbar_{x_1=z_1,\; x_2=z_2}}\nno\\
&&=\langle e^{z_1L'(1)}w'_{(4)}, \Omega_{-1}({\cal
Y}_1)(\Omega_{-1}({\cal Y}_2)(w_{(3)}, y_0)w_{(2)},
y_2)\cdot\nn
&&\quad\quad\quad\quad\quad\quad\quad\quad\quad\cdot
w_{(1)} \rangle_{W_4}\lbar_{y^{n}_0= e^{nl_{\tilde{p}}(-z_2)},\; \log y_{0}
=l_{\tilde{p}}(-z_2),\; y^{n}_2=e^{nl_{\tilde{q}}(-z_{0})},\;
\log y_{2}=l_{\tilde{q}}(-z_{0})}\nno\\
\end{eqnarray}
for all $w_{(1)}\in W_1$, $w_{(2)}\in W_2$, $w_{(3)}\in W_3$ and
$w'_{(4)}\in W'_4$.
\end{enumerate}
\end{lemma}
\pf We prove only (\ref{(12)3-3(21)}); (\ref{1(23)-(32)1}) is proved
similarly, and at the end of the proof we discuss it briefly.  In the
first part of our proof, we shall interpret substitution notation such
as ``$x_2=e^{-i \pi}z_2$'' the same way we did in the proof of
Proposition \ref{convergence}, namely (in this instance), as the
substitution of
\[
e^{\log z_2 - \pi i}
\]
for $x_2$ (rather than as in (\ref{iterabbr}) and (\ref{prodabbr}),
where $p=0$); more precisely, for convenience we shall reverse the
occurrences of $\Omega_{0}$ and $\Omega_{-1}$ in
(\ref{nosub})--(\ref{i2p}), and correspondingly, we shall use
$x_2=e^{+i \pi}z_2$, which serves to replace $x_2$ by $e^{\log z_2 +
\pi i}$.

Using the formulas
(\ref{nosub})--(\ref{i2p}) (or more precisely, the indicated variant of
(\ref{i2p})), along with Proposition \ref{log:omega} and
(\ref{log:p1}), we have
\begin{eqnarray}\label{(12)3-3(21)-1}
\lefteqn{\langle w'_{(4)}, {\cal Y}^1({\cal
Y}^2(w_{(1)},x_0)w_{(2)},x_2)w_{(3)}\rangle_{W_4}\lbar_{x_0=z_0,\;
x_2=z_2}}\nno\\
&&=\langle e^{z_2L'(1)}w'_{(4)}, \Omega_{-1}({\cal Y}^1) (w_{(3)},
x_2)\Omega_0(\Omega_{-1}({\cal Y}^2))(w_{(1)},
x_0)w_{(2)}\rangle_{W_4}\lbar_{x_0=z_0,\; x_2=e^{\pi
i}z_2}\nno\\
&&=\langle e^{z_2L'(1)}w'_{(4)}, \Omega_{-1}({\cal Y}^1) (w_{(3)},
x_2)e^{x_0L(-1)}\Omega_{-1}({\cal Y}^2)(w_{(2)}, e^{\pi i}x_0)
w_{(1)}\rangle_{W_4}\lbar_{x_0=z_0,\; x_2=e^{\pi
i}z_2}\nno\\
&&=\langle e^{z_2L'(1)}w'_{(4)}, e^{x_0L(-1)}\Omega_{-1}({\cal Y}^1)
(w_{(3)}, x_2-x_0)\Omega_{-1}({\cal Y}^2)(w_{(2)}, e^{\pi i}x_0)
w_{(1)}\rangle_{W_4}\lbar_{x_0=z_0,\; x_2=e^{\pi i}z_2}\nno\\
&&=\langle e^{x_0L'(1)}e^{z_2L'(1)}w'_{(4)}, \Omega_{-1}({\cal Y}^1)
(w_{(3)}, x_2-x_0)\Omega_{-1}({\cal Y}^2)(w_{(2)}, e^{\pi i}x_0)
w_{(1)}\rangle_{W_4}\lbar_{x_0=z_0,\; x_2=e^{\pi i}z_2}\nno\\
&&=\langle e^{z_0L'(1)}e^{z_2L'(1)}w'_{(4)}, \Omega_{-1}({\cal Y}^1)
(w_{(3)}, x_2+x_0) \Omega_{-1}({\cal Y}^2)(w_{(2)}, x_0)
\cdot\nn
&&\quad\quad\quad\quad\quad\quad\quad\quad\quad\quad\quad\quad\cdot
w_{(1)}\rangle_{W_4}\lbar_{x_0^{n}=e^{nl_{q}(-z_{0})},\; 
\log x_{0}=l_{q}(-z_{0}),\;x^{n}_2=e^{nl_{\tilde{p}}(-z_{2})},\;
\log x_{2}=l_{\tilde{p}}(-z_{2})},\nno\\
\end{eqnarray}
for some $\tilde{p}, q\in \Z$ independent of $\Y^{1}$, $\Y^{2}$,
$w_{(1)}$, $w_{(2)}$, $w_{(3)}$ and $w'_{(4)}$ (see the discussion
before (\ref{i2p})).  Note that the left-hand side of
(\ref{(12)3-3(21)-1}), as a multisum obtained by substituting powers
of the formal variables $x_{0}$, $x_{2}$, $\log x_{0}$ and $\log
x_{2}$ by the indicated complex numbers, is absolutely convergent by
Proposition \ref{formal=proj}, and each step in (\ref{(12)3-3(21)-1})
means that the multisums on both sides are both absolutely convergent
and are equal. In particular, the right-hand side, as a multisum
obtained by substituting the powers of the formal variables $x_{0}$,
$x_{2}$, $\log x_{0}$ and $\log x_{2}$ by the indicated complex
numbers, is absolutely convergent.

Note that since $e^{\pm z_1 L'(1)}$ are linear automorphisms of
$W'_4$, $W'_4$ is spanned by homogeneous elements of the form $e^{z_1
L'(1)}w'_{(4)}$.  We need only prove (\ref{(12)3-3(21)}) for
homogeneous $w_{(1)}, w_{(2)}, w_{(3)}$ and $e^{z_1L'(1)}w'_{(4)}$,
and we assume this homogeneity.  Recalling Proposition
\ref{formal=proj} and (\ref{triple-sum}), let
\[
\Delta=-\wt e^{z_1L'(1)}w_{(4)}'+\wt w_{(1)}+\wt w_{(2)}+\wt w_{(3)} \in \R
\]
and define
\[
a_{n, j, i}=\langle e^{z_1L'(1)}w'_{(4)}, 
(w_{(3)})^{\Omega_{-1}(\mathcal{Y}^{1})}_{\Delta-n-2, j}
(w_{(2)})^{\Omega_{-1}(\mathcal{Y}^{2})}_{n, i}w_{(1)}\rangle\in \C
\] 
for $n\in \R$, $j=0, \dots, M$, $i=0, \dots, N$.
{}From this expression of $a_{n, j, i} $ and (\ref{log:ltc}),
for $\mu \in \R/\Z$, there exists $R_{\mu}\in \mu$ such that $a_{n, j, i}=0$
for any $n \in {\mu}$ with $n>R_{\mu}$.
Then since
\begin{eqnarray*}
\lefteqn{\langle e^{z_1L'(1)}w'_{(4)}, \Omega_{-1}({\cal Y}^1)
(w_{(3)}, x_2+x_0) \Omega_{-1}({\cal Y}^2)(w_{(2)}, x_0)w_{(1)}\rangle_{W_4}}\nn
&&=\sum_{n\in \R}\sum_{j=0}^{M}\sum_{i=0}^{N}a_{n, j, i} 
(x_{2}+x_{0})^{-\Delta+n+1}(\log (x_{2}+x_{0}))^{j} x_{0}^{-n-1}
(\log x_{0})^{i}\nn
&&=\sum_{m\in\R}\sum_{j=0}^{M}\sum_{i=0}^{N}\left(\sum_{k\in \N} 
a_{-m-1+k, j, i} {-\Delta-m+k\choose k} x_{2}^{-\Delta-m}
x_{0}^{m}\right)\cdot\nn
&&\quad\quad\quad\quad\quad\quad\quad\quad\quad\quad\cdot \left(\log x_{2}+
\sum_{l\in \Z_{+}}\frac{(-1)^{l-1}}{l}\frac{x_{0}^{l}}{x_{2}^{l}}\right)^{j}
(\log x_{0})^{i},
\end{eqnarray*}
the right-hand side of (\ref{(12)3-3(21)-1}) is equal to
\begin{eqnarray}\label{(12)3-3(21)-0}
\lefteqn{\sum_{m\in\R}\sum_{j=0}^{M}\sum_{i=0}^{N}\left(\sum_{k\in \N} 
a_{-m-1+k, j, i} {-\Delta-m+k\choose k} e^{(-\Delta-m)l_{\tilde{p}}(-z_{2})}
e^{ml_{q}(-z_{0})}\right)\cdot}\nn
&&\quad\quad\quad\quad\quad\quad\quad\quad\quad\quad\cdot \left(l_{\tilde{p}}(-z_{2})+
\sum_{l\in \Z_{+}}\frac{(-1)^{l-1}}{l}\frac{(-z_{0})^{l}}{(-z_{2})^{l}}\right)^{j}
l_{q}(-z_{0})^{i},
\end{eqnarray}
an absolutely convergent triple sum since $|z_{2}|>|z_{0}|>0$, with
the first of the inner sums finite (since $a_{-m-1+k, j, i}=0$ for
$k>m+1+R_{\overline{-m}}$ where $\overline{-m}$ is the congruence
class of $-m$) and the second of the inner sums absolutely convergent,
again since $|z_{2}|>|z_{0}|>0$.

Since 
\[
l_{\tilde{p}}(-z_{2})+
\sum_{l\in \Z_{+}}\frac{(-1)^{l-1}}{l}\frac{(-z_{0})^{l}}{(-z_{2})^{l}}
\]
is a value of the multivalued logarithmic function at the point $-z_{1}$, 
there exists  $p\in \Z$ such that 
\[
l_{p}(-z_{1})=l_{\tilde{p}}(-z_{2})+
\sum_{l\in \Z_{+}}\frac{(-1)^{l-1}}{l}\frac{(-z_{0})^{l}}{(-z_{2})^{l}}.
\]
Note that $p$ is independent of $\Y^{1}$, $\Y^{2}$,
$w_{(1)}$, $w_{(2)}$, $w_{(3)}$ and
$w'_{(4)}$.
Then since $|z_{1}|>|z_{0}|>0$,
Proposition \ref{formal=proj} and (\ref{triple-sum}) give that
the right-hand side of (\ref{(12)3-3(21)}) with $p$ and $q$ as above
is equal to the absolutely convergent triple sum
\begin{eqnarray}\label{(12)3-3(21)-0.2}
\lefteqn{
\sum_{n\in \R}\sum_{j=0}^{M}\sum_{i=0}^{N}a_{n, j, i} 
e^{(-\Delta+n+1)l_{p}(-z_{1})}
l_{p}(-z_{1})^{j}e^{(-n-1)l_{q}(-z_{0})}l_{q}(-z_{0})^{i}}\nn
&&=\sum_{n\in \R}\sum_{j=0}^{M}\sum_{i=0}^{N}
\left(\sum_{k\in \N} 
a_{n, j, i} {-\Delta+n+1\choose k} e^{(-\Delta+n-k+1)l_{\tilde{p}}(-z_{2})}
e^{(-n-1+k) l_{q}(-z_{0})}\right)\cdot\nn
&&\quad\quad\quad\quad\quad\quad\quad\quad\quad\quad\cdot \left(l_{\tilde{p}}(-z_{2})+
\sum_{l\in \Z_{+}}\frac{(-1)^{l-1}}{l}\frac{(-z_{0})^{l}}{(-z_{2})^{l}}\right)^{j}
l_{q}(-z_{0})^{i},
\end{eqnarray}
with the inner sums absolutely convergent binomial and logarithmic series
since $|z_{2}|>|z_{0}|>0$.

We now consider complex variables $z_{1}'$, $z_{2}'$ and
$z'_{0}=z_{1}'-z_{2}'$. We view $z'_{2}$ and $z_{0}'$ as independent
variables.  Let $U$ be any open subset of the region
$|z_{2}'+z_{0}'|>|z_{0}|$ and $|z_{2}'|>|z_{0}'|$ of $\C^{2}$ such
that its projection $U_{2}$ to the $z_{2}'$ coordinate is simply
connected and let $l(-z_{2}')$ be any single-valued analytic branch of
the logarithmic function of $-z_{2}'$ defined for $z_{2}'\in
U_{2}$. Then
\[
\tilde{l}(-z_{1}')=l(-z_{2}')+
\sum_{l\in \Z_{+}}\frac{(-1)^{l-1}}{l}\frac{(-z_{0}')^{l}}{(-z_{2}')^{l}}
\]
is a single-valued analytic branch of the logarithmic function of $-z_{1}'$
for  $z_{1}'\in z_{0}'+U_{2}$.
By Proposition \ref{formal=proj}, (\ref{triple-sum}) and 
Proposition \ref{analytic}, 
\begin{eqnarray}\label{(12)3-3(21)-0.25}
\lefteqn{\langle e^{z_1L'(1)}w'_{(4)}, \Omega_{-1}({\cal Y}^1) (w_{(3)},
y_1)\Omega_{-1}({\cal Y}^2)(w_{(2)}, y_2)
\cdot}\nn
&&\quad\quad\quad\quad\quad\quad\quad\quad\cdot
w_{(1)}\rangle_{W_4}\lbar_{y_1^{n}=e^{n\tilde{l}(-z_1')},\; 
\log y_{1}=\tilde{l}(-z_1'),\;y^{n}_2=e^{nl_{q}(-z_{0})},\;
\log y_{2}=l_{q}(-z_{0})}\nn
&&=
\sum_{n\in \R}\sum_{j=0}^{M}\sum_{i=0}^{N}a_{n, j, i} 
e^{(-\Delta+n+1)\tilde{l}(-z_{1}')}
\tilde{l}(-z_{1}')^{j}e^{(-n-1)l_{q}(-z_{0})}l_{q}(-z_{0})^{i}
\end{eqnarray}
and the corresponding series of its derivatives are absolutely
convergent as triple sums when $(z_{2}', z_{0}')\in U$. By Proposition
\ref{log-coeff-conv<=>iterate-conv} (see also Corollary
\ref{double-conv<=>iterate-conv} and its proof), for each $j=0, \dots,
M$ and $i=0, \dots, N$,
\begin{equation}\label{(12)3-3(21)-0.3}
\sum_{n\in \R}a_{n, j, i} 
e^{(-\Delta+n+1)\tilde{l}(-z_{1}')}
e^{(-n-1)l_{q}(-z_{0})}
\end{equation}
is absolutely convergent and by Lemma \ref{po-ser-an} is analytic in
$z_{1}'$ for $z_{1}'\in z_{0}'+U_{2}$.  Expanding
\[
e^{(-\Delta+n+1)\tilde{l}(-z_{1}')}=e^{(-\Delta+n+1)\tilde{l}(-z_{2}'-z_{0}')}
\]
for $n\in \R$ as a power series in $z_{0}'$ (as we did above with
$z_{1}$, $z_{2}$ and $z_{0}$ in (\ref{(12)3-3(21)-0.2})), and
recalling that for $\mu \in \R/\Z$, there exists $R_{\mu}\in \mu$ such
that $a_{n, j, i}=0$ for any $n \in {\mu}$ with $n>R_{\mu}$, we see
that (\ref{(12)3-3(21)-0.3}) is equal to the absolutely convergent
double sum
\begin{eqnarray}\label{(12)3-3(21)-0.4}
\lefteqn{\sum_{\mu\in \R/\Z}\sum_{\tilde{n}\in 
-\N}
\Biggl(\sum_{k\in \N} 
a_{\tilde{n}+R_{\mu}, j, i} {-\Delta+\tilde{n}+R_{\mu}+1\choose k} \cdot}\nn
&&\quad\quad\quad\quad\quad\quad\quad\quad\quad\quad\quad\quad\quad\cdot
e^{(-\Delta+\tilde{n}+R_{\mu}-k+1)l(-z_{2}')}(-z_{0}')^{k}
e^{(-\tilde{n}-R_{\mu}-1)l_{q}(-z_{0})}\Biggr),\nn
\end{eqnarray}
for $(z_{2}', z_{0}')\in U$, with the inner sum absolutely convergent
since $|z_{2}'|>|z_{0}'|$. In particular, as the quotient by
$e^{(-\Delta+R_{\mu}+1)l(-z_{2}')}e^{(-R_{\mu}-1)l_{q}(-z_{0})}$ of a
subsum of (\ref{(12)3-3(21)-0.4}), the series
\begin{equation}\label{(12)3-3(21)-0.5}
\sum_{\tilde{n}\in 
-\N}
\left(\sum_{k\in \N} 
a_{\tilde{n}+R_{\mu}, j, i} {-\Delta+\tilde{n}+R_{\mu}+1\choose k} 
(-z_{2}')^{\tilde{n}-k}(-z'_{0})^{k}
(-z_{0})^{-\tilde{n}}\right)
\end{equation}
for each $\mu\in \R/\Z$, $j=0, \dots, M$ and $i=0, \dots, N$ is
absolutely convergent for $(z_{2}', z_{0}')\in U$.  The sum of
(\ref{(12)3-3(21)-0.4}) as the composition of the analytic functions
(\ref{(12)3-3(21)-0.3}) and $-z_{1}'=-z_{2}'-z_{0}'$ is analytic in
each of $z_{2}'$ and $z_{0}'$ for $(z_{2}', z_{0}')\in U$, and, by
Lemma \ref{po-ser-an}, its derivatives with respect to $z_{2}'$ and
$z_{0}'$ are sums of absolutely convergent series obtained by taking
the derivatives term by term.  In particular, for each $\mu\in \R/\Z$,
$j=0, \dots, M$ and $i=0, \dots, N$, since
$e^{(-\Delta+R_{\mu}+1)l(-z_{2}')}e^{(-R_{\mu}-1)l_{q}(-z_{0})}$ is
analytic in $z_{2}'$ for $z_{2}'\in U_{2}$, the sum of
(\ref{(12)3-3(21)-0.5}) as the quotient by
$e^{(-\Delta+R_{\mu}+1)l(-z_{2}')}e^{(-R_{\mu}-1)l_{q}(-z_{0})}$ of a
subsum of (\ref{(12)3-3(21)-0.4}) is analytic in each of $z_{2}'$ and
$z_{0}'$ for $(z_{2}', z_{0}')\in U$ and its derivatives are sums of
absolutely convergent series obtained by taking the derivatives term
by term. Since $U$ is an arbitrary open subset of the region given by
$|z_{2}'+z_{0}'|>|z_{0}|$ and $|z_{2}'|>|z_{0}'|$, the sum of
(\ref{(12)3-3(21)-0.5}) is analytic with respect to each of $z_{2}'$
and $z_{0}'$ in the region $|z_{2}'+z_{0}'|>|z_{0}|$ and
$|z_{2}'|>|z_{0}'|$.

We now view (\ref{(12)3-3(21)-0.5}) as an analytic function of
$(z_{2}')^{-1}$ and $z_{0}'$.  The function (\ref{(12)3-3(21)-0.5}) is
equal to
\begin{eqnarray*}
\lefteqn{\sum_{\tilde{n}\in 
-\N}a_{\tilde{n}+R_{\mu}, j, i}(-z_{2}')^{\tilde{n}}
(1+(-z'_{0})(-z_{2}')^{-1})^{-\Delta+\tilde{n}+R_{\mu}+1}
(-z_{0})^{-\tilde{n}}}\nn
&&=(1+(-z'_{0})(-z_{2}')^{-1})^{-\Delta+R_{\mu}+1}\sum_{\tilde{n}\in 
-\N}a_{\tilde{n}+R_{\mu}, j, i}((-z_{2}')^{-1}
(1+(-z'_{0})(-z_{2}')^{-1})^{-1}(-z_{0}))^{-\tilde{n}}.
\end{eqnarray*}
Since the left-hand side is absolutely convergent when $|z_{2}'+z_{0}'|>|z_{0}|$
and $|z_{2}'|>|z_{0}'|$, the series 
\begin{equation}\label{(12)3-3(21)-0.5-2}
\sum_{\tilde{n}\in 
-\N}a_{\tilde{n}+R_{\mu}, j, i}((-z_{2}')^{-1}
(1+(-z'_{0})(-z_{2}')^{-1})^{-1}(-z_{0}))^{-\tilde{n}}
\end{equation}
is also absolutely convergent in the same region.  Consider the power
series $\sum_{\tilde{n}\in -\N}a_{\tilde{n}+R_{\mu}, j,
i}z^{-\tilde{n}}$.  {}From the discussion above, its radius of
convergence is not $0$.  In particular,
\[
\lim_{z\to 0}\sum_{\tilde{n}\in 
-\N}a_{\tilde{n}+R_{\mu}, j, i}z^{-\tilde{n}}
\]
exists and is equal to $a_{R_{\mu}, j, i}$.  Since the limit of
$(-z_{2}')^{-1} (1+(-z'_{0})(-z_{2}')^{-1})^{-1}(-z_{0})$ as
$(z_{2}')^{-1}$ approaches $0$ is $0$, the limit of
(\ref{(12)3-3(21)-0.5-2}) as $(z_{2}')^{-1}$ approaches $0$ is
$a_{R_{\mu}, j, i}$.  Thus for fixed $z_{0}'\in \C$, since the limit
of $(1+(-z'_{0})(-z_{2}')^{-1})^{-\Delta+R_{\mu}+1}$ as
$(z_{2}')^{-1}$ approaches $0$ is $1$, the limit of
(\ref{(12)3-3(21)-0.5}) as $(z_{2}')^{-1}$ approaches $0$ is
$a_{R_{\mu}, j, i}$. Hence for fixed $z_{0}'\in \C$, the singularity
$(z_{2}')^{-1}=0$ in (\ref{(12)3-3(21)-0.5}) is removable. We know
that (\ref{(12)3-3(21)-0.5}) is analytic in $z_{0}'$ for fixed
$(z_{2}')^{-1}\ne 0$, in our region.  Since the limit of the function
(\ref{(12)3-3(21)-0.5}) as $(z_{2}')^{-1}$ approaches $0$ is
$a_{R_{\mu}, j, i}$, this function is also analytic in $z_{0}'$ when
$(z_{2}')^{-1}=0$.  Hence by Hartogs' theorem (see, for example, page
8 of \cite{Sh}), this function is analytic as a function of the two
variables $(z_{2}')^{-1}$ and $z_{0}'$ in the neighborhood of $(0, 0)$
given by $|1+z_{0}'(z_{2}')^{-1}|>|z_{0}(z_{2}')^{-1}|$ and
$1>|z_{0}'(z_{2}')^{-1}|$.

Let $r$ be a real number satisfying $r>2|z_{0}|$. Then for
$(z_{2}')^{-1}$ and $z_{0}'$ satisfying $|(z_{2}')^{-1}|<r^{-1}$ and
$|z_{0}'|<r-|z_{0}|$, we have
\[
|z_{0}'(z_{2}')^{-1}|<(r-|z_{0}|)r^{-1}=1-|z_{0}|r^{-1}<1
\]
and 
\[|1+z_{0}'(z_{2}')^{-1}|\ge 1-|z_{0}'(z_{2}')^{-1}|>1-(r-|z_{0}|)r^{-1}
=|z_{0}|r^{-1}>|z_{0}(z_{2}')^{-1}|.
\]
Thus the polydisk given by $|(z_{2}')^{-1}|<r^{-1}$ and
$|z_{0}'|<r-|z_{0}|$ is in the region given by
$|1+z_{0}'(z_{2}')^{-1}|>|z_{0}(z_{2}')^{-1}|$ and
$1>|z_{0}'(z_{2}')^{-1}|$.  In particular, our function has a power
series expansion in $(z_{2}')^{-1}$ and $z_{0}'$ and the power series
is doubly absolutely convergent in the polydisk.  Since the
derivatives of this analytic function are obtained by taking the
derivatives of the series (\ref{(12)3-3(21)-0.5}) term by term, we see
that the power series expansion of this analytic function is the
double series
\begin{equation}\label{(12)3-3(21)-0.6}
\sum_{\tilde{n}\in 
-\N}
\sum_{k\in \N} 
a_{\tilde{n}+R_{\mu}, j, i} {-\Delta+\tilde{n}+R_{\mu}+1\choose k} 
(-z_{2}')^{\tilde{n}-k}(-z'_{0})^{k}
(-z_{0})^{-\tilde{n}}.
\end{equation}
Thus the two iterated series, (\ref{(12)3-3(21)-0.5})
and 
\begin{equation}\label{(12)3-3(21)-0.7}
\sum_{k\in \N} 
\left(\sum_{\tilde{n}\in 
-\N}
a_{\tilde{n}+R_{\mu}, j, i} {-\Delta+\tilde{n}+R_{\mu}+1\choose k} 
(-z_{2}')^{\tilde{n}-k}(-z'_{0})^{k}
(-z_{0})^{-\tilde{n}}\right),
\end{equation}
associated to (\ref{(12)3-3(21)-0.6}) are also absolutely convergent
in the polydisk and their sums are equal to the double sum of
(\ref{(12)3-3(21)-0.6}) in the polydisk.

Also, $a_{-\tilde{m}+R_{\mu}-1+k, j, i}=0$ when $-\tilde{m}-1+k>0$.
Thus in the polydisk, we obtain
\begin{eqnarray*}
\lefteqn{\sum_{\tilde{n}\in 
-\N}
\left(\sum_{k\in \N} 
a_{\tilde{n}+R_{\mu}, j, i} {-\Delta+\tilde{n}+R_{\mu}+1\choose k} 
(-z_{2}')^{\tilde{n}-k}(-z'_{0})^{k}
(-z_{0})^{-\tilde{n}}\right)}\nn
&&=\sum_{k\in \N} \sum_{\tilde{n}\in 
-\N}
a_{\tilde{n}+R_{\mu}, j, i} {-\Delta+\tilde{n}+R_{\mu}+1\choose k} 
(-z_{2}')^{\tilde{n}-k}(-z'_{0})^{k}
(-z_{0})^{-\tilde{n}}\nn
&&=\sum_{k\in \N} \sum_{\tilde{m}\in 
\N-1+k}
a_{-\tilde{m}+R_{\mu}-1+k, j, i} {-\Delta-\tilde{m}+R_{\mu}+k\choose k} 
(-z_{2}')^{-\tilde{m}-1}(-z'_{0})^{k}
(-z_{0})^{\tilde{m}+1-k}\nn
&&=\sum_{k\in \N} \sum_{\tilde{m}\in 
\N-1}
a_{-\tilde{m}+R_{\mu}-1+k, j, i} {-\Delta-\tilde{m}+R_{\mu}+k\choose k} 
(-z_{2}')^{-\tilde{m}-1}(-z'_{0})^{k}
(-z_{0})^{\tilde{m}+1-k}\nn
&&= \sum_{\tilde{m}\in 
\N-1}\sum_{k\in \N}
a_{-\tilde{m}+R_{\mu}-1+k, j, i} {-\Delta-\tilde{m}+R_{\mu}+k\choose k} 
(-z_{2}')^{-\tilde{m}-1}(-z'_{0})^{k}
(-z_{0})^{\tilde{m}+1-k}\nn
&&= \sum_{\tilde{m}\in 
\N-1}\left(\sum_{k\in \N}
a_{-\tilde{m}+R_{\mu}-1+k, j, i} {-\Delta-\tilde{m}+R_{\mu}+k\choose k} 
(-z_{2}')^{-\tilde{m}-1}(-z'_{0})^{k}
(-z_{0})^{\tilde{m}+1-k}\right)\nn
\end{eqnarray*}
for $\mu\in \R/\Z$,  $j=0, \dots, M$ and $i=0, \dots, N$. Then in the polydisk, we have 
\begin{eqnarray}\label{(12)3-3(21)-0.8}
\lefteqn{\sum_{\tilde{n}\in 
-\N}
\left(\sum_{k\in \N} 
a_{\tilde{n}+R_{\mu}, j, i} {-\Delta+\tilde{n}+R_{\mu}+1\choose k} 
(-z_{2}')^{\tilde{n}-k}(-z'_{0})^{k}
(-z_{0})^{-\tilde{n}}\right)\cdot }\nn
&&\quad\quad\quad\quad\quad\quad\cdot \left(l_{\tilde{p}}(-z_{2})+
\sum_{l\in \Z_{+}}\frac{(-1)^{l-1}}{l}\frac{(-z_{0})^{l}}{(-z_{2})^{l}}\right)^{j}
l_{q}(-z_{0})^{i}\nn
&&= \sum_{\tilde{m}\in 
\N-1}\left(\sum_{k\in \N}
a_{-\tilde{m}+R_{\mu}-1+k, j, i} {-\Delta-\tilde{m}+R_{\mu}+k\choose k} 
(-z_{2}')^{-\tilde{m}-1}(-z'_{0})^{k}
(-z_{0})^{\tilde{m}+1-k}\right)\cdot\nn
&&\quad\quad\quad\quad\quad\quad\cdot \left(l_{\tilde{p}}(-z_{2})+
\sum_{l\in \Z_{+}}\frac{(-1)^{l-1}}{l}\frac{(-z_{0})^{l}}{(-z_{2})^{l}}\right)^{j}
l_{q}(-z_{0})^{i}
\end{eqnarray}
for $\mu\in \R/\Z$, $j=0, \dots, M$ and $i=0, \dots, N$.  In
particular, since when $|(z_{2}')^{-1}|<r^{-1}$, $((z_{2}')^{-1},
z_{0})$ is in the polydisk, (\ref{(12)3-3(21)-0.8}) holds for such
$z_{2}'$ and $z_{0}'=z_{0}$.  We know that the left-hand side of
(\ref{(12)3-3(21)-0.8}) with $z_{0}'=z_{0}$ is analytic in
$(z_{2}')^{-1}$ for $|1+z_{0}(z_{2}')^{-1}|>|z_{0}(z_{2}')^{-1}|$ and
$1>|z_{0}(z_{2}')^{-1}|$. In particular, the value at
$(z_{2}')^{-1}=z_{2}^{-1}$ of the left-hand side of
(\ref{(12)3-3(21)-0.8}) with $z_{0}'=z_{0}$ is determined by analytic
extension from its values on the disk $|(z_{2}')^{-1}|<r^{-1}$.  We
also know, from (\ref{(12)3-3(21)-0}), that the right-hand side of
(\ref{(12)3-3(21)-0.8}) with $z_{0}'=z_{0}$ is absolutely convergent
when $z_{2}'=z_{2}$. Thus as a power series in $(z_{2}')^{-1}$, the
right-hand side of (\ref{(12)3-3(21)-0.8}) with $z_{0}'=z_{0}$ is
absolutely convergent when $|(z_{2}')^{-1}|\le |z_{2}^{-1}|$.  Hence
the sum of the right-hand side of (\ref{(12)3-3(21)-0.8}) with
$z_{0}'=z_{0}$ is analytic in $(z_{2}')^{-1}$ for $|(z_{2}')^{-1}|<
|z_{2}^{-1}|$ and is continuous on the closed disk $|(z_{2}')^{-1}|\le
|z_{2}^{-1}|$. In particular, the value at $(z_{2}')^{-1}=z_{2}^{-1}$
of the right-hand side of (\ref{(12)3-3(21)-0.8}) with $z_{0}'=z_{0}$
is determined by analytically extending its values on the disk
$|(z_{2}')^{-1}|<r^{-1}$ to the open disk $|(z_{2}')^{-1}|<
|z_{2}^{-1}|$ and then taking the limit $(z_{2}')^{-1}\to
z_{2}^{-1}$. Since (\ref{(12)3-3(21)-0.8}) holds in the disk
$|(z_{2}')^{-1}|<r^{-1}$ and the closed segment from $(z_{2}')^{-1}=0$
to $(z_{2}')^{-1}=z_{2}^{-1}$ lies in the domain of the function given
by the left-hand side of (\ref{(12)3-3(21)-0.8}) (with
$z_{0}'=z_{0}$), the two sides of (\ref{(12)3-3(21)-0.8}) must be
equal when $(z_{2}')^{-1}=z_{2}^{-1}$, that is,
\begin{eqnarray}\label{(12)3-3(21)-0.9}
\lefteqn{\sum_{\tilde{n}\in 
-\N}
\left(\sum_{k\in \N} 
a_{\tilde{n}+R_{\mu}, j, i} {-\Delta+\tilde{n}+R_{\mu}+1\choose k} 
(-z_{2})^{\tilde{n}-k}
(-z_{0})^{-\tilde{n}+k}\right)\cdot}\nn
&&\quad\quad\quad\quad\quad\quad\cdot \left(l_{\tilde{p}}(-z_{2})+
\sum_{l\in \Z_{+}}\frac{(-1)^{l-1}}{l}\frac{(-z_{0})^{l}}{(-z_{2})^{l}}\right)^{j}
l_{q}(-z_{0})^{i}\nn
&&= \sum_{\tilde{m}\in 
\N-1}\left(\sum_{k\in \N}
a_{-\tilde{m}+R_{\mu}-1+k, j, i} {-\Delta-\tilde{m}+R_{\mu}+k\choose k} 
(-z_{2})^{-\tilde{m}-1}
(-z_{0})^{\tilde{m}+1}\right)\cdot\nn
&&\quad\quad\quad\quad\quad\quad\cdot \left(l_{\tilde{p}}(-z_{2})+
\sum_{l\in \Z_{+}}\frac{(-1)^{l-1}}{l}\frac{(-z_{0})^{l}}{(-z_{2})^{l}}\right)^{j}
l_{q}(-z_{0})^{i}
\end{eqnarray}
for $\mu\in \R/\Z$,  $j=0, \dots, M$ and $i=0, \dots, N$.
Hence the right-hand side of (\ref{(12)3-3(21)-0.2}) is equal to 
\begin{eqnarray}\label{(12)3-3(21)-0.10}
\lefteqn{\sum_{\mu\in \R/\Z}\sum_{j=0}^{M}\sum_{i=0}^{N}\sum_{\tilde{n}\in 
-\N}
\left(\sum_{k\in \N} 
a_{\tilde{n}+R_{\mu}, j, i} {-\Delta+\tilde{n}+R_{\mu}+1\choose k} 
e^{(-\Delta+\tilde{n}+R_{\mu}-k+1)l_{\tilde{p}}(-z_{2})}
e^{(-\tilde{n}-R_{\mu}-1+k)l_{q}(-z_{0})}\right)\cdot}\nn
&&\quad\quad\quad\quad\quad\quad\quad\quad\quad\quad\cdot \left(l_{\tilde{p}}(-z_{2})+
\sum_{l\in \Z_{+}}\frac{(-1)^{l-1}}{l}\frac{(-z_{0})^{l}}{(-z_{2})^{l}}\right)^{j}
l_{q}(-z_{0})^{i}\nn
&&=\sum_{\mu\in \R/\Z}\sum_{j=0}^{M}\sum_{i=0}^{N}\Biggl(\sum_{\tilde{m}\in 
\N-1}\Biggl(\sum_{k\in \N}
a_{-\tilde{m}+R_{\mu}-1+k, j, i} {-\Delta-\tilde{m}+R_{\mu}+k\choose k} 
\cdot\nn
&&\quad\quad\quad\quad\quad\cdot 
e^{(-\Delta-\tilde{m}+R_{\mu})l_{\tilde{p}}(-z_{2})}
e^{(\tilde{m}-R_{\mu})l_{q}(-z_{0})}\Biggr)\Biggr)\left(l_{\tilde{p}}(-z_{2})+
\sum_{l\in \Z_{+}}\frac{(-1)^{l-1}}{l}\frac{(-z_{0})^{l}}{(-z_{2})^{l}}\right)^{j}
l_{q}(-z_{0})^{i}\nn
&&=\sum_{m\in\R}\sum_{j=0}^{M}\sum_{i=0}^{N}\left(\sum_{k\in \N}
a_{-m-1+k, j, i} {-\Delta-m+k\choose k} 
e^{(-\Delta-m)l_{\tilde{p}}(-z_{2})}
e^{ml_{q}(-z_{0})}\right)\cdot\nn
&&\quad\quad\quad\quad\quad\quad\quad\quad\quad\quad\cdot \left(l_{\tilde{p}}(-z_{2})+
\sum_{l\in \Z_{+}}\frac{(-1)^{l-1}}{l}\frac{(-z_{0})^{l}}{(-z_{2})^{l}}\right)^{j}
l_{q}(-z_{0})^{i}.\nn
\end{eqnarray}
Since the right-hand side of (\ref{(12)3-3(21)-0.10}) is exactly (\ref{(12)3-3(21)-0}),
which in turn is 
equal to 
the right-hand side of (\ref{(12)3-3(21)-1}), and 
the left-hand side of (\ref{(12)3-3(21)-0.10}) is equal to 
the right-hand side of (\ref{(12)3-3(21)}), (\ref{(12)3-3(21)}) holds.

For (\ref{1(23)-(32)1}), we have the following analogue of
(\ref{(12)3-3(21)-1}), using (\ref{nosub2})--(\ref{p2i}):
\begin{eqnarray}\label{1(23)-(32)1-1}
\lefteqn{\langle w'_{(4)}, {\cal Y}_1(w_{(1)},x_1){\cal
Y}_2(w_{(2)},x_2)w_{(3)}\rangle_{W_4}\lbar_{x_1=z_1,\; x_2=z_2}}\nno\\
&&=\langle e^{z_1L'(1)}w'_{(4)}, \Omega_{-1}({\cal Y}_1)
(\Omega_0(\Omega_{-1}({\cal Y}_2))(w_{(2)},
x_2)w_{(3)}, x_1)w_{(1)} \rangle_{W_4}\lbar_{x_1= e^{\pi i}z_1,\;
x_2=z_2}\nno\\
&&=\langle e^{z_1L'(1)}w'_{(4)}, \Omega_{-1}({\cal
Y}_1)(e^{x_2L(-1)}\Omega_{-1}({\cal Y}_2)(w_{(3)}, e^{\pi
i}x_2)w_{(2)}, x_1)w_{(1)} \rangle_{W_4}\lbar_{x_1= e^{\pi i}z_1,\;
x_2=z_2}\nno\\
&&=\langle e^{z_1L'(1)}w'_{(4)}, \Omega_{-1}({\cal
Y}_1)(\Omega_{-1}({\cal Y}_2)(w_{(3)}, e^{\pi
i}x_2)w_{(2)}, x_1+x_2)w_{(1)} \rangle_{W_4}\lbar_{x_1= e^{\pi i}z_1,\;
x_2=z_2}.\nno\\
\end{eqnarray}
This exhibits the format of (\ref{1(23)-(32)1}), and arguments similar
to those in the proof of (\ref{(12)3-3(21)}) above prove
(\ref{1(23)-(32)1}).  \epfv

In the following consequence of Lemma \ref{123=321}, we assert that
two conditions are equivalent; note that the appropriate hypotheses
about the generalized modules and the complex numbers are part of the
conditions:

\begin{theo}\label{asso-io}
Assume that the convergence condition for intertwining maps in ${\cal
C}$ holds. Then the following two conditions are equivalent:

\begin{enumerate}

\item For any objects $W_1$, $W_2$, $W_3$, $W_4$ and $M_1$ of ${\cal
C}$, any nonzero complex numbers $z_1$ and $z_2$ satisfying
$|z_1|>|z_2|>|z_0|>0$, and any logarithmic intertwining operators
(ordinary intertwining operators in the case that $\mathcal{C}$ is in
$\mathcal{M}_{sg}$) ${\cal Y}_1$ and ${\cal Y}_2$ of types
${W_4}\choose {W_1 M_1}$ and ${M_1}\choose {W_2 W_3}$, respectively,
there exist an object $M_2$ of ${\cal C}$ and logarithmic intertwining
operators (ordinary intertwining operators in the case that
$\mathcal{C}$ is in $\mathcal{M}_{sg}$) ${\cal Y}^1$ and ${\cal Y}^2$
of types ${W_4}\choose {M_2 W_3}$ and ${M_2}\choose {W_1 W_2}$,
respectively, such that
\begin{eqnarray}\label{extcnd1}
\lefteqn{\langle w'_{(4)}, {\cal Y}_1(w_{(1)},x_1){\cal Y}_2(w_{(2)},
x_2)w_{(3)}\rangle\lbar_{x_1=z_1,\,x_2=z_2}}\nno\\
&&=\langle w'_{(4)}, {\cal Y}^1({\cal Y}^2(w_{(1)},
x_0)w_{(2)},x_2)w_{(3)}\rangle\lbar_{x_0=z_0,\,x_2=z_2}
\end{eqnarray}
for all $w_{(1)}\in W_1$, $w_{(2)}\in W_2$, $w_{(3)}\in W_3$ and
$w'_{(4)}\in W'_4$.  (Here the substitution notation is as indicated
in (\ref{iterabbr}) and (\ref{prodabbr})).

\item For any objects $W_1$, $W_2$, $W_3$, $W_4$ and $M_2$ of ${\cal
C}$, any nonzero complex numbers $z_1$ and $z_2$ satisfying
$|z_1|>|z_2|>|z_0|>0$, and any logarithmic intertwining operators
(ordinary intertwining operators in the case that $\mathcal{C}$ is in
$\mathcal{M}_{sg}$) ${\cal Y}^1$ and ${\cal Y}^2$ of types
${W_4}\choose {M_2 W_3}$ and ${M_2}\choose {W_1 W_2}$, respectively,
there exist an object $M_1$ of ${\cal C}$ and logarithmic intertwining
operators (ordinary intertwining operators in the case that
$\mathcal{C}$ is in $\mathcal{M}_{sg}$) ${\cal Y}_1$ and ${\cal Y}_2$
of types ${W_4}\choose {W_1 M_1}$ and ${M_1}\choose {W_2 W_3}$,
respectively, such that
\begin{eqnarray}\label{extcnd2}
\lefteqn{\langle w'_{(4)}, {\cal Y}^1({\cal Y}^2(w_{(1)},
x_0)w_{(2)},x_2)w_{(3)}\rangle\lbar_{x_0=z_0,\,x_2=z_2}}\nno\\
&&=\langle w'_{(4)}, {\cal Y}_1(w_{(1)},x_1){\cal Y}_2(w_{(2)},
x_2)w_{(3)}\rangle\lbar_{x_1=z_1,\,x_2=z_2}
\end{eqnarray}
for all $w_{(1)}\in W_1$, $w_{(2)}\in W_2$, $w_{(3)}\in W_3$ and
$w'_{(4)}\in W'_4$.
\end{enumerate}
\end{theo}
\pf First we note that if there exist $M_2$, ${\cal Y}^1$ and ${\cal
Y}^2$, or $M_1$, ${\cal Y}_1$ and ${\cal Y}_2$ such that Condition 1
or Condition 2, respectively, holds for some particular $z_{1},
z_{2}\in \C$ satisfying $|z_1|>|z_2|>|z_0|>0$ and for all $w_{(1)}\in
W_1$, $w_{(2)}\in W_2$, $w_{(3)}\in W_3$ and $w'_{(4)}\in W'_4$, then
with the same $M_2$, ${\cal Y}^1$ and ${\cal Y}^2$, or $M_1$, ${\cal
Y}_1$ and ${\cal Y}_2$, Condition 1 or Condition 2, respectively,
holds for all $z_{1}, z_{2}\in \C$ satisfying $|z_1|>|z_2|>|z_0|>0$
and for all $w_{(1)}\in W_1$, $w_{(2)}\in W_2$, $w_{(3)}\in W_3$ and
$w'_{(4)}\in W'_4$.  This follows from the analyticity (Proposition
\ref{analytic}), the $L(-1)$-derivative property for logarithmic
intertwining operators, and the general fact that if two analytic
functions and their derivatives are equal at a particular point, then
they are equal on the intersection of their domains by Taylor's
theorem and analytic extension, assuming the intersection is
connected.  In fact, if Condition 1 holds for some particular $z_{1},
z_{2}\in \C$ satisfying $|z_1|>|z_2|>|z_0|>0$ and for all $w_{(1)}\in
W_1$, $w_{(2)}\in W_2$, $w_{(3)}\in W_3$ and $w'_{(4)}\in W'_4$, then
by the $L(-1)$-derivative property and the $L(-1)$-bracket relation,
\begin{eqnarray*}
\lefteqn{\Biggl(\frac{\partial^{k}}{\partial (z_{1}')^{k}}
\frac{\partial^{l}}{\partial (z_{2}')^{l}}
\langle w'_{(4)}, {\cal Y}_1(w_{(1)},x_1){\cal Y}_2(w_{(2)},
x_2)w_{(3)}\rangle\lbar_{x_1=z_1',\,x_2=z_2'}\Biggr)
\Bigg|_{z_{1}'=z_{1},\; z_{2}'=z_{2}}}\nno\\
&&=\Biggl(\frac{\partial^{k}}{\partial (z_{1}')^{k}}
\frac{\partial^{l}}{\partial (z_{2}')^{l}}\langle w'_{(4)}, {\cal Y}^1({\cal Y}^2(w_{(1)},
x_0)w_{(2)},x_2)w_{(3)}\rangle\lbar_{x_0=z_0',\,x_2=z_2'}
\Biggr)
\Bigg|_{z_{1}'=z_{1},\; z_{2}'=z_{2}}
\end{eqnarray*}
for all $w_{(1)}\in W_1$, $w_{(2)}\in W_2$, $w_{(3)}\in W_3$ and
$w'_{(4)}\in W'_4$, where $z_{1}'$, $z_{2}'$ are complex variables and
$z_{0}'=z_{1}'-z_{2}'$; if $z'=z_{1}'$, $z_{2}'$ or $z_{0}'$ is a
positive real number, then we compute the derivatives using the
branches with $\arg z' \ge 0$.  Then by Taylor's theorem, there is an
open subset of the region $|z_1'|>|z_2'|>|z_0'|>0$ whose closure
contains $(z_{1}, z_{2})$ such that on the closure of this open
subset,
\begin{eqnarray}\label{extcnd1-1}
\lefteqn{\langle w'_{(4)}, {\cal Y}_1(w_{(1)},x_1){\cal Y}_2(w_{(2)},
x_2)w_{(3)}\rangle\lbar_{x_1=z_1',\,x_2=z_2'}}\nno\\
&&=\langle w'_{(4)}, {\cal Y}^1({\cal Y}^2(w_{(1)},
x_0)w_{(2)},x_2)w_{(3)}\rangle\lbar_{x_0=z_0',\,x_2=z_2'}
\end{eqnarray}
for all $w_{(1)}\in W_1$, $w_{(2)}\in W_2$, $w_{(3)}\in W_3$ and
$w'_{(4)}\in W'_4$.  The two sides of (\ref{extcnd1-1}) are analytic
in $z_{1}'$ and $z_{2}'$ on the regions $|z_1'|>|z_2'|>0$, $\arg
z_{1}', \arg z_{2}'> 0$, and $|z_2'|>|z_0'|>0$, $\arg z_{2}', \arg
z_{0}'> 0$, respectively, and thus are equal on their intersection,
$|z_1'|>|z_2'|>|z_0'|>0$, $\arg z_{1}', \arg z_{2}', \arg z_{0}'> 0$.
Hence (\ref{extcnd1-1}) holds on the region $|z_1'|>|z_2'|>|z_0'|>0$.
The argument for Condition 2 is similar.

We shall prove only that Condition 1 implies Condition 2, the other
direction being similar.

Suppose that Condition 1 holds. Then for $z_1$ and $z_2$ as in the
statement of Condition 2, by the first part of Lemma \ref{123=321},
there exist $p, q\in \Z$ such that for any logarithmic intertwining
operators ${\cal Y}^1$ and ${\cal Y}^2$ as in the statement of
Condition 2,
\begin{eqnarray*}
\lefteqn{\langle w'_{(4)}, {\cal Y}^1({\cal
Y}^2(w_{(1)},x_0)w_{(2)},x_2)w_{(3)}\rangle_{W_4}\lbar_{x_0=z_0,\;
x_2=z_2}}\nno\\
&&=\langle e^{z_1L'(1)}w'_{(4)}, \Omega_{-1}({\cal Y}^1) (w_{(3)},
y_1) \Omega_{-1}({\cal Y}^2)(w_{(2)}, y_2)
\cdot\nn
&&\quad\quad\quad\quad\quad\quad\quad\quad\cdot
w_{(1)}\rangle_{W_4}\lbar_{y_1^{n}=e^{nl_{p}(-z_1)},\; 
\log y_{1}=l_{p}(-z_1),\;y^{n}_2=e^{nl_{q}(-z_{0})},\;
\log y_{2}=l_{q}(-z_{0})}
\end{eqnarray*}
for all $w_{(1)}\in W_1$, $w_{(2)}\in W_2$, $w_{(3)}\in W_3$ and
$w'_{(4)}\in W'_4$. The same argument as in the proof of
(\ref{extcnd1-1}) above gives
\begin{eqnarray}\label{step1}
\lefteqn{\langle w'_{(4)}, {\cal Y}^1({\cal
Y}^2(w_{(1)},x_0)w_{(2)},x_2)w_{(3)}\rangle_{W_4}\lbar_{x_0=z_0',\;
x_2=z_2'}}\nno\\
&&=\langle e^{z_1L'(1)}w'_{(4)}, \Omega_{-1}({\cal Y}^1) (w_{(3)},
y_1) \Omega_{-1}({\cal Y}^2)(w_{(2)}, y_2)
\cdot\nn
&&\quad\quad\quad\quad\quad\quad\quad\quad\cdot
w_{(1)}\rangle_{W_4}\lbar_{y_1^{n}=e^{nl_{p}(-z_1')},\; 
\log y_{1}=l_{p}(-z_1'),\;y^{n}_2=e^{nl_{q}(-z_{0}')},\;
\log y_{2}=l_{q}(-z_{0}')}
\end{eqnarray}
for $|z_{1}'|>|z_{0}'|>0$, $|z_{2}'|>|z_{0}'|>0$ and for
$w_{(1)}\in W_1$, $w_{(2)}\in W_2$, $w_{(3)}\in W_3$ and
$w'_{(4)}\in W'_4$.

By Remark \ref{formalinvariance}, there exist logarithmic intertwining
operators $\tilde{\cal Y}^1$ and $\tilde{\cal Y}^2$ of types
${W_4}\choose {W_3 M_2}$ and ${M_2}\choose {W_2 W_1}$, respectively,
such that
\begin{eqnarray}\label{step1.5}
\lefteqn{\langle e^{z_1L'(1)}w'_{(4)}, \Omega_{-1}({\cal Y}^1) (w_{(3)},
y_1) \Omega_{-1}({\cal Y}^2)(w_{(2)}, y_2)
\cdot}\nn
&&\quad\quad\quad\quad\quad\quad\quad\quad\cdot
w_{(1)}\rangle_{W_4}\lbar_{y_1^{n}=e^{nl_{p}(-z_1')},\; 
\log y_{1}=l_{p}(-z_1'),\;y^{n}_2=e^{nl_{q}(-z_{0}')},\;
\log y_{2}=l_{q}(-z_{0}')}\nn
&&=\langle e^{z_1L'(1)}w'_{(4)}, \tilde{\cal Y}^1 (w_{(3)},
y_1) \tilde{\cal Y}^2(w_{(2)}, y_2)
w_{(1)}\rangle_{W_4}\lbar_{y_1=-z_1',\; y_2=-z_0'}
\end{eqnarray}
for $|z_{1}'|>|z_{0}'|>0$ and for $w_{(1)}\in W_1$, $w_{(2)}\in W_2$,
$w_{(3)}\in W_3$ and $w'_{(4)}\in W'_4$.  Since the last expression is
of the same form as the left-hand side of (\ref{extcnd1-1}), we have
{}from Condition 1 and (\ref{extcnd1-1}) that there exist an object
$M_3$ of ${\cal C}$ and logarithmic intertwining operators
$\tilde{\cal Y}^3$ and $\tilde{\cal Y}^4$ of types ${W_4\choose M_3
W_1}$ and ${M_3\choose W_3 W_2}$, respectively, such that
\begin{eqnarray}\label{step2}
\lefteqn{\langle e^{z_1L'(1)}w'_{(4)}, \tilde{\cal Y}^1 (w_{(3)},
y_1) \tilde{\cal Y}^2(w_{(2)}, y_2)
w_{(1)}\rangle_{W_4}\lbar_{y_1=-z'_1,\; y_2=-z'_0}}\nno\\
&&=\langle e^{z_1L'(1)}w'_{(4)}, \tilde{\cal Y}^3(\tilde{\cal
Y}^4(w_{(3)},y_0)w_{(2)},y_2)w_{(1)}\rangle\lbar_{y_0=-z'_2,\;
y_2=-z'_0}
\end{eqnarray}
for $|z'_1|>|z'_0|>|z'_2|>0$ and for $w_{(1)}\in W_1$, $w_{(2)}\in
W_2$, $w_{(3)}\in W_3$ and $w'_{(4)}\in W'_4$.  (Note that the
inequality $|z'_0|>|z'_2|>0$ fails for $z'_0 = z_0$ and $z'_2 = z_2$.)
Again by Remark \ref{formalinvariance}, for any $\tilde{p},
\tilde{q}\in \Z$, there exist logarithmic intertwining operators
${\cal Y}^3$ and ${\cal Y}^4$ of types ${W_4}\choose {M_3 W_1}$ and
${M_3}\choose {W_3 W_2}$, respectively, such that
\begin{eqnarray}\label{step2.5}
\lefteqn{\langle e^{z_1L'(1)}w'_{(4)}, \tilde{\cal Y}^3(\tilde{\cal
Y}^4(w_{(3)},y_0)w_{(2)},y_2)w_{(1)}\rangle\lbar_{y_0=-z'_2,\;
y_2=-z'_0}}\nn
&&=\langle e^{z_1L'(1)}w'_{(4)}, {\cal Y}^3({\cal
Y}^4(w_{(3)},y_0)w_{(2)},y_2)\cdot\nn
&&\quad\quad\quad\quad\quad\quad\quad\quad\cdot
w_{(1)}\rangle\lbar_{y^{n}_0= e^{nl_{\tilde{p}}(-z'_2)},\; \log y_{0}
=l_{\tilde{p}}(-z'_2)\; y^{n}_2=e^{nl_{\tilde{q}}(-z'_{0})},
\;\log y_{2}=l_{\tilde{q}}(-z'_{0})}
\end{eqnarray}
for $|z'_0|>|z'_2|>0$ and for $w_{(1)}\in W_1$, $w_{(2)}\in W_2$,
$w_{(3)}\in W_3$ and $w'_{(4)}\in W'_4$.  Let $z_{2}^{0}$ be a fixed
complex number satisfying $|z_1|>|z_2^{0}|>0$ and
$|z_0^{0}|>|z_2^{0}|>0$, with $z_0^{0}=z_1 -z_2^{0}$.  By using the
fact that ${\cal Y}=\Omega_{-1}(\Omega_0({\cal Y}))$ and comparing the
last expression to the right-hand side of (\ref{1(23)-(32)1}), we see
{}from the second part of Lemma \ref{123=321} that there exist
$\tilde{p}, \tilde{q}\in \Z$ independent of ${\cal Y}^3$ and ${\cal
Y}^4$ such that
\begin{eqnarray*}
\lefteqn{\langle e^{z_1L'(1)}w'_{(4)}, {\cal Y}^3({\cal
Y}^4(w_{(3)},y_0)w_{(2)},y_2)\cdot}\nn
&&\quad\quad\quad\quad\quad\quad\quad\quad\cdot
w_{(1)}\rangle\lbar_{y^{n}_0= e^{nl_{\tilde{p}}(-z^{0}_2)},\; \log y_{0}
=l_{\tilde{p}}(-z^{0}_2)\; y^{n}_2=e^{nl_{\tilde{q}}(-z^{0}_{0})},
\;\log y_{2}=l_{\tilde{q}}(-z^{0}_{0})}\nno\\
&&=\langle w'_{(4)}, \Omega_0({\cal
Y}^3)(w_{(1)},x_1)\Omega_0({\cal
Y}^4)(w_{(2)},x_2)w_{(3)}\rangle_{W_4}\lbar_{x_1=z_1,\; x_2=z^{0}_2}
\end{eqnarray*}
for $w_{(1)}\in W_1$, $w_{(2)}\in W_2$, $w_{(3)}\in W_3$ and
$w'_{(4)}\in W'_4$. The same argument as in the proof of (\ref{extcnd1-1})
gives
\begin{eqnarray}\label{step3}
\lefteqn{\langle e^{z_1L'(1)}w'_{(4)}, {\cal Y}^3({\cal
Y}^4(w_{(3)},y_0)w_{(2)},y_2)\cdot}\nn
&&\quad\quad\quad\quad\quad\quad\quad\quad\cdot
w_{(1)}\rangle\lbar_{y^{n}_0= e^{nl_{\tilde{p}}(-z_2')},\; \log y_{0}
=l_{\tilde{p}}(-z_2')\; y^{n}_2=e^{nl_{\tilde{q}}(-z_{0}')},
\;\log y_{2}=l_{\tilde{q}}(-z_{0}')}\nno\\
&&=\langle w'_{(4)}, \Omega_0({\cal
Y}^3)(w_{(1)},x_1)\Omega_0({\cal
Y}^4)(w_{(2)},x_2)w_{(3)}\rangle_{W_4}\lbar_{x_1=z_1',\; x_2=z_2'}
\end{eqnarray}
for $|z'_{1}|>|z'_{2}|>0$, $|z'_0|>|z'_2|>0$ and for $w_{(1)}\in W_1$, 
$w_{(2)}\in W_2$, $w_{(3)}\in W_3$ and
$w'_{(4)}\in W'_4$.  The right-hand side is of course defined for
$|z'_{1}|>|z'_{2}|>0$.

By Proposition \ref{analytic}, we have that both sides of
(\ref{step1}), (\ref{step1.5}), (\ref{step2}), (\ref{step2.5}) and
(\ref{step3}) define analytic functions of $z_1'$ and $z_2'$ in the
indicated regions, with the cuts handled as in Propositon
\ref{analytic}.  Thus when restricted to the region
$|z_1'|>|z_2'|>|z_0'|>0$, the left-hand side of (\ref{step1}) and the
right-hand side of (\ref{step3}) are analytic extensions of each other
along loops, avoiding crossing the cuts, starting in the region
$|z_1'|>|z_2'|>|z_0'|>0$, passing through the region
$|z'_1|>|z'_0|>0$, the region $|z'_1|>|z'_0|>|z'_2|>0$, the region
$|z_0'|>|z_2'|>0$, the region $|z_1'|>|z_2'|>0$, and coming back to
the region $|z_1'|>|z_2'|>|z_0'|>0$ again.  We take $z_{1}, z_{2},
z_{2}^{0}\in \R_{+}$ satisfying $z_1>z_2>z_0>0$ and
$z_{1}>z_0^{0}=z_{1}-z_{2}^{0}>z_2^{0}>0$ and consider the path
$\gamma$ from $(z_{1}, z_{2})$ to $(z_{1}, z_{2}^{0})$ given by
$\gamma(t)=(z_{1}, (1-t)z_{2}+tz_{2}^{0})$ for $t\in [0, 1]$.  Then
the product $\gamma^{-1}\circ \gamma$ is a loop starting at the point
$(z_{1}, z_{2})$ in the region $|z_1'|>|z_2'|>|z_0'|>0$, passing
through the region $|z'_1|>|z'_0|>0$, the region
$|z'_1|>|z'_0|>|z'_2|>0$, the region $|z_0'|>|z_2'|>0$, the region
$|z_1'|>|z_2'|>0$, and coming back to the same point $(z_{1}, z_{2})$
in the region $|z_1'|>|z_2'|>|z_0'|>0$ again.  Thus the value of the
right-hand side of (\ref{step3}) at $(z_{1}, z_{2})$ is the analytic
extension of the left-hand side of (\ref{step1}) along the loop
$\gamma^{-1}\circ \gamma$, which is homotopic to the trivial loop, and
so the analytic extension must give the same value. Thus if we take
$\Y_{1}$ and $\Y_{2}$ to be $\Omega_0({\cal Y}^3)$ and $\Omega_0({\cal
Y}^4)$, respectively, (\ref{extcnd2}) holds at this particular point
$(z_{1}, z_{2})$. By the discussion in the beginning of this proof,
(\ref{extcnd2}) holds for all $z_{1}, z_{2}$ satisfying
$|z_1|>|z_2|>|z_0|>0$ and hence Condition 2 holds.

In the case that $\mathcal{C}$ is in $\mathcal{M}_{sg}$, the same
arguments still hold except that all the logarithmic intertwining
operators involved are ordinary intertwining operators.  \epfv

Using Theorem \ref{asso-io}, we now prove that under the global
assumptions in Theorem \ref{lgr=>asso}, the assumptions in Part 1 of
Theorem \ref{lgr=>asso} (or equivalently, of its reformulation,
Corollary \ref{lgr=>asso-op}) and the assumptions in Part 2 of Theorem
\ref{lgr=>asso} (or of Corollary \ref{lgr=>asso-op}) are equivalent.
These sets of assumptions, which are stated as Conditions 1 and 2 in
the theorem below, are the two (equivalent) statements of what we will
call the expansion condition below. {}From Proposition \ref{4.19}, the
statement that $\mathcal{C}$ is closed under the $P(z)$-tensor product
operation for {\it some} $z\in \C^{\times}$ is equivalent to the
statement that $\mathcal{C}$ is closed under the $P(z)$-tensor product
operation for {\it every} $z\in \C^{\times}$.  In particular, in the
following results, instead of assuming that
$W_{1}\boxtimes_{P(z_{0})}W_{2}$ and $W_{2}\boxtimes_{P(z_{2})} W_{3}$
exist in $\mathcal{C}$ for all objects $W_{1}$, $W_{2}$ and $W_{3}$ of
$\mathcal{C}$, we assume that $\mathcal{C}$ is closed under the
$P(z)$-tensor product operation for some $z\in \C^{\times}$.  Since
Conditions 1 and 2 below are about to be used as the two equivalent
formulations of the expansion condition, we include the hypotheses on
the generalized modules and the complex numbers in Conditions 1 and 2
(as we did in Theorem \ref{asso-io}).

\begin{theo}\label{expansion}
Assume that $\mathcal{C}$ is closed under images and under the
$P(z)$-tensor product operation for some $z\in \C^{\times}$, and that
the convergence condition for intertwining maps in ${\cal C}$
holds. Then the following two conditions are equivalent:
\begin{enumerate}
\item For any objects $W_1$, $W_2$, $W_3$, $W_4$ and $M_1$ of ${\cal
C}$, any nonzero complex numbers $z_1$ and $z_2$ satisfying
$|z_1|>|z_2|>|z_0|>0$, any $P(z_1)$-intertwining map $I_1$ of type
${W_4}\choose {W_1 M_1}$ and $P(z_2)$-intertwining map $I_2$ of type
${M_1}\choose {W_2W_3}$, and any $w'_{(4)}\in W'_4$, 
\[
(I_1\circ (1_{W_1}\otimes I_2))'(w'_{(4)})\in (W_1\otimes W_2\otimes
W_3)^{*}
\]
satisfies the $P^{(2)}(z_0)$-local grading restriction condition (or
the $L(0)$-semisimple $P^{(2)}(z_0)$-local grading restriction
condition when $\mathcal{C}$ is in $\mathcal{M}_{sg}$).  Moreover, for
any $w_{(3)}\in W_{3}$ and $n\in \R$, the smallest doubly graded
subspace of $W^{(2)}_{(I_1\circ (1_{W_1}\otimes I_2))'(w'_{(4)}),
w_{(3)}}$ containing the term $\lambda_{n}^{(2)}$ of the (unique)
series $\sum_{n\in \R}\lambda_{n}^{(2)}$ weakly absolutely convergent
to $\mu^{(2)}_{(I_1\circ (1_{W_1}\otimes I_2))'(w'_{(4)}), w_{(3)}}$
as indicated in the $P^{(2)}(z_0)$-grading condition (or the
$L(0)$-semisimple $P^{(2)}(z_0)$-grading condition) and stable under
the action of $V$ and of $\mathfrak{sl}(2)$ (which is a generalized
$V$-module (or a $V$-module) by Theorem \ref{9.7-1}) is a generalized
V-submodule (or a $V$-submodule) of some object of $\mathcal{C}$
included in $(W_{1}\otimes W_{2})^{*}$.

\item For any objects $W_1$, $W_2$, $W_3$, $W_4$ and $M_2$ of ${\cal
C}$, any nonzero complex numbers $z_1$ and $z_2$ satisfying
$|z_1|>|z_2|>|z_0|>0$, any $P(z_2)$-intertwining map $I^1$ of type
${W_4}\choose {M_2 W_3}$ and $P(z_0)$-intertwining map $I^2$ of type
${M_2}\choose {W_1W_2}$, and any $w'_{(4)}\in W'_4$,
\[
(I^1\circ (I^2\otimes 1_{W_3}))'(w'_{(4)})\in (W_1\otimes W_2\otimes
W_3)^{*}
\]
satisfies the $P^{(1)}(z_2)$-local grading restriction condition (or
the $L(0)$-semisimple $P^{(1)}(z_2)$-local grading restriction
condition when $\mathcal{C}$ is in $\mathcal{M}_{sg}$).  Moreover, for
any $w_{(1)}\in W_{1}$ and $n\in \R$, the smallest doubly graded
subspace of $W^{(1)}_{(I^1\circ (I^2\otimes 1_{W_3}))'(w'_{(4)}),
w_{(1)}}$ containing the term $\lambda_{n}^{(1)}$ of the (unique)
series $\sum_{n\in \R}\lambda_{n}^{(1)}$ weakly absolutely convergent
to $\mu^{(1)}_{(I^1\circ (I^2\otimes 1_{W_3}))'(w'_{(4)}), w_{(1)}}$
as indicated in the $P^{(1)}(z_2)$-grading condition (or the
$L(0)$-semisimple $P^{(1)}(z_2)$-grading condition) and stable under
the action of $V$ and of $\mathfrak{sl}(2)$ (which is a generalized
$V$-module (or a $V$-module) by Theorem \ref{9.7-1}) is a generalized
V-submodule (or a $V$-submodule) of some object of $\mathcal{C}$
included in $(W_{2}\otimes W_{3})^{*}$.
\end{enumerate}
\end{theo}
\pf By Propositions \ref{im:correspond} and \ref{9.7}, together with
Proposition \ref{unique-lambda-n}, Condition 1 (respectively,
Condition 2) in Theorem \ref{asso-io} implies Condition 1
(respectively, Condition 2) in the present theorem.  Conversely, by
Proposition \ref{im:correspond} and Theorem \ref{lgr=>asso} (also
recall the formulation in Corollary \ref{lgr=>asso-op}), Condition 1
(respectively, Condition 2) in the present theorem implies Condition 1
(respectively, Condition 2) in Theorem \ref{asso-io}.  Thus the
present theorem follows immediately {}from Theorem \ref{asso-io}.
\epfv

We are finally ready to define formally, in the following precise
sense, the main concept whose theory has been developed in this
section:

\begin{defi}\label{expansion-conditions}
{\rm Assume that $\mathcal{C}$ is closed under images and under the
$P(z)$-tensor product operation for some $z\in \C^{\times}$, and that
the convergence condition for intertwining maps in ${\cal C}$ holds.
We call either of the two equivalent conditions in Theorem
\ref{expansion} the {\it expansion condition for intertwining maps in
the category ${\cal C}$}.}
\end{defi}

Then Theorem \ref{lgr=>asso} can be reformulated as the following
result, stating that the convergence and expansion conditions together
with certain ``minor'' conditions imply both versions of associativity
of intertwining maps:

\begin{theo}\label{conv-exp=>asso}
Assume that $\mathcal{C}$ is closed under images and under the
$P(z)$-tensor product operation for some $z\in \C^{\times}$, and that
the convergence condition and the expansion condition for intertwining
maps in the category ${\cal C}$ hold, and assume that
\[
|z_1|>|z_2|>|z_{0}|>0.
\]
Let $W_{1}$, $W_{2}$, $W_{3}$,
$W_{4}$, $M_{1}$ and $M_{2}$ be objects of $\mathcal{C}$. 
\begin{enumerate}

\item Let $I_1$ and $I_2$ be $P(z_1)$- and $P(z_2)$-intertwining maps
of types ${W_4 \choose W_1\, M_1}$ and ${M_1 \choose W_2\, W_3}$,
respectively.  Then there exists a unique $P(z_2)$-intertwining map $I^{1}$ of
type ${W_4\choose W_1\boxtimes_{P(z_0)} W_2\,\,W_3}$ such that
\[
\langle w'_{(4)},I_1(w_{(1)} \otimes I_2(w_{(2)} \otimes w_{(3)}))\rangle
=\langle w'_{(4)}, I^{1}((w_{(1)}\boxtimes_{P(z_0)} w_{(2)})\otimes 
w_{(3)})\rangle
\]
for all $w_{(1)}\in W_1$, $w_{(2)}\in W_2$, $w_{(3)}\in W_3$ and
$w'_{(4)}\in W'_4$.  

\item Analogously, let $I^1$ and $I^2$ be $P(z_2)$- and
$P(z_0)$-intertwining maps of types ${W_4 \choose M_2\, W_3}$ and
${M_2 \choose W_1\, W_2}$, respectively.  Then there exists a unique
$P(z_1)$-intertwining map $I_{1}$ of type ${W_4\choose
W_1\,\,W_2\boxtimes_{P(z_2)} W_3}$ such that
\[
\langle w'_{(4)}, I^1(I^2(w_{(1)}\otimes w_{(2)})\otimes w_{(3)})\rangle
=\langle w'_{(4)}, I_{1}(w_{(1)}\otimes (w_{(2)}\boxtimes_{P(z_2)}w_{(3)}))
\rangle
\]
for all $w_{(1)}\in W_1$, $w_{(2)}\in W_2$, $w_{(3)}\in W_3$ and
$w'_{(4)}\in W'_4$. \epf

\end{enumerate}
\end{theo}

We also have the corresponding reformulation of Corollary
\ref{lgr=>asso-op}, asserting the associativity of logarithmic and of
ordinary intertwining operators, under our global conditions:

\begin{corol}\label{conv-exp=>asso-op}
Assume that $\mathcal{C}$ is closed under images and under the
$P(z)$-tensor product operation for some $z\in \C^{\times}$, and that
the convergence condition and the expansion condition for intertwining
maps in the category ${\cal C}$ hold, and assume that
\[
|z_1|>|z_2|>|z_{0}|>0.
\]
Let $W_{1}$, $W_{2}$, $W_{3}$,
$W_{4}$, $M_{1}$ and $M_{2}$ be objects of $\mathcal{C}$. 
\begin{enumerate}

\item Let $\Y_{1}$ and $\Y_{2}$ be logarithmic intertwining operators
(ordinary intertwining operators in the case that $\mathcal{C}$ is in
$\mathcal{M}_{sg}$) of types ${W_4}\choose {W_1M_1}$ and ${M_1}\choose
{W_2W_3}$, respectively.  Then there exists a unique logarithmic
intertwining operator (a unique ordinary intertwining operator in the
case that $\mathcal{C}$ is in $\mathcal{M}_{sg}$) $\Y^{1}$ of type
${W_4\choose W_1\boxtimes_{P(z_0)} W_2\,\,W_3}$ such that
\begin{eqnarray*}
\lefteqn{\langle w'_{(4)},\Y_1(w_{(1)}, x_{1}) \Y_2(w_{(2)}, x_{2}) w_{(3)}\rangle
\lbar_{x_{1}=z_{1},\;x_{2}=z_{2}}}\nn
&&=\langle w'_{(4)}, \Y^{1}(\Y_{\boxtimes_{P(z_{0})}, 0}
(w_{(1)}, x_{0})w_{(2)}, x_{2})w_{(3)})\rangle
\lbar_{x_{0}=z_{0},\;x_{2}=z_{2}}
\end{eqnarray*}
(recalling (\ref{recover}), (\ref{iterabbr}) and (\ref{prodabbr})) for
all $w_{(1)}\in W_1$, $w_{(2)}\in W_2$, $w_{(3)}\in W_3$ and
$w'_{(4)}\in W'_4$.  In particular, the product of the logarithmic
intertwining operators (ordinary intertwining operators in the case
that $\mathcal{C}$ is in $\mathcal{M}_{sg}$) $\Y_{1}$ and $\Y_{2}$
evaluated at $z_{1}$ and $z_{2}$, respectively, can be expressed as an
iterate (with the intermediate generalized $V$-module
$W_1\boxtimes_{P(z_0)} W_2$) of logarithmic intertwining operators
(ordinary intertwining operators in the case that $\mathcal{C}$ is in
$\mathcal{M}_{sg}$) evaluated at $z_{2}$ and $z_{0}$.

\item Analogously, let $\Y^1$ and $\Y^2$ be logarithmic intertwining
operators (ordinary intertwining operators in the case that
$\mathcal{C}$ is in $\mathcal{M}_{sg}$) of types ${W_4}\choose
{M_2W_3}$ and ${M_2}\choose {W_1W_2}$, respectively. Then there exists
a unique logarithmic intertwining operator (a unique ordinary
intertwining operator in the case that $\mathcal{C}$ is in
$\mathcal{M}_{sg}$) $\Y_{1}$ of type ${W_4\choose
W_1\,\,W_2\boxtimes_{P(z_2)} W_3}$ such that
\begin{eqnarray*}
\lefteqn{\langle w'_{(4)}, \Y^1(\Y^2(w_{(1)}, x_{0})w_{(2)}, x_{2})w_{(3)}\rangle
\lbar_{x_{0}=z_{0},\;x_{2}=z_{2}}}\nn
&&=\langle w'_{(4)}, \Y_{1}(w_{(1)}, x_{1})\Y_{\boxtimes_{P(z_2)}, 0}(w_{(2)}, x_{2})w_{(3)}
\rangle\lbar_{x_{1}=z_{1},\;x_{2}=z_{2}}
\end{eqnarray*}
(again recalling (\ref{recover}), (\ref{iterabbr}) and
(\ref{prodabbr})) for all $w_{(1)}\in W_1$, $w_{(2)}\in W_2$,
$w_{(3)}\in W_3$ and $w'_{(4)}\in W'_4$.  In particular, the iterate
of the logarithmic intertwining operators (ordinary intertwining
operators in the case that $\mathcal{C}$ is in $\mathcal{M}_{sg}$)
$\Y^{1}$ and $\Y^{2}$ evaluated at $z_{2}$ and $z_{0}$, respectively,
can be expressed as a product (with the intermediate generalized
$V$-module $W_2\boxtimes_{P(z_2)} W_3$) of logarithmic intertwining
operators (ordinary intertwining operators in the case that
$\mathcal{C}$ is in $\mathcal{M}_{sg}$) evaluated at $z_{1}$ and
$z_{2}$.\epf

\end{enumerate}
\end{corol}

\begin{rema}
{\rm Theorem \ref{conv-exp=>asso} or, respectively, Corollary
\ref{conv-exp=>asso-op}, in fact says that the product of $I_{1}$ and
$I_{2}$ or, respectively, the product of $\Y_{1}$ and $\Y_{2}$,
uniquely ``factors through'' $W_{1}\boxtimes_{P(z_{0})} W_{2}$, and
analogously, that the iterate of $I^{1}$ and $I^{2}$ or, respectively,
the iterate of $\Y^{1}$ and $\Y^{2}$, uniquely ``factors through''
$W_{2}\boxtimes_{P(z_{2})} W_{3}$; cf. Remark \ref{factor-thr}.}
\end{rema}

\newpage

\setcounter{equation}{0}
\setcounter{rema}{0}

\section{The associativity isomorphisms}

We are now in a position to construct our associativity isomorphisms,
assuming the convergence and expansion conditions for intertwining
maps. The strategy and steps in our construction in this section are
essentially the same as those in \cite{tensor4} in the finitely
reductive case but in place of the corresponding results in
\cite{tensor1}, \cite{tensor2}, \cite{tensor3} and \cite{tensor4}, we
have to use virtually all the constructions and results that we have
obtained so far in this work. We remark that the construction
presented here will make the proofs of the coherence and other
properties in our construction of braided tensor category structure
straightforward.  At the end of this section, we show that the
validity of the expansion condition is forced by the assumption of the
existence of natural associativity maps, no matter how they may be
constructed; this exhibits the naturality of the expansion condition.

In the remainder of this work, in addition to Assumptions \ref{assum},
\ref{assum-c} and \ref{assum-exp-set}, we shall also assume that our
category $\mathcal{C}$ is closed under images and that for some $z\in
\C^{\times}$, $\mathcal{C}$ is closed under $P(z)$-tensor products,
that is, the $P(z)$-tensor product of $W_{1}, W_{2}\in \ob
\mathcal{C}$ exists (in $\mathcal{C}$).  For the reader's convenience,
we combine all these assumptions as follows:

\begin{assum}\label{assum-assoc}
Throughout the remainder of this work, we shall assume the following,
unless other assumptions are explicitly made: 

\begin{enumerate}

\item $A$ is an abelian group
and $\tilde{A}$ is an abelian group containing $A$ as a subgroup.

\item $V$
is a strongly $A$-graded M\"{o}bius or conformal vertex algebra.

\item All
$V$-modules and generalized $V$-modules considered are strongly
$\tilde{A}$-graded.

\item All intertwining operators and logarithmic
intertwining operators considered are grading-compatible.

\item $\mathcal{C}$ is a full subcategory of the category
$\mathcal{M}_{sg}$ or $\mathcal{G}\mathcal{M}_{sg}$ (recall Notation
\ref{MGM}).

\item For any object of ${\cal C}$, the (generalized) weights are real
numbers and in addition there exists $K\in \Z_{+}$ such that
$(L(0)-L(0)_{s})^{K}=0$ on the generalized module (when $\mathcal{C}$
is in $\mathcal{M}_{sg}$, the latter assertion holds vacuously).

\item $\mathcal{C}$ is closed under images, under the contragredient
functor, under taking finite direct sums, and under $P(z)$-tensor
products for some $z\in \C^{\times}$.

\end{enumerate}
\end{assum}

\begin{rema}\label{boxtensorchoice}
{\rm {}From Proposition \ref{4.19}, for {\it every} $z\in
\C^{\times}$, $\mathcal{C}$ is closed under $P(z)$-tensor products.
Also, by Proposition \ref{tensor1-13.7}, the assumption that
$\mathcal{C}$ is closed under $P(z)$-tensor products for some $z\in
\C^{\times}$ is equivalent to the assumption that for any $W_{1},
W_{2}\in \ob \mathcal{C}$, $W_{1}\hboxtr_{P(z)}W_{2}$ is an object of
$\mathcal{C}$, and in this case,
\[
W_{1}\boxtimes_{P(z)}W_{2} = (W_{1}\hboxtr_{P(z)}W_{2})'.
\]
{\it From now on, for $z\in \C^{\times}$ we shall take our tensor
product bifunctor $\boxtimes_{P(z)}$ to be}
\[
\boxtimes_{P(z)} = \hboxtr_{P(z)}'.
\]
}
\end{rema}

We shall construct our associativity isomorphisms using the next
theorem.  The proof of this theorem is analogous to the proof of the
corresponding statement in Theorem 14.10 in \cite{tensor4}, but the
results used in the proof below are those developed in the present
work from Section 2 through Section 9.  We shall be using the usual
notation
\[
\overline{\eta}:\overline{W_1} \to \overline{W_2}
\]
to denote the natural extension of a map $\eta:W_1 \to W_2$ of
generalized modules to the formal completions.

We shall be constructing a natural isomorphism between the two
functors from ${\cal C} \times {\cal C} \times {\cal C}$ to ${\cal C}$
in (\ref{naturaliso}) below.  We first determine how the functor
\[
\boxtimes_{P(z_{1})}\circ (1 \times \boxtimes_{P(z_2)})
\]
acts on maps and elements when the convergence condition holds and
when $|z_1|>|z_2|>0$: Consider maps
\[
\sigma_1 : W_1 \rightarrow W_4,
\]
\[
\sigma_2 : W_2 \rightarrow W_5,
\]
\[
\sigma_3 : W_3 \rightarrow W_6
\]
between objects of ${\cal C}$.  Recall from Remark \ref{bifunctor}
that for $z \in \C^{\times}$ the functor $\boxtimes_{P(z)}$ acts on
maps and elements by:
\begin{equation}\label{sigma1sigma2}
\overline{\sigma_1\boxtimes_{P(z)}\sigma_2}(w_{(1)}\boxtimes_{P(z)} w_{(2)})
=\sigma_1 (w_{(1)}) \boxtimes_{P(z)} \sigma_2 (w_{(2)}),
\end{equation}
and recall that by Proposition \ref{span}, (\ref{sigma1sigma2})
determines the $V$-module map $\sigma_1\boxtimes_{P(z)}\sigma_2$
uniquely.  We have
\[
(\boxtimes_{P(z_{1})}\circ (1 \times \boxtimes_{P(z_2)}))
(\sigma_1, \sigma_2, \sigma_3)
=\boxtimes_{P(z_1)}(\sigma_1, \sigma_2 \boxtimes_{P(z_2)}
\sigma_3)
=\sigma_1 \boxtimes_{P(z_1)} (\sigma_2 \boxtimes_{P(z_2)} \sigma_3),
\]
and the effect of this map on elements is determined as follows:
Since
\[
\overline{\sigma_2\boxtimes_{P(z_2)}\sigma_3}(w_{(2)}\boxtimes_{P(z_2)}
w_{(3)}) =\sigma_2 (w_{(2)}) \boxtimes_{P(z_2)} \sigma_3 (w_{(3)}),
\]
we have (using the projection notation $\pi_n$)
\[
(\sigma_2\boxtimes_{P(z_2)}\sigma_3)(\pi_n (w_{(2)}\boxtimes_{P(z_2)}
w_{(3)}))
=\pi_n (\sigma_2 (w_{(2)}) \boxtimes_{P(z_2)} \sigma_3 (w_{(3)}))
\]
for all $n \in \R$, so that
\begin{eqnarray*}
\lefteqn{\overline{\sigma_1 \boxtimes_{P(z_1)}(\sigma_2\boxtimes_{P(z_2)}\sigma_3)}
(w_{(1)}\boxtimes_{P(z_1)} \pi_n (w_{(2)}\boxtimes_{P(z_2)} w_{(3)}))}\nn
&&= \sigma_1 (w_{(1)})\boxtimes_{P(z_1)} \pi_n (\sigma_2 (w_{(2)})
\boxtimes_{P(z_2)} \sigma_3 (w_{(3)})).
\end{eqnarray*}
Thus for 
\[
w' \in (W_4\boxtimes_{P(z_{1})} (W_5\boxtimes_{P(z_2)} W_6))',
\]
we have
\begin{eqnarray*}
\lefteqn{\langle (\sigma_1 \boxtimes_{P(z_1)}(\sigma_2\boxtimes_{P(z_2)}\sigma_3))'(w'),
w_{(1)}\boxtimes_{P(z_1)} \pi_n (w_{(2)}\boxtimes_{P(z_2)} w_{(3)})
\rangle}\nn
&&=\langle w',
\sigma_1 (w_{(1)})\boxtimes_{P(z_1)} \pi_n (\sigma_2 (w_{(2)})
\boxtimes_{P(z_2)} \sigma_3 (w_{(3)}))\rangle,
\end{eqnarray*}
and so by the convergence condition,
\begin{eqnarray*}
\lefteqn{\langle (\sigma_1
\boxtimes_{P(z_1)}(\sigma_2\boxtimes_{P(z_2)}\sigma_3))'(w'),
w_{(1)}\boxtimes_{P(z_1)} (w_{(2)}\boxtimes_{P(z_2)} w_{(3)})
\rangle}\nn
&&=\langle w',
\sigma_1 (w_{(1)})\boxtimes_{P(z_1)} (\sigma_2 (w_{(2)})
\boxtimes_{P(z_2)} \sigma_3 (w_{(3)}))\rangle.
\end{eqnarray*}
Hence
\begin{eqnarray}\label{sigma1(23)onelements}
\lefteqn{\overline{\sigma_1 \boxtimes_{P(z_1)}(\sigma_2\boxtimes_{P(z_2)}\sigma_3)}
(w_{(1)}\boxtimes_{P(z_1)} (w_{(2)}\boxtimes_{P(z_2)} w_{(3)}))}\nn
&&=\sigma_1 (w_{(1)})\boxtimes_{P(z_1)} (\sigma_2 (w_{(2)})
\boxtimes_{P(z_2)} \sigma_3 (w_{(3)})),
\end{eqnarray}
and by Corollary \ref{prospan}, this determines the $V$-module map
$\sigma_1 \boxtimes_{P(z_1)}(\sigma_2\boxtimes_{P(z_2)}\sigma_3)$
uniquely (cf. (\ref{sigma1sigma2})).  Analogously, when the
convergence condition holds and when $|z_2|>|z_1-z_{2}|>0,$ the
functor
\[
\boxtimes_{P(z_2)}\circ (\boxtimes_{P(z_1-z_2)}\times 1)
\]
acts on elements by:
\begin{eqnarray}\label{sigma(12)3onelements}
\lefteqn{\overline{(\sigma_1 \boxtimes_{P(z_1-z_2)}\sigma_2)
\boxtimes_{P(z_2)}\sigma_3}
((w_{(1)}\boxtimes_{P(z_1-z_2)} w_{(2)})\boxtimes_{P(z_2)} w_{(3)})}\nn
&&= (\sigma_1 (w_{(1)})\boxtimes_{P(z_1-z_2)} \sigma
(w_{(2)}))\boxtimes_{P(z_2)} \sigma_3 (w_{(3)}),
\end{eqnarray}
and by Corollary \ref{iterspan}, this determines the corresponding
$V$-module map uniquely.  The formulas (\ref{sigma1(23)onelements})
and (\ref{sigma(12)3onelements}), which extend (\ref{sigma1sigma2}),
are crucial.

\begin{theo}\label{assoc-thm}
Assume that the convergence condition and the expansion condition for
intertwining maps in ${\cal C}$ (see Definitions \ref{conv-conditions}
and \ref{expansion-conditions}) both hold.  Let $z_1$, $z_2$ be
complex numbers satisfying
\[
|z_1|>|z_2|>|z_1-z_{2}|>0
\]
(so that in particular, $z_1\neq 0$, $z_2\neq 0$ and $z_1\neq z_2$).
Then there exists a unique natural isomorphism
\begin{equation}\label{naturaliso}
\mathcal{A}_{P(z_{1}), P(z_{2})}^{P(z_{1}-z_{2}), P(z_{2})}: 
\boxtimes_{P(z_{1})}\circ (1 \times \boxtimes_{P(z_2)})
\to \boxtimes_{P(z_2)}\circ (\boxtimes_{P(z_1-z_2)}\times 1)
\end{equation}
such that for all $w_{(1)}\in W_1$, $w_{(2)}\in W_2$ and
$w_{(3)}\in W_3$, with $W_j$ objects of ${\cal C}$,
\begin{equation}\label{assoc-elt-1}
\overline{\mathcal{A}_{P(z_{1}), P(z_{2})}^{P(z_{1}-z_{2}), P(z_{2})}}
(w_{(1)}\boxtimes_{P(z_1)}
(w_{(2)}\boxtimes_{P(z_2)} w_{(3)})) = (w_{(1)}\boxtimes_{P(z_1-z_2)}
w_{(2)})\boxtimes_{P(z_2)} w_{(3)},
\end{equation}
where for simplicity we use the same notation $\mathcal{A}_{P(z_{1}),
P(z_{2})}^{P(z_{1}-z_{2}), P(z_{2})}$ to denote the isomorphism of
(generalized) modules
\begin{equation}\label{assoc-iso}
\mathcal{A}_{P(z_{1}), P(z_{2})}^{P(z_{1}-z_{2}), P(z_{2})}: 
W_1\boxtimes_{P(z_{1})}
(W_2\boxtimes_{P(z_2)} W_3) \longrightarrow (W_1\boxtimes_{P(z_1-z_2)}
W_2)\boxtimes_{P(z_2)} W_3.
\end{equation}
\end{theo}
\pf The uniqueness of $\mathcal{A}_{P(z_{1}),
P(z_{2})}^{P(z_{1}-z_{2}), P(z_{2})}$ follows {}from Corollary
\ref{prospan}.

Let $W_{P(z_1, z_2)}$ be the subspace of 
\[
(W_1\otimes
W_2\otimes W_3)^{*}
\] 
consisting of the elements $\lambda$ satisfying the following
conditions:
\begin{enumerate}

\item The
$P(z_1,z_2)$-compatibility condition (see Section 8).

\item The $P(z_1, z_2)$-local grading restriction condition (the
$L(0)$-semisimple $P(z_1, z_2)$-local grading restriction condition in
the case that $\mathcal{C}$ is in $\mathcal{M}_{sg}$) (see Section 8).

\item Either one of the following conditions (see Section 9):

\begin{enumerate}

\item The $P^{(1)}(z_{2})$-local grading restriction
condition (the $L(0)$-semisimple $P^{(1)}(z_{2})$-local grading 
restriction condition in the case that  $\mathcal{C}$ is in 
$\mathcal{M}_{sg}$).

\item The $P^{(2)}(z_{1}-z_{2})$-local grading restriction
condition (the $L(0)$-semisimple $P^{(2)}(z_{1}-z_{2})$-local grading restriction
condition in the case that  $\mathcal{C}$ is in 
$\mathcal{M}_{sg}$).

\end{enumerate}

\item Either one of the following conditions, depending on which
condition is satisfied in 3 above (that is, either 3(a) and 4(a) hold
or 3(b) and 4(b) hold):

\begin{enumerate}

\item For any $w_{(1)}\in W_{1}$ and $n\in \R$, the smallest doubly
graded subspace of $W^{(1)}_{\lambda, w_{(1)}}$ containing the term
$\lambda_{n}^{(1)}$ of the (unique) series $\sum_{n\in
\R}\lambda_{n}^{(1)}$ weakly absolutely convergent to
$\mu^{(1)}_{\lambda, w_{(1)}}$ as indicated in the
$P^{(1)}(z_2)$-grading condition (or the $L(0)$-semisimple
$P^{(1)}(z_2)$-grading condition when $\mathcal{C}$ is in
$\mathcal{M}_{sg}$) and stable under the action of $V$ and of
$\mathfrak{sl}(2)$ is a generalized $V$-module (or a $V$-module) and
is in fact a generalized $V$-submodule (or a $V$-submodule) of some
object of $\mathcal{C}$ included in $(W_{2}\otimes W_{3})^{*}$.

\item For any $w_{(3)}\in W_{3}$ and $n\in \R$, the smallest doubly
graded subspace of $W^{(2)}_{\lambda, w_{(3)}}$ containing the term
$\lambda_{n}^{(2)}$ of the (unique) series $\sum_{n\in
\R}\lambda_{n}^{(2)}$ weakly absolutely convergent to
$\mu^{(2)}_{\lambda, w_{(3)}}$ as indicated in the
$P^{(2)}(z_1-z_2)$-grading condition (or the $L(0)$-semisimple
$P^{(2)}(z_1-z_2)$-grading condition when $\mathcal{C}$ is in
$\mathcal{M}_{sg}$) and stable under the action of $V$ and of
$\mathfrak{sl}(2)$ is a generalized $V$-module (or a $V$-module) and
is in fact a generalized $V$-submodule (or a $V$-submodule) of some
object of $\mathcal{C}$ included in $(W_{1}\otimes W_{2})^{*}$.

\end{enumerate}

\end{enumerate}

By Proposition \ref{productanditerateareintwmaps}, the natural map
\[
\boxtimes_{P(z_{1})}\circ (1_{W_{1}}\otimes \boxtimes_{P(z_{2})})
: W_1\otimes W_2\otimes W_3 \to
\overline{W_1\boxtimes_{P(z_1)} (W_2\boxtimes_{P(z_2)} W_3)}
\]
is a $P(z_1, z_2)$-intertwining map.  Recalling Remark
\ref{Atildecompatcorrespondence}, let
\[
\Psi^{(1)}_{P(z_1, z_2)}
=(\boxtimes_{P(z_{1})}\circ (1_{W_{1}}\otimes \boxtimes_{P(z_{2})}))',
\]
the natural map {}from
\[
W_1\hboxtr_{P(z_1)} (W_2\boxtimes_{P(z_2)} W_3)=(W_1\boxtimes_{P(z_1)}
(W_2\boxtimes_{P(z_2)} W_3))'
\]
to $(W_1\otimes W_2\otimes W_3)^*$ given by
\begin{eqnarray}\label{Psi1}
\lefteqn{(\Psi^{(1)}_{P(z_1, z_2)}(\nu))(w_{(1)}\otimes
w_{(2)}\otimes w_{(3)})}\nn
&&=\langle \nu,
w_{(1)}\boxtimes_{P(z_1)}(w_{(2)}
\boxtimes_{P(z_2)}w_{(3)})\rangle_{W_1\boxtimes_{P(z_1)}(W_2
\boxtimes_{P(z_2)}W_3)}
\end{eqnarray}
for
\[
\nu \in W_1\hboxtr_{P(z_1)}(W_2\boxtimes_{P(z_2)}W_3),
\]
$w_{(1)}\in W_1$, $w_{(2)}\in W_2$ and $w_{(3)}\in W_3$.  Then by
Proposition \ref{zzcor}, $\Psi^{(1)}_{P(z_1, z_2)}$ is an
$\tilde{A}$-compatible map and it intertwines the actions of
\[
V\otimes \iota _{+}{\mathbb C}[t,t^{-1},(z_1^{-1}-t)^{-1},(z_2^{-1}-t)^{-1}]
\]
and of $L'_{P(z_{1})}(j)$ and $L'_{P(z_1,z_2)}(j)$, $j=-1,0,1$, on
$W_1\hboxtr_{P(z_1)} (W_2\boxtimes_{P(z_2)} W_3)$ and on $(W_1\otimes
W_2\otimes W_3)^*$.  In particular,
\[
\Psi^{(1)}_{P(z_1, z_2)} \circ L'_{P(z_{1})}(0)
=L'_{P(z_{1}, z_{2})}(0) \circ \Psi^{(1)}_{P(z_1, z_2)}
\]
and
\[
\Psi^{(1)}_{P(z_1, z_2)}\circ Y'_{P(z_{1})}(u, x) 
=Y'_{P(z_{1}, z_{2})}(u, x) \circ \Psi^{(1)}_{P(z_1, z_2)}
\]
for $u\in V$.  Thus $\Psi^{(1)}_{P(z_1, z_2)}$ preserves generalized
weights, the image
\[
\Psi^{(1)}_{P(z_1, z_2)}(W_1\hboxtr_{P(z_1)} (W_2\boxtimes_{P(z_2)} W_3))
\]
of $\Psi^{(1)}_{P(z_1, z_2)}$ is a generalized module (recall the
proof of Proposition \ref{8.12}), and 
$\Psi^{(1)}_{P(z_1, z_2)}$ is a map of 
generalized modules from 
\[
W_1\hboxtr_{P(z_1)} (W_2\boxtimes_{P(z_2)} W_3)
\]
to this image.

By Propositions \ref{productanditerateareintwmaps}, \ref{8.12}
and \ref{9.7}, $\Psi^{(1)}_{P(z_1, z_2)}$
in fact maps $W_1\hboxtr_{P(z_1)} (W_2\boxtimes_{P(z_2)} W_3)$
into $W_{P(z_1, z_2)}$; the elements of the image satisfy 3(a) and
4(a).  By the expansion condition, the elements of this image also
satisfy 3(b) and 4(b).

Analogously, let
\[
\Psi^{(2)}_{P(z_1, z_2)}=
(\boxtimes_{P(z_{2})}\circ (\boxtimes_{P(z_{1}-z_{2})}\otimes 1_{W_{3}}))',
\]
the natural map {}from
\[
(W_1\boxtimes_{P(z_1-z_2)}W_2)\hboxtr_{P(z_2)}W_3
=((W_1\boxtimes_{P(z_1-z_2)}W_2)\boxtimes_{P(z_2)}W_3)'
\]
to $(W_1\otimes W_2\otimes W_3)^{*}$ given by
\begin{eqnarray}\label{Psi2}
\lefteqn{(\Psi^{(2)}_{P(z_1, z_2)}(\xi))(w_{(1)}\otimes w_{(2)}
\otimes w_{(3)})}\nn
&&=\langle \xi, (w_{(1)}\boxtimes_{P(z_1-z_2)}w_{(2)})
\boxtimes_{P(z_2)}w_{(3)}\rangle_{(W_1\boxtimes_{P(z_1-z_2)}W_2)
\boxtimes_{P(z_2)}W_3}
\end{eqnarray}
for
\begin{equation}\label{xi}
\xi \in (W_1\boxtimes_{P(z_1-z_2)}W_2)\hboxtr_{P(z_2)}W_3,
\end{equation}
$w_{(1)}\in W_1$, $w_{(2)}\in W_2$ and $w_{(3)}\in W_3$.  Again by
Propositions \ref{productanditerateareintwmaps} and \ref{zzcor},
$\Psi^{(2)}_{P(z_1, z_2)}$ is an $\tilde{A}$-compatible map and it
intertwines the actions of
\[
V\otimes \iota _{+}{\mathbb C}[t,t^{-1},(z_1^{-1}-t)^{-1},(z_2^{-1}-t)^{-1}]
\]
and of $L'_{P(z_{2})}(j)$ and $L'_{P(z_1,z_2)}(j)$, $j=-1,0,1$, on
$(W_1\boxtimes_{P(z_1-z_2)}W_2)\hboxtr_{P(z_2)}W_3$ and on
$(W_1\otimes W_2\otimes W_3)^*$, and in particular,
\[
\Psi^{(2)}_{P(z_1, z_2)} \circ L'_{P(z_{2})}(0)
=L'_{P(z_{1}, z_{2})}(0) \circ \Psi^{(2)}_{P(z_1, z_2)}
\]
and
\[
\Psi^{(2)}_{P(z_1, z_2)}\circ Y'_{P(z_{2})}(u, x) 
=Y'_{P(z_{1}, z_{2})}(u, x) \circ \Psi^{(2)}_{P(z_1, z_2)}
\]
for $u\in V$.  Thus $\Psi^{(2)}_{P(z_1, z_2)}$ preserves generalized
weights, the image
\[
\Psi^{(2)}_{P(z_1, z_2)}((W_1\boxtimes_{P(z_1-z_2)}W_2)\hboxtr_{P(z_2)}W_3)
\]
of $\Psi^{(2)}_{P(z_1, z_2)}$ is a generalized module, and 
$\Psi^{(2)}_{P(z_1, z_2)}$ is a map of 
generalized modules from 
\[
(W_1\boxtimes_{P(z_1-z_2)}W_2)\hboxtr_{P(z_2)}W_3
\]
to this image.  Again by Propositions
\ref{productanditerateareintwmaps}, \ref{8.12} and
\ref{9.7}, $\Psi^{(2)}_{P(z_1, z_2)}$ maps
$(W_1\boxtimes_{P(z_1-z_2)}W_2)\hboxtr_{P(z_2)}W_3$ into $W_{P(z_1,
z_2)}$; the elements of the image satisfy 3(b) and
4(b).  By the expansion condition, the elements of this image also
satisfy 3(a) and 4(a).

We next show that both $\Psi^{(1)}_{P(z_1,z_2)}$ and
$\Psi^{(2)}_{P(z_1, z_2)}$ are injective.  Let
\[
\nu \in
W_1\hboxtr_{P(z_1)}(W_2\boxtimes_{P(z_2)}W_3)
\]
be such that
\[
\Psi^{(1)}_{P(z_1, z_2)}(\nu)=0,
\]
that is,
\begin{equation}\label{assoc-pf-2}
\langle \nu, w_{(1)}\boxtimes_{P(z_1)}(w_{(2)}
\boxtimes_{P(z_2)}w_{(3)})\rangle_{W_1\boxtimes_{P(z_1)}
(W_2\boxtimes_{P(z_2)}W_3)}=0
\end{equation}
for $w_{(1)}\in W_1$, $w_{(2)}\in W_2$ and $w_{(3)}\in W_3$.  Since
$\Psi^{(1)}_{P(z_1,z_2)}$ preserves generalized weights, we can assume
that $\nu$ is homogeneous.  Then (\ref{assoc-pf-2}) implies that
for all $n \in \R$,
\[
\langle \nu, \pi_{n}(w_{(1)}\boxtimes_{P(z_1)}(w_{(2)}
\boxtimes_{P(z_2)}w_{(3)}))\rangle_{W_1\boxtimes_{P(z_1)}
(W_2\boxtimes_{P(z_2)}W_3)}=0
\]
for $w_{(1)}\in W_1$, $w_{(2)}\in W_2$ and $w_{(3)}\in W_3$.  But by
Corollary \ref{prospan}, the elements
\[
\pi_{n}(w_{(1)}\boxtimes_{P(z_1)}(w_{(2)} \boxtimes_{P(z_2)}w_{(3)}))
\]
for $n\in\R$, $w_{(1)}\in W_1$, $w_{(2)}\in W_2$ and $w_{(3)}\in W_3$
span the space $W_1\boxtimes_{P(z_1)} (W_2\boxtimes_{P(z_2)} W_3)$, so
that $\nu=0$, and we have the injectivity of $\Psi^{(1)}_{P(z_1,
z_2)}$.  The proof of the injectivity of $\Psi^{(2)}_{P(z_1, z_2)}$ is
completely analogous.

Now we want to prove that the images of our two maps are equal:
\[
\Psi^{(1)}_{P(z_1,z_2)}(W_1\hboxtr_{P(z_1)}(W_2\boxtimes_{P(z_2)}W_3))
=\Psi^{(2)}_{P(z_1, z_2)}((W_1\boxtimes_{P(z_1-z_2)}W_2)\hboxtr_{P(z_2)}W_3).
\]
Let
\[
\lambda\in
\Psi^{(1)}_{P(z_1,z_2)}(W_1\hboxtr_{P(z_1)}(W_2\boxtimes_{P(z_2)}W_3))
\]
and take $\nu\in
W_1\hboxtr_{P(z_1)}(W_2\boxtimes_{P(z_2)}W_3)$ so that
$\lambda=\Psi^{(1)}_{P(z_1,z_2)}(\nu)$.  Then by
Theorem \ref{lgr=>asso}, there exists a unique
$P(z_{2})$-intertwining map $I^{1}$ of type
\[
{W_1\boxtimes_{P(z_1)}(W_2\boxtimes_{P(z_2)}W_3) \choose
W_1\boxtimes_{P(z_1-z_2)}W_2 \;\;\;\;\;\;\;\; W_3}
\]
such that
\begin{eqnarray}\label{lambda-I-1}
\lambda(w_{(1)}\otimes w_{(2)}\otimes
w_{(3)})&=&(\Psi^{(1)}_{P(z_1,z_2)}(\nu))(w_{(1)}\otimes w_{(2)}\otimes
w_{(3)})\nn
&=&\langle \nu, w_{(1)}\boxtimes_{P(z_1)}(w_{(2)}\boxtimes_{P(z_2)}
w_{(3)})\rangle\nn
&=&\langle \nu, I^{1}((w_{(1)}\boxtimes_{P(z_{1}-z_{2})} w_{(2)})
\otimes w_{(3)})\rangle
\end{eqnarray}
for $w_{(1)}\in W_{1}$, $w_{(2)}\in W_{2}$ and
$w_{(3)}\in W_{3}$. 

By the definition of the $P(z_{2})$-tensor product of
$W_{1}\boxtimes_{P(z_{1}-z_{2})} W_{2}$ and $W_{3}$ in $\mathcal{C}$, there
exists a unique map of generalized modules
\[
\eta: (W_1\boxtimes_{P(z_1-z_2)}W_2)\boxtimes_{P(z_2)}W_3
\to W_1\boxtimes_{P(z_1)}(W_2\boxtimes_{P(z_2)}W_3)
\]
such that 
\[
I^{1}=\overline{\eta}\circ \boxtimes_{P(z_2)}.
\] 
Then by (\ref{lambda-I-1}), we obtain
\begin{eqnarray*}
\lambda(w_{(1)}\otimes w_{(2)}\otimes
w_{(3)})&=&
\langle \nu, I^{1}
((w_{(1)}\boxtimes_{P(z_{1}-z_{2})} w_{(2)})
\otimes w_{(3)})\rangle\nn
&=&
\langle \nu, (\overline{\eta}\circ \boxtimes_{P(z_2)})
((w_{(1)}\boxtimes_{P(z_{1}-z_{2})} w_{(2)})
\otimes w_{(3)})\rangle\nn
&=&\langle \eta'(\nu),  
(w_{(1)}\boxtimes_{P(z_{1}-z_{2})} w_{(2)})
\boxtimes_{P(z_2)} w_{(3)}\rangle\nn
&=&(\Psi^{(2)}_{P(z_1,z_2)}(\eta'(\nu)))(w_{(1)}\otimes w_{(2)}\otimes
w_{(3)})
\end{eqnarray*}
for $w_{(1)}\in W_{1}$, $w_{(2)}\in W_{2}$ and
$w_{(3)}\in W_{3}$. Thus 
$\lambda=\Psi^{(2)}_{P(z_1,z_2)}(\eta'(\nu))$, 
proving that $\lambda$ lies in $\Psi^{(2)}_{P(z_1,z_2)}
((W_1\boxtimes_{P(z_1-z_2)}W_2)\hboxtr_{P(z_2)}W_3)$.
Thus
\[
\Psi^{(1)}_{P(z_1,z_2)}(W_1\hboxtr_{P(z_1)}(W_2\boxtimes_{P(z_2)}W_3))
\subset \Psi^{(2)}_{P(z_1, z_2)}((W_1\boxtimes_{P(z_1-z_2)}W_2)
\hboxtr_{P(z_2)}W_3).
\]
The proof of the opposite inclusion,
\[
\Psi^{(1)}_{P(z_1,z_2)}(W_1\hboxtr_{P(z_1)}(W_2\boxtimes_{P(z_2)}W_3))
\supset \Psi^{(2)}_{P(z_1, z_2)}((W_1\boxtimes_{P(z_1-z_2)}W_2)
\hboxtr_{P(z_2)}W_3),
\]
is completely analogous.

Since
\[
\Psi^{(1)}_{P(z_1,z_2)}:
W_1\hboxtr_{P(z_1)}(W_2\boxtimes_{P(z_2)}W_3)
\to \Psi^{(1)}_{P(z_1,z_2)}(W_1\hboxtr_{P(z_1)}(W_2\boxtimes_{P(z_2)}W_3))
\]
is injective, we have the natural map
\[
(\Psi^{(1)}_{P(z_1,z_2)})^{-1}:
\Psi^{(1)}_{P(z_1,z_2)}(W_1\hboxtr_{P(z_1)}(W_2\boxtimes_{P(z_2)}W_3))
\to W_1\hboxtr_{P(z_1)}(W_2\boxtimes_{P(z_2)}W_3)
\]
of generalized modules.  Thus we have a natural isomorphism
\[
(\Psi^{(1)}_{P(z_1,z_2)})^{-1} \circ \Psi^{(2)}_{P(z_1, z_2)}:
(W_1\boxtimes_{P(z_1-z_2)}W_2)\hboxtr_{P(z_2)}W_3
\to W_1\hboxtr_{P(z_1)} (W_2\boxtimes_{P(z_2)} W_3),
\]
of generalized modules.  Its contragredient map $\mathcal{A}$
hence gives a natural isomorphism (\ref{assoc-iso}), and
(\ref{assoc-elt-1}) indeed holds for the map $\mathcal{A}$.  In fact,
for $\xi$ as in (\ref{xi}), (\ref{Psi2}) holds, but on the other hand,
for
\[
\nu = (\Psi^{(1)}_{P(z_1,z_2)})^{-1}\Psi^{(2)}_{P(z_1, z_2)}(\xi),
\]
(\ref{Psi1}) gives
\[
(\Psi^{(2)}_{P(z_1, z_2)})(\xi)(w_{(1)}\otimes w_{(2)}\otimes
w_{(3)})=
\langle (\Psi^{(1)}_{P(z_1,z_2)})^{-1}\Psi^{(2)}_{P(z_1, z_2)}(\xi),
w_{(1)}\boxtimes_{P(z_1)}(w_{(2)}\boxtimes_{P(z_2)}w_{(3)})\rangle,
\]
and this equals
\[
\langle \xi, \overline{\mathcal{A}}
(w_{(1)}\boxtimes_{P(z_1)}(w_{(2)}\boxtimes_{P(z_2)}w_{(3)}))\rangle.
\]

The formulas (\ref{sigma1(23)onelements}),
(\ref{sigma(12)3onelements}) and (\ref{assoc-elt-1}) exhibit the
naturality of (\ref{naturaliso}).  \epfv

\begin{defi}
{\rm For  $z_1, z_2\in \C$ satisfying
\[
|z_1|>|z_2|>|z_1-z_{2}|>0
\]
and objects $W_1$, $W_2$, and $W_3$ of ${\cal C}$, the {\it
associativity isomorphism from}
\[
W_1\boxtimes_{P(z_{1})} (W_2\boxtimes_{P(z_2)} W_3)
\]
{\it to}
\[
(W_1\boxtimes_{P(z_1-z_2)} W_2)\boxtimes_{P(z_2)} W_3
\]
is the natural isomorphism $\mathcal{A}_{P(z_{1}),
P(z_{2})}^{P(z_{1}-z_{2}), P(z_{2})}$ given in Theorem
\ref{assoc-iso}. We also have the natural {\it inverse associativity
isomorphism}
\begin{equation}\label{assoc-iso-inv}
\alpha_{P(z_{1}), P(z_{2})}^{P(z_{1}-z_{2}), P(z_{2})}: 
(W_1\boxtimes_{P(z_1-z_2)}
W_2)\boxtimes_{P(z_2)} W_3 \to W_1\boxtimes_{P(z_1)}
(W_2\boxtimes_{P(z_2)} W_3).
\end{equation}}
\end{defi}

\begin{rema}{\rm
The inverse associativity isomorphism $\alpha_{P(z_{1}),
P(z_{2})}^{P(z_{1}-z_{2}), P(z_{2})}$ satisfies
\begin{equation}\label{assoc-elt-2}
\overline{\alpha_{P(z_{1}), P(z_{2})}^{P(z_{1}-z_{2}), P(z_{2})}}
((w_{(1)}\boxtimes_{P(z_1-z_2)}
w_{(2)})\boxtimes_{P(z_2)} w_{(3)}) = w_{(1)}\boxtimes_{P(z_1)}
(w_{(2)}\boxtimes_{P(z_2)} w_{(3)})
\end{equation}
for $w_{(1)}\in W_1$, $w_{(2)}\in W_2$ and $w_{(3)}\in W_3$, and
(\ref{assoc-elt-2}) determines (\ref{assoc-iso-inv}) uniquely.  }
\end{rema}

As in the setting of \cite{tensor4}, the existence of such
associativity isomorphisms implies the expansion condition and hence
that products can be expressed as iterates and vice versa; the
following converse of Theorem \ref{assoc-thm} essentially says, then,
that the expansion condition is equivalent to the existence of natural
associativity isomorphisms:

\begin{theo}
Assume that the convergence condition for intertwining maps in ${\cal
C}$ (see Definition \ref{conv-conditions}) holds.  Suppose that for
any complex numbers $z_1$, $z_2$ satisfying
$|z_1|>|z_2|>|z_1-z_{2}|>0$ and any objects $W_1$, $W_2$ and $W_3$ of
$\mathcal{C}$, there exists a map $\mathcal{A}_{P(z_{1}),
P(z_{2})}^{P(z_{1}-z_{2}), P(z_{2})}$ of (generalized) modules of the
form (\ref{assoc-iso}) such that (\ref{assoc-elt-1}) holds for
$w_{(1)}\in W_1$, $w_{(2)}\in W_2$ and $w_{(3)}\in W_3$. Then
$\mathcal{A}_{P(z_{1}), P(z_{2})}^{P(z_{1}-z_{2}), P(z_{2})}$ for
$W_1$, $W_2$ and $W_3$ is uniquely determined and is a module
isomorphism; these maps define a natural isomorphism of functors of
the form (\ref{naturaliso}); and furthermore, the expansion condition
holds, and in particular, products of intertwining maps can be
expressed as iterates and conversely, as in Theorem
\ref{conv-exp=>asso} and Corollary \ref{conv-exp=>asso-op}.
Analogously, suppose that for any complex numbers $z_1$, $z_2$
satisfying $|z_1|>|z_2|>|z_1-z_{2}|>0$ and for any objects $W_1$,
$W_2$ and $W_3$ of $\mathcal{C}$, there exists a map of (generalized)
modules of the form (\ref{assoc-iso-inv}) such that
(\ref{assoc-elt-2}) holds for all $w_{(1)}\in W_1$, $w_{(2)}\in W_2$
and $w_{(3)}\in W_3$.  Then the analogous conclusions hold; in
particular, the expansion condition again holds, and products of
intertwining maps can be expressed as iterates and conversely.
\end{theo}
\pf We shall prove the first half; the second half is proved
analogously. We need only prove that given the maps $\mathcal{A}$ as
indicated, the expansion condition holds; all the other conclusions
are either clear or immediate consequences, using Theorem
\ref{assoc-thm} and its proof.

To prove the expansion condition, we shall prove Condition 2 in
Theorem \ref{expansion}.  With $z_1$ and $z_2$ as indicated, for any
objects $W_1$, $W_2$, $W_3$, $W_4$ and $M_2$ of ${\cal C}$ and any
$P(z_2)$-intertwining map $I^1$ of type ${W_4}\choose {M_2 W_3}$ and
$P(z_1-z_2)$-intertwining map $I^2$ of type ${M_2}\choose {W_1W_2}$,
by Proposition \ref{intermediate} there exists a unique
\[
\widetilde{I}^{1}\in
\mathcal{M}[P(z_{2})]_{(W_{1}\boxtimes_{P(z_1-z_2)}W_{2})\;W_{3}}^{W_{4}}
\]
such that
\[
I^{1}\circ (I^{2}\otimes 1_{W_{3}})
=\widetilde{I}^{1}\circ (\boxtimes_{P(z_1-z_2)} \otimes
1_{W_{3}}).
\]
By the definition of $(W_{1}
\boxtimes_{P(z_1-z_2)}W_{2})\boxtimes_{P(z_{2})}W_{3}$
(cf. Proposition \ref{pz-iso}), there exists a unique module map
\[
\eta: (W_{1} \boxtimes_{P(z_1-z_2)}W_{2})\boxtimes_{P(z_{2})}W_{3}
\to W_{4}
\]
such that 
\[
\widetilde{I}^{1}=\overline{\eta}\circ \boxtimes_{P(z_{2})}.
\]
Then we have
\[
I^{1}\circ (I^{2}\otimes 1_{W_{3}})
=\widetilde{I}^{1}\circ (\boxtimes_{P(z_1-z_2)} \otimes 1_{W_{3}})
=\overline{\eta}\circ (\boxtimes_{P(z_{2})}\circ (\boxtimes_{P(z_1-z_2)} \otimes
1_{W_{3}})).
\]
Thus for $w_{(1)}\in W_{1}$, $w_{(2)}\in W_{2}$ and $w_{(3)}\in
W_{3}$, by (\ref{assoc-elt-1}) we have
\begin{eqnarray}\label{conv-1}
\lefteqn{I^{1}(I^{2}(w_{(1)}\otimes w_{(2)})\otimes w_{(3)})}\nn
&&=\overline{\eta}((w_{(1)}\boxtimes_{P(z_1-z_2)}w_{(2)})
\boxtimes_{P(z_{2})}w_{(3)})\nn
&&=\overline{\eta}(\overline{{\cal A}_{P(z_{1}),
P(z_{2})}^{P(z_{1}-z_{2}), P(z_{2})}}(w_{(1)}\boxtimes_{P(z_{1})}(w_{(2)}
\boxtimes_{P(z_{2})}w_{(3)}))).
\end{eqnarray}
By Proposition \ref{9.7}, for $w'\in W_{1}\hboxtr_{P(z_{1})}
(W_{2}\boxtimes_{P(z_{2})}W_{3}))$,
\[
(\boxtimes_{P(z_{1})}\circ (1_{W_{1}}\otimes
\boxtimes_{P(z_{2})}))'(w') \in (W_1\otimes W_2\otimes W_3)^*
\]
satisfies the $P^{(1)}(z_{2})$-local grading restriction condition (or
the $L(0)$-semisimple $P^{(1)}(z_{2})$-local grading restriction
condition when $\mathcal{C}$ is in $\mathcal{M}_{sg}$) and the other
condition in Condition 2 in Theorem \ref{expansion}.  Since $\eta$ and
${\cal A}_{P(z_{1}), P(z_{2})}^{P(z_{1}-z_{2}), P(z_{2})}$ are module
maps, for $w_{(4)}'\in W_{4}'$ we have that
\begin{eqnarray*}
\lefteqn{(\overline{\eta}\circ \overline{{\cal A}_{P(z_{1}),
P(z_{2})}^{P(z_{1}-z_{2}), P(z_{2})}}\circ
(\boxtimes_{P(z_{1})}\circ (1_{W_{1}}\otimes \boxtimes_{P(z_{2})})))'(w_{(4)}')}\nn
&&=(\boxtimes_{P(z_{1})}\circ (1_{W_{1}}\otimes \boxtimes_{P(z_{2})}))'
(({\cal A}_{P(z_{1}),
P(z_{2})}^{P(z_{1}-z_{2}), P(z_{2})})'(\eta'(w_{(4)}')))
\end{eqnarray*}
also satisfies Condition 2 in Theorem \ref{expansion}.  Thus by
(\ref{conv-1}), $(I^{1}\circ (I^{2}\otimes 1_{W_{3}}))'(w_{(4)}')$
satisfies this condition, that is, the expansion condition holds.
\epfv


\bigskip

\noindent {\small \sc Department of Mathematics, Rutgers University,
Piscataway, NJ 08854 (permanent address)}

\noindent {\it and}

\noindent {\small \sc Beijing International Center for Mathematical Research,
Peking University, Beijing, China}

\noindent {\em E-mail address}: yzhuang@math.rutgers.edu

\vspace{1em}

\noindent {\small \sc Department of Mathematics, Rutgers University,
Piscataway, NJ 08854}

\noindent {\em E-mail address}: lepowsky@math.rutgers.edu

\vspace{1em}

\noindent {\small \sc Department of Mathematics, Rutgers University,
Piscataway, NJ 08854}

\noindent {\em E-mail address}: linzhang@math.rutgers.edu

\end{document}